\documentclass[graybox]{svmult}

\usepackage{mathptmx}       
\usepackage{helvet}         
\usepackage{courier}        
\usepackage{type1cm}        
\usepackage{makeidx}         
\usepackage{graphicx}        
\usepackage{multicol}        
\usepackage[bottom]{footmisc}

\usepackage{amsmath}
\usepackage{amssymb}
\usepackage{amscd}
\usepackage{indentfirst}

\usepackage{cite}
\usepackage[ruled,linesnumbered]{algorithm2e}
\usepackage{benstyle_springer}

\makeindex             

\begin{document}

\title*{Towards optimal sampling for learning sparse approximations in high dimensions}
\titlerunning{Towards optimal sampling for learning sparse approximations in high dimensions}
\author{Ben Adcock, Juan M.\ Cardenas, Nick Dexter and Sebastian Moraga}
\authorrunning{B.\ Adcock, J.\ M.\ Cardenas, N.\ Dexter and S.\ Moraga}
\institute{B.\ Adcock, J.\ M.\ Cardenas, N.\ Dexter and S.\ Moraga \at Department of Mathematics, Simon Fraser University\\ \email{ben_adcock@sfu.ca, jcardena@sfu.ca, nicholas_dexter@sfu.ca, sebastian_moraga_scheuermann@sfu.ca} }

\maketitle
\setcounter{page}{1}

\abstract{In this chapter, we discuss recent work on learning sparse approximations to high-dimensional functions on data, where the target functions may be scalar-, vector- or even Hilbert space-valued. Our main objective is to study how the sampling strategy affects the sample complexity -- that is, the number of samples that suffice for accurate and stable recovery -- and to use this insight to obtain optimal or near-optimal sampling procedures. We consider two settings. First, when a target sparse representation is known, in which case we present a near-complete answer based on drawing independent random samples from carefully-designed probability measures. Second, we consider the more challenging scenario when such representation is unknown. In this case, while not giving a full answer, we describe a general construction of sampling measures that improves over standard Monte Carlo sampling. We present examples using algebraic and trigonometric polynomials, and for the former, we also introduce a new procedure for function approximation on irregular (i.e., nontensorial) domains. The effectiveness of this procedure is shown through numerical examples. Finally, we discuss a number of structured sparsity models, and how they may lead to better approximations.
}


\section{Introduction}

Learning an accurate approximation to an unknown function from data is a fundamental problem at the heart of many key tasks in applied mathematics and computer science. This problem is rendered challenging by the famous \textit{curse of dimensionality}. In many relevant applications, the domain of the function is a high-dimensional space, thus standard algorithms (those well suited in lower dimensions) often suffer from an exponential blow-up in \textit{sample complexity} (the number of samples required to obtain an accurate approximation). This is particularly problematic in many practical settings, since the amount of data available is often highly limited.

Fortunately, it is well known that functions arising in practice often possess low-dimensional structure. Specifically, they admit approximately \textit{sparse} representations, meaning that they can be efficiently approximated using a relatively small number $s$ of functions from a particular dictionary. With this in mind, the aim of this chapter is to address the following fundamental question: \textit{supposing a function has an approximately sparse representation, how many samples (of a given type) suffice to learn such an approximation from data, and how can it be computed?}

\subsection{Main problem}

Let $(D,\cD,\rho)$ be a probability space. Here $D$ is typically a subset of $\bbR^d$, where $d \gg 1$, is the dimension of the problem. We consider approximating functions defined over $D$. In many applications, such a function takes scalar values. However, other applications call for the approximation of functions that are vector- or function-space valued. To this end, in this work we let $\bbV$ be a separable Hilbert space over the field $\bbC$ and consider a function of the form
\bes{
f : D \rightarrow \bbV.
}
Note that $\bbV$ may be taken as $(\bbC,\abs{\cdot})$ in the case of scalar-valued function approximation or $(\bbC^k,\nm{\cdot}_{\ell^2})$ in the case of vector-valued function approximation. Alternatively, it may be an infinite-dimensional Hilbert space of functions. We discuss several motivations for studying this case later. Note also that we consider vector spaces over complex fields. Doing so presents a number of additional challenges over considering the real case only.

Let $L^2_{\rho}(D)$ be the Lebesgue space of complex scalar-valued, square-integrable functions on $D$. We now consider a known dictionary of functions
\bes{
\Phi = \{ \phi_{\iota} : \iota \in \cI \} \subset L^2_{\rho}(D),
}
which may be finite, countable or uncountable, and we assume that $f$ has an approximate $s$-sparse representation in $\Phi$. That is to say, there exists a set $S \subset \cI$ of size $|S| \leq s$ for which
\be{
\label{sparse_rep_intro}
f \approx f_{S} : = \sum_{\iota \in S} c_{\iota} \phi_{\iota},
}
where the coefficients $c_{\iota} \in \bbV$ are elements of the Hilbert space $\bbV$.

Motivated by the high-dimensionality of the domain $D$, our primary focus in this work is on random sampling schemes. To this end, we assume that there are probability measures $\mu_1,\ldots,\mu_m$ on $D$, we draw $m$ independent samples $y_1,\ldots,y_m$ with $y_i \sim \mu_i$, $i = 1,\ldots,m$, and we assume that the data takes the form
\be{
\label{f_data_intro}
(y_i , f(y_i) + n_i),\qquad i = 1,\ldots,m,
}
where $n_i \in \bbV$ is measurement noise. With this in hand, we may reformulate the main question stated previously as follows: \textit{how should one choose the number of samples $m$, the sampling measures $\mu_1,\ldots,\mu_m$, and the learning procedure so that an approximation to $f$ yielding an error close to that of the sparse representation $f_{S}$ can be computed from the data \R{f_data_intro}? Furthermore, is this approximation stable to measurement noise?}

Lacking any further insight, the standard random sampling strategy involves drawing samples in a Monte Carlo fashion from the underlying measure $\rho$; in other words, we let $\mu_1 = \ldots = \mu_m = \rho$. We consider this strategy the starting point for the discussion. To this end, we also consider the related question: \textit{what is the sample complexity of Monte Carlo sampling, and to what extent can this be improved by changing the sampling measures $\mu_i$?}

Note that the focus of this work is on approximations that can be computed (potentially up to some tolerance) in finite time. When $f$ takes values in an infinite-dimensional Hilbert space $\bbV$, this presents an issue. To address it, we assume the $\bbV$ can be discretized via a finite-dimensional space $\bbV_h$ (here $h > 0$ is a discretization parameter), and then proceed to perform computations in $\bbV_h$, as opposed to $\bbV$. Thus, another important question we discuss in this chapter is: \textit{what is the effect of this discretization on the ensuing approximation to $f$?}

\subsection{Overview}

The purpose of this chapter is to survey a recent body of work that has sought to answer these questions. See \S \ref{ss:literature} for a detailed summary of relevant literature. We divide our discussion into two main cases:

First, we consider the case where the target set $S$ in the sparse representation \R{sparse_rep_intro} is known. This is by far the simpler situation, yet it can indeed occur in certain problems arising in practice. For example, $S$ may be obtained by \textit{a priori} regularity estimates on $f$. Moreover, even though it may not be applicable in general, examining this case helps provide insight into what can possibly be achieved in the second setting, where $S$ is unknown. 

We provide an almost complete set of answers to the above questions in this first setting. The approximation is learned through a simple (weighted) least-squares fit, which is readily shown to provide accurate and stable approximations. We also obtain a general condition on the sampling measures $\mu_i$, as well as explicit examples satisfying such a condition, for which only
\be{
\label{LS-near-opt-samp-comp-intro}
m \gtrsim s \cdot \log(2s/\epsilon)
}
such samples suffice for recovery, with probability at least $1-\epsilon$ for some $\epsilon > 0$. This condition is optimal up to the constant implied by the $\gtrsim$ symbol and the log factor. As we also discuss, the near-optimal sample complexity bound \R{LS-near-opt-samp-comp-intro} typically does not hold in the case of Monte Carlo sampling. We discuss examples where the corresponding bound for Monte Carlo sampling can be arbitrarily large. \S \ref{s:wLS} is devoted to weighted least-squares approximation. 

Unfortunately, the first case is rather rare in practice. It is more common to encounter the situation where $S$ is unknown \textit{a priori}. To overcome this, one may seek to estimate employ adaptive sampling while building $S$ in an iterative manner, typically via a greedy scheme. While such procedures can sometimes work well in practice -- especially when $s$ is relatively small -- they often lack theoretical guarantees. 
Instead, we pursue a different approach using tools from sparse regularization, in which we seek to promote the sparsity of $f$ in the dictionary $\Phi$ via $\ell^1$-minimization-type techniques. Analysis of this case can then be performed using tools from compressed sensing theory.  

In this case, we assume that $\Phi$ is a finite set of linearly-independent elements. Let $n = | \Phi |$. Our main result on sample complexity in this case demonstrates that there exist choices of sampling measures $\mu_1,\ldots,\mu_m$ for which
\bes{
m \gtrsim (b/a) \cdot (\theta^2 /a) \cdot s \cdot \left (\log(\E n) \cdot \log^2(\E (b/a) (\theta^2 / a)  s) \ + \log(2/\epsilon) \right )
}
sample suffice for recovery, where $a,b > 0$ are the Riesz basis constants of $\Phi$ (see \R{Riesz-bounds}). Here $\theta$ is an explicit constant, given by
\bes{
\theta^2 = \int_{D} \max_{\iota \in \cI} | \phi_{\iota}(y) |^2 \D \rho(y) .
}
We also present a sample complexity bound for Monte Carlo sampling, which takes the form
\bes{
m \gtrsim (b/a) \cdot (\Theta^2 /a) \cdot s \cdot \left (\log(\E n) \cdot \log^2(\E (b/a) (\theta^2 / a)  s) \ + \log(2/\epsilon) \right ),
}
where
\bes{
\Theta^2 = \max_{\iota \in \cI} \nm{\phi_{\iota}}^2_{L^{\infty}_{\rho}(D)}.
}
Notice that $\theta \leq \Theta$. Hence the former strategy is always at least as good as Monte Carlo sampling. We present examples where $\theta = \Theta = 1$ (in which case, Monte Carlo sampling is sufficient) and where $\theta \ll \Theta$ (in which case, the former strategy is strictly better).

\subsection{Additional contributions}

In tandem with the various sample complexity bounds, we also present error bounds for the learned approximations. These show that such approximations are \textit{accurate} -- i.e.\ the error is bounded by the best approximation error $f - f_S$, measured in some norm -- and \textit{stable} to noise, i.e.\ the error scales linearly with the noise values $n_i$. In the Hilbert-valued setting, we also determine stability to discretization error, in the sense that the error involves an additional term that is proportional to the orthogonal projection onto $\bbV_h$.

Several of our examples consider function approximation on tensor-product domains, such as the {symmetric hypercube} $D = [-1,1]^d$ in $d$ dimensions. However, certain practical applications result in approximation problems on irregular domains. Another contribution of this chapter is to introduce a new approach for function approximation on irregular domains via sparse regularization. We demonstrate the efficacy of this new approach both through theoretical guarantees and numerical examples.

Finally, we also discuss settings where $f$ admits a structured sparse approximation in the dictionary $\Phi$. Such representations arise frequently in practice, and can lead to tangible benefits in accuracy. We consider two such models, \textit{weighted sparsity} and \textit{lower set sparsity}, and briefly describe the extension of the main results to these settings. Focusing on the irregular domain case, we also showcase the benefits of such structured sparsity models via numerical examples.

\subsection{Related literature}\label{ss:literature}

This work is motivated in great part by applications arising in \textit{parametric models}. Here, one seeks to understand how the parameters in a physical model --  a weather or climate model, a chemical or biological process, a fluid flow model such as groundwater flow, a nuclear reactor, an aircraft engine, etc --  affect its output. Parametric models are ubiquitous in engineering and the physical sciences. Approximating the input-output map of a parametric model is a problem that lies at the heart of many key tasks in parametric modelling, such as performing uncertainty quantification, parameter optimization or solving parametric inverse problems. See \cite{ghanem2017handbook,le2010spectral,smith2013uncertainty,sullivan2015introduction} for detailed introductions to this topic.

Parametric models are often formulated as (systems of) DEs. In such problems, the function $f$ is the solution $u$ of a PDE system of the form
\be{
\label{UQ_model_problem}
\mathcal{L}_{x} (u, y) = 0,
}
defined over a physical domain $\Omega$ and subject to suitable boundary conditions, where $\mathcal{L}_x(\cdot,y)$ denotes a differential operator in the physical variable $x$ which depends on $y$. Therefore, the solution $u$ is also a function defined over $\Omega \times D$, and for each fixed $y\in D$ the solution $u(\cdot, y)$ is an element of a function space $\bbV$. Note here that $\bbV$ may be a Hilbert or a Banach space, depending on the particular form of \eqref{UQ_model_problem}. The typical goal in such settings is then to compute a quantity of interest (QoI) $Q: \bbV \to \bbR$ depending on $u$, e.g., the expectation or variance of $u$ with respect to $y$ at certain points $(x_i)$ or the integral of $u$ with respect to the physical variable $x$ as a function of $y$. Depending on the task at hand, a number of quantities of interest may be required, and in such scenarios computation of a fast surrogate of the full parameter-to-solution map $y\mapsto u(\cdot, y) \in \bbV$ is desirable. 

Generally speaking, evaluating $u(\cdot,y)$ (or some QoI of $u$) at a fixed value of $y$ is expensive. This either involves a costly physical experiment, or a computationally-intensive numerical simulation to (approximately) solve \R{UQ_model_problem}. Hence, the objective approximate $u$, or some QoI, from as few sample values
\bes{
u(\cdot,y_1),\ldots,u(\cdot,y_m),
}
as possible.
There are many different approaches to effect such an approximation, many of which seek to exploit low-dimensional structure of the solution $u$, typically in the form of sparsity with respect to a dictionary. Amongst the most popular methods are those which use a basis of algebraic polynomials (termed \textit{polynomial chaos} expansions in uncertainty quantification), which are motivated by the fact that solutions of many parametric DEs \R{UQ_model_problem} are smooth functions of their parameters. But there are also techniques based on multiscale or hierarchical bases, radial basis functions, trigonometric polynomials, and various others. Furthermore, there are adaptive or learned bases methods, such as, most recently, techniques involving deep neural networks.

The systematic study of least-squares approximation in general finite-dimensional subspaces from Monte Carlo samples began with the work of \cite{cohen2013stability}, with a focus on spaces of algebraic polynomials. Other early works on algebraic polynomials include \cite{migliorati2014analysis,migliorati2013polynomial,chkifa2015discrete}. It was observed that Monte Carlo sampling can lead to large sample complexities or poor approximations, which in turn led to a series of investigations into the design of improved sampling strategies. See \cite{adcock2021sparse,tang2014discrete,narayan2017christoffel,hampton2015coherence,zhou2015weighted,zhou2014multivariate,migliorati2015analysis,fajraoui2017sequential,hadigol2018least,shin2016nonadaptive,zein2013efficient,dolbeault2020optimal} and references therein. The matter of optimal sampling was theoretically resolved in \cite{hampton2015coherence} for specific polynomial subspaces, and later  \cite{cohen2017optimal} for general spaces. However, drawing samples from the resulting measures may not always be straightforward in practice. The measures are also {nonadaptive.} This led to various further extensions, including adaptive strategies \cite{arras2019sequential,migliorati2019adaptive}, more practical approaches based on discrete measures \cite{dolbeault2020optimal,adcock2020nearoptimal,migliorati2021multivariate} and recent work on boosting \cite{haberstich2019boosted,dolbeault2020optimal}. For other reviews of this topic, see \cite{adcock2021sparse,cohen2018multivariate,hadigol2018least,guo2020constructing}.

The application of $\ell^1$-minimization for computing sparse polynomial approximations of functions was first considered in \cite{blatman2011adaptive,rauhut2012sparse,doostan2011nonadapted,mathelin2012compressed,yan2012stochastic}. This led to substantial amounts of subsequent research, including \cite{tran2018analysis,yang2013reweighted,rauhut2017compressive,tsilifis2019compressive,yang2018sliced,yang2016enhancing,yang2019general,jakeman2015enhancing,alemazkoor2017divide,hampton2018basis,luthen2021sparsesolvers,guo2017sparse,yan2017sparse,tran2019class,xu2020analysis,choi2021sparse,choi2021sparse,tang2013methods,peng2016polynomial,adcock2019compressive,guo2017gradient,shin2016correcting,adcock2019correcting,adcock2018compressed2,ho2020recovery,bouchot2017multilevel,ng2012multifidelity}. Specific extensions to weighted and lower sparsity models were developed in \cite{adcock2021sparse,adcock2020sparse,adcock2017infinite,adcock2018infinite,adcock2019correcting,chkifa2018polynomial,peng2014weighted,rauhut2016interpolation,yang2013reweighted}. The generalization to Hilbert-valued functions was considered in \cite{dexter2019mixed}.
As in the case of least squares, Monte Carlo sampling can lead to poor sample complexity bounds. Thus, a series of works considered improved sampling strategies \cite{hampton2015compressive,alemazkoor2018near-optimal,xu2014sparse,tang2014subsampled,jakeman2017generalized,guo2017stochastic,liu2016stochastic,diaz2018sparse}. Weighted and lower set sparsity were developed in series of works \cite{adcock2018infinite,adcock2018compressed,chkifa2018polynomial,rauhut2016interpolation}.
For additional reviews of this topic, see \cite{adcock2021sparse,hampton2017compressive,adcock2018compressed,narayan2015stochastic,kougioumtzoglou2020sparse,luthen2021sparsesolvers,luthen2021sparseliterature}.

\subsection{Outline}

The remainder of this chapter surveys the topic of constructing sparse approximations to scalar- or Hilbert-valued functions from sample values via least squares or $\ell^1$-minimization. Our focus is on the question of sampling, and, in particular, whether or not optimal sampling can be achieved. We combines ideas from many of the aforementioned works, which are generally specific to polynomial approximations, and describe them in the setting of general dictionaries of functions. 

The outline of the remainder of this chapter is as follows. First, in \S \ref{s:preliminaries} we introduce various preliminary {concepts and notation}. We then formalize the main problem and three main questions, and introduce the main examples considered later to highlight the main results. Next, in \S \ref{s:wLS} we consider least-squares approximation. We provide definitive answers to all three main questions, and present several numerical examples. In \S \ref{s:CS} we consider $\ell^1$-minimization. We present a series of theoretical results and then describe the extent to which they resolve the three main questions. In \S \ref{s:lower-sets} we consider the extension to weighted and lower set sparsity models. Finally, we end in \S \ref{s:conclusions} with some conclusions and open problems.

\section{Preliminaries}\label{s:preliminaries}

In this section, we first provide some key notation, then we describe the setup and main problems in further detail.

\subsection{Notation}

As noted, throughout $(D,\cD,\rho)$ is a probability space and $\bbV$ is separable Hilbert space over the field $\bbC$ with inner product $\langle \cdot, \cdot \rangle_\bbV$ and corresponding norm $\| v\|_\bbV = \sqrt{\langle v , v \rangle_\bbV}$ for $v\in \bbV$. We write $L^2_{\rho}(D;\bbV)$ for the Lebesgue--Bochner space of functions $f : D \rightarrow \bbV$ for which the norm
\bes{
\nm{f}_{L^2_{\rho}(D ; \bbV)} : = \left ( \int_{D} \nm{f(y)}^2_{\bbV} \D \rho(y) \right )^{1/2} < \infty.
}
We also write $L^{\infty}(D ; \bbV)$ for the Lebesgue--Bochner space of functions $f : D \rightarrow \bbV$ for which the norm
\bes{
\nm{f}_{L^{\infty}_{\rho}(D;\bbV)} : = \esssup_{y \in D} \nm{f(y)}_{\bbV} < \infty.
}
Note that we also denote the classical Lebesgue spaces of complex-valued functions $f : D \rightarrow \bbC$ as $L^2_{\rho}(D)$ and $L^{\infty}_{\rho}(D)$. We write $\nm{\cdot}_{L^2_{\rho}(D)}$ and $\nm{\cdot}_{L^{\infty}_{\rho}(D)}$ for their norms, respectively. These coincide with the Lebesgue--Bochner spaces whenever $\bbV$ is taken as $\bbC$ with the obvious inner product.

Given an index set $\cI$ that is at most countable, we write $\ell^p(\cI ; \bbV)$ for the $\ell^p$ space of $\bbV$-valued sequences $(v_{\iota})_{\iota \in \cI}$ with finite $\ell^p$-norm, defined by
\eas{
\nm{v}_{\ell^p(\cI ; \bbV)} &= \left ( \sum_{\iota \in \cI} \nm{v_i}^p_{\bbV} \right )^{1/p},\quad 1 \leq p < \infty,
\\
\nm{v}_{\ell^{p}(\cI ; \bbV)} &= \sup_{\iota \in \cI} \nm{v_{\iota}}_{\bbV},\quad p = \infty.
}
When $p = 2$, we also write $\ip{\cdot}{\cdot}_{\ell^2(\cI;\bbV)}$ for its inner product. 
Note that when $\bbV = \bbC$, we write $\ell^p(\cI)$ and $\nm{\cdot}_{\ell^p(\cI)}$, or simply $\nm{\cdot}_{\ell^p}$ when the choice of $\cI$ is clear. Likewise, for $p = 2$ and $\bbV = \bbC$, we write $\ip{\cdot}{\cdot}$ for the $\ell^2$-inner product on $\ell^2(\cI)$.

As discussed above, the space $\bbV$ may be infinite dimensional. Hence, performing computations in $\bbV$ directly is often not possible. To this end, we introduce a finite-dimensional \textit{discretization} of $\bbV$, denoted by $\bbV_h$, where $h > 0$ is a discretization parameter. {We assume that $\bbV_h \subseteq \bbV$ is a subspace of $\bbV$ and} write  
\bes{
\cP_{h} : \bbV \rightarrow \bbV_{h}
}
for the orthogonal projection onto this subspace. Further, given $f \in L^2_{\varrho}(D ; \bbV)$, we write $\cP_{h} f \in L^2_{\varrho}(D ; \bbV_{h})$ for the almost everywhere defined function given by
\bes{
(\cP_{h} f)(y) = \cP_{h}(f(y)),\quad y \in D.
}
When necessary, we also employ a (not necessarily orthonormal basis) of $\bbV_h$. We write $\{ \psi_i \}^{k}_{i=1}$ for such a basis, where $k = \dim(\bbV_h)$.

Finally, we require a few additional pieces of notation. For convenience, we write $[n] : = \{1,\ldots,n\}$ for $n\in \bbN$. We also use the notation $A \lesssim B$ to mean that there exists a numerical constant $c > 0$ such that $A \leq c B$, and likewise for $A \gtrsim B$. Further, we write $A \lesssim_x B$ if $A \leq c_x B$ for some constant $c_x >0$ depending on a variable $x$, and likewise for $A \gtrsim_x B$.

\subsection{Problem and key questions}

As above, we let $\Phi = \{ \phi_{\iota} : \iota \in \cI \} \subset L^2_{\rho}(D)$ be a dictionary and $f \in L^2_{\rho}(D ; \bbV)$ be the function we seek to learn. We consider $s$-sparse representations of $f$ of the form
\be{
\label{sparse_rep}
f_{S} : = \sum_{\iota \in S} c_{\iota} \phi_{\iota},\quad c_{\iota} \in \bbV,
}
where $S \subseteq \cI$, $|S| \leq s$ is a subset of $s$ indices.
We now formalize the two main settings we consider in this work:

\begin{problem}[Sparsity in a known subset]
\label{ass:sparsity_known}
The function $f$ has an approximate $s$-sparse representation of the form \R{sparse_rep} for some known set $S$.
\end{problem}

\begin{problem}
[Sparsity in an unknown subset]
\label{ass:sparsity_unknown}
The function $f$ has an approximate $s$-sparse representation of the form \R{sparse_rep} for some unknown set $S$.
\end{problem}

As discussed above, we consider sample points drawn randomly according to probability measures $\mu_1,\ldots,\mu_m$. We term these the \textit{sampling measures}. We make the following assumption:

\assum{
[Absolute continuity and positivity]
\label{ass:abs-cont-pos}
The additive mixture
\bes{
\mu : = \frac{1}{m} \sum^{m}_{i=1} \mu_i,
}
is absolutely continuous with respect to $\rho$ and moreover its Radon--Nikodym derivative is strictly positive almost everywhere on $\mathrm{supp}(\mu)$.
}

This means that we can write
\be{
\label{mu_weight_fn}
\frac1m \sum^{m}_{i=1} \D \mu_i(y) = \frac{1}{w(y)} \D \rho(y),
} 
where $w : D \rightarrow \bbR$ is finite almost everywhere on $\mathrm{supp}(\mu)$.
We refer to $w$ as the \textit{weight} function. Note that it satisfies
\be{
\label{w_normalization}
\int_D \frac{1}{w(y)} \D \rho(y) = 1.
}
Given such sampling measures, we now draw samples $y_i \sim \mu_i$, $i = 1,\ldots,m$, independently from these measures and consider noisy data of the form  
\be{
\label{f_meas}
(y_i,f(y_i) + n_i )\in D \times \bbV_{h},\quad i = 1,\ldots,m.
}
Here, the $n_i \in \bbV$ are terms that capture the measurement error. We focus on the case where these terms are small in norm, but we do not assume they follow a specific distribution (e.g.\ Gaussian noise in the scalar or vector-valued case).
Note that we assume the measurements $f(y_i) + n_i $ are elements of the finite-dimensional space $\bbV_{h}$. Our motivation for doing so is the following. Since $f$ is $\bbV$-valued, the noiseless sample $f(y_i)$ is an element of the (potentially) infinite-dimensional Hilbert space $\bbV$. In general, this quantity cannot be stored, let alone used as the input to an algorithm for learning an approximation to $f$. Hence, we assume that the measurements are elements of the finite-dimensional subspace $\bbV_h$, which means they can be both stored -- for example, by storing their coefficients with respect to the basis $\{ \psi_i\}^{k}_{i=1}$ for $\bbV_h$ -- and used as input to a learning algorithm. Note that the quantity $n_i \in \bbV$ accounts for both the discretization error in approximating the true sample $f(y_i)$ by an element of $\bbV_h$, as well as any other errors that arise in the measurement process (e.g.\ noise, numerical error, and so forth). Further, we do not specify how measurements are processed to give elements of $\bbV_h$; we simply assume that this process yields an error that can be captured by the generic noise term $n_i$. In other words, the model \R{f_meas} is very general, and therefore sufficient for many applications.

We {now} formalize the three main questions considered in this work:

\begin{question}\label{q:acc-stab}
\textit{Suppose $f$ satisfies either Problem \ref{ass:sparsity_known} or \ref{ass:sparsity_unknown}. How does one learn an approximation to $f$ from the data \R{f_meas} that is \textit{accurate} -- i.e.\ the approximation error is bounded by the errors $f - f_S$ and, in the case of Hilbert-valued functions, $f - \cP_{h}(f)$, measured in suitable norms -- and \textit{stable} -- i.e.\ the error depends moderately on the noise values $n_i$, measured in a suitable norm?}
\end{question}

\begin{question}\label{q:num-samps}
\textit{In the setting of either Problem \ref{ass:sparsity_known} or \ref{ass:sparsity_unknown}, how many samples $m$ are sufficient to obtain such an approximation, and how does the choice of sampling measures affect this bound?}
\end{question}

\begin{question}\label{q:how-good-MC}
\textit{In the setting of either Problem \ref{ass:sparsity_known} or \ref{ass:sparsity_unknown}, does Monte Carlo sampling lead to \textit{near-optimal} sample complexity -- i.e.\ linear in $s$, up to constants and log factors? If not, is there a near-optimal choice of sampling measures?}
\end{question}

\subsection{Examples}\label{ss:examples}

We now introduce the main examples considered in this work.

\begin{example}[Trigonometric polynomial approximation on the $d$-torus]
\label{ex:trig-poly}
Let $d \geq 1$, {$D = \bbT^d$} is the unit torus in $d$ dimensions and  $\D \rho(y) = \D y$ be the uniform measure. In this example, we set $\cI = \bbZ^d$ and consider the set of functions
\bes{
\phi_{\iota}(y) = \exp(2 \pi \I \iota \cdot y),\quad \iota = (\iota_1,\ldots,\iota_d) \in \cI,\ y = (y_1,\ldots,y_d) \in D.
}
Observe that the dictionary $\Phi = \{ \phi_{\iota} : \iota \in \cI \}$ forms an orthonormal basis of $L^2_{\rho}(D)$.
\end{example}

Trigonometric polynomial approximation of smooth and periodic functions in high dimensions is a classical topic \cite{plonka2018numerical,dung2018hyperbolic,temlyakov2018multivariate}. The  nature of the Fourier basis makes it a relatively straightforward case to study, and as we see later, this also yields clear answers to Questions \ref{q:acc-stab}--\ref{q:how-good-MC}.

Unfortunately, many problems -- in particular parametric model problems -- do not involve periodic functions. Since such functions are often smooth, however, this motivates the study of algebraic polynomial approximations:

\examp{
[Algebraic polynomial approximation {in the symmetric hypercube}]
\label{ex:alg-poly}
Let $D = [-1,1]^d$ be the symmetric hypercube in $d \geq 1$ dimensions  of side length $2$, and $\D\rho(y) = 2^{-d} \D y$ be the uniform measure. We set $\cI = \bbN^d_0$ and consider
\be{
\label{tensor-Leg-def}
\phi_{\iota}(y) = p_{\iota_1}(y_1) \cdots p_{\iota_d}(y_d),\quad \iota = (\iota_1,\ldots,\iota_d) \in \cI,\ y = (y_1,\ldots,y_d) \in D,
}
where, on the right-hand side, $p_{\iota}$ denotes the one-dimensional Legendre polynomial of degree, $\iota \in \bbN_0$, normalized with respect to the one-dimensional uniform measure. Note that $p_{\iota}(y) = \sqrt{2 \iota + 1} P_{\iota}(y)$, where $P_{\iota}(y)$ is the classical Legendre polynomial with normalization $P_{\iota}(1) = 1$.
}

As observed, algebraic polynomial approximation is used widely in parametric model problems. This is motivated by the observation that many classes of parametric differential equations are holomorphic (analytic) functions of their parameters (see \cite{adcock2021sparse,hansen2013analytic,hoang2012regularity,cohen2011analytic,cohen2015approximation,chkifa2014high,chkifa2015breaking,tran2017analysis} and references therein). This means the polynomial coefficients decay rapidly, yielding approximately sparse representations in a given polynomial basis. Note that Problem \ref{ass:sparsity_known} is naturally motivated by such problems. For certain classes of parametric differential equations, one can use \textit{a priori} analysis to determine coefficient estimates, and using these obtain a candidate set $S$. However, this is not feasible for more complicated parametric differential equations, or problems where $f$ is given as a black box. In this setting we resort to Problem \ref{ass:sparsity_unknown}.

Note that it is also common to consider other systems of polynomials. Common examples include Chebyshev polynomials on $[-1,1]^d$, Laguerre or Hermite polynomials on $[0,\infty)^d$ and $\bbR^d$, respectively, or nonorthogonal polynomials such as Taylor polynomials. For succinctness we consider Legendre polynomials only, although our analysis readily extends to more general settings.

Many polynomial approximation problems are naturally formulated on compact hyperrectangles. Using a change of variables, these can all be reduced to the setting of Example \ref{ex:alg-poly}. In parametric models, this is inspired by the notion that the parameters are independent, with each one varying between a finite upper and lower value. However, this assumption can fail in practice. In parametric models, for example, there may often be dependencies between the parameters \cite{soize_physical_2004, le2010spectral, ernst_convergence_2012,jakeman2019polynomial}. 
This leads to polynomial approximation problems where the domain $D$, while still compact, is no longer a hypercube but an irregular-shaped domain. This poses a number of challenges, which we shall review later in this chapter. It motivates our third and final example:

\examp{
[Polynomial approximation on general compact domains] 
\label{ex:alg-poly-general}
Let $D \subset \bbR^d$ be a measurable set with nonzero measure. 
We assume without loss of generality that $D \subseteq [-1,1]^d$ is contained the symmetric, $d$-dimensional hypercube with side length $2$. Following a well-known approach, studied in detail in \cite{adcock2020approximating}, we then construct a polynomial dictionary by restricting the orthonormal basis of Example \ref{ex:alg-poly} to $D$. For consistency of notation, we continue to denote this dictionary as $\Phi = \{ \phi_{\iota} : \iota \in \cI \}$. We also let $\rho$ denote the uniform measure on $D$, i.e.\ $\D \rho(y) = |D|^{-1} \D y$.
An important observation in this case is that, in contrast to the previous two examples, this dictionary is not a basis of $L^2_{\rho}(D)$. Rather, it forms a \textit{frame} \cite{christensen2016introduction}. On the other hand, every finite set of elements from $\Phi$ is linearly independent -- indeed, no finite linear combination of polynomials can vanish on a set of nonzero measure --  and therefore a \textit{Riesz basis} for its span. These two observations will be particularly important later when we consider $\ell^1$-minimization techniques in the setting of Problem \ref{ass:sparsity_unknown}.
}

\subsection{Multi-index sets}

Notice that the dictionary $\Phi$ in Example \ref{ex:trig-poly}, \ref{ex:alg-poly} or \ref{ex:alg-poly-general} is indexed over a multi-index $\iota \in \cI =  \bbZ^d$ or $\iota \in \cI = \bbN^d_0$. It is useful to define a number of standard choices for finite subsets of $\cI$. In the case of Problem \ref{ass:sparsity_known} such an index set could be used as a potential choice for $S$. Whereas in Problem \ref{ass:sparsity_unknown} we see later that it is important to truncate the infinite set of multi-indices $\bbN^d_0$ or $\bbZ^d$ to some finite, but large subset in which we expect the indices of the sparse representation to belong.

Several standard subsets are the \textit{tensor product} index set
\be{
\label{TP-index}
\cI^{\mathrm{TP}}_{t} = \left\{ \iota = (\iota_k)^{d}_{k=1} : \iota_k \leq t,\ k = 1,\ldots,d \right \},
}
of \textit{order} $t \in \bbN_0$, the \textit{total degree} index set
\be{
\label{TD-index}
\cI^{\mathrm{TD}}_{t} = \left\{ \iota = (\iota_k)^{d}_{k=1} : \sum^{d}_{k=1} \iota_k \leq t \right \},
}
of order $t$, and the \textit{hyperbolic cross} index set
\be{
\label{HC-index}
\cI^{\mathrm{HC}}_{t} = \left\{ \iota = (\iota_k)^{d}_{k=1} :\prod^{d}_{k=1} (\iota_k +1) \leq t+1\right \}, 
}
of order $t$. Note that these index sets are subsets of $\bbN^d_0$, and therefore suitable for Examples \ref{ex:alg-poly} and \ref{ex:alg-poly-general}. We define analogous subsets of $\bbZ^d$ for Example \ref{ex:trig-poly} simply by replacing $\iota_k$ by its absolute value $|\iota_k|$ in \R{TP-index}--\R{HC-index}.

The choices \R{TP-index}--\R{HC-index} are commonly used in low to moderate dimensional problems when selecting the index set $S$ in the setting of Problem \ref{ass:sparsity_known}. However, their respective cardinalities grow rapidly with dimension; this is in particular true of \R{TP-index}, whose cardinality is $(n+1)^d$. This makes their applicability limited in higher dimensions, as, for a fixed maximum cardinality $s$, it may be impossible to achieve high orders, which are generally necessary and correspond to better accuracy. Indeed, in higher dimensions, it is often important to incorporate \textit{anisotropy} into the index set $S$ to take into account different rates of variation of the function in different coordinate directions. By contrast, the index sets \R{TP-index}--\R{HC-index} are \textit{isotropic}; indices in $S$ remain in $S$ when their entries are permuted. While it is possible to define anisotropic versions of each of these index sets (see, for example, \cite{back2011stochastic}), the challenge becomes to set the anisotropy parameters in an \textit{a priori} manner without knowledge of the underlying function $f$. Instead, we adopt the setting of Problem \ref{ass:sparsity_unknown} and suppose $f$ has a sparse representation in some unknown index set $S$ contained within a larger, but finite index set of the above form -- the goal then being to compute an approximation achieving a similar error as that of the sparse representation, without necessarily computing $S$ itself.

\section{Sparse approximation via (weighted) least squares}\label{s:wLS}

We first suppose that Problem \ref{ass:sparsity_known} holds and also that $m \geq s$. Let $S \subseteq \cI$, $|S| \leq s$, be the corresponding subset, and define the resulting subspace
\bes{
P_{S;\bbV} = \left \{ \sum_{\iota \in S} c_{\iota} \phi_{\iota} : c_{\iota} \in \bbV \right \} \subset L^2_{\rho}(D ; \bbV).
}
Note that if $\bbV = \bbC$, we simply write $P_S$ for subspace of complex-valued functions $P_{S} = P_{S,\bbC} = \left \{ \sum_{\iota \in S} c_{\iota} \phi_{\iota} : c_{\iota}\in \bbC \right \} \subset L^2_{\rho}(D)$. Next, we recall the discretized subspace $\bbV_{h}$ of $\bbV$ and the noisy samples \R{f_meas}. With this in hand, we follow a similar approach of \cite{cohen2017optimal} (which considers only the real scalar-valued case) and define the \textit{weighted least-squares} approximation to $f$ as:
\be{
\label{wLS_approx}
\hat{f} \in  \argmin{p \in P_{S ; \bbV_{h}}} \left \{ \frac1m \sum^{m}_{i=1} w(y_i) \nm{ f(y_i) + n_i - p(y_i) }^2_{\bbV} \right \}.
}
Here $w$ is the weight function specified in \R{mu_weight_fn}. Notice that we form an approximation in the subspace $P_{S ; \bbV_h}$, as opposed to $P_{S ; \bbV}$, since we generally cannot perform computations over the infinite-dimensional Hilbert space $\bbV$. In the scalar-valued case, we simply have $\bbV_h = \bbV = (\bbC,\abs{\cdot})$.

\subsection{Computation of the least-squares approximation}

We first describe the computation of the approximation \R{wLS_approx}.
Since any $p \in P_{S ; \bbV_{h}}$ can be expressed  $p = \sum_{\iota \in S} c_{\iota} \phi_{\iota}$ with $c_{\iota} \in \bbV_h$, we can rewrite \R{wLS_approx} as
\be{
\label{wLS_alg1}
\hat{f} = \sum_{\iota \in S} \hat{c}_{\iota} \phi_{\iota},\quad \hat{c} = (\hat{c}_{\iota})_{\iota \in S} \in \argmin{c \in \bbV^s_h} \nm{A c - {v}}_{\ell^2([m];\bbV)},
}
where 
{
\be{
\label{A-b-def-LS}
\begin{split}
A &= \frac{1}{\sqrt{m}} \left ( \sqrt{w(y_i)} \phi_{\iota_j}(y_i) \right )_{i  \in [m], j \in [s]} \in \bbC^{m \times s},\\
v &= \frac{1}{\sqrt{m}} \left ( \sqrt{w(y_i)} (f(y_i)+n_i) \right )_{i \in [m]} \in \bbV^m_h,
\end{split}
}
}
and $\{\iota_1,\ldots,\iota_s\}$ is an enumeration of the indices in $S$.
Note that we consider $A$ both as an $m \times s$ matrix and as a mapping $\bbV^s \rightarrow \bbV^m$ defined in the obvious way, i.e.\ $A c = \left ( \sum_{j \in [s]} A_{ij} c_j \right )_{i \in [m]} \in \bbV^m$ for $c = (c_j)_{j \in [s]} \in \bbV^s$. In the case of scalar-valued function approximation (i.e.\ $\bbV = \bbV_h = \bbC$), the problem \R{wLS_alg1} is a standard algebraic least-squares problem.  

Moreover, in the general Hilbert-valued case, it is a straightforward exercise to show that a solution of \R{wLS_alg1} is given by
\be{
\label{hatc-pseudoinv}
\hat{c} = A^{\dag} {v}.
}
Here $A^{\dag}$ is the pseudoinverse of $A$, or more precisely, its extension in the above manner to a mapping $\bbV^m \rightarrow \bbV^s$. In particular, if $A$ is full rank, then this is the unique solution of \R{wLS_alg1}. Now recall the basis $\{ \psi_i \}^{k}_{i=1}$ for $\bbV_h$. Observe that we can write the $i$th component ${v_i}$ of the $\bbV_h$-valued vector ${v}$ defined above as
\bes{
{v_i} = \sum_{j \in [k]} {V_{ij}} \psi_j,\quad i \in [m],
}
for scalar coefficients ${V_{ij}} \in \bbC$. Likewise, we can also write 
\bes{
\hat{c} = (\hat{c}_i)_{i \in [s]},\qquad \hat{c}_i =  \sum_{j \in [k]} \hat{C}_{ij},\quad i \in [s],
}
for scalar coefficients $\hat{C}_{ij} \in \bbC$, so that $\hat{f}$ can be expressed as
\bes{
\hat{f} = \sum_{i \in S} \sum_{j \in [k]} \hat{C}_{ij} \phi_i \otimes \psi_j.
}
Letting $\hat{C} = (C_{ij})_{i \in [s],j \in [k]}$ and ${V} = ({V_{ij}})_{i \in [m],j \in [k]}$ and 
using \R{hatc-pseudoinv}, we see that
\bes{
\hat{C} = A^{\dag} {V}.
}
Hence, the coefficients $\hat{C}$ can be computed by first computing the pseudoinverse $A^{\dag}$ and then performing the above matrix-matrix multiplication, for a total of $\ord{m s^2 + m s k}$ floating point operations. Alternatively, one could solve $k$ standard algebraic least-squares problems for the columns of ${V}$. If conjugate gradients are used, for example, the cost of obtaining a residual error of size $\eta$ is $\ord{\mathrm{cond}(A) m s k \log(\eta^{-1})}$, where $\mathrm{cond}(A)$ is the condition number of $A$. This may be more efficient in the case where $k \ll s$; in particular, the scalar-valued case, where $k = 1$.

\subsection{Accuracy, stability and sample complexity}\label{ss:acc-stab-samp-LS}

In this and the next several subsections, we investigate Questions \ref{q:acc-stab}--\ref{q:how-good-MC}. We commence with Question \ref{q:acc-stab}.
Accuracy and stability of the approximation \R{wLS_approx} is governed by the existence of a \textit{norm equivalence} over $P_S$. Specifically, we assume that
\be{
\label{RIP_over_S}
\alpha \nm{p}^2_{L^2_{\varrho}(D)} \leq \frac{1}{m} \sum^{m}_{i=1} w(y_i) | p(y_i)|^2 \leq \beta \nm{p}^2_{L^2_{\varrho}(D)},\quad \forall p \in P_{S},
}
for constants $0 \leq \alpha \leq \beta < \infty$.  In other words, the functional $p \mapsto \sqrt{\frac{1}{m} \sum^{m}_{i=1} w(y_i) | p(y_i)|^2}$ is an equivalent norm over $P_{S}$ to the $L^2_{\varrho}$-norm.
We remark also that \R{RIP_over_S} is a condition for the space $P_{S} = P_{S;\bbC}$ consisting of scalar-valued functions. As the next theorem shows, however, such a condition also determines accuracy and stability for the approximation of Hilbert-valued functions in the space $P_{S;\bbV_h}$. With this in hand, we now also define the discrete semi-inner product
\bes{
\ip{g}{h}_{\mathrm{disc}} = \frac1m \sum^{m}_{i=1} w(y_i) \ip{g(y_i)}{h(y_i)}_{\bbV},\quad g , h \in L^2_{\rho}(D;\bbV),
}
and corresponding discrete semi-norm $\nm{g}_{\mathrm{disc}} = \sqrt{\ip{g}{g}_{\mathrm{disc}}}$, $g \in L^2_{\rho}(D;\bbV)$.

\thm{
[Accuracy and stability of weighted least squares]
\label{t:LS_acc_stab}
Let $f \in L^2_{\rho}(D)$, $0 < \alpha \leq \beta < \infty$, $\{y_i\}^{m}_{i=1} \subseteq D$ , $w : D \rightarrow \bbR$ be such that $w(y_i)$ is well defined for all $i$, and suppose that  \R{RIP_over_S} holds.
Then the approximation $\hat{f}$ in \R{wLS_approx} is unique, and satisfies
\eas{
\nmu{f - \hat{f}}_{L^2_{\varrho}(D;\bbV)} \leq & \inf_{p \in P_{S;\bbV}} \left \{ \nm{f - p}_{L^2_{\rho}(D;\bbV)} + \left ( 1 + \frac{1}{\sqrt{\alpha}} \right ) \nm{f-p}_{\mathrm{disc}} \right \}
\\
& + \nmu{f - \cP_h(f)}_{L^2_{\rho}(D ; \bbV)}  + \frac{1}{\sqrt{\alpha}} \nm{e}_{\ell^2([m];\bbV)}
}
where $e = \frac{1}{\sqrt{m}} \left ( \sqrt{w(y_i)} n_i \right )^{m}_{i=1} \in \bbV^m$.
}

This result (see \S \ref{ss:proof-LS} for its proof) asserts that the error for the learned approximation $\hat{f}$ splits into three quantities. First, a best approximation error term in the subspace $P_{S;\bbV}$. Second, a \textit{space discretization} error, which accounts for the fact that the least-squares problem is formulated over $\bbV_h$ as opposed to $\bbV$, and is equal to the projection error $f - \cP_{h}(f)$. And third, a term depending on the measurement noise values $n_i$. Note that this theorem does not require the points $\{ y_i\}^{m}_{i=1}$ to be random. It holds for any fixed set of sample points whenever \R{RIP_over_S} also holds.

\begin{remark}
\label{r:LS-error}
This result has several disadvantages. First, the noise terms $n_i$ are multiplied by the weight factors $\sqrt{w(y_i)}$, meaning that noise terms corresponding to large values of $w$ are weighted more heavily. Second, the best approximation error mixes the $L^2_{\rho}(D;\bbV)$-norm (which is the norm in which the error $f - \hat{f}$ is measured) with the discrete norm $\nm{f - p}_{\mathrm{disc}}$.
When the sample points $y_i$ are random variables (as they will be below), one can use this fact to slightly modify the approximation $\hat{f}$ in a way in which error bounds involving only $\nm{f - p}_{L^2_{\rho}(D;\bbV)}$ can be obtained. We omit the details. See \cite{cohen2013stability,cohen2017optimal,cohen2018multivariate} for further information in the scalar-valued case.
\end{remark}

\begin{remark}
The reader will notice that Theorem \ref{t:LS_acc_stab} does not involve the upper constant $\beta$ in \R{RIP_over_S}. While not strictly needed for this theorem, this constant plays a role in the computation of the least-squares approximation. Indeed, it is straightforward to show that the condition number $\mathrm{cond}(A)$ is bounded by $\sqrt{\beta/\alpha}$ whenever $\{ \phi_{\iota} : \iota \in S \}$ forms an orthonormal basis for $P_{S}$. Hence, when the ratio $\beta/\alpha$ is small, the least-squares system can be solved more efficiently (when employing conjugate gradients) and its output is less affected by floating point errors.

This {property is relevant to} Examples \ref{ex:trig-poly} and \ref{ex:alg-poly}, {since they involve orthonormal bases}. On the other hand, the least-squares matrix $A$ will be poorly conditioned whenever the system $\{ \phi_{\iota} : \iota \in S \}$ is near-linear dependent. This occurs notably in Example \ref{ex:alg-poly-general} \cite{adcock2020approximating}. Perhaps counter-intuitively, this does not necessarily lead to substantial errors in the resulting least-squares approximation. In fact, whenever the infinite system of functions forms a frame (as it does Example \ref{ex:alg-poly-general}), this property endows the problem with sufficient structure to ensure accurate and stable (regularized) least-squares approximations. See \cite{adcock2020approximating} for further discussion.
\end{remark}

We now progress to the matter of sample complexity, which will lead to answers to Questions \ref{q:num-samps} and \ref{q:how-good-MC}.
As shown in \cite{cohen2017optimal} (see also \cite{adcock2020nearoptimal}) Sample complexity of the least-squares scheme is determined by the existence of a so-called \textit{weighted Nikolskii-type} inequality over $P_S$. Specifically, let $\cN(P_S,w)$ be the smallest constant such that
\be{
\label{Nikolskii-type}
\nm{p}_{L^{\infty}_{\rho}(D)} \leq \cN(P_S,w) \nm{p}_{L^2_{\rho}(D)},\quad \forall p \in P_S.
}
Again, we observe that this inequality is formulated for the space $P_{S} = P_{S;\bbC}$ of scalar-valued functions.
We remark also that $\cN(P_S,w)$ is related to the \textit{Christoffel} function $K(P_S)$ of the subspace $P_S$. Specifically,
\be{
\label{Nikolskii_Christoffel}
\cN(P_S , w) = \nm{\sqrt{w(\cdot) K(P_S)(\cdot)}}_{L^{\infty}_{\rho}(D)},
}
where $K(P_S)$ is the reciprocal of the Christoffel function of $P_S$. Let $\{ \upsilon_{\iota} \}_{\iota \in S} \subset L^2_{\rho}(D)$ be any orthonormal basis for $P_S$. Then this function has the explicit expression
\be{
\label{Christoffel_fn}
K(P_S)(y) = \sum_{\iota \in S} | \upsilon_{\iota}(y) |^2.
}

\thm{
[Sample complexity of weighted least squares]
\label{t:LS_samp_comp}
Let $0 < \epsilon < 1$, $0 < \delta < 1$, $\mu_1,\ldots,\mu_m$ be probability measures on $D$ satisfying Assumption \ref{ass:abs-cont-pos} and $y_1,\ldots,y_m$ be independent with $y_i \sim \mu_i$ for $i = 1,\ldots,m$. Suppose that
\be{
\label{LS_samp_comp_bd}
m \geq c_{\delta} \cdot \left ( \cN(P_S,w) \right )^2 \cdot \log(2s/\epsilon),\qquad c_{\delta} = ((1-\delta) \log(1-\delta) + \delta)^{-1},
}
where $w$ is the weight function specified in \R{mu_weight_fn}.
Then \R{RIP_over_S} holds with $1-\delta \leq \alpha \leq \beta \leq 1+\delta$, with probability at least $1-\epsilon$.
}

This theorem (see \S \ref{ss:proof-LS} for the proof) states that the sample complexity is dominated by the behaviour of the weighted Nikolskii constant $\cN(P_S,w)$. Observe that 
\be{
\label{Nikolskii-lower-bound}
(\cN(P_S , w))^2 \geq s,
}
for any choice of $w$. Indeed, $(\cN(P_S,w))^2 \geq w(y) K(P_S)(y)$ 
for almost every $y$, and therefore
\bes{
(\cN(P_S,w))^2 \int_{D} \frac{1}{w(y)} \D \rho(y) \geq \int_{D} K(P_S)(y) \D \rho(y).
}
The left-hand side is equal to $(\cN(P_S,w))^2$ due to \R{w_normalization}, and the right hand side is equal to $s$, due to the relation \R{Christoffel_fn} and the fact that the $\upsilon_{\iota}$'s are orthonormal.

\subsection{Monte Carlo sampling}\label{MC-sampling-LS}

We are now ready to discuss the first part of Question \ref{q:how-good-MC} in the context of the examples introduced in \S \ref{ss:examples}.
Recall that Monte Carlo sampling corresponds to setting
\bes{
\mu_1 = \ldots = \mu_m = \rho.
}
In this case, it follows from \R{mu_weight_fn} that the function $w(y) \equiv 1$. Hence, $\hat{f}$ is a standard unweighted least-squares approximation. As shown by Theorem \ref{t:LS_samp_comp}, the sample complexity
\be{
\label{MC-LS-sample-complexity}
m \geq c_{\delta} \cdot (\cN(P_S))^2 \cdot \log(2 s /\epsilon),
}
is governed by the unweighted Nikolskii constant
\be{
\label{unweighted-Nikolskii}
\cN(P_S) = \nm{\sqrt{K(P_S)(\cdot)}}_{L^{\infty}_{\rho}(D)}.
}
We are interested in the behaviour of $(\cN(P_S))^2$ in relation to $s = |S|$. Clearly, there are instances where $(\cN(P_S))^2$ attains the optimal value $(\cN(P_S))^2 = s$ (recall \R{Nikolskii-lower-bound}).  Indeed, the functions $\phi_{\iota}$ of Example \ref{ex:trig-poly} are orthonormal and equal to one in absolute value. Hence $K(P_S) \equiv s$ by \R{Christoffel_fn}, and \R{MC-LS-sample-complexity} yields the sample estimate $m \gtrsim s \cdot \log(2 s/\epsilon)$, which is optimal up to constants and log factors.

Unfortunately, this desirable property does not hold in general. As the next result attests, the constant $\cN(P_S)$ can generally be arbitrarily large in comparison to $s$:

\lem{
There exists a probability space $(D,\cD,\rho)$ such that following holds. For every $s \in \bbN$ and $C>0$ there exists a subspace $P \subset L^2_{\rho}(D)$ of dimension $s$ such that $\cN(P) \geq C$.
}
\prf{
We consider Example \ref{ex:alg-poly} in the case $d = 1$. The classical Legendre polynomial $P_{\iota}$ attains its maximum value at $y = 1$ and takes value $P_{\iota}(1) = 1$. Hence,
\be{
\label{1D-Leg-bound}
\nm{p_{\iota}}_{L^{\infty}_{\rho}([-1,1])} = p_{\iota}(1) = \sqrt{2 \iota + 1}.
}
It follows that for any subspace $P = P_{S}$, we have
\bes{
(\cN(P_S))^2 = \nm{K(P_S)(\cdot)}_{L^{\infty}_{\rho}(D)} = K(P_S)(1) = \sum_{\iota \in S} (2 \iota + 1).
}
Since $S \subset \bbN^d_0$, $|S| = s$ can be arbitrary, we now choose it so that the right-hand side exceeds $C$.
}

This lemma and its proof suggest that Monte Carlo sampling may be highly suboptimal in the setting of Example \ref{ex:alg-poly} (and therefore Example \ref{ex:alg-poly-general} as well) if the indices in the target set $S$ are allowed to become arbitrarily large. One way to mitigate this is to impose additional structure on $S$. A common structure is that of \textit{lower sets}:

\defn{
\label{lower-set}
A multi-index set $\cI\subseteq \mathbb{N}_0^d$ is lower if, whenever $\iota \in \cI$ and $\kappa \leq \iota$ (this inequality is understood componentwise), then $\kappa \in \cI$.
}

Note that many common index sets used in polynomial approximation are lower. For example, the sets \R{TP-index}--\R{HC-index} are all lower.
In general, lower sets are known to be good candidates for the support sets of polynomial coefficients of smooth functions in high dimensions \cite{adcock2021sparse,adcock2018infinite,adcock2018compressed,chkifa2018polynomial,chkifa2015discrete,chkifa2018polynomial,cohen2018multivariate}. In particular, this is true for solutions to parametric PDEs, where the lower set sparsity has been studied and variously exploited to construct effective polynomial approximations \cite{adcock2021sparse,cohen2018multivariate,cohen2015approximation,chkifa2015breaking,chkifa2015discrete,chkifa2014high,chkifa2013sparse}. Motivated by Example \ref{ex:trig-poly}, we observe that is also straightforward to define lower subsets of $\bbZ^d$. In this case, we replace the inequality by $| \kappa | \leq | \iota |$, where, for a multi-index $\iota = (\iota_k)^{d}_{k=1} \in \bbZ^d$, $| \iota | = ( | \iota_k |)^{d}_{k=1}$ is the multi-index of its absolute values.

In the case of Example \ref{ex:alg-poly}, it is known that when $S \subset \bbN^d_0$, $|S| \leq s$, is a lower set, one has
\bes{
(\cN(P_S))^2 \leq s^2.
}
See \cite{chkifa2014high,chkifa2015discrete}.
Furthermore, this bound is sharp, in the sense that there exists a lower set $S$ of size $s$ --  specifically, the set $S = \{ (\iota,0,\ldots,0) : \nu = 0,\ldots,s-1 \}$ -- for which $(\cN(P_S))^2 = s^2$. Hence, imposing a lower set structure reduces the sample complexity for Monte Carlo sampling to at worst quadratic in $s$, up to log factors.

\rem{
In view of Example \ref{ex:alg-poly-general}, we remark in passing that this quadratic bound also holds for arbitrary lower sets and large classes of irregular domains \cite{adcock2020approximating}, up to a domain-dependent constant. Moreover, this also holds for any Lipschitz domain in the case where $S = \cI^{\mathrm{TD}}_n$ is the total degree index set \R{TD-index} \cite{dolbeault2020optimal}. On the other hand, for domains with $C^2$ boundary and  $S = \cI^{\mathrm{TD}}_n$, one has a better scaling in higher dimensions; namely, $(\cN(P_S))^2 \leq c_D s^{1+1/d}$, where $c_D > 0$ is a constant depending on the domain $D$ only \cite{dolbeault2020optimal}.
}

\subsection{Optimal sampling}\label{ss:optimal-LS}

We now answer the second part of Question \ref{q:how-good-MC} in the {affirmative}. Our aim is to choose the weight function $w$ to minimize $\cN(P_S,w)$, and then choose the measures satisying Assumption \ref{ass:abs-cont-pos}. To do this, we appeal to \R{Nikolskii_Christoffel} and, keeping in mind the normalization \R{w_normalization}, set
\bes{
w(y) = \left ( \frac1s K(P_S)(y) \right )^{-1} = \left ( \frac1s \sum_{\iota \in S} | \upsilon_{\iota}(y) |^2 \right )^{-1}.
}
Notice that this yields, via \R{Nikolskii_Christoffel}, the optimal Nikolskii constant
\bes{
\cN(P_S,w) = \sqrt{s}.
}
In particular, the sample complexity estimate \R{LS_samp_comp_bd} becomes
\bes{
m \geq c_{\delta} \cdot s \cdot \log(2 s / \epsilon),
}
which is optimal up to the log factor.

Having chosen $w$, we now choose the measures $\mu_i$ so that \R{mu_weight_fn} holds. We consider two possibilities. The first we term \textit{nonhierarchical}, and is given simply by
\be{
\label{nonhierarchical_sampling}
\mu_1 = \ldots = \mu_m = \mu,\qquad \D \mu(y) = \frac{1}{w(y)} \D \rho(y) = \frac1s \sum_{\iota \in S} | \upsilon_{\iota}(y) |^2 \D \rho(y).
}
Clearly, \R{mu_weight_fn} holds in this case.
The second scheme is \textit{hierarchical}. In this scheme, we suppose that $m = k s$ for some $k \in \bbN$. Then we define
\be{
\label{hierarchical_sampling}
\mu_i = | \upsilon_{\iota_j}(y) |^2 \D \rho(y),\quad (j-1) k < i \leq j k,\quad j = 1,\ldots, s,
}
where $\{ \iota_1,\ldots,\iota_s\}$ is an enumeration of the indices in $S$. Notice that
\bes{
\frac1m \sum^{m}_{i=1} \D \mu_i(y) = \frac{k}{m} \sum^{s}_{j=1} | \upsilon_{\iota_j}(y) |^2 \D \rho(y) = \frac{1}{s} \sum_{\iota \in S} | \upsilon_{\iota}(y) |^2 \D \rho(y).
}
Therefore \R{mu_weight_fn} also holds in this case.

The nonhierarchical scheme \R{nonhierarchical_sampling} was introduced in \cite{cohen2017optimal} and is suitable for learning an approximation in a fixed subspace $S$. However, as discussed in \cite{arras2019sequential,migliorati2019adaptive}, it is not well suited to the problem where one seeks to learn a sequence of approximations $\hat{f}_1,\hat{f}_2,\ldots$ in a hierarchy of nested subspaces $S_1 \subseteq S_2 \subseteq \cdots$. The issue is that as the subspace $S = S_i$ changes, the measure $\mu$ defined in \R{nonhierarchical_sampling} changes, hence the existing samples are effectively drawn from the wrong distribution for the purposes of constructing an approximation in the new subspace $S_{i+1}$. The hierarchical scheme \R{hierarchical_sampling}, introduced in  \cite{migliorati2019adaptive}, overcomes this problem; see also \cite{arras2019sequential} for a different approach. 
We refer to \cite{adcock2020nearoptimal,migliorati2019adaptive} for further information.

\subsection{Practical optimal sampling via discrete measures}\label{ss:practical-optimal-samp}

Unfortunately, generating samples from either the measure \R{nonhierarchical_sampling} or the measures \R{hierarchical_sampling} may not be straightforward, since it requires an orthonormal basis $\{ \upsilon_{\iota} \}_{\iota \in S}$ of $P_S$. This may not be available in practice, and even it is, drawing samples from the resulting measures may be computationally challenging. See  \cite{cohen2017optimal,adcock2020nearoptimal,arras2019sequential,migliorati2021multivariate} for further information on this issue, as well as \cite{hampton2015coherence,narayan2018computation} for the specific case of tensor-product polynomial approximation.

A remedy to this situation was proposed in \cite{adcock2020nearoptimal,migliorati2021multivariate}. The idea is to replace $\rho$ (which is typically a continuous measure) by a discrete measure, supported on a finite grid, so that both constructing an orthonormal basis and sampling from the corresponding measures are automatically straightforward. Let $Z = \{ z_i \}^{k}_{i=1} \subset D$ be a finite grid. We consider the discrete uniform measure given by
\be{
\label{discrete-tau}
\tau = \frac1k \sum^{k}_{i=1} \delta_{z_i}.
}
The idea is now to replace $\rho$ by $\tau$ throughout. Consider the nonhierarchical scheme for simplicity. Then, doing so, we deduce that if
\be{
\label{samp-cond-disc-LS}
m \gtrsim s \cdot \log(2 s/\epsilon),
}
then the error bound
\ea{
\nmu{f - \hat{f}}_{L^2_{\tau}(D;\bbV)} \lesssim & \inf_{p \in P_{S,\bbV}} \left \{ \nm{f - p}_{L^2_{\tau}(D;\bbV)} + \nm{f-p}_{\mathrm{disc}} \right \} \nn
\\
& + \nmu{f - \cP_h(f)}_{L^2_{\tau}(D ; \bbV)}+ \nm{e}_{\ell^2([m];\bbV)}, \label{disc-err-bd-LS}
}
holds probability at least $1-\epsilon$, where
\bes{
w(y) = \left ( \frac1s \sum_{\iota \in S} | \upsilon_{\iota}(y) |^2 \right )^{-1},\qquad y \in \supp(\tau) = Z, 
}
\be{
\label{mu-disc-def}
\D \mu(y) = \frac1s \sum_{\iota \in S} | \upsilon_{\iota}(y) |^2 \D \tau(y),\qquad y \in \supp(\tau),
}
and $\{ \upsilon_{\iota} \}_{\iota \in S} \subset L^2_{\tau}(D)$ is an orthonormal basis for $P_S$ with respect to $\tau$.

Since $\tau$ is a discrete measure, this orthonormal basis can be constructed via straightforward linear algebra. Indeed, define the matrix
\bes{
B = \left ( \phi_{\iota_j}(z_i) / \sqrt{k} \right )_{i \in [k],j \in [s]} \in \bbC^{k \times s},
}
and suppose that it has the QR-factorization $B = Q R$, where $Q \in \bbC^{k \times s}$ has orthonormal columns and $R \in \bbC^{s \times s}$ is upper triangular. Then the orthonormal basis $\{ \upsilon_{\iota} \}_{\iota \in S}$ is given by
\bes{
\upsilon_{\iota_i}(y) = \sum^{i}_{j=1} (R^{-{\top}})_{ij} \phi_{\iota_j}(y),\quad i \in [s].
}
In particular, its values on the grid $Z$ are precisely
\bes{
\upsilon_{\iota_j}(z_i) = \sqrt{k} Q_{ij},\quad i \in [k],\ j \in [s].
}
Substituting this into \R{mu-disc-def} and recalling the definition of $\tau$, we see that the discrete measure $\mu$ is given by
\bes{
\D \mu(y) = \sum_{i \in [k]} \left ( \frac1s \sum_{j \in [s]} |Q_{ij} |^2 \right ) \D \delta_{z_i} (y) .
}
Hence, sampling from $\mu$ is now trivial. Indeed, $y \sim \mu$ if
\bes{
\bbP(y = z_i ) =  \frac1s \sum_{j \in [s]} |Q_{ij} |^2 ,\quad i \in [k].
}
The reader will have no doubt noticed that the error bound \R{disc-err-bd-LS} is with respect to the discrete measure $\tau$. It is often preferable to also have an error bound over the original measure $\rho$. Such an error bound is guaranteed whenever the $L^2_{\tau}$-norm is equivalent to the $L^2_{\rho}$-norm over $P_{S}$, i.e.
\be{
\label{grid-fine-cond}
\alpha' \nm{p}^2_{L^2_{\rho}(D)} \leq \nm{p}^2_{L^2_{\tau}(D)} \leq \beta' \nm{p}^2_{L^2_{\rho}(D)},\qquad \forall p \in P_{S}.
}
Indeed, recall that the sampling condition \R{samp-cond-disc-LS} and the choices of $\mu$ and $w$ imply a norm equivalence between the $L^2_{\tau}$-norm and the discrete norm over the sample points, i.e.
\bes{
(1-\delta) \nm{p}^2_{L^2_{\tau}(D)} \leq  \frac{1}{m} \sum^{m}_{i=1} w(y_i) | p(y_i)|^2 \leq (1+\delta) \nm{p}^2_{L^2_{\tau}(D)},\quad \forall p \in P_S. 
 }
Hence, we deduce that 
\bes{
(1-\delta) \alpha' \nm{p}^2_{L^2_{\varrho}(D)} \leq \frac{1}{m} \sum^{m}_{i=1} w(y_i) | p(y_i)|^2 \leq (1+\delta) \beta' \nm{p}^2_{L^2_{\varrho}(D)}  ,\quad \forall p \in P_S.
}
Therefore the norm equivalence \R{RIP_over_S} with respect to the original $L^2_{\rho}$-norm also holds, meaning that an error bound in this norm follows immediately from Theorem \ref{t:LS_acc_stab} (with constant $\alpha = (1-\delta) \alpha'$).

\begin{remark}\label{r:MC-grid}
A simple means to ensure \R{grid-fine-cond} is to construct $Z = \{z_i \}^{k}_{i=1}$ as a random Monte Carlo grid (independently of the sample points $y_i$). That is, we let the $z_i$ be independent random variables drawn according to the measure $\rho$. Observe that \R{grid-fine-cond} is precisely 
\bes{
\alpha' \nm{p}^2_{L^2_{\rho}(D)} \leq \frac1k \sum^{k}_{i=1} |p(z_i)|^2 \leq \beta' \nm{p}^2_{L^2_{\rho}(D)},\qquad \forall p \in P_{S}.
}
This is nothing more than the special case of \R{RIP_over_S} for the grid $Z$ (recall that $w \equiv 1$ for Monte Carlo sampling). Hence, \R{grid-fine-cond} is ensured by Theorem \ref{t:LS_samp_comp}. In particular, it holds with $1-\delta \leq \alpha' \leq \beta' \leq 1+\delta$, provided
\bes{
k \geq c_{\delta} \cdot (\cN(P_S))^2 \cdot \log(2 s/\epsilon),
}
where, as in \S \ref{MC-sampling-LS}, $\cN(P_S) = \nm{K(P_S)(\cdot)}_{L^{\infty}_{\rho}(D)}$ is the unweighted Nikolskii constant. Of course, this constant may be very large depending on the choice of $P_S$; recall the in discussion \S \ref{MC-sampling-LS}. Yet, this grid is only used to define the optimal sampling measure $\mu$. {Therefore, the number of grid points $k$ only affects the computational cost for generating the sample points}. It does not affect the sample complexity of the weighted least-squares approximation, which is $m \gtrsim s \log(2 s / \epsilon)$ for the optimal measure. 

Having said this, a practical problem is that estimates for $\cN(P_S)$ may not be available, or if they are, they may not be particularly tight, thus leading to overly large grids. In \cite{dolbeault2020optimal} an empirical strategy is described to mitigate this issue, based on independently drawing an auxiliary grid that is used to test the quality of the grid $Z$.
\end{remark}

\subsection{Numerical examples}\label{ss:examples-LS}

We conclude this discussion on least-squares approximation with several numerical examples. In these and other examples considered later in this chapter, we consider the scalar-valued functions
\be{
\label{functions-define}
\begin{split}
f_1(y)& = \exp\left(-\frac{1}{d}\sum_{k=1}^d y_k\right),
\\
f_2(y) &= \frac{\prod_{\lceil k=d/2\rceil + 1}^d\cos(16y_k/2^k)}{\prod_{k=1}^{\lceil d/2 \rceil}(1-y_k/4^k)},
\\
f_3(y) & = \prod^{d}_{i=1} \frac{d/4}{d/4 + (y_i + (-1)^{i+1} / (i+1))^2},
\\
f_4(y) & = \frac{1}{\sum_{i=1}^d \sqrt{|y_i|}}.
\end{split}
}
Since our goal is to compare Monte Carlo sampling with the optimal sampling procedures described above, we focus on Examples \ref{ex:alg-poly} and \ref{ex:alg-poly-general} (recall that Monte Carlo sampling is optimal, up to the log term, in the case of Example \ref{ex:trig-poly}).
To this end, we consider the domains
\be{
\label{domains-define}
\begin{split}
D_1 &= [-1,1]^d,
\\
D_2 &= \{y\in\bbR^d: 1/4\leq y^2_1 + \ldots + y^2_d \leq 1\},
\\
D_3 &= \{y\in [-1,1]^d: y_1+\ldots+y_d\leq 1\}.
\end{split}
}
We follow the approach of \S \ref{ss:practical-optimal-samp} and, in particular, Remark \ref{r:MC-grid}, to generate a Monte Carlo grid $Z$ and corresponding discrete measure $\tau$ as in \R{discrete-tau}. Here, $k = 30 s_{\max}$, where $s_{\max}$ is the maximum size of $s$ used in the given experiment. For the error, we compute the relative $L^{\infty}_{\tau}$-norm error, i.e.
\be{
\label{relative-error}
\frac{\nmu{f - \hat{f}}_{L^{\infty}_{\tau}(D)}}{\nm{f}_{L^{\infty}_{\tau}(D)}}.
}
We perform a total of $T = 50$ trials. In the Monte Carlo and optimal nonhierarchical schemes, each trial corresponds to a single draw of the sample points $y_1,\ldots,y_m$ at each value of $m$ considered. For the optimal hierarchical scheme, a single trial is a full set of points $y_1,\ldots,y_{m_{\max}}$, where $m_{\max}$ is the maximum value of $m$ considered. In all cases, we report the log-average of the error \R{relative-error} over these trials, with the shaded regions corresponding to one log-standard deviation (see \cite[App.\ A]{adcock2021sparse} for further information).

\begin{figure}[t]
	\begin{center}
	\begin{small}
 \begin{tabular}{@{\hspace{0pt}}c@{\hspace{-0.5pc}}c@{\hspace{-0.5pc}}c@{\hspace{0pt}}}
 			\includegraphics[width = 0.35\textwidth]{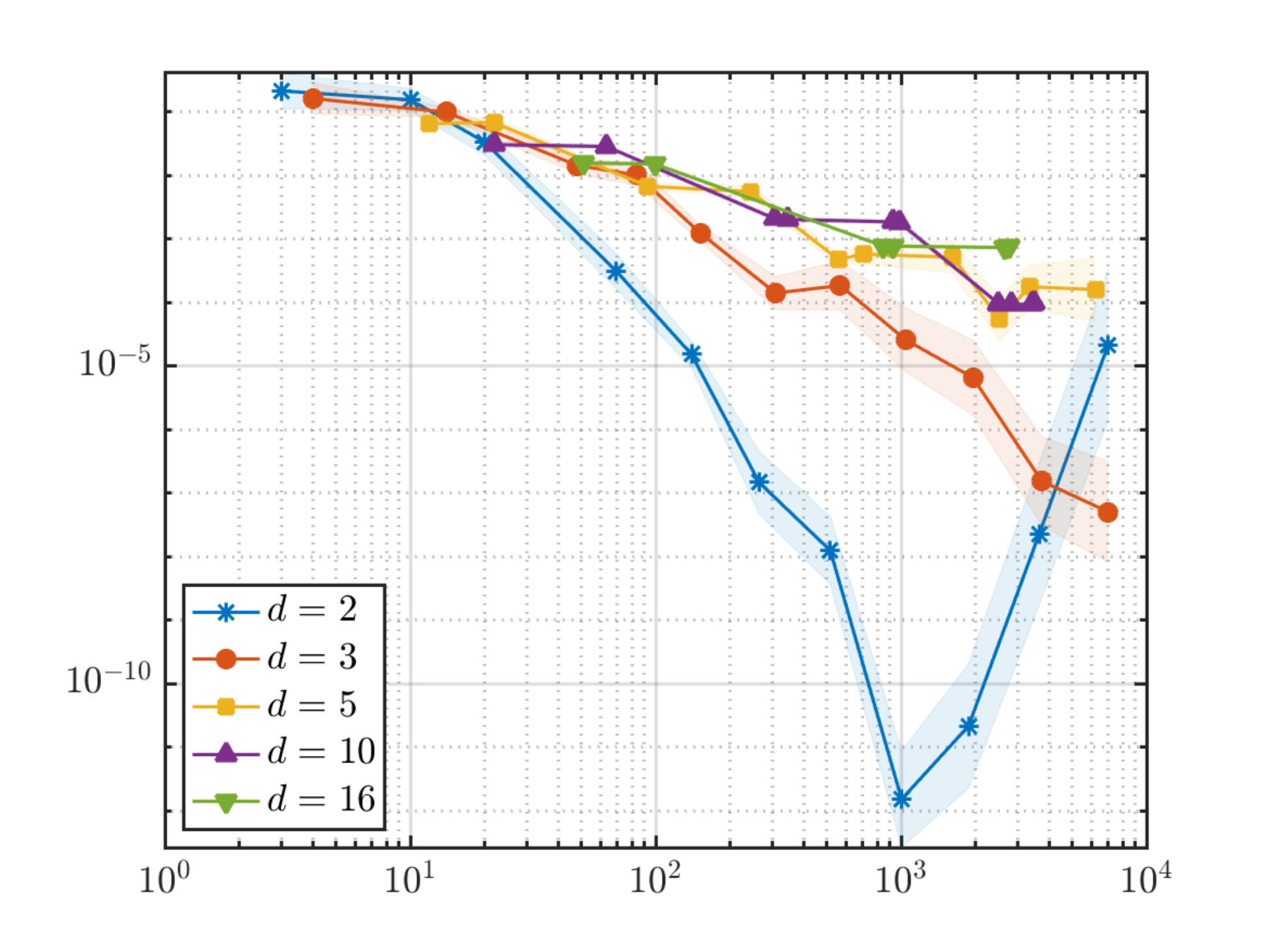} &
 			\includegraphics[width = 0.35\textwidth]{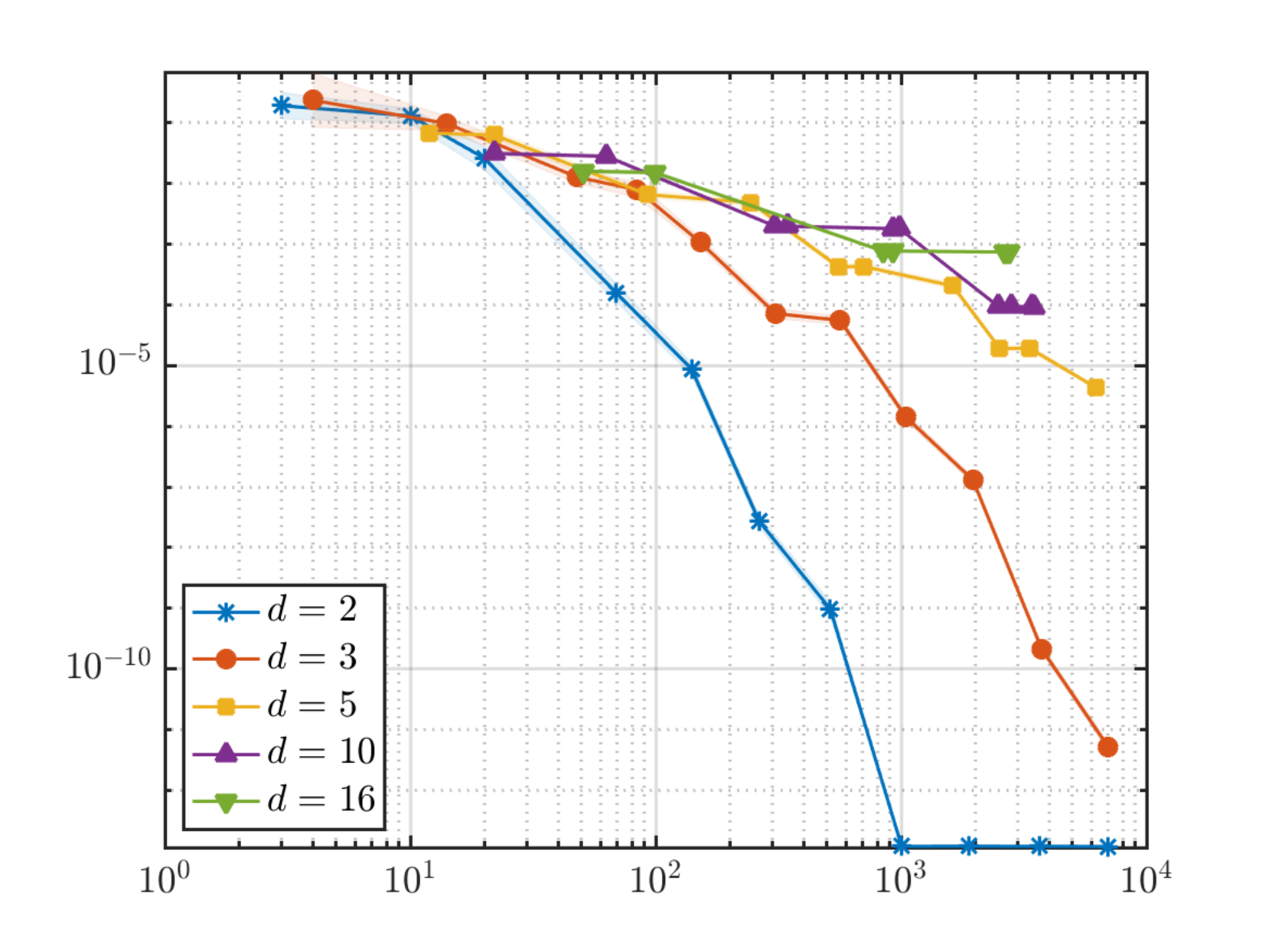} &
 			\includegraphics[width = 0.35\textwidth]{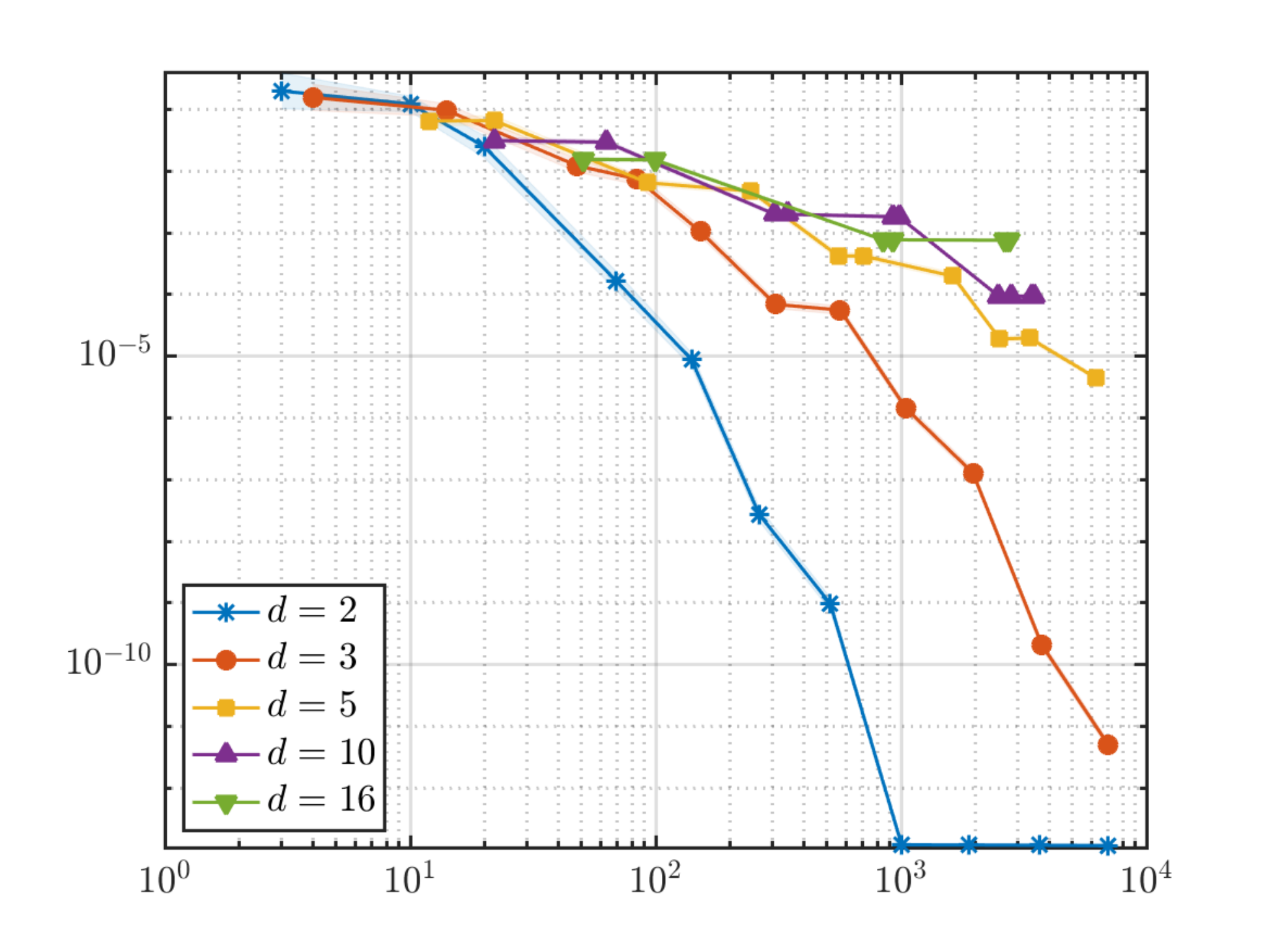}\\
 			\includegraphics[width = 0.35\textwidth]{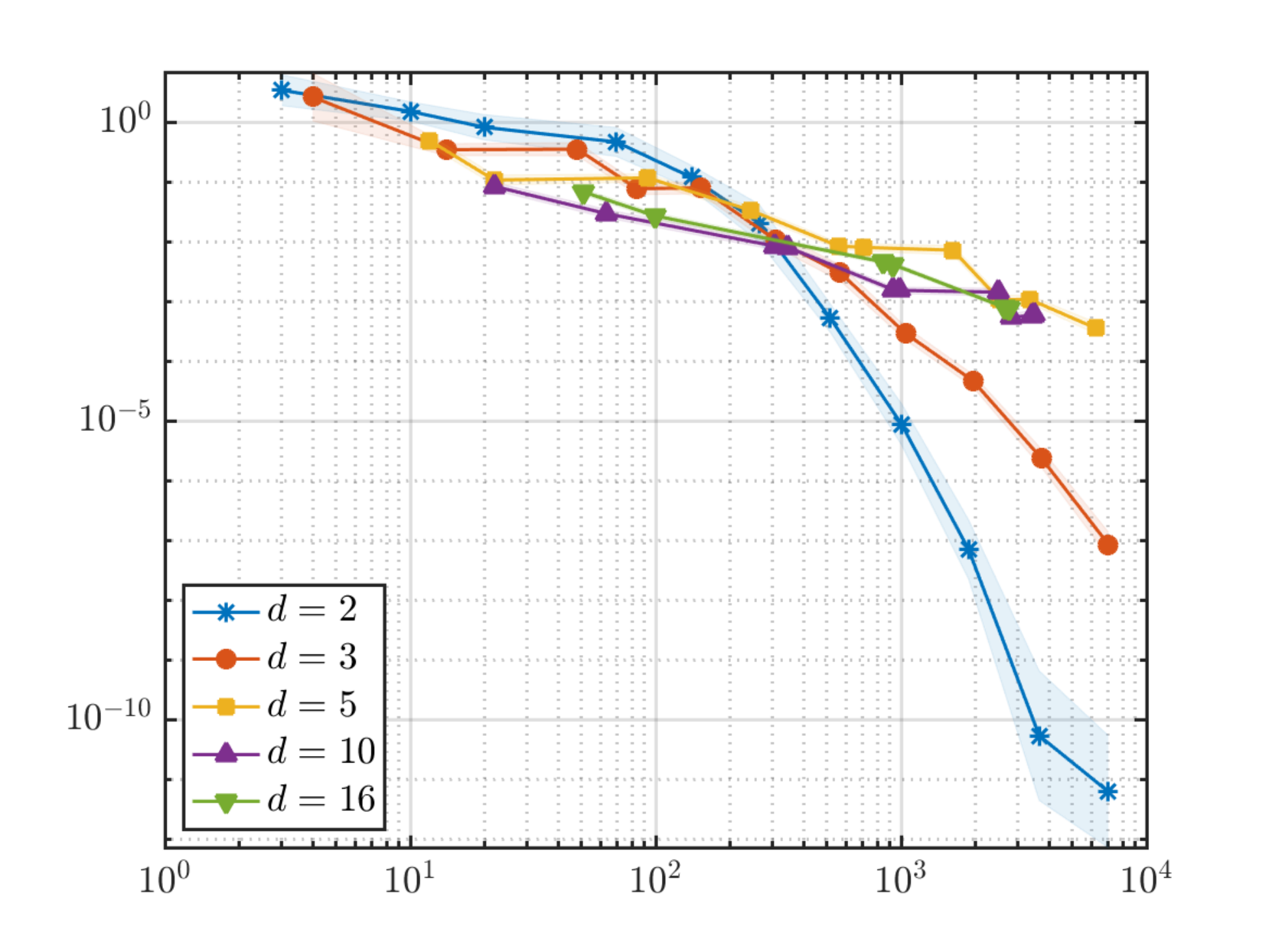} &
 			\includegraphics[width = 0.35\textwidth]{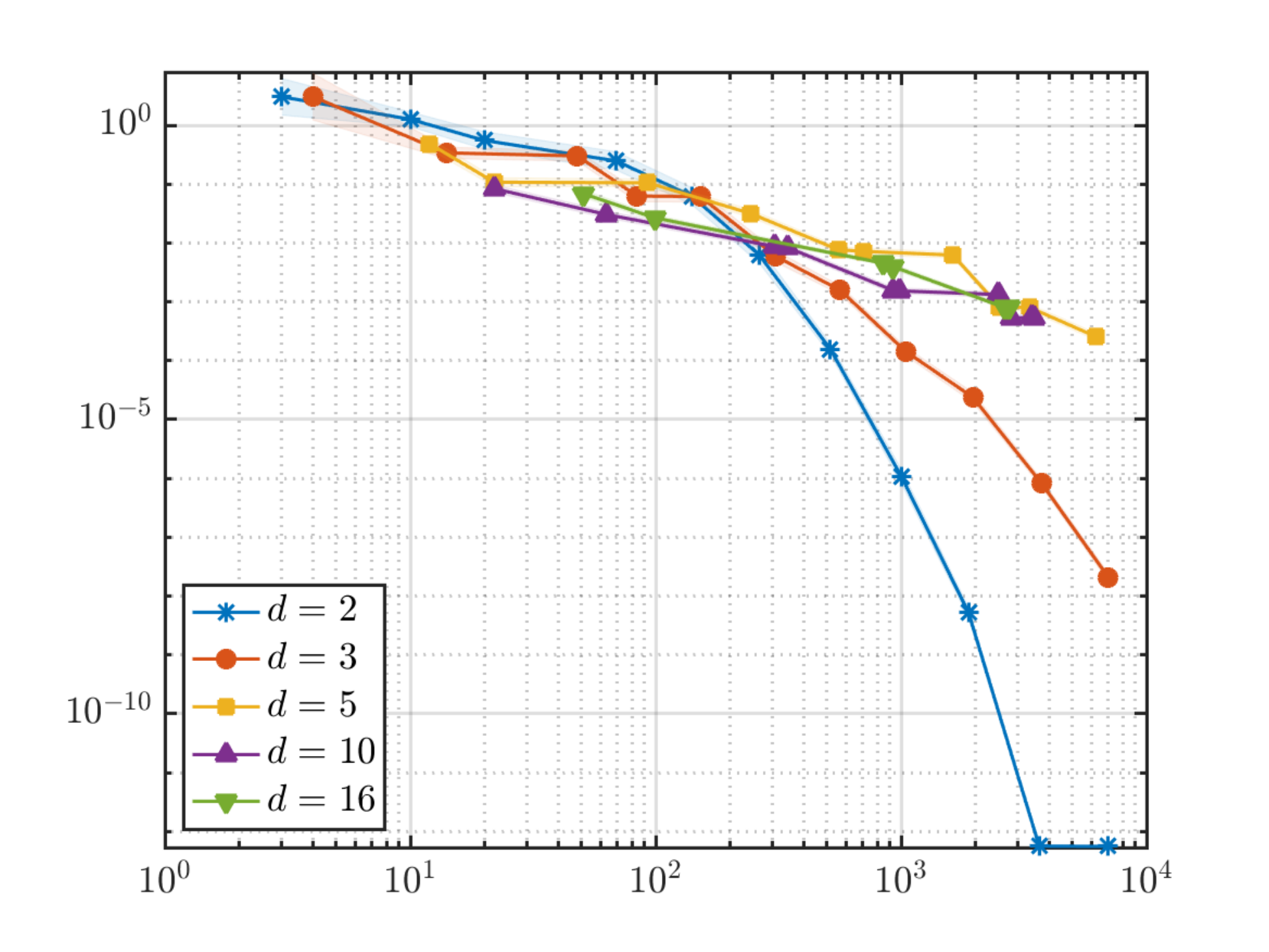} &
 			\includegraphics[width = 0.35\textwidth]{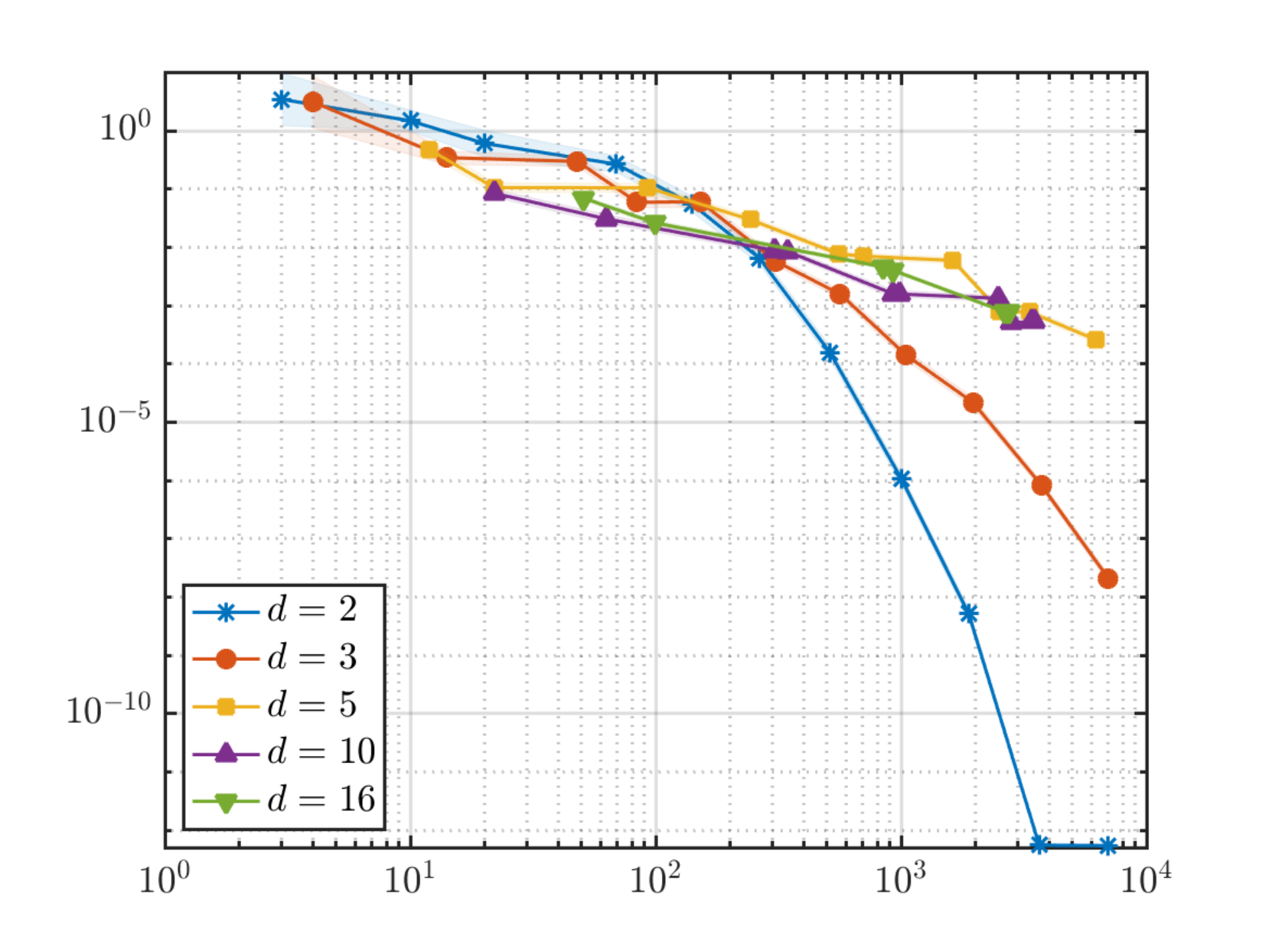} \\ 
 			\includegraphics[width = 0.35\textwidth]{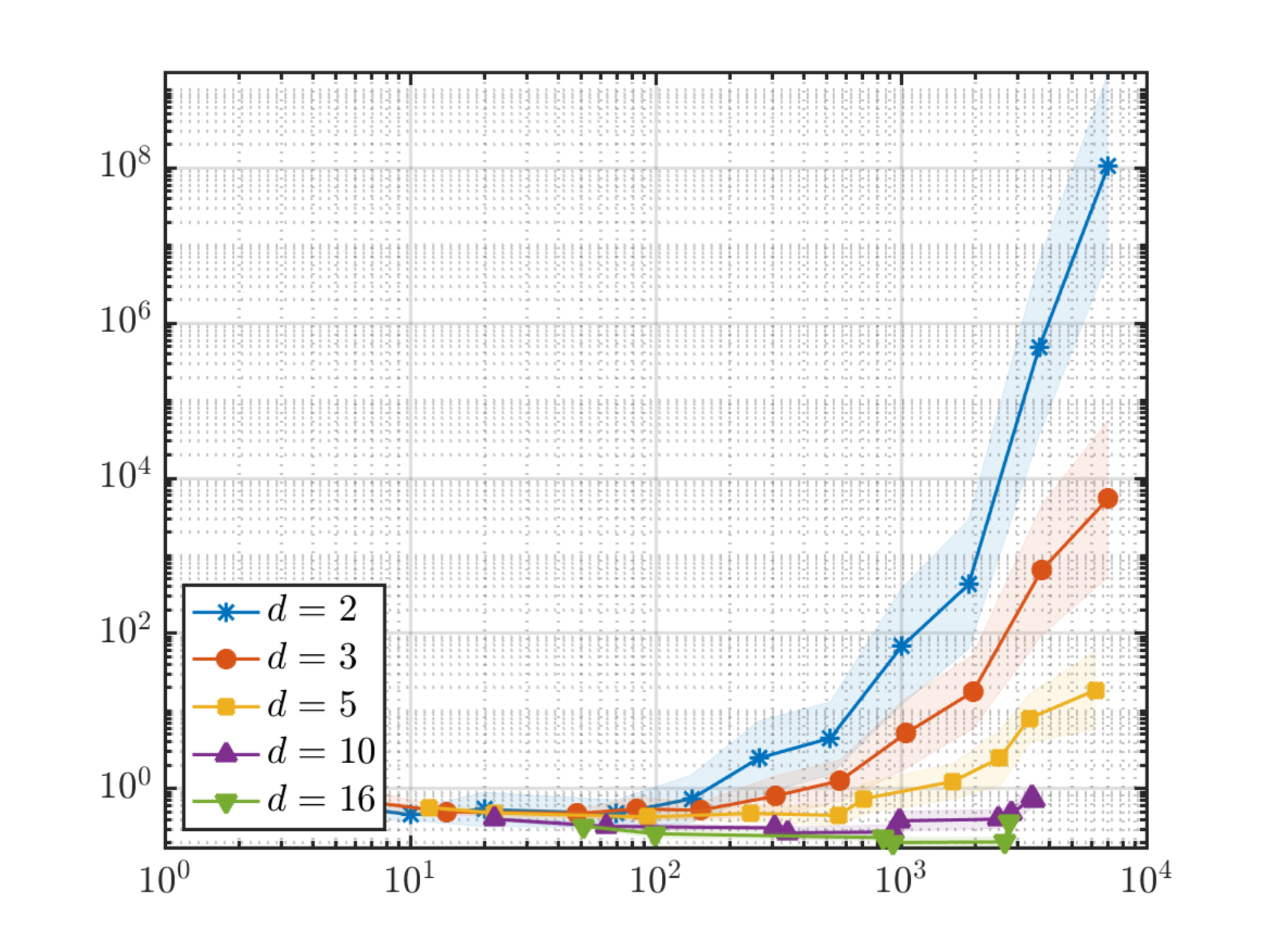} &
 			\includegraphics[width = 0.35\textwidth]{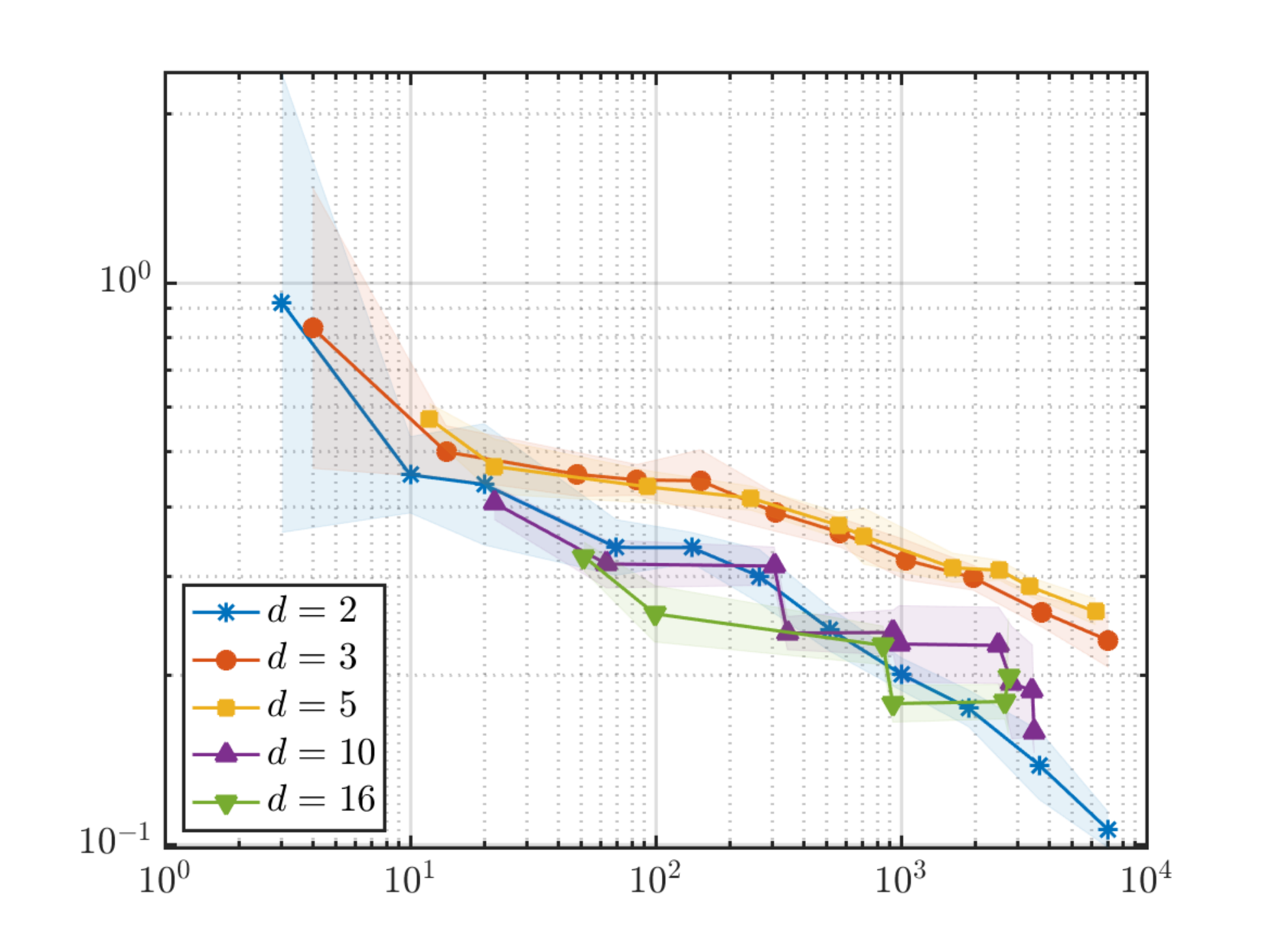} &
 			\includegraphics[width = 0.35\textwidth]{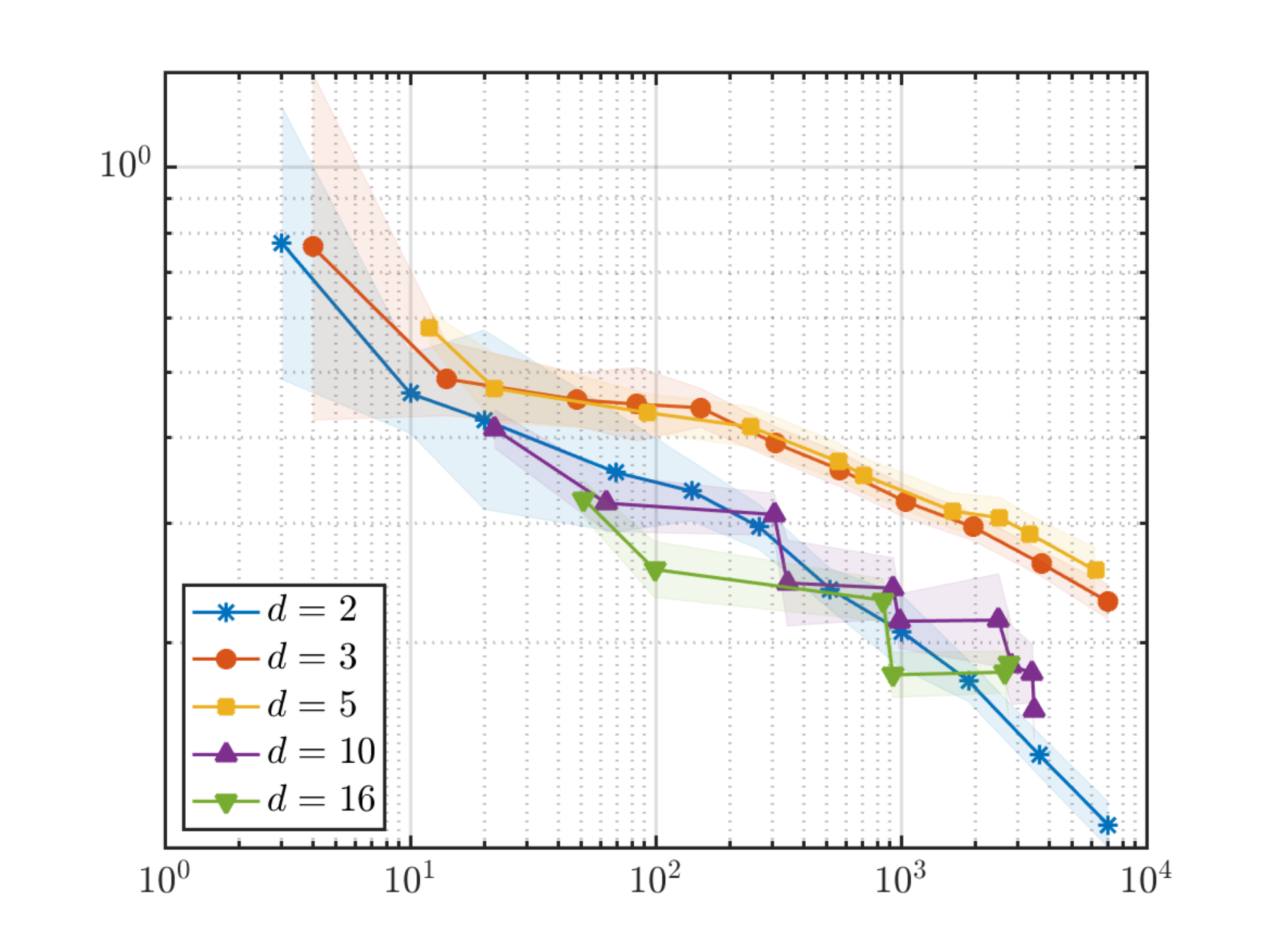} \\
		\end{tabular}
	\end{small}
	\end{center}
\caption{The relative error \R{relative-error} versus $m$ for (weighted) least squares in the case of Examples \ref{ex:alg-poly} and \ref{ex:alg-poly-general}. We compare Monte Carlo sampling (left), discrete optimal nonhierarchical sampling (middle) and discrete optimal hierarchical sampling (right) for $f = f_1$, $D = D_2$ (top row), $f = f_2$, $D = D_3$ (middle row), and $f=f_4$, $D=D_2$ (bottom row). In all experiments, $S = \cI^{\mathrm{HC}}_{t-1}$ is the hyperbolic cross index set \R{HC-index} and, for each value of $t$, $m$ is chosen as the smallest integer such that $m \geq s \log(s)$, where $s = |\cI^{\mathrm{HC}}_{t-1}|$.
} 
\label{fig:LS-1}
\end{figure}

\begin{figure}[t]
	\begin{center}
	\begin{small}
 		\begin{tabular}{ccc}
 			\includegraphics[width=0.3\textwidth]{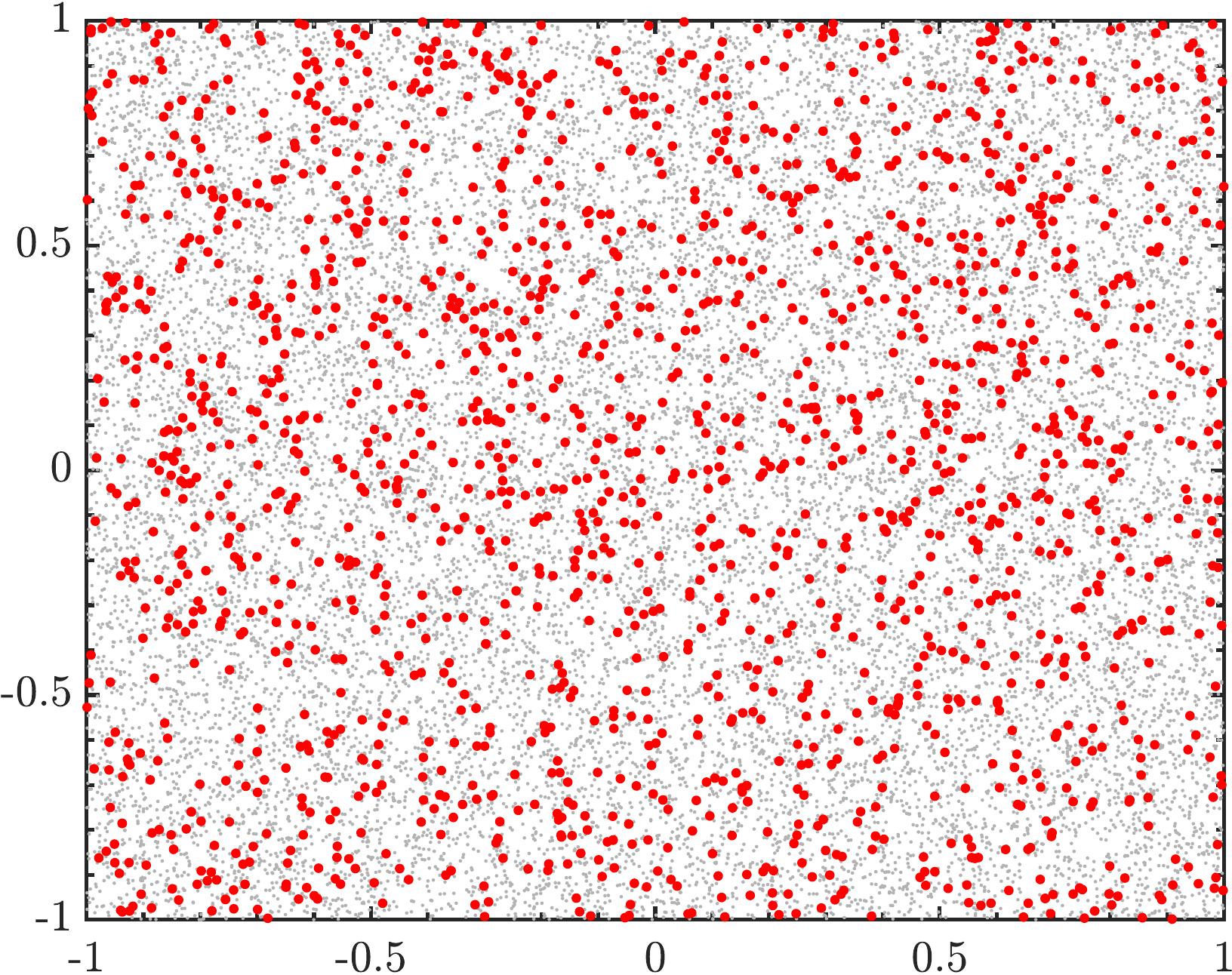} &
 			\includegraphics[width=0.3\textwidth]{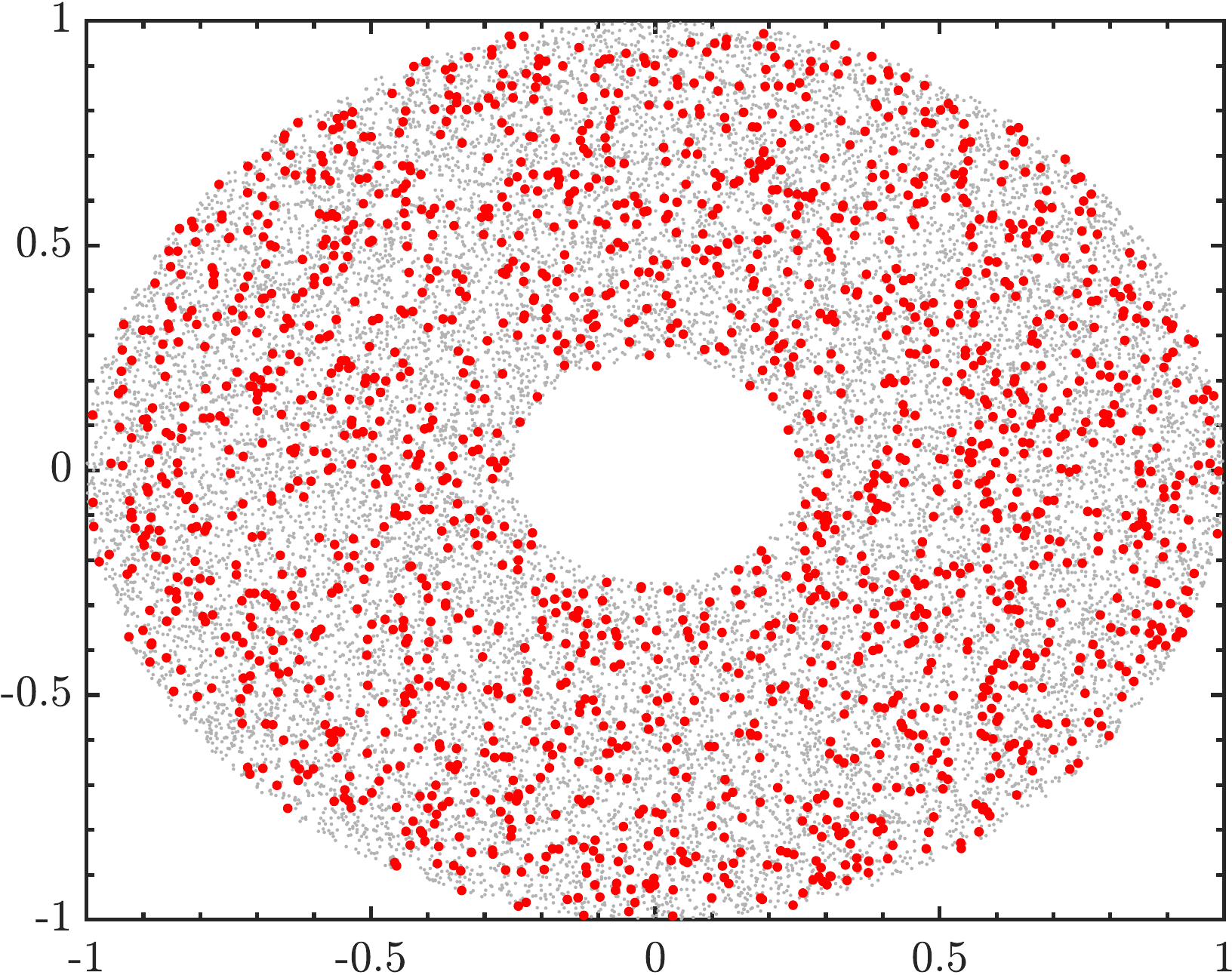} &
 			\includegraphics[width=0.3\textwidth]{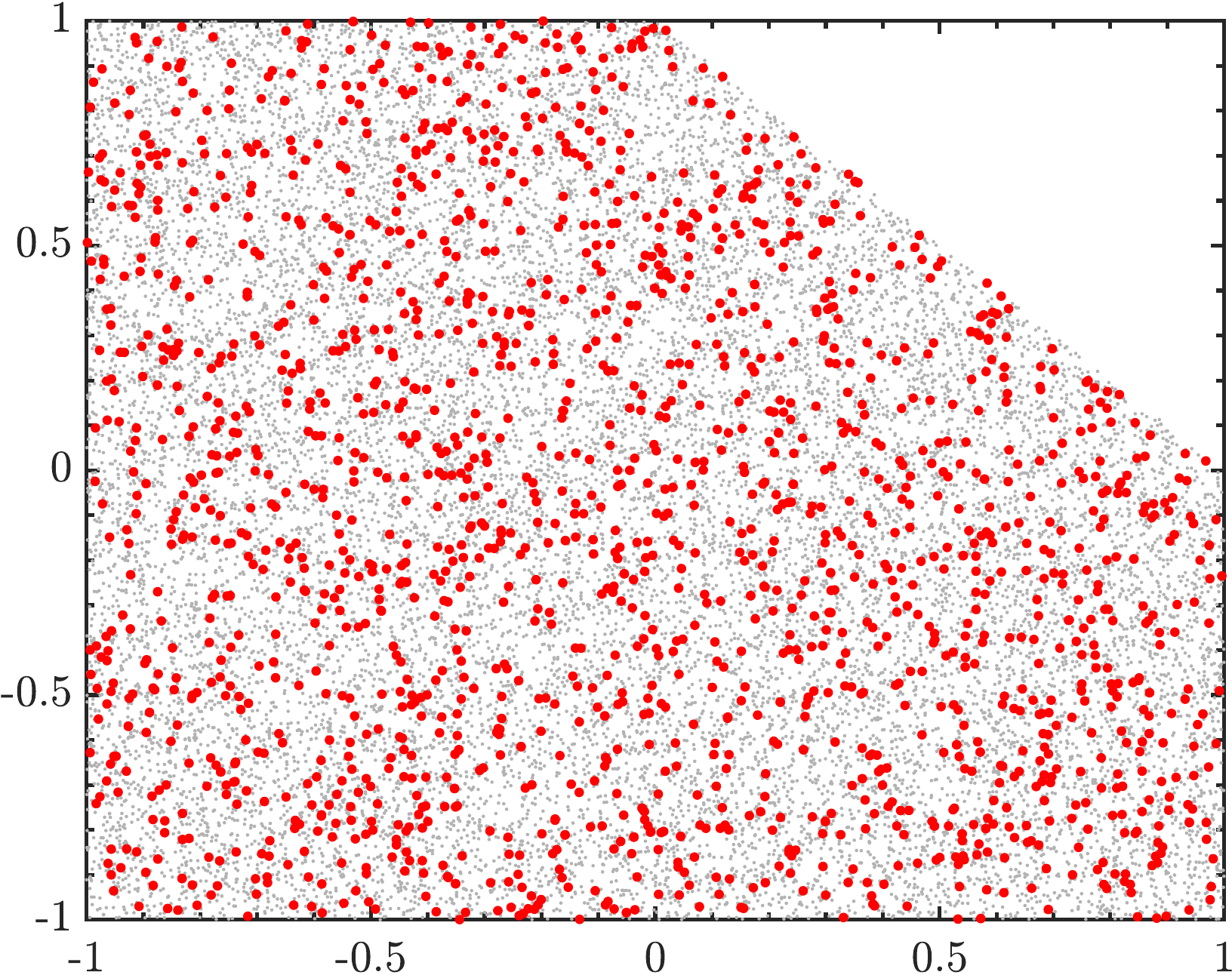} \\
 			\includegraphics[width=0.3\textwidth]{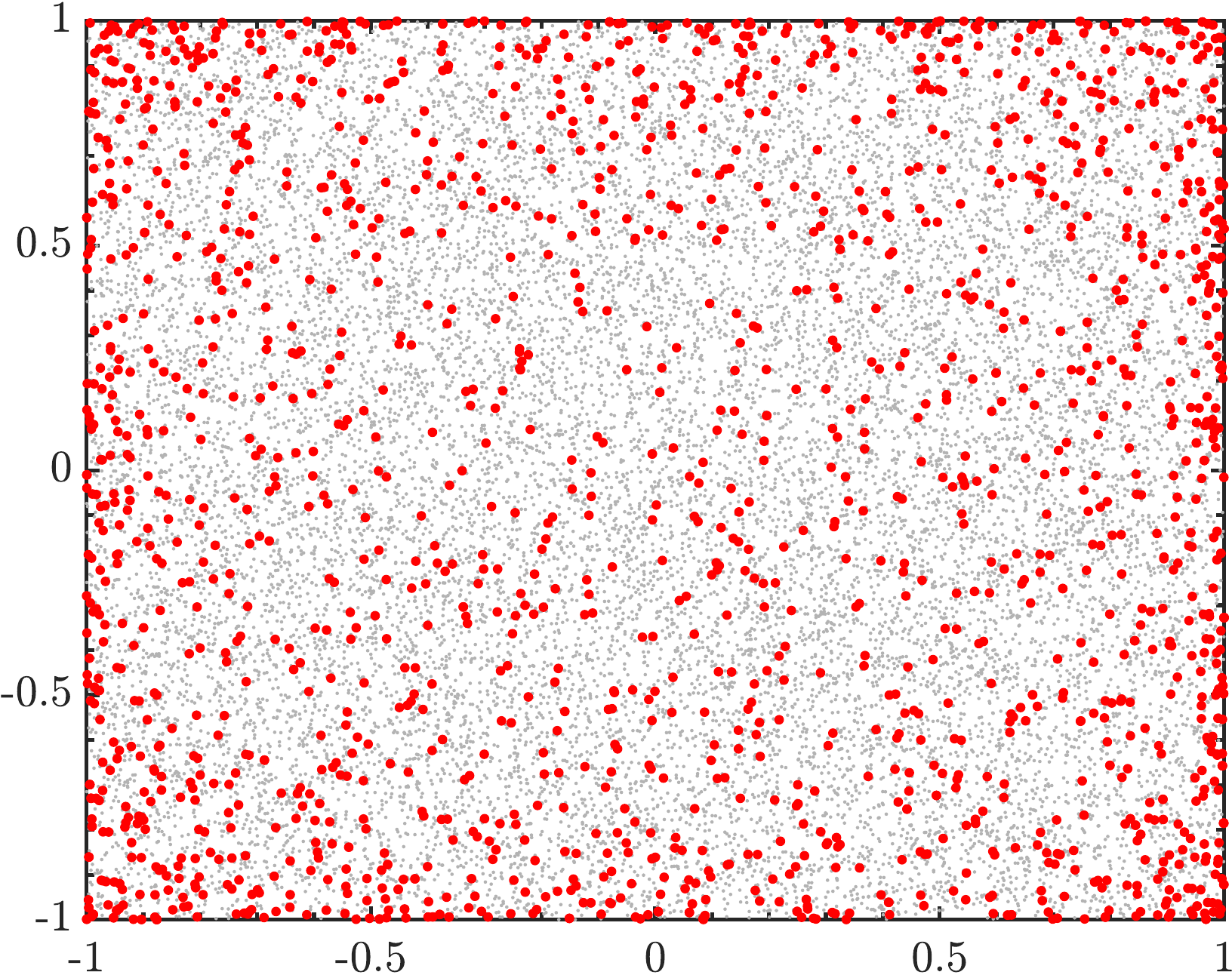} &
 			\includegraphics[width=0.3\textwidth]{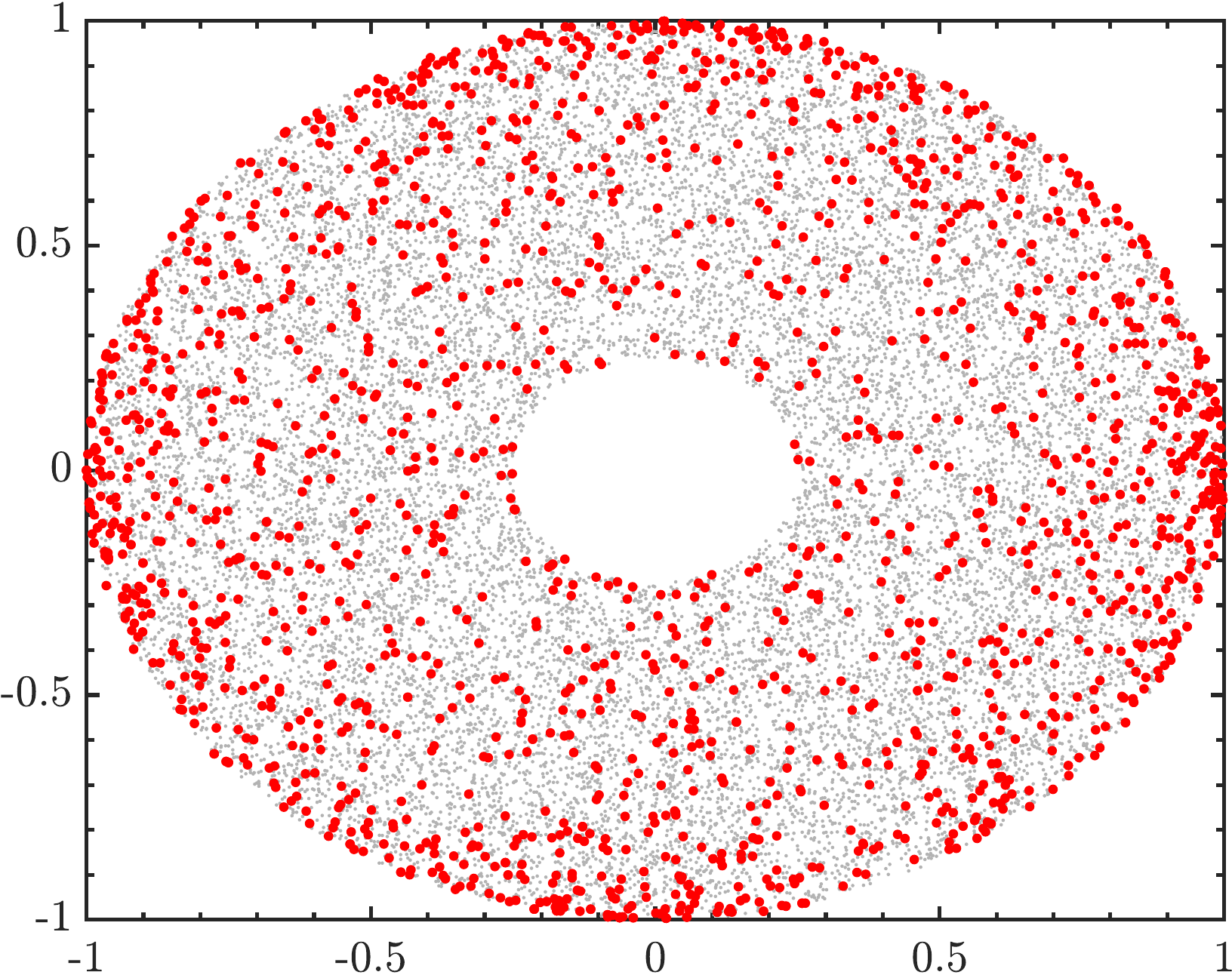} &
 			\includegraphics[width=0.3\textwidth]{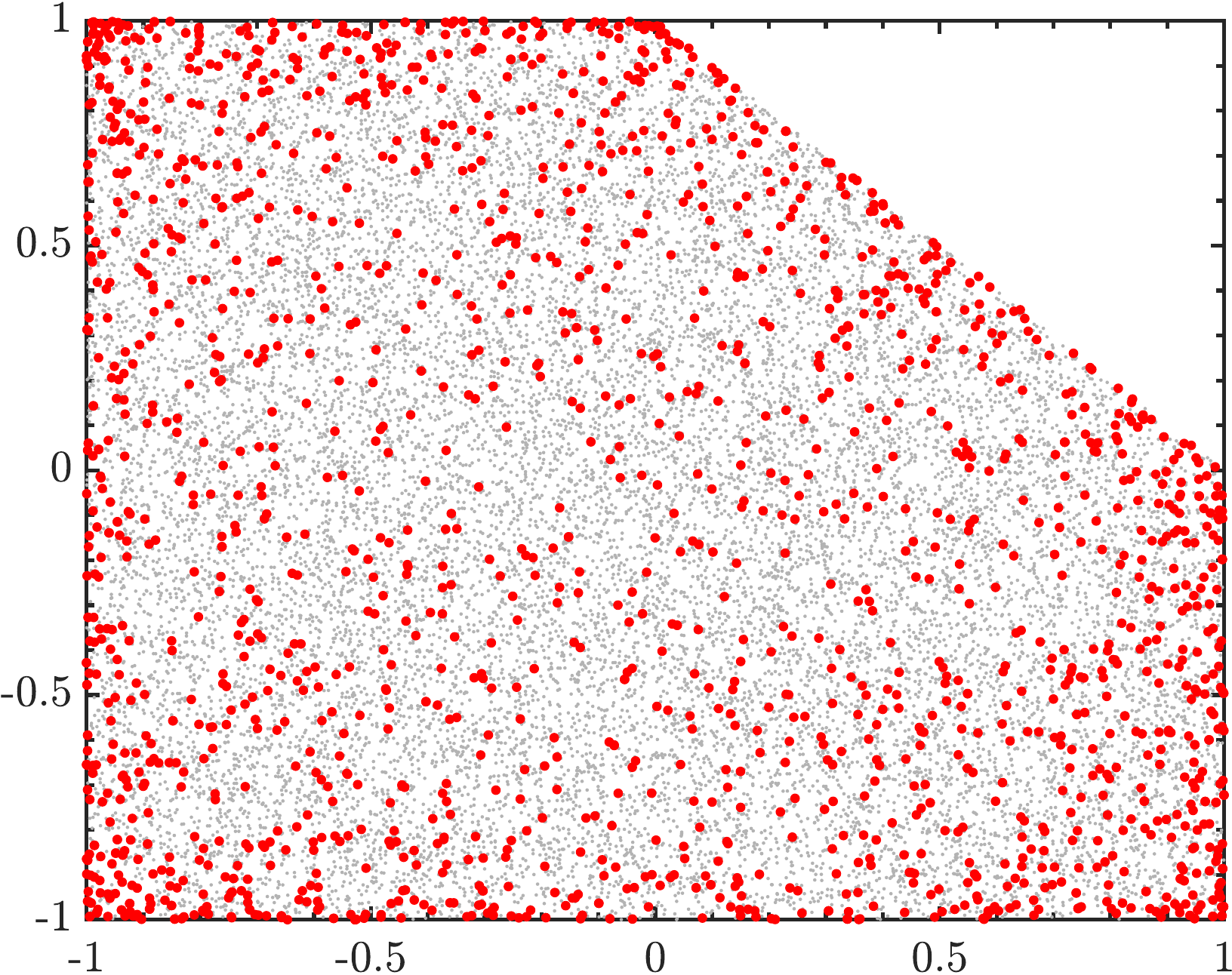} \\  			
		\end{tabular}
	\end{small}
	\end{center}
\caption{
The domains $D_1$, $D_2$ and $D_3$ (left to right) in $d = 2$ dimensions. The grey dots are the finite Monte Carlo grid $Z = \{z_i\}^{k}_{i=1}$ with $k = 20,000$. The top row shows Monte Carlo sampling with $m = 1,824$ points. The second row shows the same number of points generated from the discrete optimal nonhierarchical sampling measure based on $\cI = \cI^{\mathrm{HC}}_{t-1}$, where $t = 68$ and $s=308$.
} 
\label{fig:LS-2}
\end{figure}

\begin{figure}[t]
	\begin{center}
	\begin{small}
 \begin{tabular}{@{\hspace{0pt}}c@{\hspace{-0.5pc}}c@{\hspace{-0.5pc}}c@{\hspace{0pt}}}
 			\includegraphics[width=0.35\textwidth]{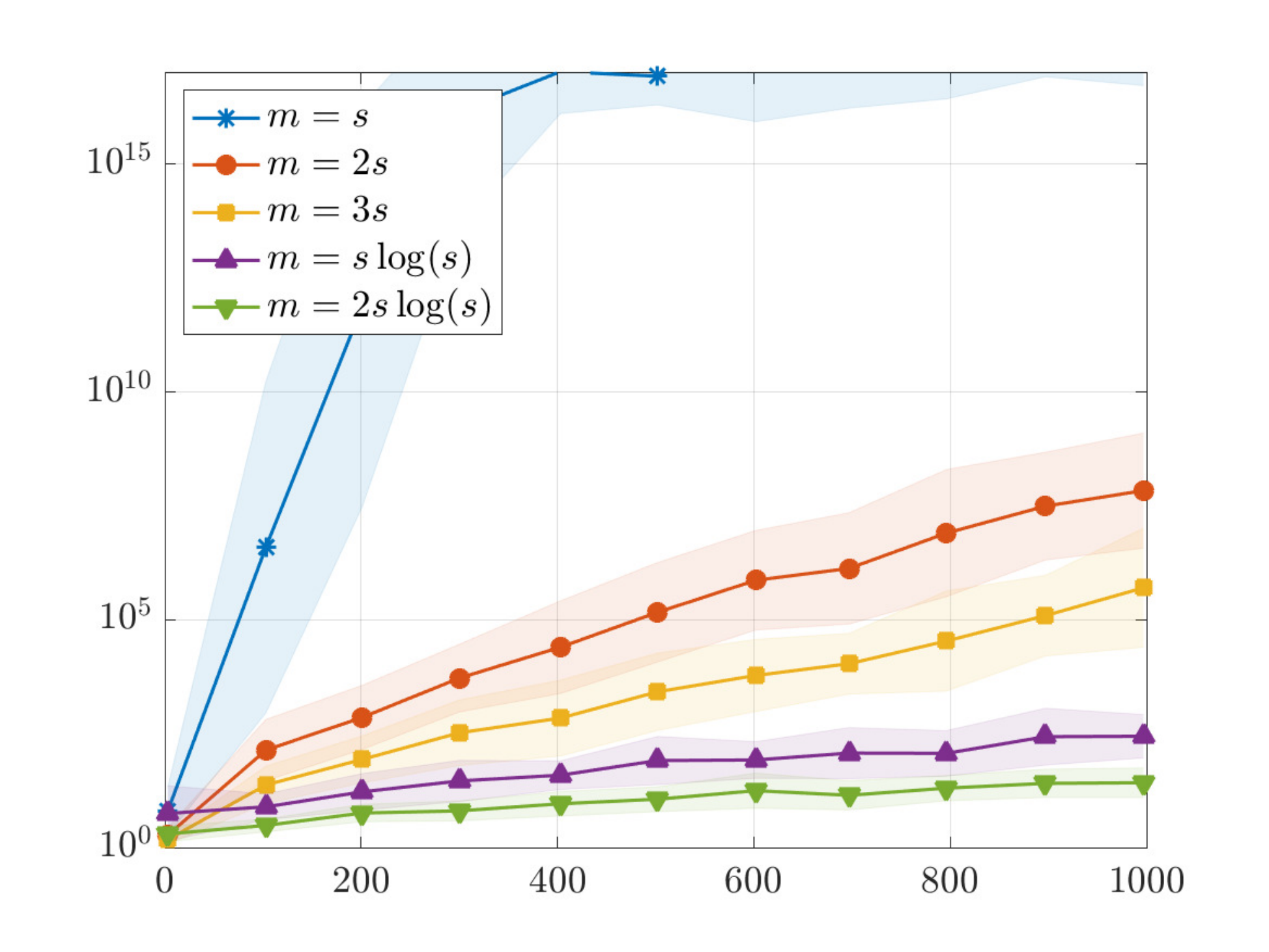} &
 			\includegraphics[width=0.35\textwidth]{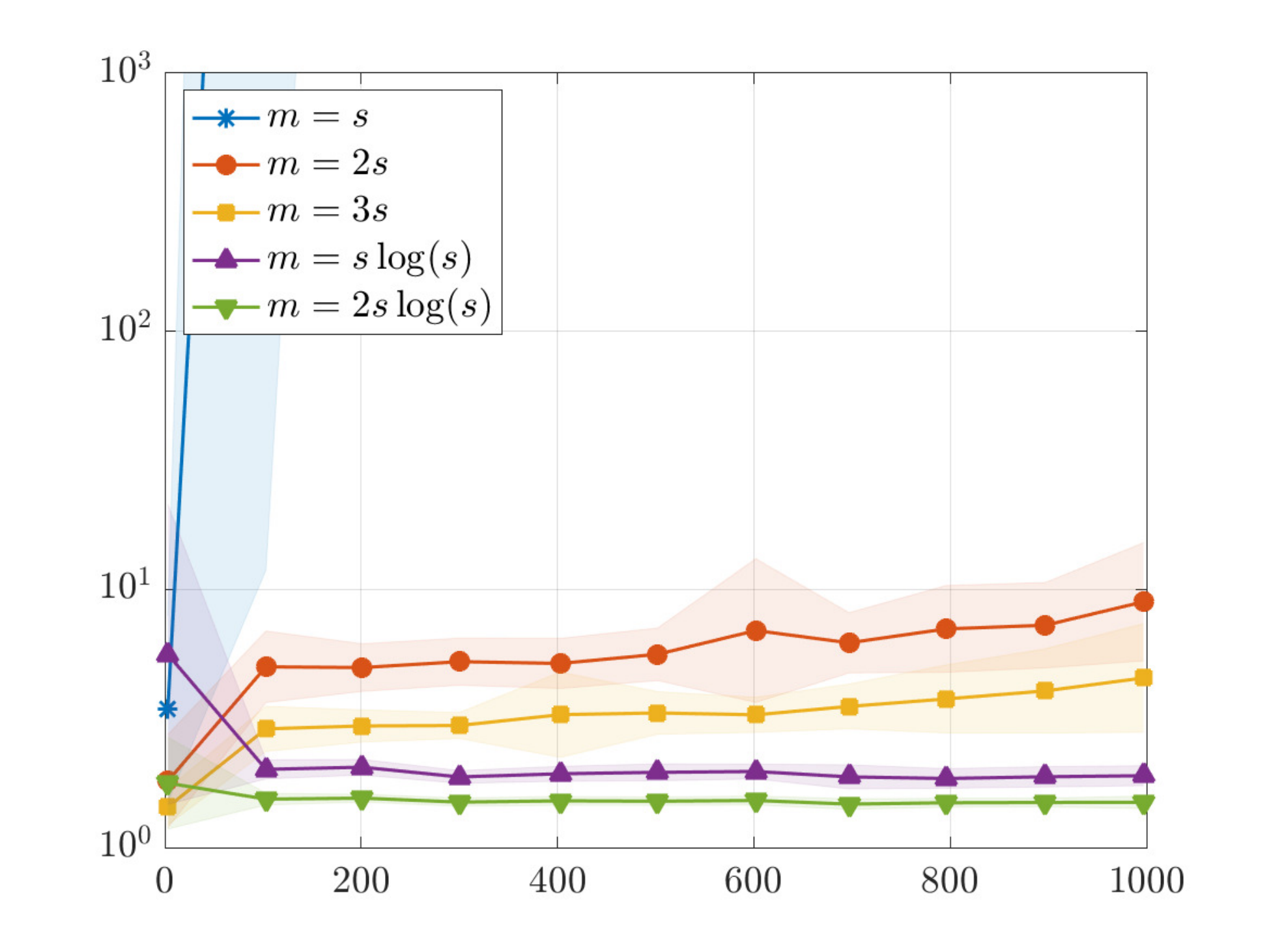} &
 			\includegraphics[width=0.35\textwidth]{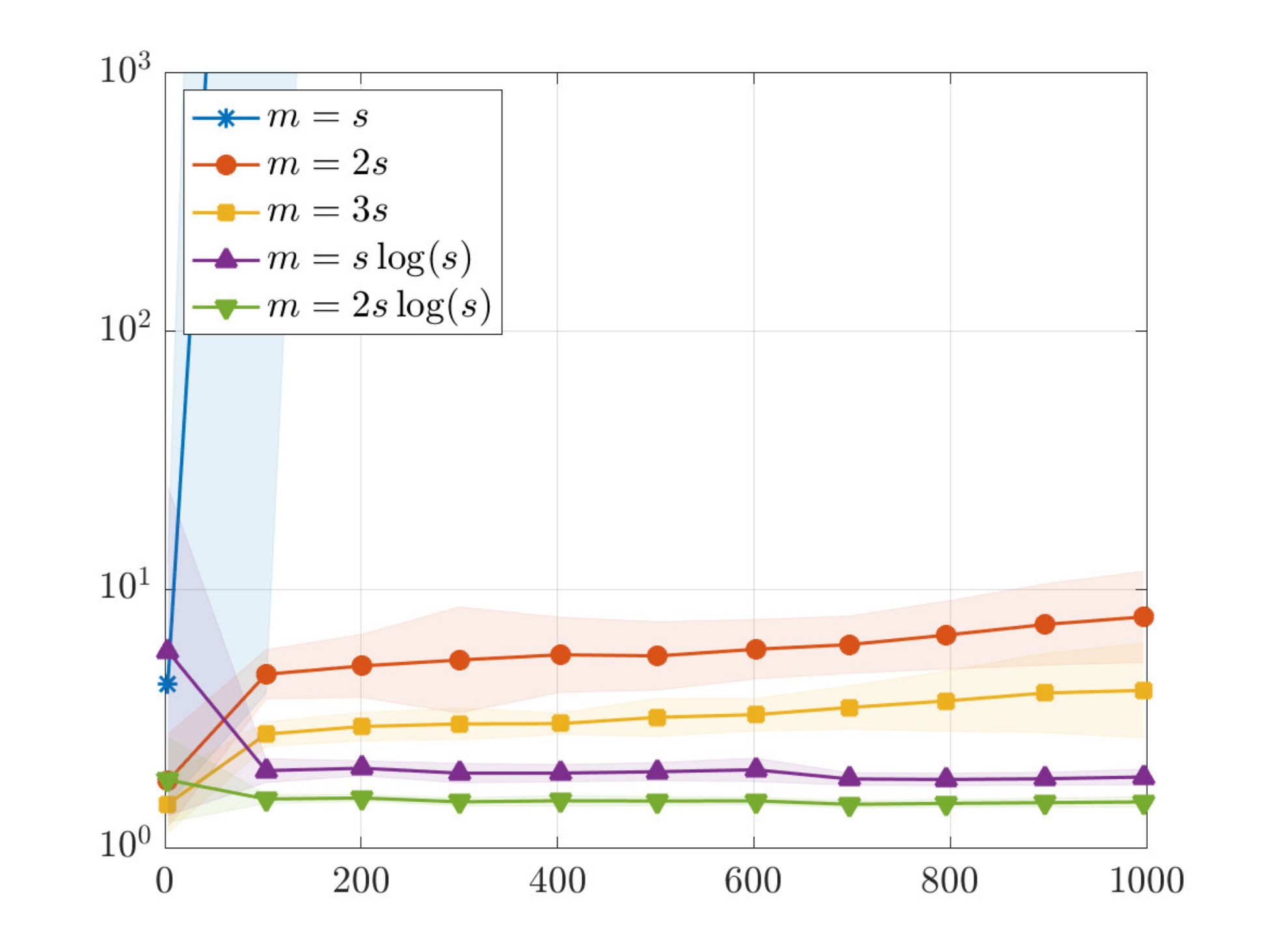}\\
 			\includegraphics[width=0.35\textwidth]{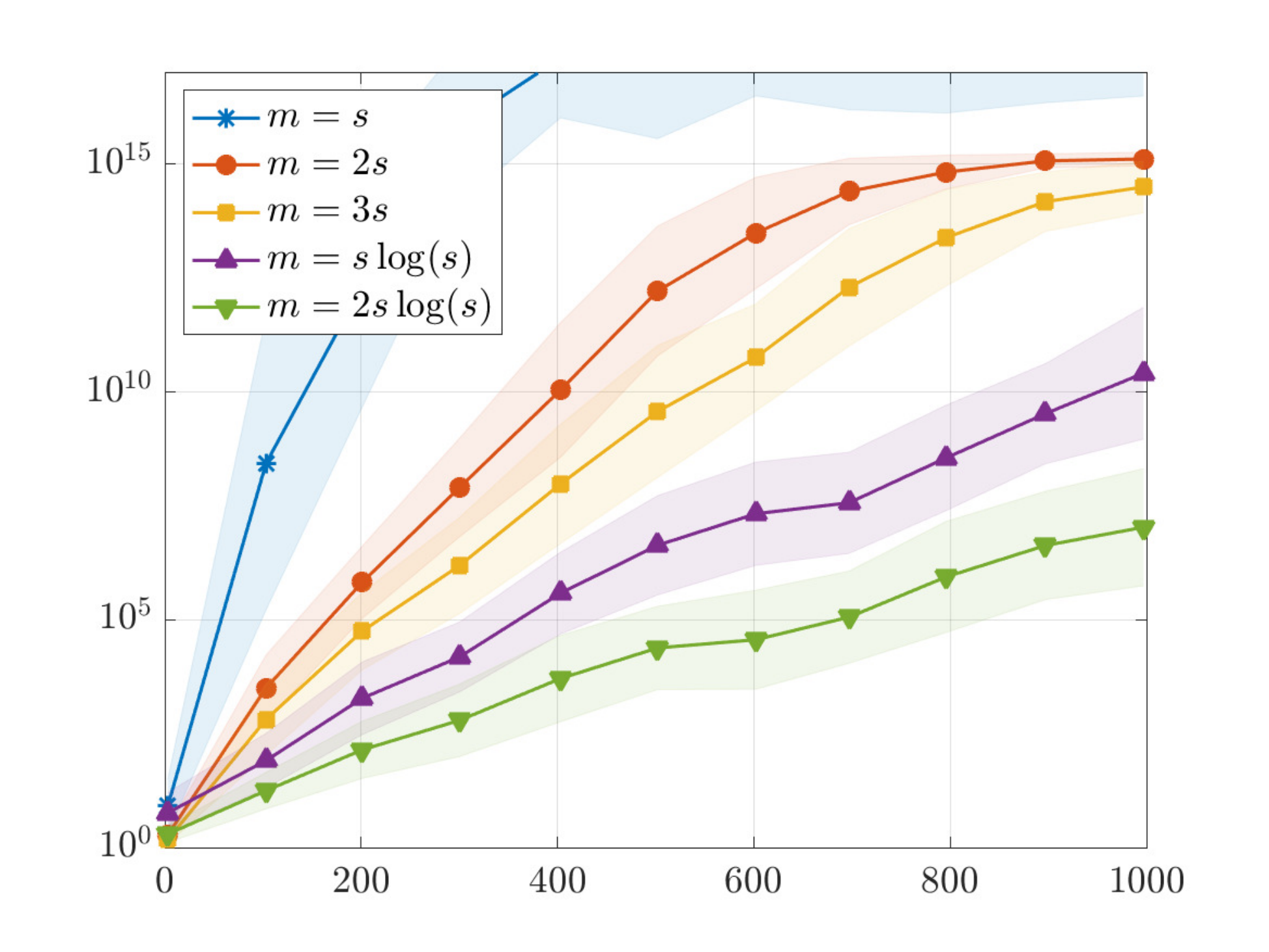} &
 			\includegraphics[width=0.35\textwidth]{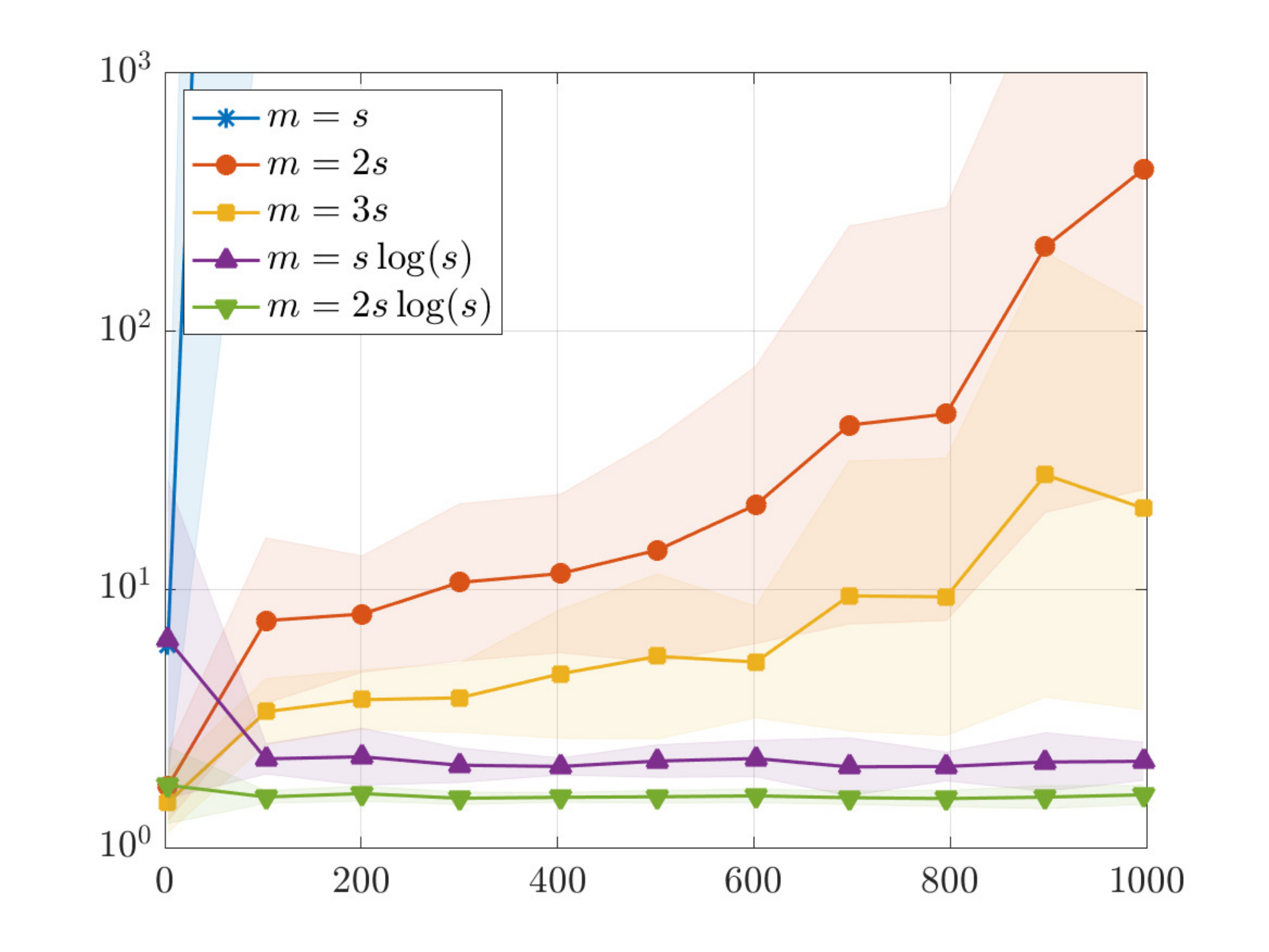} &
 			\includegraphics[width=0.35\textwidth]{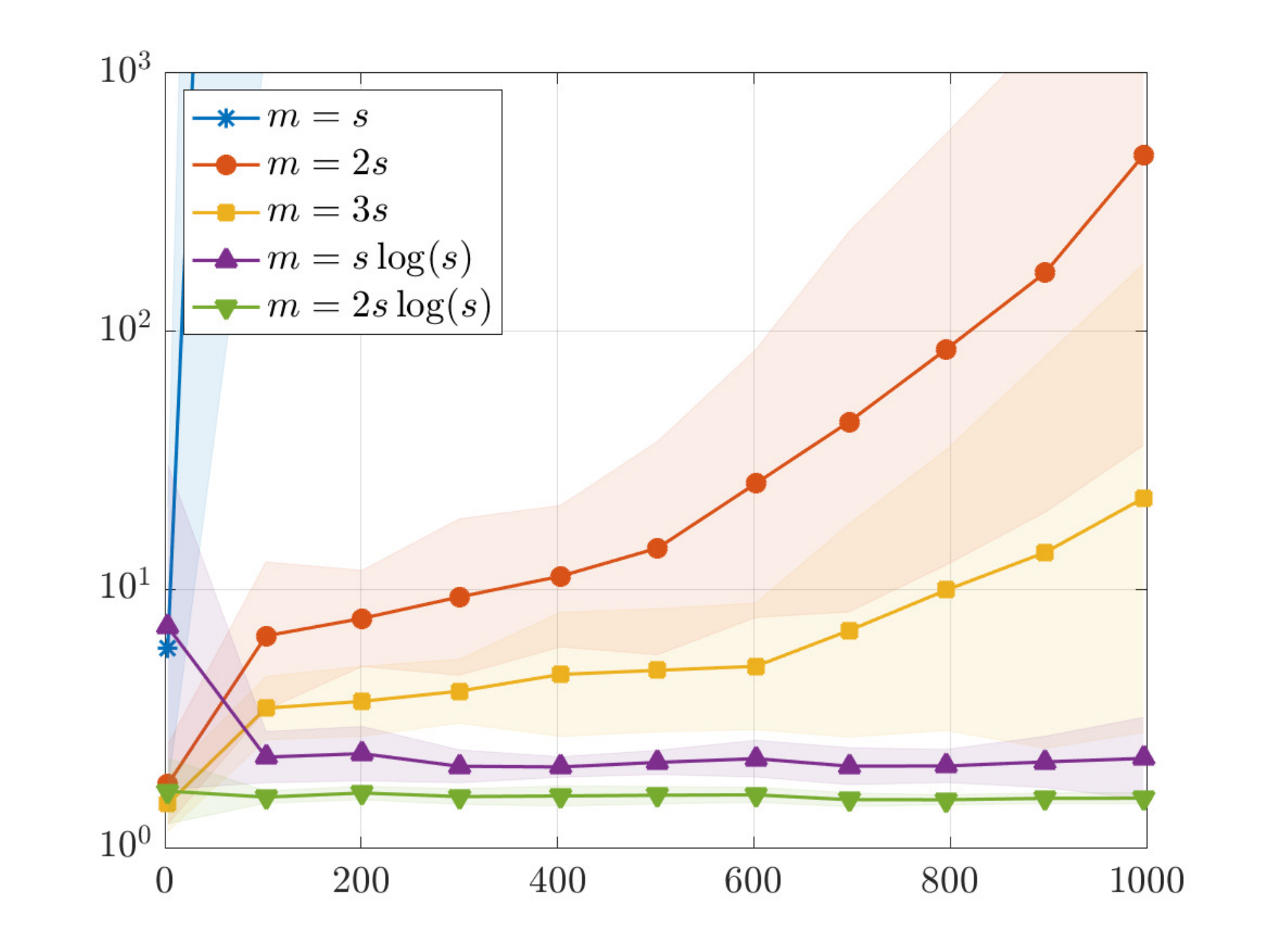}		
		\end{tabular}
	\end{small}
	\end{center}
\caption{The constant $1/\sqrt{\alpha}$ versus $s$ in $d=2$ dimensions for Monte Carlo sampling (left), discrete optimal nonhierarchical sampling (middle) and discrete optimal hierarchical sampling (right) and for $D_1$ (top row) and $D_2$ (bottom row), where $\cI = \cI^{\mathrm{HC}}_{t-1}$.} 
\label{fig:LS-3}
\end{figure}

In Fig.\ \ref{fig:LS-1} we compare Monte Carlo with both the hierarchical and nonhierarchical optimal sampling schemes. In two dimensions, typical sample points generated by these schemes for different domains are shown in Fig.\ \ref{fig:LS-2}. As we see from Fig.\ \ref{fig:LS-1}, Monte Carlo sampling leads to worse performance compared to both optimal sampling schemes, especially in lower-dimensional problems. It is notable that Monte Carlo sampling also leads to an increasing approximation error in several cases, since the number of samples is chosen to scale log-linearly with $s$, rather than log-quadratically (recall the discussion in \S \ref{MC-sampling-LS}). This is corroborated in Fig.\ \ref{fig:LS-3}, where we plot the constant $\alpha$ for the different sampling schemes.  On the other hand, we observe that the relative performance of Monte Carlo sampling improves in higher dimensions, where it offers similar approximation errors to the optimal schemes. We see this effect consistently throughout this work. Finally, we remark in passing that there is virtually no difference between the nonhierarchical and hierarchical versions of the optimal sampling scheme.

\subsection{Proofs of Theorems \ref{t:LS_acc_stab} and \ref{t:LS_samp_comp}}\label{ss:proof-LS}

The proofs of Theorems \ref{t:LS_acc_stab} and \ref{t:LS_samp_comp} follow ideas that are now well established in the literature (see, for example, \cite{cohen2013stability,cohen2017optimal,adcock2020approximating,adcock2020nearoptimal,migliorati2021multivariate}). They are included for completeness.
We commence with Theorem \ref{t:LS_acc_stab}.
We now observe the following:

\lem{
\label{l:RIP_over_S_C_V}
Suppose that \R{RIP_over_S} holds. Then
\be{
\label{RIP_over_S_V}
\alpha \nm{p}^2_{L^2_{\rho}(D;\bbV)} \leq \nm{p}^2_{\mathrm{disc}} \leq \beta \nm{p}^2_{L^2_{\rho}(D;\bbV)},\quad \forall p \in P_{S ; \bbV_h},
}
where $\nm{\cdot}_{\mathrm{disc}}$ is as in \S \ref{ss:acc-stab-samp-LS}.
}

Simply put, this lemma states that if there is a norm equivalence in the scalar-valued case over $P_S = P_{S;\bbC}$ between the continuous $L^2_{\rho}$-norm and the discrete norm defined by the sample points, then there is also the same norm equivalence in the Hilbert-valued case over $P_{S;\bbV}$.

\prf{
First, let $\{ \psi_i \}_{i=1}^K$ be an orthonormal basis of $\bbV_h \subseteq \bbV$. Let $p \in P_{S;\bbV}$, and observe that it has the unique expression
\bes{
p = \sum_{i \in [s]} \sum_{j \in [K]} c_{ij} \phi_{\iota_i} \otimes \psi_j, \qquad c_{ij} \in \bbC,\ i \in [s],\ j \in [K],
}
where $\{ \iota_1,\ldots, \iota_s \}$ is an enumeration of the indices in $S$. Let $q_j = \sum_{i \in [s]} c_{i j} \phi_{\iota_i}$ and observe that $q_j \in P_S$. Notice also that
\bes{
\nm{p(y)}^2_{\bbV} = \sum_{j \in [K]} \left | \sum_{i \in [s]} c_{i j}\phi_{\iota_j}(y) \right |^2 =  \sum_{j \in [K]}  \left |  q_j(y) \right |^2 , \quad \forall y \in D.
}
In particular, this implies that
\be{
\label{p-qj-sum-norm}
\nm{p}^2_{L^2_{\rho}(D;\bbV)} =\int_D  \nm{p(y)}_{\bbV}^2 d\rho(y) = \int_D  \sum_{j \in [K]} \left |  q_j(y)\right |^2 d\rho(y) =  \sum_{j \in [K]}  \nm{q_j}^2_{L^2_{\rho}(D)},
}
and also that
\bes{
\nm{p}^2_{\mathrm{disc}} = \frac1m \sum_{i \in [m]} w(y_i) \nm{p(y_i)}^2_{\bbV} =  \sum_{j \in [K]} \left ( \frac1m \sum_{i \in [m]} w(y_i) |q_j(y)|^2 \right ).
}
Hence, by the scalar-valued norm equivalence \R{RIP_over_S}, we deduce that
\bes{
\alpha \sum_{j \in [K]} \nm{q_j}^2_{L^2_{\rho}(D)} \leq \nm{p}^2_{\mathrm{disc}}  \leq \beta \sum_{j \in [K]} \nm{q_j}^2_{L^2_{\rho}(D)}.
}
The result now follows immediately from \R{p-qj-sum-norm}.
}

\prf{[Proof of Theorem \ref{t:LS_acc_stab}]
Since $\hat{f}$ is a solution of the least-squares problem, it is also a solution of the variational equations
\bes{
\mbox{find $\hat{f} \in P_{S;\bbV_h}$ such that}\ \ip{\hat{f}}{q}_{\mathrm{disc}} = \ip{f}{q}_{\mathrm{disc}} + \frac1m \sum_{i \in [m]} w(y_i) \ip{n_i}{q(y_i)}_{\bbV},\ \forall q \in P_{S;\bbV_h}.
}
Uniqueness of $\hat{f}$ now follows immediately from Lemma \ref{l:RIP_over_S_C_V}, since $\ip{\cdot}{\cdot}_{\mathrm{disc}}$ forms an inner product on $P_{S;\bbV}$, and therefore $P_{S;\bbV_h}$.

We now derive the desired error bound. First, observe that since $q \in P_{S ; \bbV_h}$, these equations are equivalent to
\bes{
\mbox{find $\hat{f} \in P_{S;\bbV_h}$ such that}\ \ip{\hat{f}}{q}_{\mathrm{disc}} = \ip{\cP_h(f)}{q}_{\mathrm{disc}}+ \frac1m \sum_{i \in [m]} w(y_i) \ip{n_i}{q(y_i)}_{\bbV},\ \forall q \in P_{S;\bbV_h}.
} 
Now let $q \in P_{S ; \bbV_h}$ be arbitrary. Then these equations give
\bes{
\nmu{\hat{f} - q}^2_{\mathrm{disc}} = \ip{\cP_h(f) - q}{\hat{f} - q}_{\mathrm{disc}}+ \frac1m \sum_{i \in [m]} w(y_i) \ip{n_i}{\hat{f}(y_i) - q(y_i)}_{\bbV},
}
and applying the Cauchy--Schwarz inequality several times to the right-hand side we deduce that
\bes{
\nmu{\hat{f} - q}_{\mathrm{disc}} \leq \nmu{\cP_h(f) - q}_{\mathrm{disc}} + \nm{e}_{\ell^2([m];\bbV)}. 
}
Here, we also recall the definition of $e$.
Further, using \R{RIP_over_S_V} and the fact that $q \in P_{S ; \bbV_h}$ we get
\bes{
 \nmu{\hat{f} - q}_{L^2_{\varrho}(D;\bbV)} \leq  \frac{1}{\sqrt{\alpha}} \left ( \nm{\cP_h(f)-q}_{\mathrm{disc}} + \nm{e}_{\ell^2([m];\bbV)} \right ).
}
Now let $p \in P_{S ; \bbV}$ be arbitrary and write $q = \cP_h(p)$. Therefore,
\eas{
 \nmu{\hat{f} - f}_{L^2_{\varrho}(D;\bbV)} \leq &  \nmu{\hat{f} - \cP_h(p)}_{L^2_{\varrho}(D;\bbV)} + \nmu{\cP_h(p) - \cP_h(f)}_{L^2_{\varrho}(D;\bbV)} + \nmu{\cP_h(f) - f}_{L^2_{\varrho}(D;\bbV)} 
 \\
\leq &  \frac{1}{\sqrt{\alpha}} \left ( \nm{\cP_h(f)-\cP_h(p)}_{\mathrm{disc}} + \nm{e}_{\ell^2([m];\bbV)} \right )
\\
& + \nmu{\cP_h(f-p)}_{L^2_{\varrho}(D;\bbV)}+ \nmu{\cP_h(f) - f}_{L^2_{\varrho}(D;\bbV)} 
\\
\leq &  \frac{1}{\sqrt{\alpha}} \left ( \nm{f - p}_{\mathrm{disc}} + \nm{e}_{\ell^2([m];\bbV)} \right )
\\
& +  \nmu{f-p}_{L^2_{\varrho}(D;\bbV)}+ \nmu{\cP_h(f) - f}_{L^2_{\varrho}(D;\bbV)} .
}
Here, in the final step, we used the fact that $\nm{\cP_h(g)}_{\mathrm{disc}} \leq \nm{g}_{\mathrm{disc}}$ for all $g \in L^2_{\rho}(D;\bbV)$, since $\cP_h$ is an orthogonal projection, and likewise $\nm{\cP_h(g)}_{L^2_{\varrho}(D;\bbV)} \leq \nm{g}_{L^2_{\varrho}(D;\bbV)}$. 
This completes the proof.
}

We now consider Theorem \ref{t:LS_samp_comp}. This is most commonly established using the matrix Chernoff bound \cite[Thm.\ 1.1]{tropp2012user}, which we restate here for convenience:

\thm{[Matrix Chernoff bound]
\label{t:matrixChernoff}
Let $X_1,\ldots,X_m$ be independent, self-adjoint random matrices of {size $s \times s$}. Assume that {$X_i$ is positive semidefinite and $\lambda_{\max}(X_i) \leq R$ almost surely for each $i=1,\ldots,m$,}
and define
\bes{
\mu_{\min} = \lambda_{\min} \left ( \sum^{m}_{i=1} \bbE X_i \right ) ,\quad \mu_{\max} = \lambda_{\max} \left ( \sum^{m}_{i=1} \bbE X_i \right ).
}
Then, for $0 \leq \delta \leq 1$,
\bes{
\bbP \left ( \lambda_{\min}  \left ( \sum^{m}_{i=1} X_i \right ) \leq \left ( 1 - \delta \right ) \mu_{\min} \right ) \leq s \cdot \exp \left ( -\frac{\mu_{\min} ( (1-\delta) \log(1-\delta) + \delta )}{R} \right ),
}
and, for $\delta \geq 0$,
\bes{
\bbP \left ( \lambda_{\max}  \left ( \sum^{m}_{i=1} X_i \right ) \geq \left ( 1 + \delta \right ) \mu_{\max} \right ) \leq s \cdot \exp \left ( -\frac{\mu_{\max} ( (1+\delta) \log(1+\delta)-\delta)}{R} \right ).
}
}

\prf{[Proof of Theorem \ref{t:LS_samp_comp}]
Let $\{\upsilon_{i}\}_{i \in [s]}$ be an orthonormal basis of $P_S \subset L^2_{\rho}(D)$ with respect to $\rho$, Let $p\in P_S$, $p \neq 0$, be arbitrary, and write $p=\sum_{i \in [s]}c_{i}\upsilon_{i}$ and $c = (c_{i})_{i \in [s]}$. Then
\bes{
\nm{p}^2_{L^2_\rho(D)} = \int_{D}\left|\sum_{i\in S}c_{i} \upsilon_{i}(y)\right|^2 \D\rho(y) = \sum_{i \in [s]}|c_{i}|^2 = \nm{c}_2^2, 
}
and
\bes{
\frac{1}{m}\sum_{i=1}^m w(y_i)|p(y_i)|^2 = c^*Gc, 
}
where $G = A^*A \in \bbC^{s\times s}$ is the self adjoint matrix with entries $G_{j,k}=\frac{1}{m}\sum_{i=1}^m w(y_i)\ip{\upsilon_j}{\upsilon_k}_{L^2_{\rho}(D)}$. It therefore suffices to show that $\lambda_{\min}(G) \geq 1-\delta$ and $\lambda_{\max}(G) \leq 1+\delta$. Write 
\bes{
G=\sum_{i=1}^m X_i,\quad X_i=\left\lbrace \frac{1}{m}w(y_i)\upsilon_j(y_i)\overline{\upsilon_k(y_i)}\right\rbrace_{j,k=1}^m.
}
By construction, these matrices are independent and {positive semidefinite}. Also, 
\bes{
\left(\bbE(X_i)\right)_{j,k} = \int_{D}\upsilon_j(y)\overline{\upsilon_k(y)}w(y)\frac{1}{m}d\mu_i(y),
}
which gives 
\bes{
\left(\sum_{i=1}^m\bbE(X_i)\right)_{j,k} 
= \int_D \upsilon_j(y)\overline{\upsilon_k(y)}w(y)\frac{1}{m}\sum_{i=1}^md\mu_i(y)
= \int_D\upsilon_j(y)\overline{\upsilon_k(y)}d\rho(y) 
= \delta_{j,k}.
}
Hence $\sum_{i=1}^m\bbE(X_i)=I$ is the identity matrix. Moreover, for any $c\in \bbC^s$ we have
\bes{
cX_i c = \frac{1}{m}\left|\sum_{j\in S}c_j\sqrt{w(y_i)}\upsilon_j(y_i)\right|^2
\leq \frac{\left(\mathcal{N}(P_S,w)\right)^2}{m} \nm{\sum_{j\in S}c_j\upsilon_j}^2_{L^2_\rho(D)}
= \frac{\left(\mathcal{N}(P_S,w)\right)^2}{m}\nm{c}^2_2
} 
Since these matrices are self adjoint and {positive semidefinite}, we deduce that 
\bes{
\lambda_{\max}(X_i)\leq \frac{\left(\mathcal{N}(P_S,w)\right)^2}{m}.
}
We now apply the matrix Chernoff bound \eqref{t:matrixChernoff} with $s$, $R=\left(\mathcal{N}(P_S,w)\right)^2/m$ and 
\bes{
\mu_{\min} = \lambda_{\min}\left(\sum_{i=1}^m\bbE(X_i)\right)=\lambda_{\min}(I)=1,
}
and likewise $\mu_{\max}=1$. This gives
\eas{
\bbP&\left( \lambda_{\min}(G) \leq 1-\delta\ \mbox{or}\ \lambda_{\max}(G) \geq 1+\delta \right) 
\\
&\leq  \bbP\left(\lambda_{\min}(G)\leq(1+\delta)\right) + \bbP\left(\lambda_{\max}(G)\geq (1-\delta)\right)
\\
& \leq s\cdot \left(\exp\left(-\frac{(1-\delta)\log(1-\delta)+\delta}{m^{-1}\left(\mathcal{N}(P_S,w)\right)^2}\right)
+\exp\left(-\frac{(1+\delta)\log(1+\delta)-\delta}{m^{-1}\left(\mathcal{N}(P_S,w)\right)^2}\right)\right)
}
Note that $(1+\delta)\log(1+\delta)-\delta\leq (1-\delta)\log(1-\delta)+\delta$ for $0<\delta<1$. Hence
\bes{
\bbP\left( \lambda_{\min}(G) \leq 1-\delta\ \mbox{or}\ \lambda_{\max}(G) \geq 1+\delta \right) 
\leq 2s \cdot \exp\left(\frac{-c_{\delta}}{m^{-1}\left(\mathcal{N}(P_S,w)\right)^{2}}\right) \leq \epsilon,
}
where in the last step we use the condition on $m$. This completes the proof.
}

\section{Sparse approximation via $\ell^1$-minimization}\label{s:CS}

Having discussed the case of Problem \ref{ass:sparsity_known}, we now consider the substantially more challenging setting of Problem \ref{ass:sparsity_unknown}. In order to facilitate its solution, we now also make an additional assumption on the dictionary $\Phi$: namely, the index set $\cI$ is finite, and the functions $\phi_{\iota}$, $\iota \in \cI$, are linearly independent. In what follows we write $n = |\cI|$. Note that, typically, $n \gg m$, where $m$ is the number of measurements. It is notable that the examples described in \S \ref{ss:examples} correspond to cases where the index set is countable. In this case, one may define $\cI$ as a large, but finite truncated index set in which the target set $S$ in the sparse representation \R{sparse_rep} is expected to lie. We shall return to this matter briefly in \S \ref{s:lower-sets} (see Remark \ref{r:choice-I}).

This aside, we now also assume that $\mu_1 = \ldots = \mu_m = \mu$ in Assumption \ref{ass:abs-cont-pos}, i.e.\ $\mu$ is a probability measure that is absolutely continuous with respect to $\rho$ and for which the Radon--Nikodym derivative is strictly positive almost everywhere. In this case, the corresponding weight function $w$ satisfies
\be{
\label{mu_weight_fn_single}
\D \mu(y) = \frac{1}{w(y)} \D \rho(y).
} 
This is done to simplify several of the arguments later. However, it is also possible to consider distinct measures $\mu_1,\ldots,\mu_m$ as in the previous section.

\subsection{Formulation}

Given $\Phi = \{ \phi_{\iota} : \iota \in \cI \}$, we strive to exploit the fact that $f$ is assumed to have an approximate sparse representation in $\Phi$. In this section, we do this via minimizing the $\ell^1(\cI;\bbV)$-norm of the coefficients, while also promoting fidelity of the resulting approximation to the measurements \R{f_data_intro}.  There are various ways to do this, including the ($\bbV$-valued) \textit{Quadratically-Constrained Basis Pursuit (QCBP)}
\be{
\label{qcbp}
\hat{f} \in \argmin{p \in P_{\cI ; \bbV_h}} \left \{ \lambda \nm{c}_{\ell^1(\cI ; \bbV)} : \frac1m \sum^{m}_{i=1} w(y_i) \nm{f(y_i) + n_i - p(y_i) }^2_{\bbV} \leq \eta^2 \right \}, 
}
the LASSO
\be{
\label{lasso}
\hat{f} \in \argmin{p \in P_{\cI ; \bbV_h}} \left \{ \lambda \nm{c}_{\ell^1(\cI ; \bbV)} + \frac1m \sum^{m}_{i=1} w(y_i) \nm{f(y_i) + n_i - p(y_i) }^2_{\bbV} \right \}, 
}
or the \textit{Square-Root LASSO (SR-LASSO)}
\be{
\label{sr-lasso}
\hat{f} \in \argmin{p \in P_{\cI ; \bbV_h}} \left \{ \lambda \nm{c}_{\ell^1(\cI ; \bbV)} + \sqrt{\frac1m \sum^{m}_{i=1} w(y_i) \nm{f(y_i) + n_i - p(y_i) }^2_{\bbV}} \right \}.
}
Here, in all cases, $c = (c_{\iota})_{\iota \in \cI} \in \bbV^n_h$ are the coefficients of $p = \sum_{\iota \in \cI} c_{\iota} \phi_{\iota} \in P_{\cI ; \bbV_h}$. The focus of this work is not the choice \R{qcbp}, \R{lasso} or \R{sr-lasso}. We remark that \R{sr-lasso} enjoys a known advantage over the other problems in that the theoretically-optimal value of the tuning parameter $\lambda$ is independent of the noise (see also Theorem \ref{t:CS_acc_stab}), which in this case also includes the typically unknown error $f - f_{\cI}$. For further background and in-depth comparison of these optimization problems, see \cite{adcock2019correcting}.

Now let
\ea{
A &= \frac{1}{\sqrt{m}} \left ( \sqrt{w(y_i)} \phi_{\iota_j}(y_i) \right )_{i  \in [m]; j \in [n]} \in \bbC^{m \times n}, \label{A-CS-def}
\\
{v} &= \frac{1}{\sqrt{m}} \left ( \sqrt{w(y_i)} (f(y_i)+n_i) \right )_{i \in [m]} \in \bbV^m_h,
\label{b-CS-def}
}
where $\{\iota_1,\ldots,\iota_n\}$ is an enumeration of the indices in $\cI$. Notice that we use the same notation for this matrix as in the previous section (see \R{A-b-def-LS}). However, it is important to note that this matrix is generally {short (fat)}, since $n \geq m$, whereas the matrix \R{A-b-def-LS} is $m \times s$, and therefore tall. Then $\hat{f} = \sum_{\iota \in \cI } \hat{c}_{\iota} \phi_{\iota}$ is a solution of \R{qcbp}, \R{lasso} or \R{sr-lasso} if and only if $\hat{c} = (\hat{c}_{\iota})_{\iota \in \cI} \in \bbV^n_h$ is given by
\bes{
\hat{c} \in \argmin{z \in \bbV^n_h}\left \{ \lambda \nm{z}_{\ell^1([n];\bbV)} + \nm{A z - {v}}^2_{\ell^2([m];\bbV)} \right \},
}
\bes{
\hat{c} \in \argmin{z \in \bbV^n_h}\left \{ \lambda \nm{z}_{\ell^1([n];\bbV)} + \nm{A z - {v}}^2_{\ell^2([m];\bbV)} \right \},
}
or
\bes{
\hat{c} \in \argmin{z \in \bbV^n_h}\left \{ \nm{z}_{\ell^1([n];\bbV)} : \nm{A z - {v}}_{\ell^2([m];\bbV)} \leq \eta \right \},
}
respectively.

\begin{remark}[Algorithms for solving \R{qcbp}--\R{sr-lasso}]
The $\bbV$-valued versions of the QCBP, LASSO, and SR-LASSO problems can be solved by considering reformulations of standard methods for solving their real and complex-valued counterparts. For example, in \cite{dexter2019mixed} the LASSO problem was solved by extending Bregman iterations and forward-backward iterations to the $\bbV$-valued case, while in \cite{adcock2021algorithmsHilbertvalued} the $\bbV$-valued SR-LASSO problem is solved via primal-dual iterations. We shall not describe algorithms for solving \R{qcbp}--\R{sr-lasso} in any further detail, and refer the interested reader to \cite{dexter2019mixed,adcock2021algorithmsHilbertvalued}.
\end{remark}

\subsection{Accuracy, stability and sample complexity}

As in \S \ref{s:wLS}, our main assumption will be a condition of the form \R{RIP_over_S}, but with two differences. First, since the target set from which the sparse representation of $f$ is obtained is unknown, we require this to hold for all subsets, not just a fixed subset. Second, as we see in the theorem below, we also require for it to hold for some value $t \geq s$. The precise condition is as follows:
\be{
\label{RIP_over_all_S}
\alpha \nm{p}^2_{L^2_{\varrho}(D)} \leq \frac{1}{m} \sum^{m}_{i=1} w(y_i) | p(y_i)|^2 \leq \beta \nm{p}^2_{L^2_{\varrho}(D)},\quad \forall p \in P_{T},\ T \subseteq \cI,\ |T|\leq t.
}
In addition to this, we also recall that $\Phi$ is a finite dictionary consisting of linearly independent elements. Therefore it is a Riesz basis, meaning that
\be{
\label{Riesz-bounds}
a \nm{c}^2_{\ell^2(\cI)} \leq \nm{\sum_{\iota \in \cI} c_{\iota} \phi_{\iota} }^2_{L^2_{\rho}(D)} \leq b \nm{c}^2_{\ell^2(\cI)},\quad \forall c = (c_{\iota})_{\iota \in \cI} \in \bbC^n.
}
for constants $0 < a \leq b < \infty$. Finally, before stating the main result, we need some additional notation. We write
\be{
\label{best-s-term}
\sigma_s(x)_{\ell^1(\cI;\bbV)} = \inf \left \{ \nm{x-z}_{\ell^1(\cI;\bbV)} : \mbox{$z \in \bbV^{n}$ is $s$-sparse} \right \}, \quad x  \in \bbV^{n}.
}
Here we recall that a Hilbert-valued vector $z = (z_i)_{i \in [n]} \in \bbV^n$ is \textit{$s$-sparse} if it has at most $s$ nonzero entries, i.e.\ $| \{ i : z_i \neq 0 \} | \leq s$.

\thm{[Accuracy and stability of $\ell^1$-minimization]
\label{t:CS_acc_stab}
Let $\Phi = \{ \phi_{\iota} : \iota \in \cI \} \subset L^2_{\rho}(D)$ be a finite dictionary consisting of $n$ linearly-independent functions, with bounds $a,b > 0$ as in \R{Riesz-bounds}. Let $1 \leq s \leq n $, $0 < \alpha \leq \beta < \infty$, $\{y_i\}^{m}_{i=1} \subseteq D$, $w : D \rightarrow [0,\infty)$ be such that $w(y_i)$ is well defined for all $i$, and suppose that \R{RIP_over_all_S} holds with $t = \min \{ n ,2 \lceil 4 s \frac{b\beta}{a\alpha} \rceil  \}$. Let $f \in L^2_{\rho}(D;\bbV)$ with measurements \R{f_meas} and consider the problem \R{sr-lasso} with $\lambda \leq \frac{15\sqrt{a \alpha}}{26\sqrt{s}}$. Then any solution $\hat{f}$ of \R{sr-lasso} satisfies
\eas{
\nm{f -\hat{f}}_{L^2_{\varrho}(D ; \bbV)} 
\leq & c_1 \frac{\sigma_{s}(c)_{\ell^1([n];\bbV)}}{\sqrt{s}}  
\\
& + c_2 \left( \nm{f - f_{\cI}}_{\mathrm{disc}} + \nm{f - \cP_h(f)}_{\mathrm{disc}} + \nm{ {e}}_{\ell^2([m];\bbV)} \right),
}
where $c = (c_{\iota})_{\iota \in \cI}$, $f_{\cI} = \sum_{\iota \in \cI} c_{\iota} \phi_{\iota}$ is the orthogonal projection (best approximation) of $f$ in $P_{\cI ; \bbV}$, $e = \frac{1}{\sqrt{m}} \left ( \sqrt{w(y_i)} n_i \right )^{m}_{i=1} \in \bbV^m$
and
\bes{
c_1 = 8 \sqrt{b},\qquad c_2=8 \sqrt{b} \left ( \frac{1}{2 \lambda \sqrt{s}} + \frac{1}{\sqrt{a \alpha}}\right ) + 1.
}
}

This result (see \S \ref{ss:CSproofs} for its proof) shows stable and accurate recovery for the solution $\hat{f}$ of \R{sr-lasso} (similar results can also be shown for \R{qcbp} and \R{lasso} -- see \cite{adcock2019correcting} and \cite[Chpt.\ 6]{adcock2021compressive}). Specifically, the error is bounded by a multiple of \R{best-s-term} and $\nm{f - f_{\cI}}_{\mathrm{disc}}$, which together measure how well $f$ can be represented by an $s$-sparse representation in $\Phi$ (observe that these terms vanish when $f$ has an exact $s$-sparse representation). The other terms are the space discretization error and the noise error. As in the case of weighted least squares (see Remark \ref{r:LS-error}), it is also possible to replace the $\nm{\cdot}_{\mathrm{disc}}$-norm by the $L^2_{\rho}(D)$-norm when the sample points are random variables. We also remark in passing that the factor $8$ in the constants is somewhat arbitrary. Other numerical values could also be used, subject to changing the numerical values in the definition of $t$ and $\lambda$.

We next consider sample complexity. The following result is analogous to Theorem \ref{t:LS_samp_comp} for the case of compressed sensing.

\thm{[Sample complexity for \R{RIP_over_all_S}]
\label{t:BOS_RIP}
Let $\Phi = \{ \phi_{\iota} : \iota \in \cI \} \subset L^2_{\rho}(D)$ be a finite dictionary consisting of $n$ linearly-independent functions, with bounds $a,b > 0$ as in \R{Riesz-bounds}. Let $\mu$ be a probability measure satisfying Assumption \ref{ass:abs-cont-pos}, $1 \leq t \leq n$, $0 <  \delta < \delta^*$ for some universal constant $0 < \delta^*< 1$, $0 < \epsilon < 1$, $0 < \alpha \leq \beta <\infty$ and $y_1,\ldots,y_m$ be independent with $y_i \sim \mu$ for $i = 1,\ldots,m$. Define
\be{
\label{Gamma-def}
\Gamma = \Gamma(\Phi,w) = \nm{ \max_{\iota \in \cI}\left \{ \sqrt{w(\cdot)} | \phi_{\iota}(\cdot) | \right \} }_{L^{\infty}_{\rho}(D)},
}
where $w$ is the weight function specified in \eqref{mu_weight_fn_single},
and suppose that 
\bes{
m \geq C \cdot \delta^{-2} \cdot (\Gamma^2/a) \cdot t \cdot \left ( \log(\E n) \cdot \log^2(\E (\Gamma^2/a) t / \delta )  + \log(2/\epsilon) \right ),
}
for some universal constant $C >0$.
Then \R{RIP_over_all_S} holds with $1-\delta \leq \alpha \leq \beta \leq 1+\delta$, with probability at least $1-\epsilon$.
}

As in the least-squares case, this result reduces the question of sample complexity to the matter of estimating a certain constant $\Gamma(\Phi,w)$ depending on the system $\Phi$ and the weight function $w$. Observe that $\Gamma^2 \geq a$ for any $\Phi$ and $w$. Indeed, $\Gamma^2 \geq w(y) | \phi_{\iota}(y) |^2$ almost everywhere, and therefore \R{w_normalization} and \R{Riesz-bounds} give $\Gamma^2 \geq \nm{\phi_{\iota}}^2_{L^2_{\rho}(D)} \geq a$.

Combining this with Theorem \ref{t:CS_acc_stab} we deduce the following:

\cor{[Sample complexity of $\ell^1$-minimization]
\label{c:samp-comp-CS}
Let $\Phi = \{ \phi_{\iota} : \iota \in \cI \} \subset L^2_{\rho}(D)$ be a finite dictionary consisting of $n$ linearly-independent functions, with bounds $a,b > 0$ as in \R{Riesz-bounds}. 
Let $1 \leq s \leq n$, $0 <  \epsilon < 1$, $\mu$ be a probability measure satisfying Assumption \ref{ass:abs-cont-pos}, $y_1,\ldots,y_m$ be independent with $y_i \sim \mu$ for $i = 1,\ldots,m$ and $\Gamma$ be as in \R{Gamma-def}. Suppose that
\bes{
m \gtrsim  (b/a) \cdot (\Gamma^2/a) \cdot s \cdot \left ( \log(\E n) \cdot \log^2(\E (b/a)(\Gamma^2/a)  s) + \log(2/\epsilon) \right ).
}
Then the following holds with probability at least $1-\epsilon$. Let $f \in L^2_{\rho}(D;\bbV)$ with measurements \R{f_meas} and consider the problem \R{sr-lasso} with $A$ and $b$ as in \R{A-CS-def} and $\lambda = c \sqrt{a/s}$ for some $0 < c < C$, where $C>0$ is a universal constant. Then any solution $\hat{f}$ of \R{sr-lasso} satisfies
\eas{
\nm{f -\hat{f}}_{L^2_{\varrho}(D ; \bbV)} \lesssim_c & \sqrt{b} \frac{\sigma_{s}(c)_{\ell^1([n];\bbV)}}{\sqrt{s}}  
\\
& + \sqrt{b/a} \left(\nm{f - f_{\cI}}_{\mathrm{disc}} + \nm{f - \cP_h(f)}_{\mathrm{disc}} + \nm{ {e}}_{\ell^2([m];\bbV)} \right).
} 
}
\prf{
Theorem \ref{t:BOS_RIP} and the condition on $m$ show that \R{RIP_over_all_S} holds with $\delta = \delta^*/2$, where $t = \min \{ n , 2 \lceil 4 s \frac{b (1+\delta)}{a(1-\delta)} \rceil \}$. We now apply Theorem \ref{t:CS_acc_stab}, noting that $c_1 \lesssim \sqrt{b}$ and $c_{2} \lesssim (1/c+1) \sqrt{b/a} + 1 \lesssim_c \sqrt{b/a}$ in this case.
}

\rem{[The Riesz basis constants $a,b$]
\label{r:Riesz-constants}
On closer inspection of the proofs, it is evident that it is possible to somewhat relax the assumption \R{Riesz-bounds} by requiring only sparse subsets of $\Phi$ to form Riesz bases (with the same bounds). In Theorem \ref{t:CS_acc_stab}, for example, it is possible to replace \R{Riesz-bounds} with the weaker condition
\be{
\label{Riesz-bounds-sparse}
a \nm{c}^2_{\ell^2(\cI)} \leq \nm{\sum_{\iota \in \cI} c_{\iota} \phi_{\iota} }^2_{L^2_{\rho}(D)} \leq b \nm{c}^2_{\ell^2(\cI)},\quad \mbox{$c = (c_{\iota})_{\iota \in \cI} \in \bbC^n$ is $t$-sparse},
}
where $t = \min \{ n , 2 \lceil 4 s \frac{b \beta}{a \alpha} \rceil \}$. Even when $\Phi$ forms a Riesz basis, when $s \ll n$ the corresponding constants in \R{Riesz-bounds} may be significantly better behaved than the Riesz basis constants in \R{Riesz-bounds}. In \S \ref{s:sparse-irregular-new} we see an example where the lower constant $a$ in \R{Riesz-bounds} is extremely small, yet recovery is still possible from a reasonable number of measurements. This suggests that it may be important to use \R{Riesz-bounds-sparse} instead of \R{Riesz-bounds} in some scenarios.
}

\subsection{Monte Carlo sampling}\label{s:MC-CS}

We now discuss the case of Monte Carlo sampling, which corresponds to the choice $\mu = \rho$, i.e.\ $w(y) \equiv 1$. Corollary \ref{c:samp-comp-CS} shows that the sample complexity of $\ell^1$-minimization is determined by the constant $\Gamma$ defined in \R{Gamma-def}. In this case, we have
\bes{
\Gamma = \max_{\iota \in \cI} \nm{\phi_{\iota}}_{L^{\infty}_{\rho}(D)} : = \Theta = \Theta(\Phi),
}
which leads to the sample complexity bound
\be{
\label{MC-CS-estimate}
m \gtrsim (b/a) \cdot (\Theta^2/a) \cdot s \cdot \left ( \log(\E n) \cdot \log^2((b/a) \E (\Theta^2/a) s) + \log(2/\epsilon) \right ).
}
It is worth comparing this bound with the least squares bound discussed in \S \ref{MC-sampling-LS}. Let $S \subseteq \cI$, $|S| \leq s$ and $p \in P_S$ be arbitrary. Write $p = \sum_{\iota \in \cI} c_{\iota} \phi_{\iota}$. Then
\bes{
\nm{p}_{L^{\infty}_{\rho}(D)} \leq \sum_{\iota \in \cI} |c_{\iota} | \nm{\phi}_{L^{\infty}_{\rho}(D)} \leq \Theta \sqrt{s} \sqrt{\sum_{\iota \in \cI} |c_{\iota} |^2 } \leq \Theta \sqrt{s} \nm{p}_{L^2_{\rho}(D)} / \sqrt{a}.
}
Hence, \R{Nikolskii-type} gives that
\bes{
(\cN(P_S) )^2 \leq \Theta^2 s / a.
}
Therefore, and unsurprisingly, the sample complexity for $\ell^1$-minimization in the setting of Problem \ref{ass:sparsity_unknown} is always at least as large as least squares in the setting of Problem \ref{ass:sparsity_known}.

On the other hand, there are clearly instances where both sample complexities are the same, at least up to log terms. Recall that the functions $\phi_{\iota}$ of Example \ref{ex:trig-poly} are equal to one in absolute value. Therefore, $\Theta = 1$ and, as discussed previously, $(\cN(P_S))^2 = s$. Since this is an orthonormal basis, we also have $a = b = 1$ in this case. Hence, \R{MC-CS-estimate} reads
\bes{
m \gtrsim s \cdot \left ( \log(\E n) \cdot \log^2(\E s) + \log(2/\epsilon) \right ).
}
We conclude that Monte Carlo sampling in combination with least squares (in the setting of Problem \ref{ass:sparsity_known}) or $\ell^1$-minimization (in the setting of Problem \ref{ass:sparsity_known}) is near-optimal for sparse approximation via trigonometric polynomials.

By contrast, in the case of Example \ref{ex:alg-poly} the size of $\Theta$ depends on the choice of the finite index set $\cI$. Using \R{tensor-Leg-def} and \R{1D-Leg-bound} we see that 
\bes{
\nm{\phi_{\iota}}^2_{L^{\infty}_{\rho}(D)} = \prod^{d}_{k=1} (2 \iota_k+1),\qquad \iota = (\iota_k)^{d}_{k=1},
}
and therefore
\be{
\label{Theta-Legendre}
\Theta^2 = \max_{\iota \in \cI} \left \{ \prod^{d}_{k=1} (2 \iota_k+1) \right \}.
}
Hence, if, for example,
\bes{
\cI = \cI^{\mathrm{TP}}_{s},
}
is the tensor product index set of order $s$ (see \R{TP-index}) then it follows immediately that
\bes{
\Theta^2 = (2s+1)^d.
}
Thus, the sample complexity bound behaves like $s (2s+1)^d$, up to log factors. This grows exponentially with $d$, and always substantially exceeds $s$. Further, since $n = |\cI| = (s+1)^d$ in this case, this means that the sample complexity bound actually exceeds $n$; a situation that is, naturally, undesirable.

This situation can be ameliorated by choosing a truncated set $\cI$ with fewer high-order polynomial indices, at the potential cost that important terms may be missed in the truncation. For example, let
\bes{
\cI = \cI^{\mathrm{HC}}_{s-1} 
}
be the hyperbolic cross index set of order $s-1$ (recall \R{HC-index}). Then \R{Theta-Legendre} gives
\bes{
\Theta^2 \leq 2^d \max_{\iota \in \cI} \left \{ \prod^{d}_{k=1} (\iota_k + 1) \right \} \leq 2^d s.
}
Hence, the sample complexity behaves like $2^d s$, up to log terms -- in other words, substantially better than in the case of the tensor-product index set, but still exponentially large in $d$. Note that this bound is well suited when $d$ is comparatively small in relation to $s$. In the setting where $d$ is large, one can also show that
\bes{
s^{\log(3)/\log(2)}/3 \leq \Theta^2 \leq s^{\log(3)/\log(2)},\quad 1 \leq s \leq 2^d.
}
See \cite[Lem.\ 3.5]{chkifa2018polynomial}. Thus, for large $d$, the same complexity bound scales like $s^{\log(3)/\log(2)+1} \approx s^{2.58}$, up to log terms -- in other words, polynomial in $s$, independently of $d$, albeit with a scaling that is substantially bigger than the optimal linear in $s$ scaling.

\subsection{`Optimal' sampling}

With this in mind, we now consider how to choose the sampling measure to obtain a smaller sample complexity. Following ideas of \cite{hampton2015compressive}, our aim is to minimize the quantity $\Gamma$ defined in \R{Gamma-def}. This is achieved by setting
\be{
\label{w-theta-opt-CS}
(w(y))^{-1} = \frac{\max_{\iota \in \cI} | \phi_{\iota}(y) |^2}{\theta^2},\qquad \theta = \theta(\Phi) : = \left ( \int_{D} \max_{\iota \in \cI} | \phi_{\iota}(y) |^2 \D \rho(y) \right )^{1/2}.
}
Notice that $\Gamma = \theta$ for this choice of $w$,  which gives the sample complexity bound
\be{
\label{opt-CS-estimate}
m \gtrsim (b/a) \cdot (\theta^2/a) \cdot s \cdot \left ( \log(\E n) \cdot \log^2(\E (\theta^2/a) (b/a) s) + \log(2/\epsilon) \right ),
}
provided $\mu$ satisfies \R{mu_weight_fn_single}, i.e.\
\be{
\label{mu-optimal-CS}
\D \mu(y) = \frac{\max_{\iota \in \cI} | \phi_{\iota}(y) |^2}{\theta^2} \D \rho(y),\quad y \in \supp(\rho).
}
Observe that the constant $\theta$ is always no larger than the constant $\Theta$ that appears in the Monte Carlo sampling estimate \R{MC-CS-estimate}. Hence, we expect this choice of measures to be no worse than Monte Carlo sampling. 

\rem{[The gap between Problem \ref{ass:sparsity_known} and Problem \ref{ass:sparsity_unknown}]
In the setting of Problem \ref{ass:sparsity_known}, we obtained a sampling measure in \S \ref{ss:optimal-LS} that lead to near-optimal sample complexity, scaling linearly in $s$ for any $\Phi$ and subset $S$ of size $|S| = s$. Critically, this measure depended on the known, target set $S$. Conversely,  in the setting of Problem \ref{ass:sparsity_unknown} the measure defined above does not, in general, lead to near-optimal sample complexity. Indeed, it is not generally the case that $\theta = 1$. This constitutes a key gap between the two settings. That it exists should come of little surprise. The sampling measure used in the former setting depends completely on the target set. Yet, the whole purpose of the latter setting is to compute sparse approximations in the absence of this assumption. Hence, it is not unexpected that the sample measure defined above (which depends on $\cI$ but not $S$) is not generally optimal.
}

\rem{[Arbitrarily-large improvements are possible]
On the other hand, there are cases where the above sampling measure leads to substantial theoretical improvements. For example, let $D = [0,1]$, $\D \rho(y) = \D y$ and
\bes{
\phi_{\iota}(y) = y^{-\alpha} \exp(2 \pi \I \iota y),
}
be the trigonometric polynomials scaled by a weight factor $y^{-\alpha}$ for some $0 < \alpha < 1/2$. Notice that $\{ \phi_{\iota} : \iota \in \cI \} \subset L^2_{\rho}(D)$ forms a Riesz basis for any finite $\cI$. Clearly, in this case one has $\Theta = \infty$. On the other hand,
\bes{
\theta^2 = \int^{1}_{0} y^{-2 \alpha} \D y = \frac{1}{1-2\alpha},
}
is bounded, for any choice of $\cI$. The reason for this is that the functions $\phi_{\iota}$ are all singular, yet their singularity occurs at the same place $y = 0$. The measure $\mu$, which has the form
\bes{
\D \mu(y) = \frac{\max_{\iota \in \cI} | \phi_{\iota}(y) |^2}{\theta^2} \D \rho(y) = \frac{y^{-2 \alpha}}{1-2 \alpha} \D y,
}
samples more densely near $y = 0$, thereby capturing the common singularity more efficiently. Note that a similar scenario occurs in the setting of algebraic polynomial approximation on the real line via Hermite polynomials. See \cite{hampton2015compressive,jakeman2017generalized}. 
}

\subsection{`Optimal' sampling and discrete measures}\label{ss:CS-opt-disc}

As in the context of least squares, drawing samples from the measure \R{mu-optimal-CS} can be challenging. Fortunately, we can overcome this issue in the same way by introducing a finite grid. Let $Z = \{ z_i \}^{k}_{i=1}$ be such a grid and
\be{
\label{tau-def-CS}
\tau = \frac1k \sum^{k}_{i=1} \delta_{z_i},
}
be the discrete uniform measure supported on it. We then define the corresponding discrete measure $\mu$ 
as
\bes{
\D \mu(y) = \frac{\max_{\iota \in \cI} | \phi_{\iota}(y) |^2}{\theta^2} \D \tau(y),
}
where
\be{
\label{theta-discrete-opt}
\theta^2 = \int_{D} \max_{\iota \in \cI} | \phi_{\iota}(y) |^2 \D \tau(y) = \frac1k \sum^{k}_{i=1} \max_{\iota \in \cI} | \phi_{\iota}(z_i) |^2.
}
In other words,
\be{
\label{mu-optimal-CS-disc}
\mu = \sum^{k}_{i=1} \frac{\max_{\iota \in \cI} | \phi_{\iota}(z_i) |^2}{\sum^{k}_{j=1} \max_{\iota \in \cI} | \phi_{\iota}(z_j) |^2 } \delta_{z_i}.
}
Sampling from this measure is achieved as follows. Define the matrix
\bes{
B = \left ( \phi_{\iota_j}(z_i) / \sqrt{k} \right )_{i \in [k], j \in [n]} \in \bbC^{k \times n},
}
and notice that
\bes{
\theta^2 = \sum^{k}_{i=1} \max_{j \in [n]} | B_{ij} |^2.
}
Hence, $y \sim \mu$ if
\bes{
\bbP(y = z_i) = \frac{\max_{j \in [n]} | B_{ij} |^2 }{\sum^{k}_{i=1} \max_{j \in [n]} | B_{ij} |^2},\quad i \in [k].
}
As in the least-squares setting, sampling with respect to this measure, following the sample complexity bound \R{opt-CS-estimate} (with $\theta$ as in \R{theta-discrete-opt}), is sufficient to ensure an estimate with respect to the discrete measure $\tau$. From this, one can also obtain an error bound with respect to the original measure $\rho$, whenever the grid is sufficiently fine. Indeed, suppose that
\be{
\label{cond-on-P-I}
\alpha' \nm{p}^2_{L^2_{\rho}(D)} \leq \nm{p}^2_{L^2_{\tau}(D)} \leq \beta' \nm{p}^2_{L^2_{\rho}(D)},\qquad \forall p \in P_{\cI},
}
for constants $\beta' \geq \alpha' > 0$ (note that this is analogous to \R{grid-fine-cond}, the difference being that we now require it to hold over $P_{\cI}$). Then \R{RIP_over_all_S} holds with respect to the $\rho$ measure with constants $\alpha' \alpha$ and $\beta' \beta$ whenever it holds with respect to the $\tau$ measure. The key point is that the grid is required to satisfy essentially the same condition \R{cond-on-P-I} -- in other words, it should give rise to a discrete norm on $P_{\cI}$. One can construct such a grid exactly as in Remark \ref{r:MC-grid}.

\subsection{Further discussion and numerical examples}\label{ss:examples-CS}

We now numerically examine Question \ref{q:how-good-MC} in the context of the Examples \ref{ex:trig-poly} and \ref{ex:alg-poly} (we discuss Example \ref{ex:alg-poly-general} further in \S \ref{s:sparse-irregular-new}). Since the complex exponentials have absolute value equal to one, Example \ref{ex:trig-poly} is an instance where $\theta = \Theta = 1$, and where $\mu = \rho$. In other words, much as in the setting of Problem \ref{ass:sparsity_known}, sparse approximation in the case of Problem \ref{ass:sparsity_unknown} using trigonometric polynomials is possible via Monte Carlo sampling with a number of samples that is proportional to $s$, multiplied by several log terms.

We next consider Example \ref{ex:alg-poly}. Here, we recall from \R{Theta-Legendre} that the relevant constant $\Theta^2$ for Monte Carlo sampling can become arbitrarily-large depending on the choice of $\cI$. By contrast, we now show that this situation cannot occur when sampling according to the measure \R{mu-optimal-CS}.

\prop{
Consider Example \ref{ex:alg-poly}. Let $\cI$ be an finite index set and $\Theta = \{ \phi_{\iota} : \iota \in \cI \}$, where the $\phi_{\iota}$ are as in \R{tensor-Leg-def}. Then the constant $\theta = \theta(\Phi)$ defined by \R{w-theta-opt-CS} satisfies
\bes{
\theta^2 < 2^d.
}
In particular, when sampling from the corresponding measure \R{mu-optimal-CS}, sample complexity estimate \R{opt-CS-estimate} is implied by
\be{
\label{opt-suffic-leg}
m \gtrsim 2^d \cdot s \cdot \left ( \log(\E n) \cdot \left ( d + \log(\E s ) \right )^2 + \log(2/\epsilon) \right ).
}
}
\prf{
The univariate Legendre polynomials $\phi_{\iota}$ satisfy the \textit{envelope bound}
\be{
\label{envelope-bound}
| \phi_{\iota}(y) | < \frac{2}{\sqrt{\pi} (1-y^2)^{1/4}},\quad -1 < y <1,\iota \in \bbN_0.
}
See, for example, \cite[Eqn.\ (5.3)]{adcock2018infinite}. Observe that $\int^{1}_{-1} \frac{4}{\pi \sqrt{1-y^2}} \frac{\D y}{2} = 2$. The result now follows by taking tensor products.
}

This result states that the sample complexity for the `optimal' measure is linear in $s$, up to a constant that scales at worst like $2^d$. This is a marked improvement over the sample complexity bounds for Monte Carlo sampling. Moreover, as we see in the examples below, the constant $\theta^2$ can be substantially smaller than $2^d$ for certain choices of index set $\cI$.

\rem{
[The preconditioning scheme]
The envelope bound \R{envelope-bound} suggests an alternative strategy for choosing $w$, based on the choice
\bes{
w(y) = \prod^{d}_{k=1} (\pi/2) (1-y^2_k)^{1/2}.
}
This is sometimes termed the \textit{preconditioning} technique for sparse approximation with Legendre polynomials \cite{rauhut2012sparse,jakeman2017generalized}. The corresponding measure $\mu$ is precisely the arcsine (Chebyshev) measure.
Because of \R{opt-suffic-leg}, it leads to the same sufficient sample complexity bound \R{opt-suffic-leg} as sampling via the `optimal' measure \R{w-theta-opt-CS}. 
}

We now explore this effect numerically in the setting of Example \ref{ex:alg-poly}. In order to avoid the difficulties of sampling from the continuous measure \R{mu-optimal-CS}, we instead use the discrete measure \R{mu-optimal-CS-disc} throughout. The fine grid consists of $k = 10 n$ Monte Carlo points (recall Remark \ref{r:MC-grid}).

\begin{figure}[t]
	\begin{center}
	\begin{small}
		\begin{tabular}{ccc}
Tensor product: $\cI = \cI^{\mathrm{TP}}_t$ & Total degree:  $\cI = \cI^{\mathrm{TD}}_t$ & Hyperbolic cross: $\cI = \cI^{\mathrm{HC}}_t$ \\
		\includegraphics[width=0.33\textwidth]{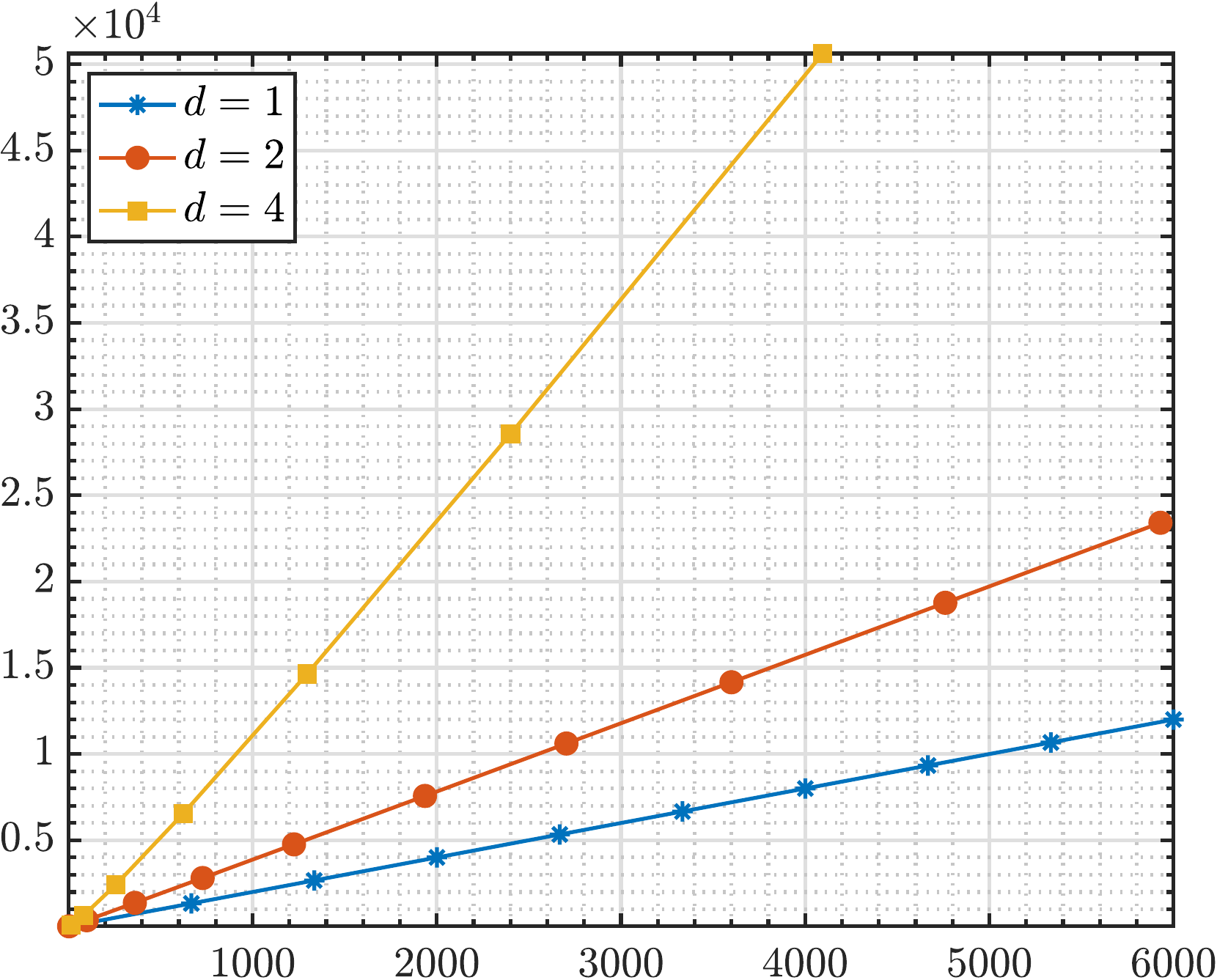} &
		\includegraphics[width=0.33\textwidth]{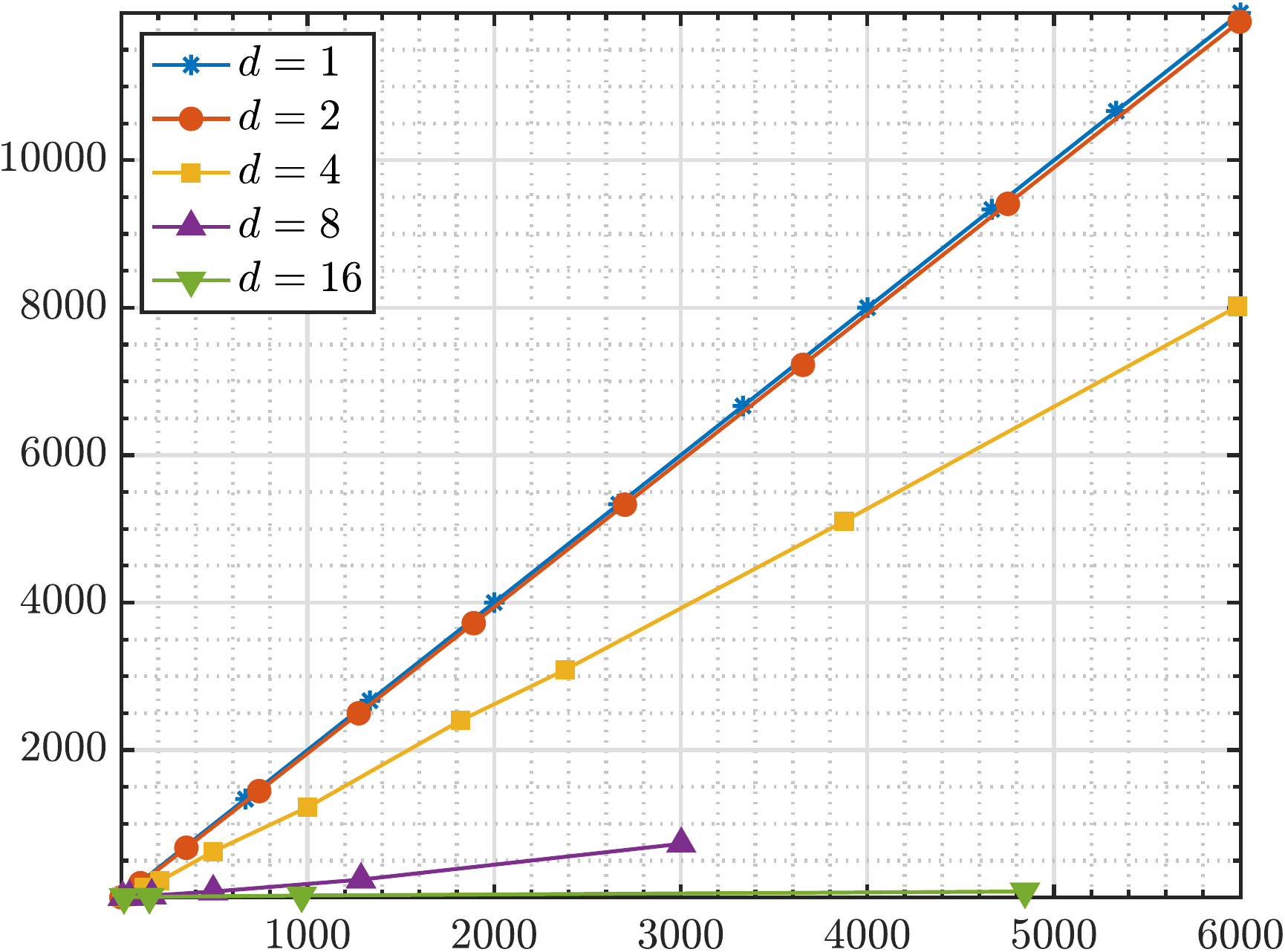} &
		\includegraphics[width=0.33\textwidth]{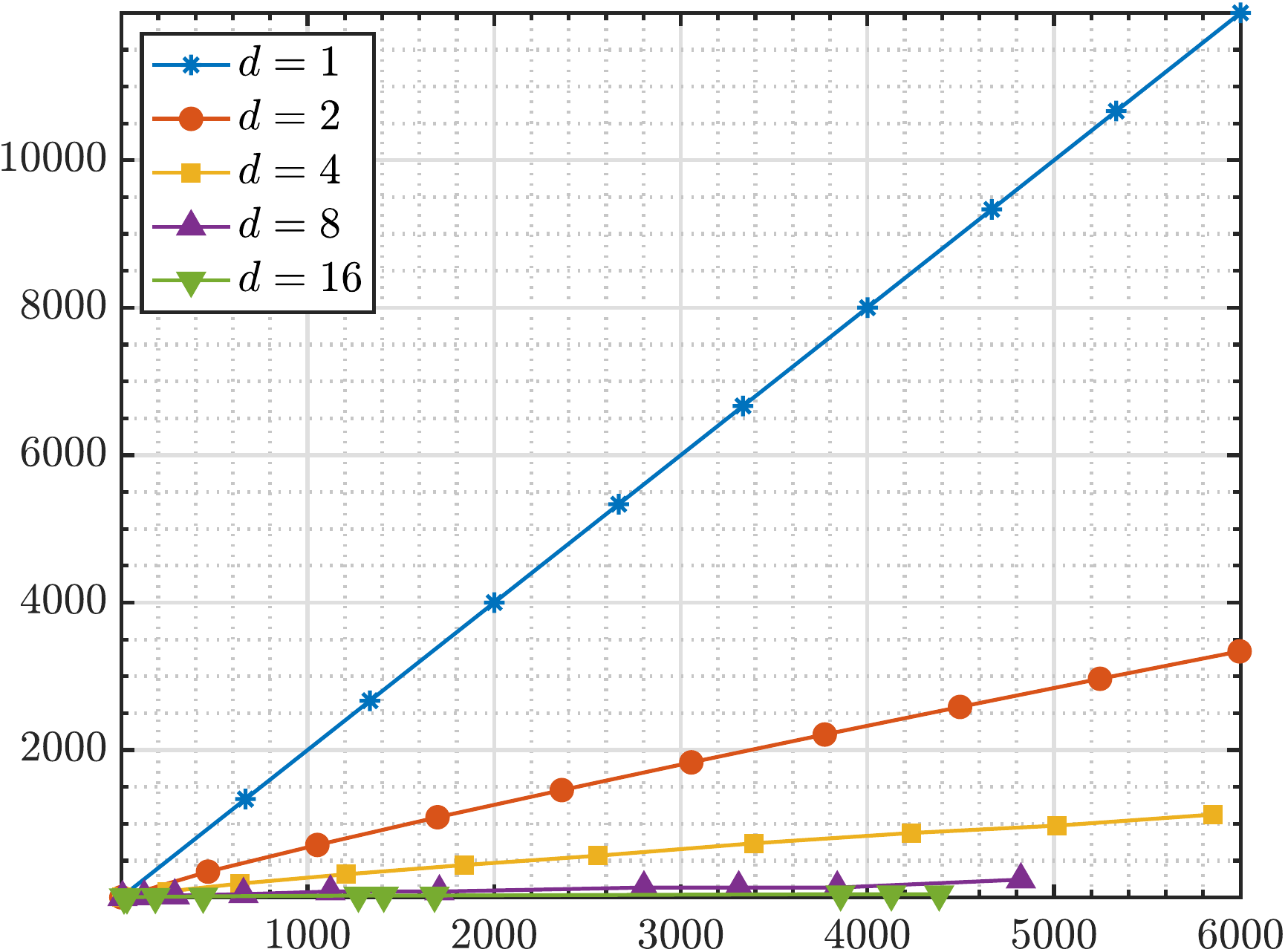}\\
		\includegraphics[width=0.33\textwidth]{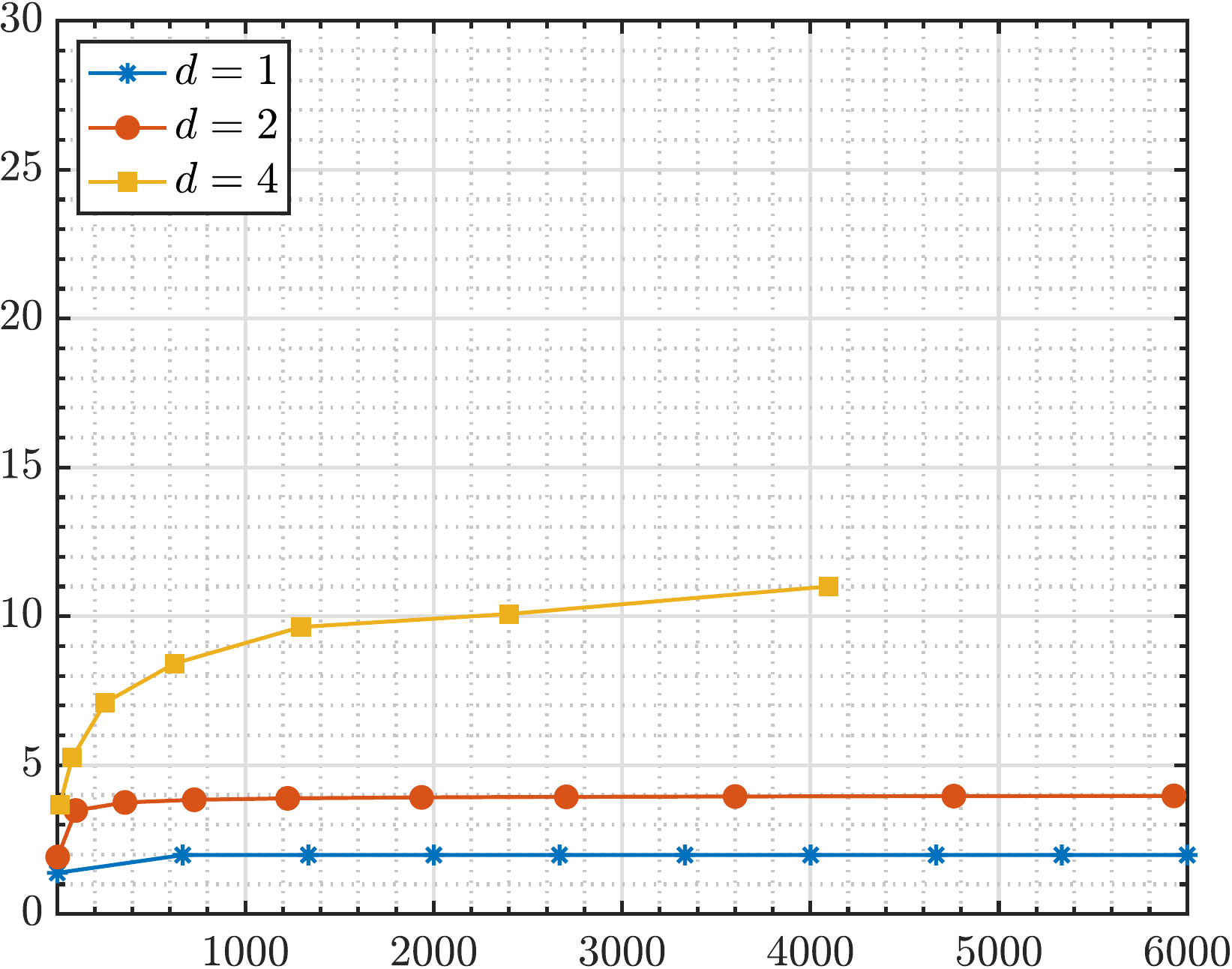} &
		\includegraphics[width=0.33\textwidth]{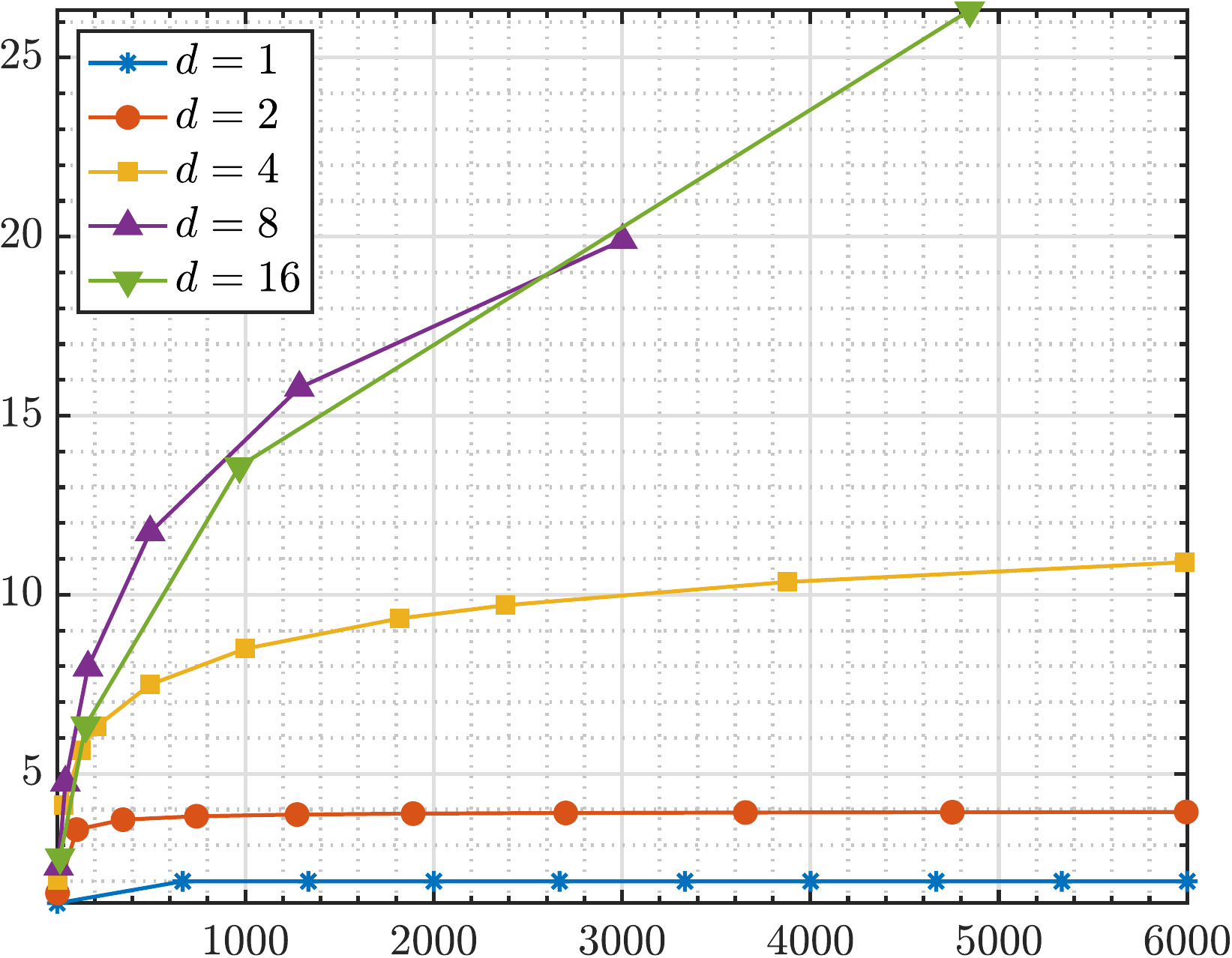} &
		\includegraphics[width=0.33\textwidth]{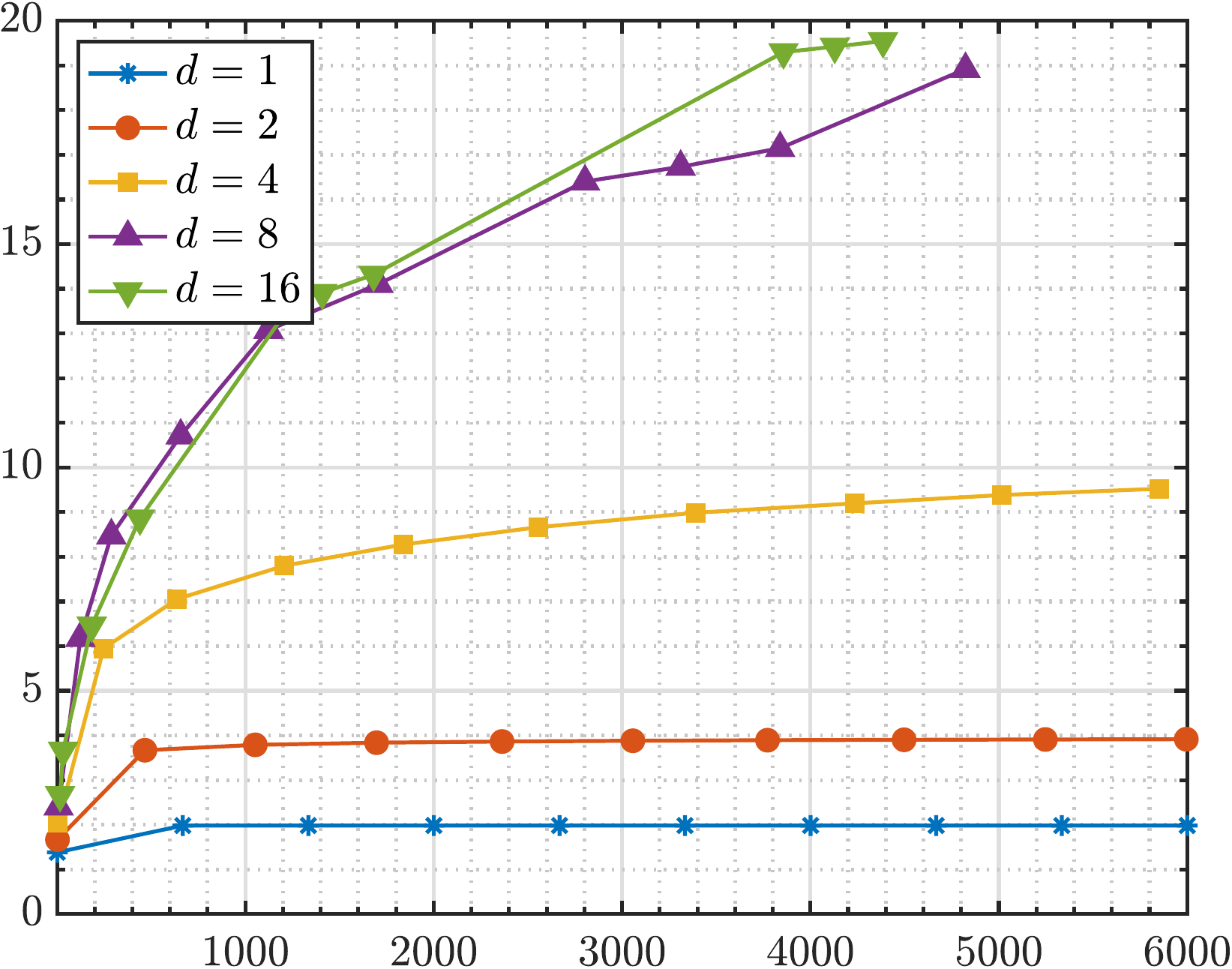} 
		\end{tabular}
	\end{small}
\end{center}
\caption{The constants $\Theta^2$ (top row), given by \R{Theta-Legendre},  and $\theta^2$ (bottom row), given by \R{theta-discrete-opt}, against $n = |\cI |$ in the case of Example \ref{ex:alg-poly} for different choices of $\cI$ and dimensions $d$.
} 
\label{fig:theta-consts}
\end{figure}

In Figure \ref{fig:theta-consts} we plot the constants $\Theta^2$ and $\theta^2$ for several different choices of $\cI$.  We notice several key effects. First, the constant $\theta^2$ is small -- in fact, no greater than $\approx 25$ in all cases. It is much smaller than the bound $2^d$ shown above (notice that $2^{16} = 65,536$), which appears to be very pessimistic in practice. We also observe that $\theta^2$ is several times smaller than $\Theta^2$, suggesting better sample complexity when sampling from \R{mu-optimal-CS-disc} instead of Monte Carlo sampling. On the other hand, the difference between the two quantities lessens in higher dimensions. This is not surprising, since the bad scaling of $\Theta$ is caused by the presence of high polynomial indices. For fixed maximum size $|\cI | = n$, the index set $\cI$ contains fewer higher-order polynomials in higher dimensions than in lower dimensions. This suggests that Monte Carlo sampling may become more acceptable in higher dimensions. We show this effect in more detail next.

\begin{figure}[t]
	\begin{center}
	\begin{small}
 \begin{tabular}{@{\hspace{0pt}}c@{\hspace{-0.5pc}}c@{\hspace{0pt}}}
		
 		\includegraphics[width=0.5\textwidth]{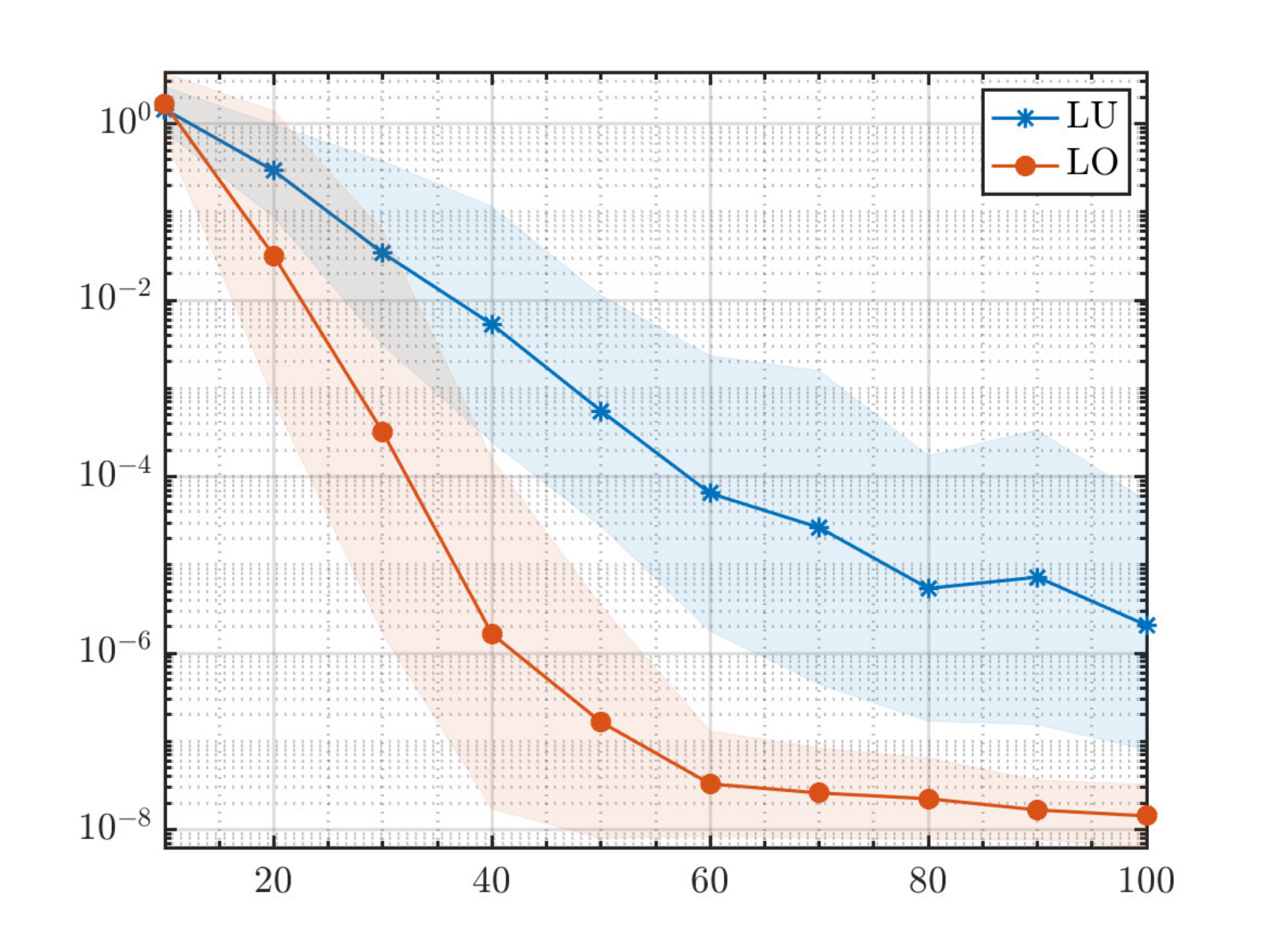} & 
 		\includegraphics[width=0.5\textwidth]{fig_5_1_1-eps-converted-to.pdf} \\
 		$(d,t,n) = (1,399,400)$ & $(d,t,n)=(2,152,796)$ \\
 		\includegraphics[width=0.5\textwidth]{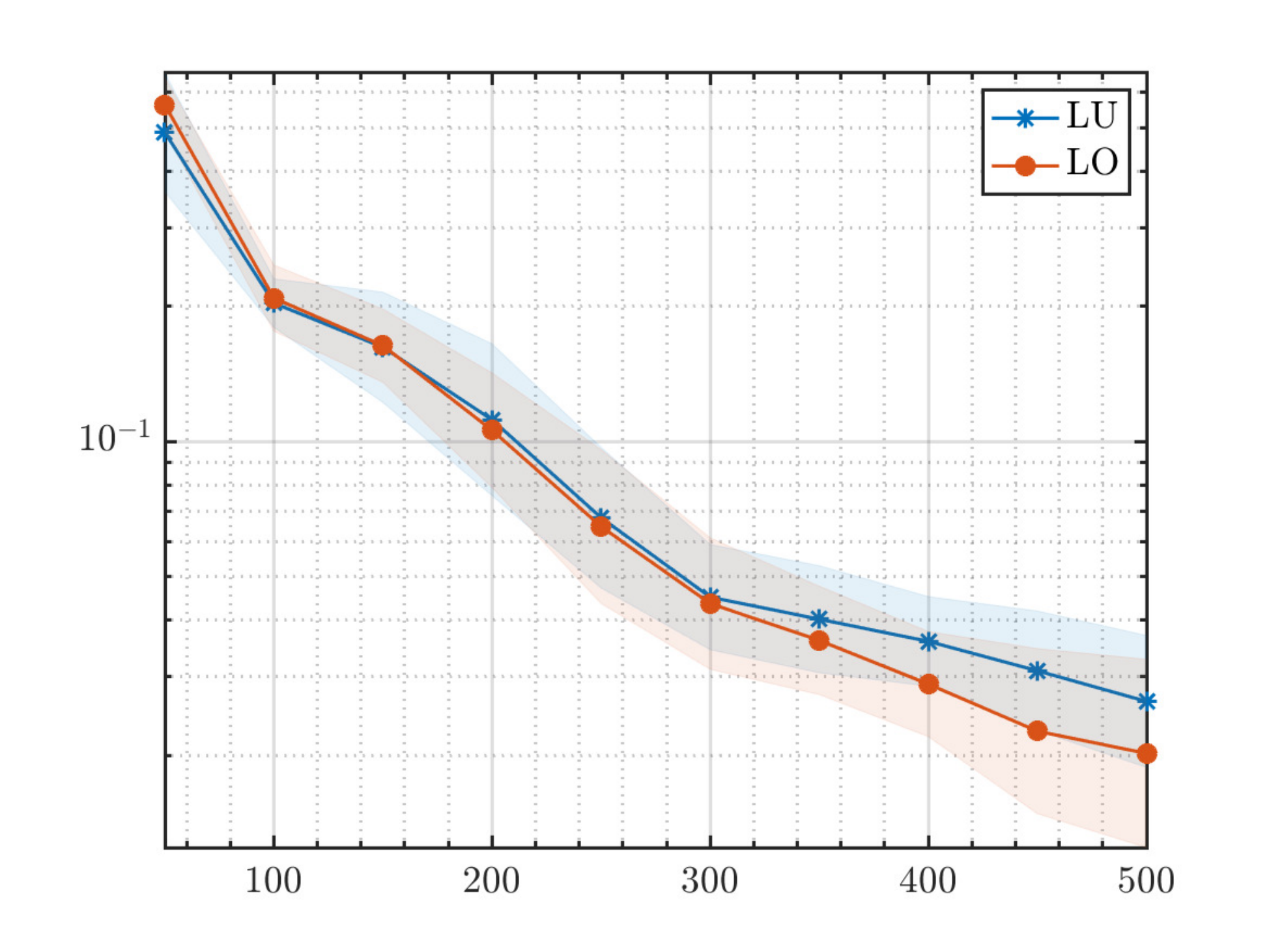} & 
 		\includegraphics[width=0.5\textwidth]{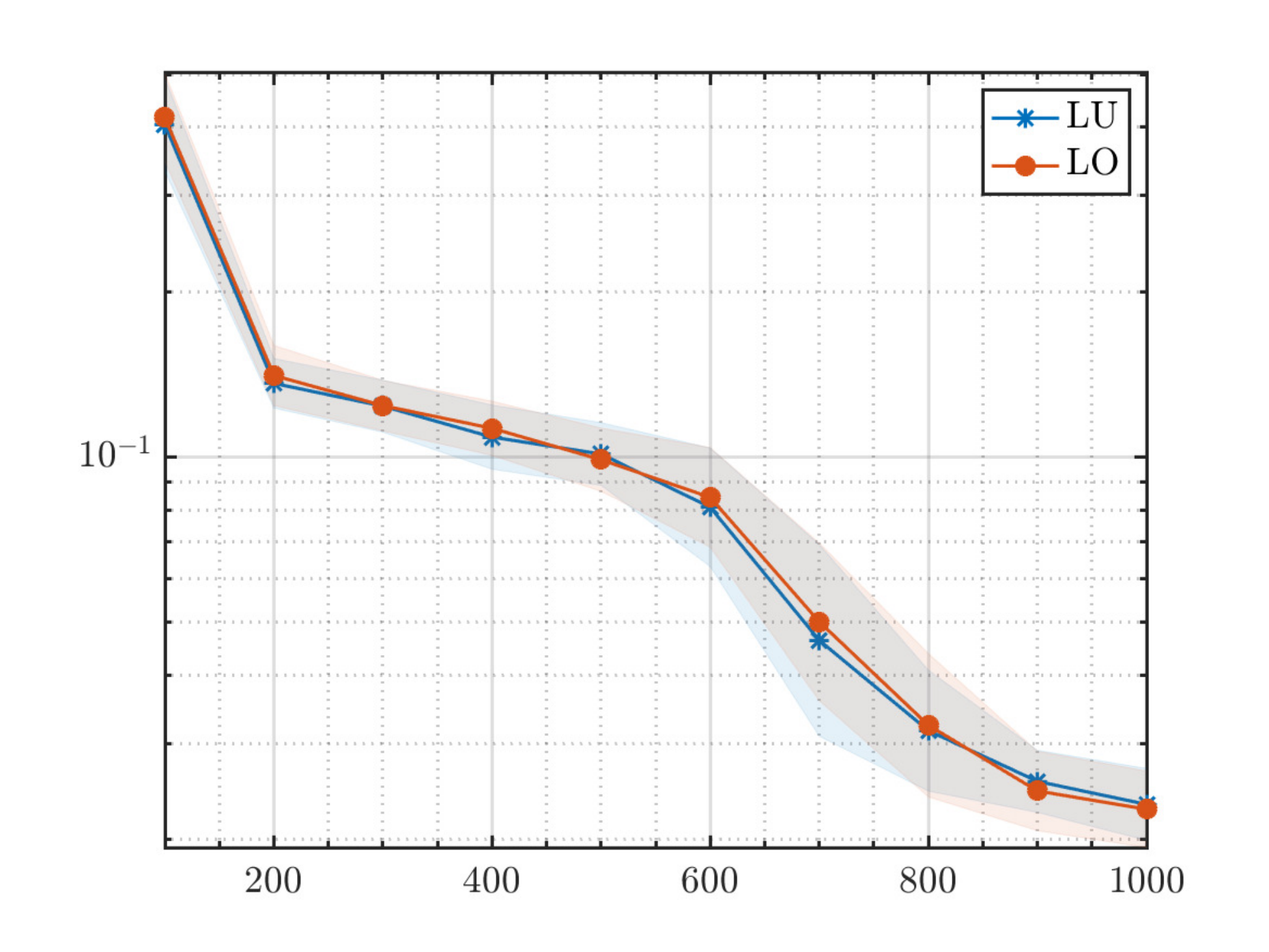} \\
 		$(d,t,n) = (8,22,1843)$ & $(d,t,n)=(16,14,4385)$ \\
	\end{tabular}
	\end{small}
\end{center}
\caption{The relative error \R{relative-error} versus $m$ for $\ell^1$-minimization in the case of Example \ref{ex:alg-poly}, where $\cI = \cI^{\mathrm{HC}}_{t-1}$ is the hyperbolic cross index set \R{HC-index} and $f = f_1$ is as in \R{functions-define}. This figure compares Monte Carlo uniform sampling (labelled `LU') and sampling from the discrete `optimal' measure \R{mu-optimal-CS-disc} (labelled `LO').}
\label{fig:cs-approx-1}
\end{figure}

\begin{figure}[t]
	\begin{center}
	\begin{small}
 \begin{tabular}{@{\hspace{0pt}}c@{\hspace{-0.5pc}}c@{\hspace{0pt}}}
		
 		\includegraphics[width=0.5\textwidth]{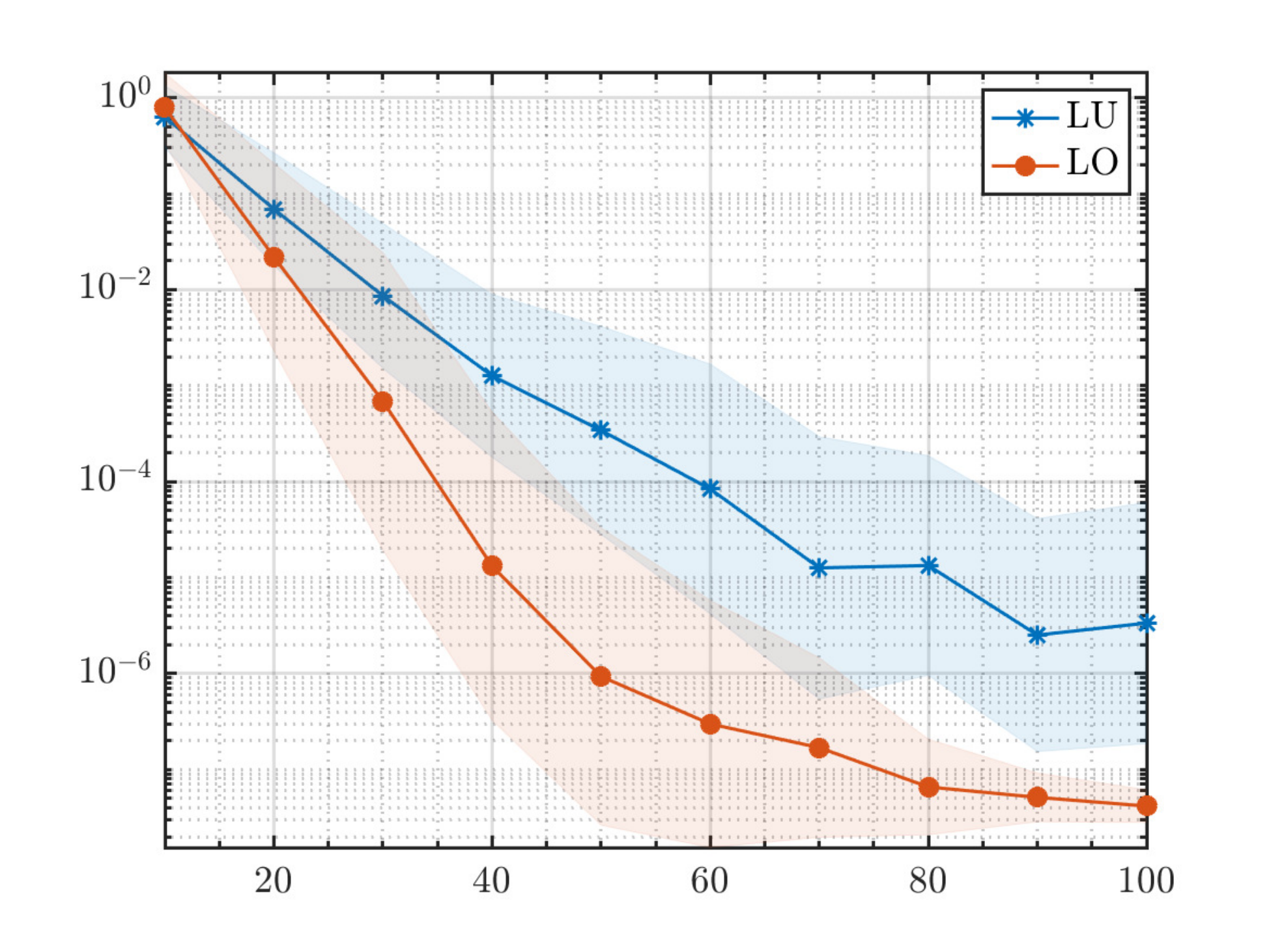} & 
 		\includegraphics[width=0.5\textwidth]{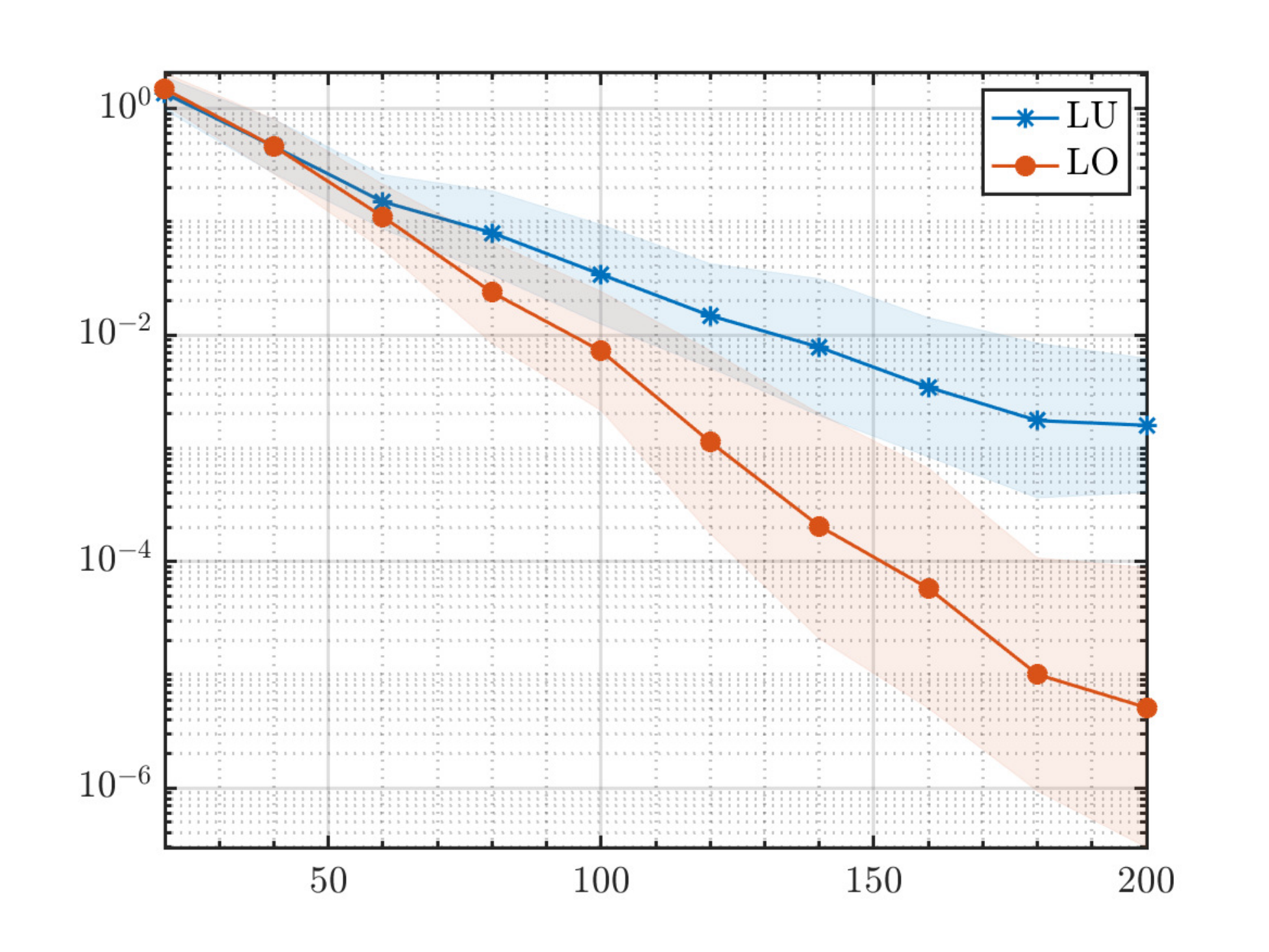} \\
 		$(d,n,N) = (1,399,400)$ & $(d,n,N)=(2,152,796)$ \\
 		\includegraphics[width=0.5\textwidth]{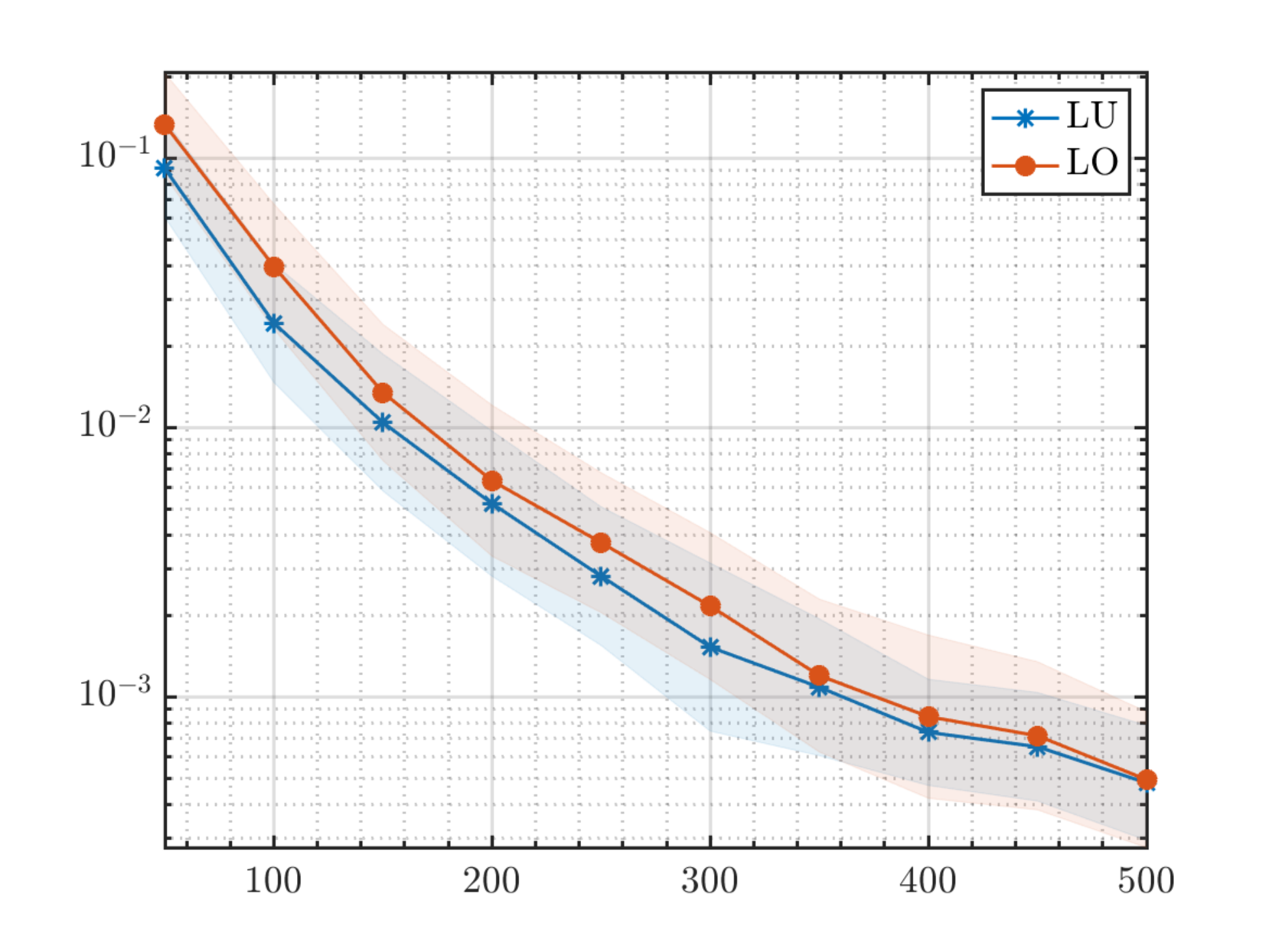} & 
 		\includegraphics[width=0.5\textwidth]{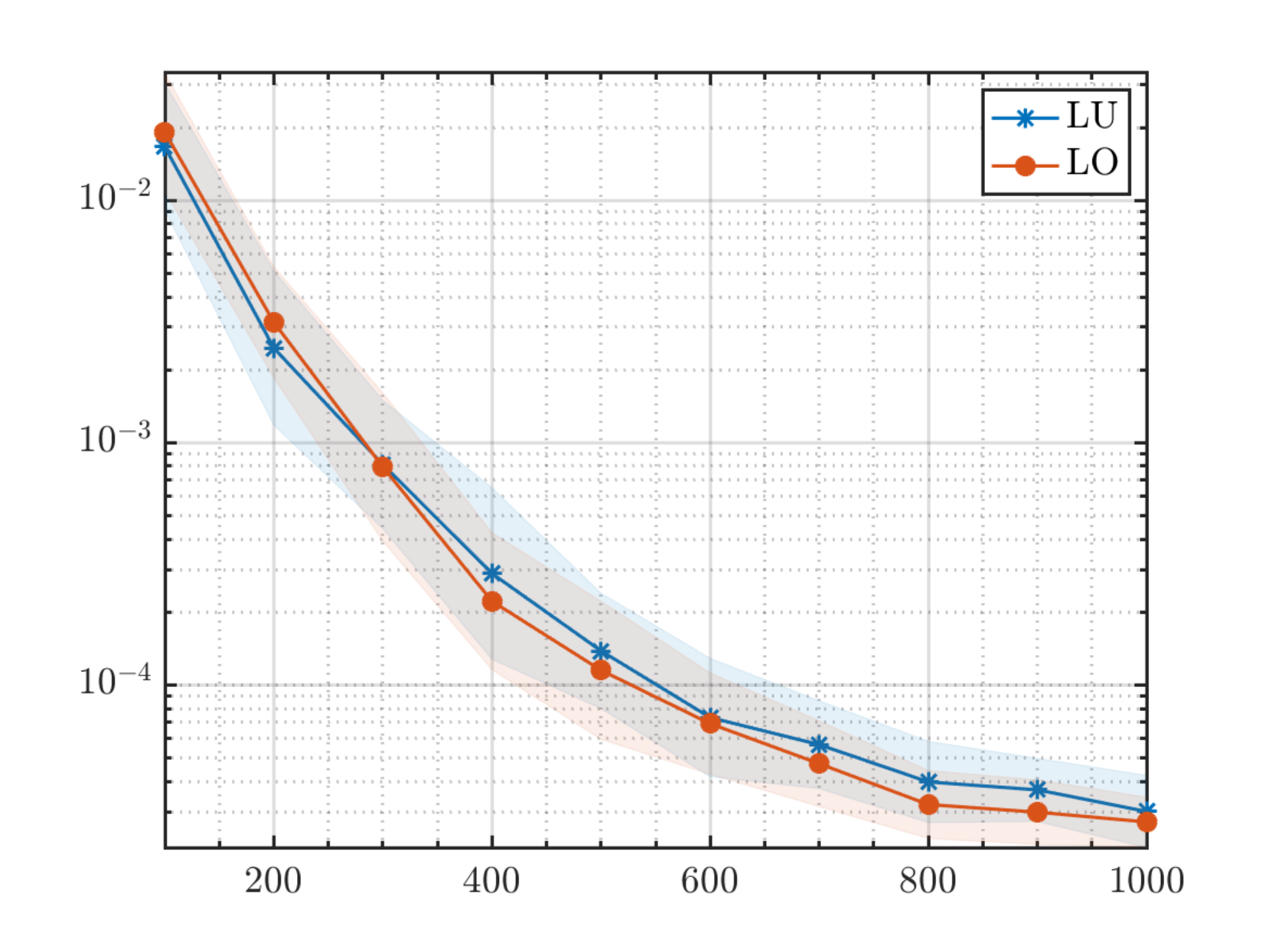} \\
 		$(d,n,N) = (8,22,1843)$ & $(d,n,N)=(16,14,4385)$ \\
	\end{tabular}
	\end{small}
\end{center}
\caption{The same as in Fig.\ \ref{fig:cs-approx-2} but with $f = f_2$ as in \R{functions-define}.}
\label{fig:cs-approx-2}
\end{figure}

\begin{figure}[t]
	\begin{center}
	\begin{small}
 \begin{tabular}{@{\hspace{0pt}}c@{\hspace{-0.5pc}}c@{\hspace{0pt}}}
		
 		\includegraphics[width=0.5\textwidth]{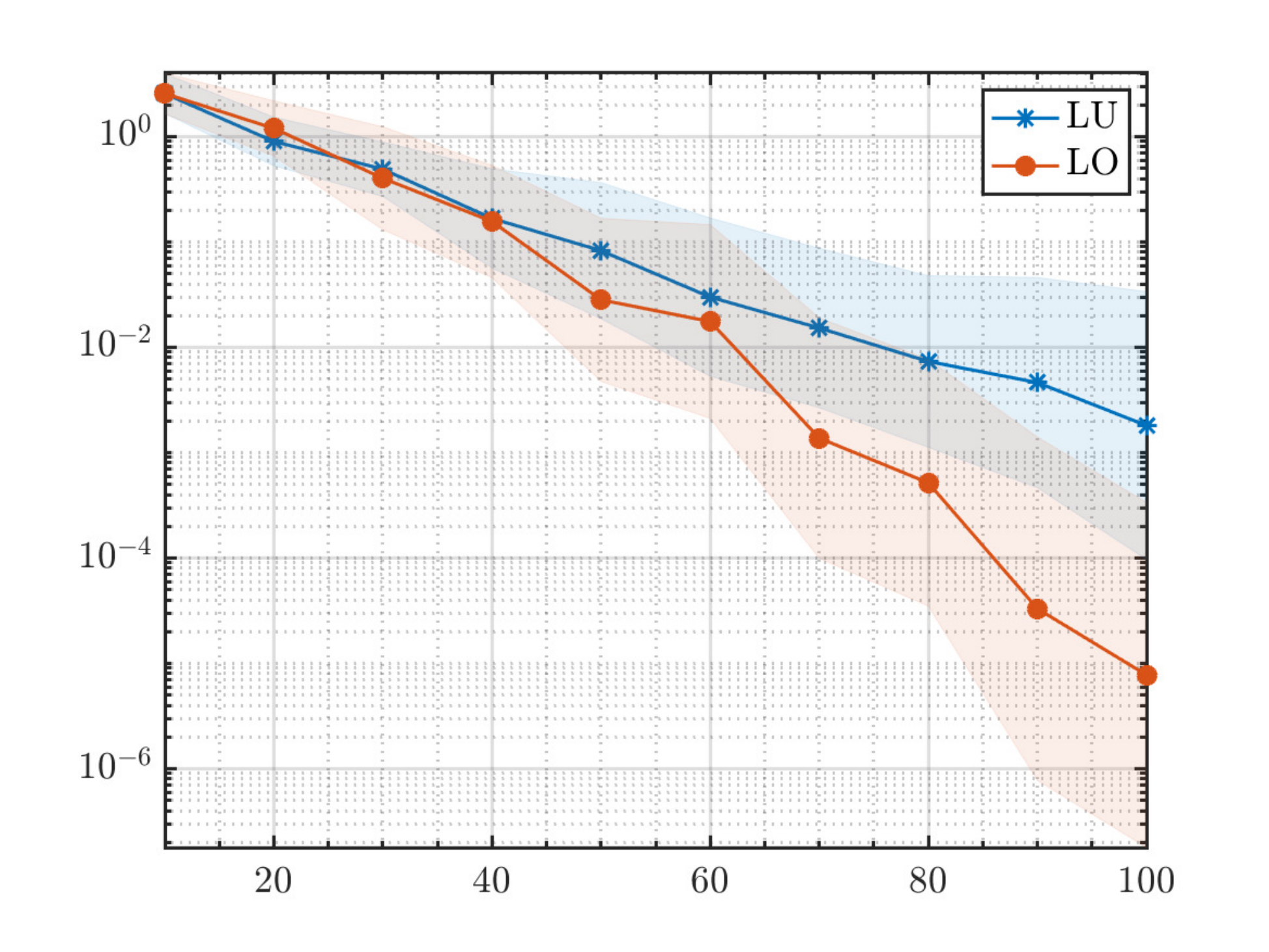} & 
 		\includegraphics[width=0.5\textwidth]{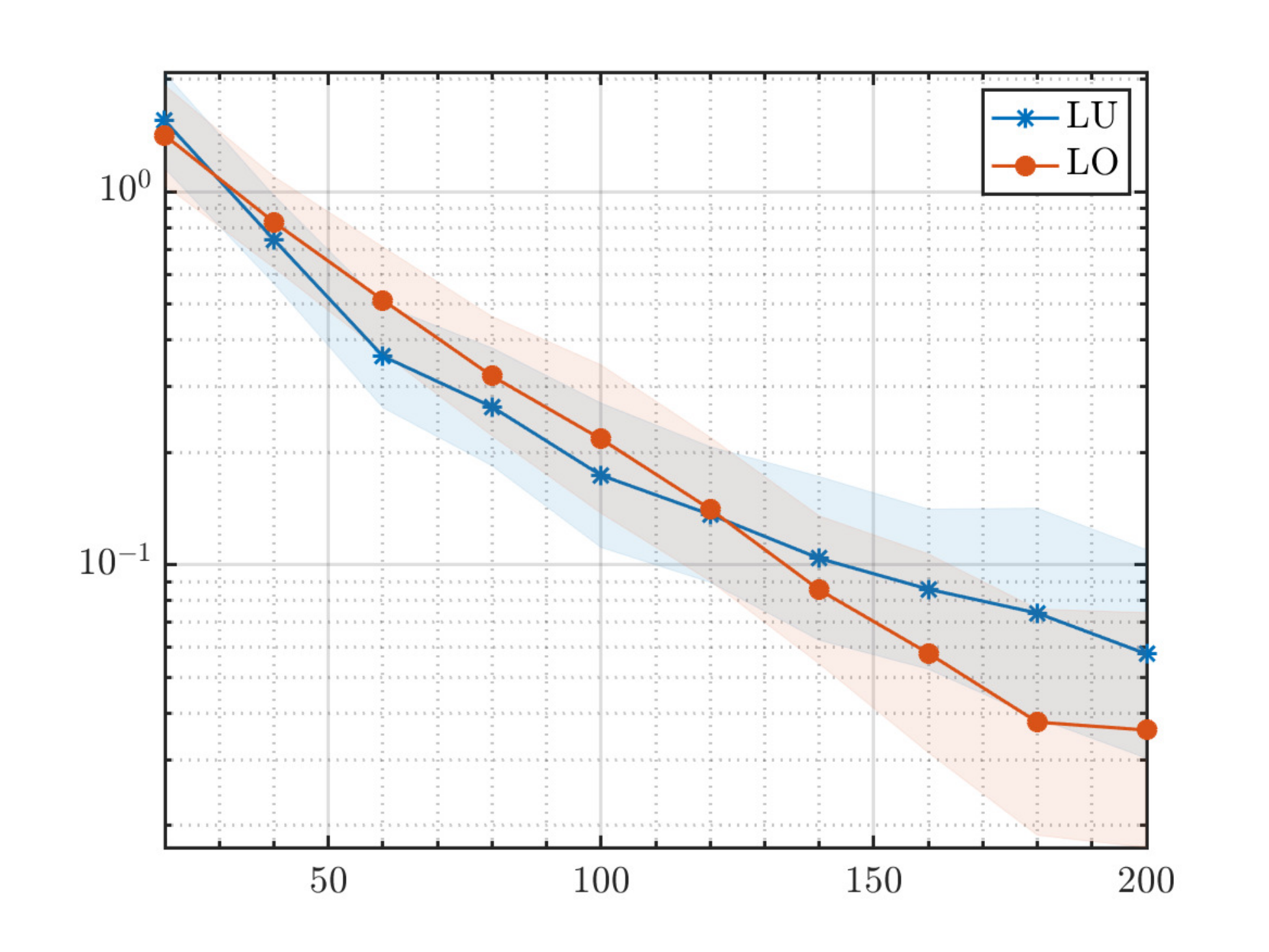} \\
 		$(d,n,N) = (1,399,400)$ & $(d,n,N)=(2,152,796)$ \\
 		\includegraphics[width=0.5\textwidth]{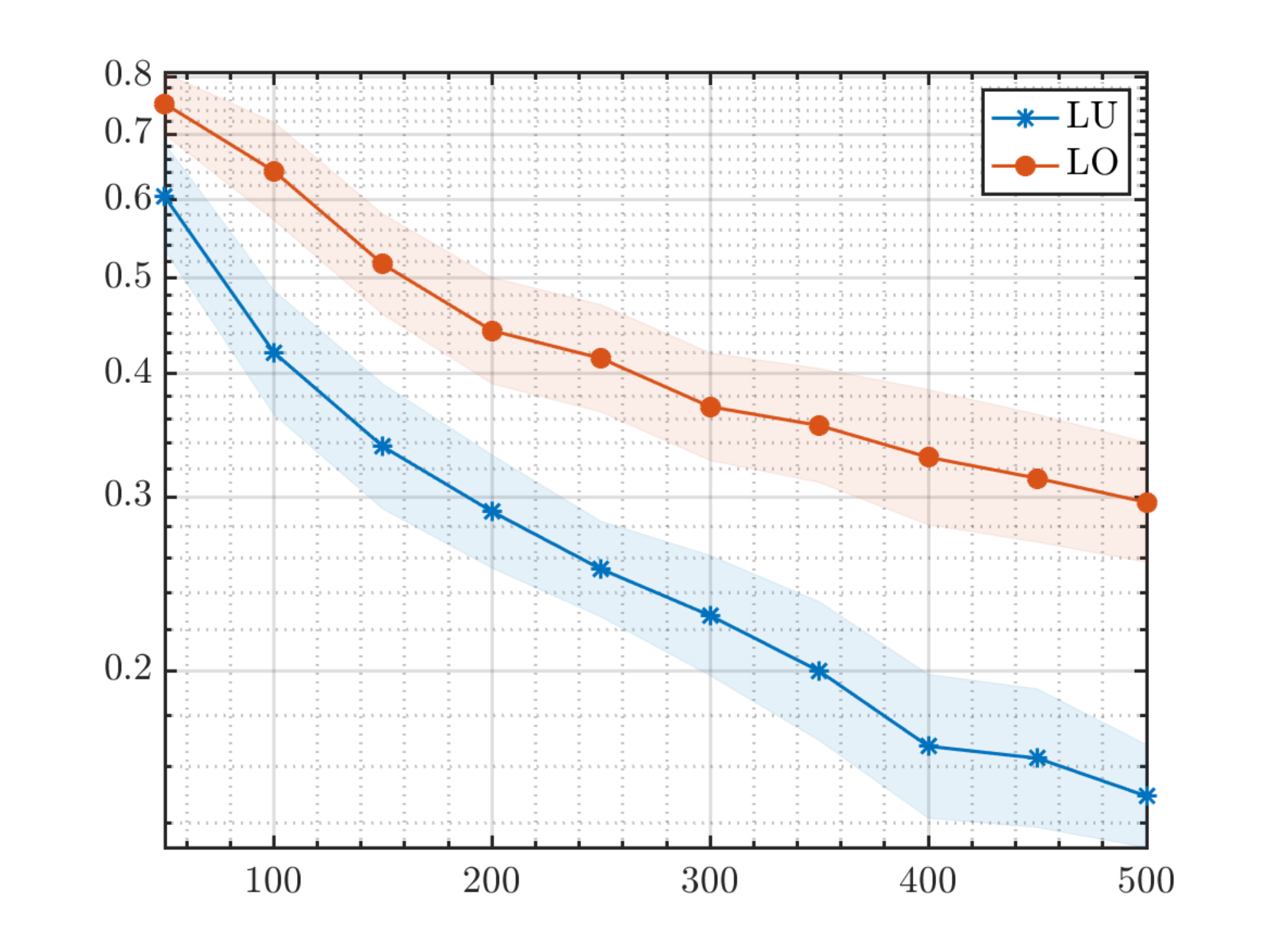} & 
 		\includegraphics[width=0.5\textwidth]{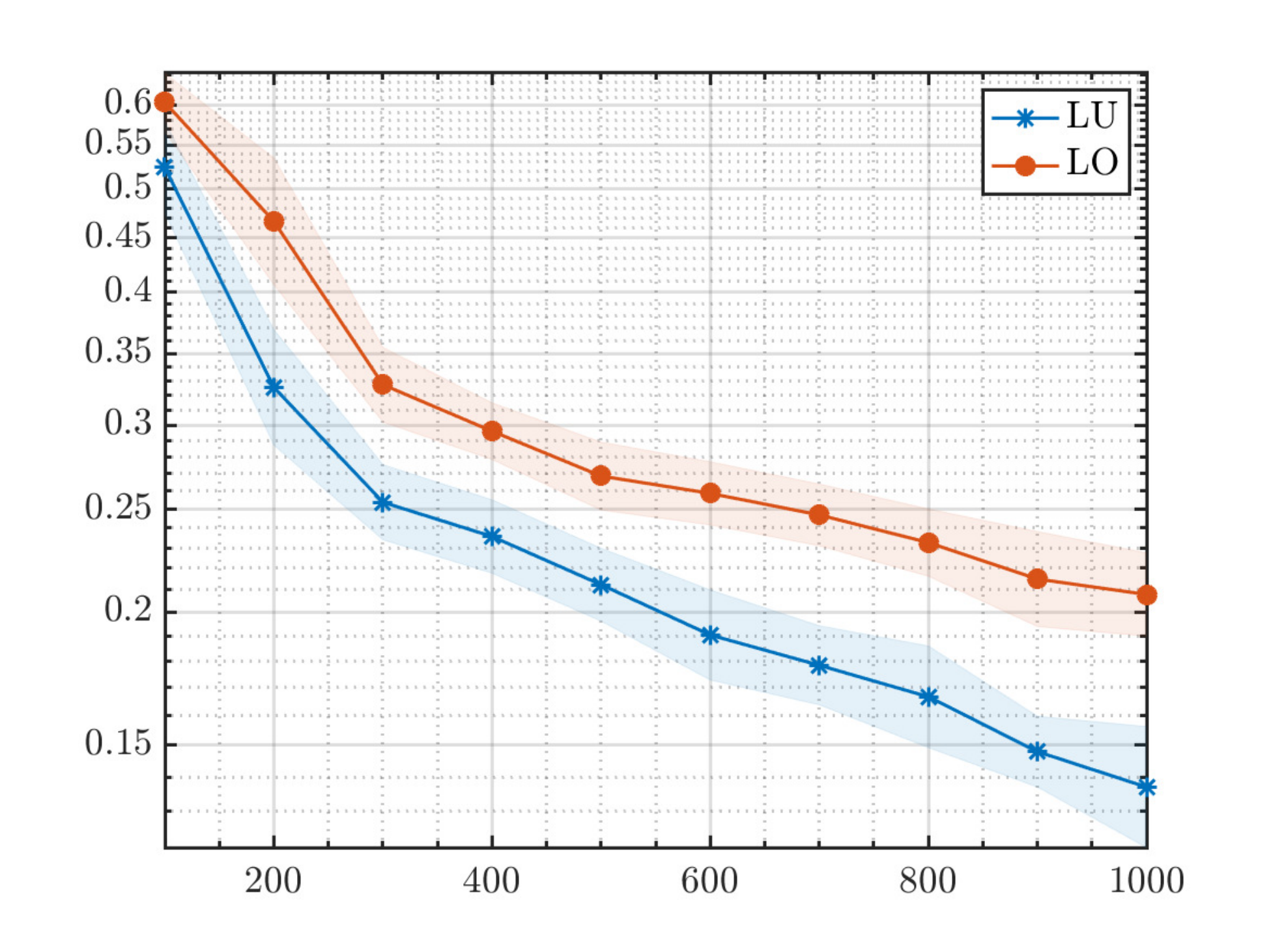} \\
 		$(d,n,N) = (8,22,1843)$ & $(d,n,N)=(16,14,4385)$ \\
	\end{tabular}
	\end{small}
\end{center}
\caption{The same as Figure \ref{fig:cs-approx-1} but with $f = f_3$ as in \R{functions-define}. 
} 
\label{fig:cs-approx-3}
\end{figure}

In Figures \ref{fig:cs-approx-1}--\ref{fig:cs-approx-3} we consider function approximation using the two sampling strategies.
As we see, the `optimal' strategy gives a nonnegligible benefit over Monte Carlo sampling in lower dimensions for $f_1$ and $f_2$. Yet, in higher dimensions this benefit lessens. On the other hand, for $f = f_3$ the `optimal' strategy yields no better performance, and actually a larger error than Monte Carlo sampling in high dimensions. The lessening benefit with increasing dimension is consistent with the results of Figure \ref{fig:theta-consts}, wherein it is shown that the difference between the constants $\Theta^2$ and $\theta^2$ decreases as $d$ increases. On the other hand, the observation that it can sometimes yield worse approximations is an important reminder that sampling strategy is designed to enhance the performance of sparse approximation in general, and may not therefore be the best strategy for any fixed function.

\rem{
To elaborate on this previous comment, notice that the `optimal' measure samples more densely near the boundary of the hypercube $[-1,1]^d$, where the Legendre polynomials are larger, and therefore less densely near the origin. Hence, any function that varies most significantly in the interior of the domain is liable to be less well approximated by sampling from the `optimal' measure. This is the case in particular for the function $f = f_3$, which is a product of one dimensional functions that are peaked around centres that get progressively closer to the origin with increasing index $i$, and which are relatively flat away from their centre.
}

\subsection{Proof of Theorems \ref{t:CS_acc_stab} and \ref{t:BOS_RIP}}\label{ss:CSproofs}

Throughout this section, if $x = (x_i)^{n}_{i=1} \in \bbV^n$ and $S \subseteq [n]$ we use the notation $x_{S} \in \bbV^n$ to denote the vector with $i$th entry equal to $x_i$ if $i \in S$ and $0$ otherwise. Note that $x_{S}$ is isomorphic to a vector in $\bbV^{|S|}$. We will sometimes consider it as an element of this space. 
We now recall the following definition and lemma, which can be found in \cite[Defn.\ 6 \& Lem.\ 7]{adcock2021deep}:

\begin{definition}
A matrix $A \in \bbC^{m \times n}$ satisfies the \textit{robust Null Space Property (rNSP)} of order $1 \leq s \leq n$ over $\bbV^n$ with constants $0 < \rho < 1$ and $\tau > 0$ if
\bes{
\nm{x_S}_{\ell^2([n];\bbV)} \leq  \frac{\rho \nmu{x_{S^c}}_{\ell^1([n];\bbV)} }{\sqrt{s}}+\tau \nm{A x}_{\ell^2([m];\bbV)} ,\quad \forall x\in \bbV^N,
}
for any $S \subseteq [n]$ with $|S| \leq s$.
\end{definition}

Here, in the second term, we recall that a matrix $A \in \bbC^{m \times n}$ extends in the obvious way to a mapping $\bbV^n \rightarrow \bbV^m$.

\begin{lemma}\label{l:rNSP_SRLASSO}
Suppose that $A \in \bbC^{m \times n}$ has the rNSP of order $1 \leq s \leq n$ with constants $0 < \rho < 1$ and $\tau > 0$. Let $x \in \bbV^n$, ${v} = A x + e\in \bbV^m$ and
\bes{
\lambda \leq \frac{C_1}{C_2 \sqrt{s}},
}
where $C_1 = \frac{(3 \rho+1)(\rho+1)}{2(1-\rho)}$ and $C_2 = \frac{(3 \rho+5) \tau}{2(1-\rho)}$. Then every minimizer $\hat{x} \in \bbV^n$ of the Hilbert-valued SR-LASSO problem
\bes{
\min_{z \in \bbV^n} \lambda \nm{z}_{\ell^1([n];\bbV)} + \nm{A z - {v}}_{\ell^2([m];\bbV)},
}
satisfies
\bes{
\nm{\hat{x} - x}_{\ell^2([n];\bbV)} \leq 2 C_1 \frac{\sigma_s(x)_{\ell^1([n];\bbV)}}{\sqrt{s}} + \left ( \frac{C_1}{\sqrt{s} \lambda} + C_2 \right ) \nm{e}_{\ell^2([m];\bbV)}.
}
\end{lemma}

\prf{[Proof of Theorem \ref{t:CS_acc_stab}]
We first claim that the matrix $A$ defined in \R{A-CS-def} satisfies
\be{\label{Double_RIP}
a \alpha \nm{c}^2_{\ell^2([n];\bbV)} \leq \nm{A c}^2_{\ell^2([m];\bbV)} \leq b \beta \nm{c}^2_{\ell^2([n];\bbV)},\quad 
}
for all vectors $c \in \bbV^n$ that are $t$-sparse. Observe that any such vector $c$ corresponds to the coefficients of an element $p \in P_{T;\bbV}$ for $T \subseteq \cI$, $|T| \leq t$, and also that $\nm{A c}^2_{\ell^2([m];\bbV)} = \frac1m \sum^{m}_{i=1} w(y_i) \nm{p(y_i)}^2_{\bbV}$. We now recall from Lemma \ref{l:RIP_over_S_C_V} that \R{RIP_over_all_S} also holds over $P_{T;\bbV}$ whenever it holds over $P_{T} \equiv P_{T;\bbC}$. Hence
\be{
\label{equiv-over-V-T}
\alpha \nm{p}^2_{L^2_{\rho}(D;\bbV)} \leq \nm{A c}^2_{\ell^2([m];\bbV)}  \leq \beta \nm{p}^2_{L^2_{\rho}(D;\bbV)} 
}
Now, using almost identical arguments to those used in the proof of Lemma \ref{l:RIP_over_S_C_V}, we find that the Riesz basis condition \R{Riesz-bounds} also extends to the $\bbV$-valued case:
\be{
\label{Riesz-over-V}
a \nm{c}^2_{\ell^2(\cI;\bbV)} \leq \nm{\sum_{\iota \in \cI} c_{\iota} \phi_{\iota} }^2_{L^2_{\rho}(D;\bbV)} \leq b \nm{c}^2_{\ell^2(\cI;\bbV)},\quad \forall c = (c_{\iota})_{\iota \in \cI} \in \bbV^n.
}
Since $p = \sum_{\iota \in \cI} c_{\iota} \phi_{\iota}$, the claim now follows from this and \R{equiv-over-V-T}.

We now show that $A$ satisfies the rNSP of order $s$ over $\bbV^n$ and derive values for the constants $\rho$ and $\tau$. The following argument is based on \cite[Lem.\ 13.8]{adcock2021compressive}.
Let $c \in \bbV^n$ and $S \subseteq [n]$ with $|S| \leq s$. Suppose first that $t < n$, so that $t= 2 \lceil 4 s \frac{b\beta}{a\alpha} \rceil $ (we consider the case $t = n$ later).
Define a partition $\Delta_1, \Delta_2 , \ldots$ of $S^c$ as follows.  
First, let $\Delta_1$ be the index set of the largest $t' = t/2$ (notice that this is an integer, due to the definition of $t$) indices of the vector $(\nmu{c_j}_{\bbV})_{j \in S^c}$,  $\Delta_2$ be the index set of the largest $t' =t/2$ indices of the vector  $(\nmu{c_j}_{\bbV})_{j \in (S \cup \Delta_1)^c}$ and so forth. This gives a partition of $S^c$ for which each set is of size $t'$, except possibly the final set (this is of no consequence to the argument). Consider the set  $S \cup \Delta_1 \supseteq S$. Since $s \leq t'$ we have $|S \cup \Delta_1| \leq s +t' \leq 2t' =t$. Hence, we may apply \R{Double_RIP} to obtain
\be{
\label{S-to-ASDelta}
\nm{P_{S} c}^2_{\ell^2([n];\bbV)} \leq \nm{P_{S \cup \Delta_1} c}^2_{\ell^2([n];\bbV)} \leq \frac{1}{ a \alpha} \nm{A P_{S \cup \Delta_1} c}^2_{\ell^2([m];\bbV)}.
}
We now write
\bes{
\nm{A P_{S \cup \Delta_1} c}^2_{\ell^2([m];\bbV)} = \ip{A P_{S \cup \Delta_1} c}{A c}_{\ell^2([m];\bbV)} - \sum_{i \geq 2} \ip{A P_{S \cup \Delta_1} c}{A P_{\Delta_i} c}_{\ell^2([n];\bbV)},
}
and then apply the Cauchy--Schwarz inequality and the \R{Double_RIP} once more to get
\bes{
\nm{A P_{S \cup \Delta_1} c}_{\ell^2([m];\bbV)} \leq \nm{A c}_{\ell^2([m];\bbV)} + \sqrt{b \beta} \sum_{i \geq 2} \nm{ P_{\Delta_i} c}_{\ell^2([n];\bbV)}.
}
Notice that, by construction of $\Delta_i$ for each $i \geq 2$, we have 
\[
\nm{P_{\Delta_i}c}_{\ell^2([n];\bbV)}^2 \leq t' \max_{j \in \Delta_i} \nm{c_j}_{\bbV}^2 \leq t' \min_{j \in \Delta_{i-1}} \nm{c_j}_{\bbV}^2.
\]
Since $\nm{ P_{\Delta_{i-1}} c}_{\ell^1([n];\bbV)}\geq t' \min_{j \in \Delta_{i-1}} \nm{c_j}_{\bbV}$, this implies the following:
\bes{
\nm{P_{\Delta_{i}} c}_{\ell^2([n];\bbV)} \leq (t')^{-1/2} \nm{P_{\Delta_{i-1}} c}_{\ell^1([n];\bbV)},\quad i \geq 2.
}
Using this and the previous expression, we deduce that
\eas{
\nm{A P_{S \cup \Delta_1} c}_{\ell^2([m];\bbV)} &\leq \nm{A c}_{\ell^2([m];\bbV)} + \sqrt{\frac{b\beta}{t'}} \sum_{i \geq 1} \nm{ P_{\Delta_i} c}_{\ell^2([n];\bbV)}
\\
& = \nm{A c}_{\ell^2([m];\bbV)} + \sqrt{\frac{b\beta}{t'}} \nm{P_{S^c} c}_{\ell^1([n];\bbV)} \\
}
Substituting this back into \R{S-to-ASDelta} and noticing that $b \beta / (a \alpha t') \leq 1/(2 \sqrt{s})$, we get
\be{
\label{Rnsp_pr}
\nm{P_{S} c}_{\ell^2([n];\bbV)} \leq \frac{1}{\sqrt{a \alpha}} \nm{A c}_{\ell^2([m];\bbV)} + \frac{1}{2 \sqrt{s}}  \nm{P_{S^c} c}_{\ell^1([n];\bbV)},
}
for the case $t < n$. Now suppose that $t = n$. Since \eqref{Double_RIP} now holds for any vector $c \in \bbV^n$, we easily see that
\bes{
\nm{P_{S} c}_{\ell^2([n];\bbV)}  \leq \frac{1}{\sqrt{a \alpha}} \nm{A   c}_{\ell^2([m];\bbV)} \leq \frac{1}{\sqrt{a \alpha}} \nm{A c}_{\ell^2([m];\bbV)} + \frac{1}{2 \sqrt{s}}  \nm{P_{S^c} c}_{\ell^1([n];\bbV)}.
}
Hence, \R{Rnsp_pr} also holds in this case as well. We deduce that $A$ satisfies the rNSP of order $s$ over $\bbV^n$ with constants $\rho=1/2$ and $\tau= 1/\sqrt{a \alpha}$.

Having shown this, we complete the proof by establishing the error bounds for $\hat{f}$. Let $\hat{ {c}}$ be the coefficients of $\hat{f}$. By the triangle inequality, we have
\eas{
&\nm{f -\hat{f}}_{L^2_{\varrho}(D ; \bbV)} 
\\
& \leq   \nm{f - \cP_h(f)}_{L^2_{\varrho}(D ; \bbV)}  + \nm{\cP_h(f) - \cP_{h}(f_{\cI}) }_{L^2_{\varrho}(D ; \bbV)} + \nm{\cP_h(f_{\cI}) - \hat{f}}_{L^2_{\varrho}(D ; \bbV)}  
}
Therefore, since $\cP_h$ is an orthogonal projection we have
\bes{
\begin{split}
& \nm{f -\hat{f}}_{L^2_{\varrho}(D ; \bbV)} 
\\
& \leq  \nm{f - \cP_h(f)}_{L^2_{\varrho}(D ; \bbV)}  + \nm{f-f_{\cI}}_{L^2_{\varrho}(D ; \bbV)} + \nm{\cP_h(f_{\cI}) - \hat{f}}_{L^2_{\varrho}(D ; \bbV)}  .    
\end{split}
}
We first bound the term $\nm{\cP_h(f_{\cI}) - \hat{f}}_{L^2_{\varrho}(D ; \bbV)}$. By \R{Riesz-over-V}, we have
\bes{
\nm{\cP_h(p) - \hat{f}}_{L^2_{\varrho}(D ; \bbV)} = \nm{\sum_{\iota \in \cI} (\cP_{h}(c_{\iota}) - \hat{c}_{\iota} )\phi_{\iota} }_{L^2_{\varrho}(D ; \bbV)} \leq \sqrt{b} \nm{\cP_{h}(c) - \hat{c} }_{\ell^2([n];\bbV)},
}
where $\cP_{h}(c)$ is the vector $(\cP_{h}(c_{\iota}))_{\iota \in \cI}$.
Using Lemma \ref{l:rNSP_SRLASSO}, we get
\be{
\label{step-b}
\begin{split}
& \nm{\cP_h(p) - \hat{f}}_{L^2_{\varrho}(D ; \bbV)}  
\\
& \leq 8 \sqrt{b} \left ( \frac{\sigma_{s}(\cP_{h}(c))_{\ell^1([n];\bbV)}}{\sqrt{s}} + \left ( \frac{1}{2 \lambda \sqrt{s}} + \frac{1}{\sqrt{a\alpha}} \right ) \nm{A\cP_{h}( {c})-b}_{\ell^2([m];\bbV)}  \right ).
\end{split}
}
Since $\cP_h$ is an orthogonal projection, we have
\be{
\label{step-a}
 \sigma_{s}(\cP_{h}( {c}))_{\ell^1([n];\bbV)}\leq \sigma_{s}(  {c})_{\ell^1([n];\bbV)}.
}
 Moreover, using \R{A-CS-def} and \R{b-CS-def}, we see that
\eas{
&\nm{A\cP_{h}( {c})-{v}}_{\ell^2([m];\bbV)} 
\\
&=  \nm{ \frac{1}{\sqrt{m}}\left ( \sqrt{w(y_i)} \left ( \cP_h(f_{\cI}(y_i)) - f(y_i) - n_i \right ) \right )_{i \in [m]}}_{\ell^2([m];\bbV)}
\\
& \leq \nm{ \frac{1}{\sqrt{m}}\left ( \sqrt{w(y_i)} \left ( \cP_h(f_{\cI}(y_i)) - \cP_h(f)(y_i) \right ) \right )_{i \in [m]} }_{\ell^2([m];\bbV)}
 \\
&+ \nm{ \frac{1}{\sqrt{m}}\left ( \sqrt{w(y_i)} \left ( f(y_i) - \cP_h(f)(y_i) \right ) \right )_{i \in [m]} }_{\ell^2([m];\bbV)} + \nm{e}_{\ell^2([m];\bbV)}.
}
Now, observe that
\eas{
\nm{\sqrt{w(y_i)} \left( \cP_h(f_{\cI}(y_i)) - \cP_h(f)(y_i) \right ) }_{\bbV} &\leq \sqrt{w(y_i)} \nm{f_{\cI}(y_i) - f(y_i)}_{\bbV}.
}
We deduce that
\bes{
\nm{A\cP_{h}( {c})-{v}}_{\ell^2([m];\bbV)} \leq \nm{f - f_{\cI} }_{\mathrm{disc}} + \nm{f - \cP_h(f) }_{\mathrm{disc}} + \nm{e}_{\ell^2([m];\bbV)}.
}
Substituting this into \R{step-b} and then combining with \R{step-a} now completes the proof.
}

We now prove Theorem \ref{t:BOS_RIP}. For this, we use the following result, which was shown in \cite[Thm.\ 1.1]{brugiapaglia2021sparse}:

\thm{
\label{t:brugiapaglia2021sparse}
There exist absolute constants $\kappa,c_0,c_1 > 0$ such that the following holds. Let $X_1,\ldots,X_m$ be independent copies of a random vector $X \in \bbC^n$ such that $\nm{X}_{\infty} \leq K$  almost surely for some $K > 0$. Let $\cT \subseteq \{ c \in \bbC^n : \nm{c}_{\ell^1} \leq \sqrt{s} \}$, $\delta \in (0,\kappa)$, $0 < \epsilon < 1$ and suppose that
\bes{
m \geq c_0 \cdot K^2 \cdot \delta^{-2} \cdot s \cdot \left ( \log(\E n) \log^2(s K^2/\delta) + \log(2/\epsilon) \right ).
}
Then, with probability at least $1-\epsilon$, we have
\bes{
\sup_{c \in \cT} \left | \frac1m \sum^{m}_{i=1} | \ip{X_i}{c} |^2 - \bbE | \ip{X}{c} |^2 \right | \leq c_1 \delta \left ( 1 + \sup_{c \in \cT} \bbE | \ip{c}{X} |^2 \right ).
}
}
\prf{
[Proof of Theorem \ref{t:BOS_RIP}]
Note that the result holds, provided
\bes{
\left | \frac1m \sum^{m}_{i=1} w(y_i) | p(y_i) |^2 - \nm{p}^2_{L^2_{\rho}(D)} \right | \leq \delta,\quad \forall p \in P_{T}, \nm{p}_{L^2_{\rho}(D)} \leq 1,\ T \subseteq \cI,\ |T| \leq t.
}
Define the random vector $X = (\sqrt{w(y)} \phi_{\iota_j}(y))_{j \in [N]}$, where $y \sim \mu$. Let $p \in P_{T}$ for some $T \subseteq \cI$ with $|T| \leq t$, and write $p = \sum_{\iota \in T} c_{\iota} \phi_{\iota}$. Then
\bes{
\frac1m \sum^{m}_{i=1} w(y_i) | p(y_i) |^2 = \frac1m \sum^{m}_{i=1} | \ip{X_i}{c} |^2,
}
and, since $\D \mu(y) = (w(y))^{-1} \D \rho(y)$,
\bes{
\nm{p}^2_{L^2_{\rho}(D)} = \int_{D} |p(y)|^2 \D \rho(y) = \int_{D} w(y) \left | \sum_{\iota \in T} c_{\iota} \phi_{\iota}(y) \right |^2 \D \mu(y) = \bbE | \ip{X}{c} |^2.
}
Hence, it suffices to show that 
\be{
\label{what-we-need-to-show}
\left | \frac1m \sum^{m}_{i=1} | \ip{X_i}{c} |^2 - \bbE | \ip{X}{c} |^2 \right | \leq \delta,\quad \forall c \in \cT, 
}
where
\bes{
\cT = \{ c \in \bbC^n : \mbox{$c$ is $t$-sparse and $\bbE | \ip{X}{c} |^2 \leq 1$} \}.
}
Notice that if $c \in \cT$ then, by \R{Riesz-bounds},
\bes{
\nm{c}_{\ell^1} \leq \sqrt{t} \nm{c}_{2} \leq \sqrt{t/a} \nm{\sum_{\iota} c_{\iota} \phi_{\iota} }^2_{L^2_{\rho}(D)} = \sqrt{t/a} \sqrt{\bbE | \ip{X}{c} |^2} \leq \sqrt{t/a}.
}
Hence $\cT \subseteq \{ c \in \bbC^n: \nm{c}_{\ell^1} \leq \sqrt{t/a} \}$. Therefore, we may apply Theorem \ref{t:brugiapaglia2021sparse} to get that
\bes{
\left | \frac1m \sum^{m}_{i=1} | \ip{X_i}{c} |^2 - \bbE | \ip{X}{c} |^2 \right | \leq c_1 \delta' \left ( 1 + \sup_{c \in \cT} \bbE | \ip{c}{X} |^2 \right ) \leq 2 c_1 \delta',
}
for $0 < \delta' < \kappa$. This holds with probability at least $1-\epsilon$, provided
\bes{
m \geq c_0 \cdot K^2 \cdot (\delta')^{-2} \cdot (t/a) \cdot \left ( \log(\E n) \log^2(t K^2 / (a \delta')) + \log(2/\epsilon) \right ).
}
Observe that $\nm{X}_{\infty} = \max_{\iota \in \cI} \{ \sqrt{w(y)} | \phi_{\iota}(y) | \}$, and therefore we may take $K = \Gamma = \nmu{ \max_{\iota \in \cI} \{ \sqrt{w(\cdot)} | \phi_{\iota}(\cdot) | \} }_{L^{\infty}_{\rho}(D)}$. We now set $\delta' = \delta / (2 c_1)$ to deduce that \R{what-we-need-to-show} holds for $0 < \delta < \kappa / (2 c_1)$ with probability at least $1-\epsilon$, provided
\bes{
m \geq 4 c_0 c^2_1 \cdot \Gamma^2 \cdot \delta^{-2} \cdot (t/a) \cdot \left ( \log(\E n) \log^2(2 c_1 t \Gamma^2 / (a \delta)) + \log(2 / \epsilon) \right ).
}
To complete the proof, we simply notice that $\log(2 c_1 t \Gamma^2/(a\delta)) \leq \log(\E t \Gamma^2/(a\delta)) + \log(2 c_1 / \E) \lesssim \log(\E t \Gamma^2/(a\delta))$, since $t \Gamma^2 / (a\delta) > 1$.
}

\section{A novel approach for sparse polynomial approximation on irregular domains}\label{s:sparse-irregular-new}

In this section, we focus on Example \ref{ex:alg-poly-general} in the context of $\ell^1$-minimization. As we have seen in the previous section, the various measurement conditions depend on the Riesz basis constants $a,b$ of the system $\Phi = \{ \phi_{\iota} :\iota \in \cI \}$. Unfortunately, in cases such as Example \ref{ex:alg-poly-general}, these constants are often very poorly behaved, especially when the polynomial degree is large. Later, in Figs.\ \ref{fig:irregular-exp-1}--\ref{fig:irregular-exp-3} we compute $a$, $b$ for various different domains. While we note that these constants may be pessimistic (in particular, recall Remark \ref{r:Riesz-constants}), the fact remains that dealing with poorly-conditioned dictionaries may well be problematic in the setting of $\ell^1$-minimization for sparse approximation.
 
\subsection{Method}

Inspired by our earlier use of discrete measures to effect optimal sampling, in this section we consider a different approach for approximation on irregular domains in which the original Riesz basis $\Phi$ is orthogonalized over the discrete grid. Specifically, let $Z = \{ z_i\}^{k}_{i=1}$ be a finite grid and $\tau$ be as in \R{tau-def-CS}. Then we construct a new basis $\Upsilon = \{ \upsilon_{\iota} : \iota \in \cI \}$ that is orthonormal with respect to $\tau$, and subsequently use this basis in which to construct the approximation $\hat{f}$ to $f$ via $\ell^1$-minimization.

As in \S \ref{ss:practical-optimal-samp}, the orthonormal basis $\Upsilon$ is constructed via QR factorization. Let 
\bes{
B = \left ( \phi_{\iota_j}(z_i) / \sqrt{k} \right )_{i \in [k], j \in [n]} \in \bbC^{k \times n},
}
and suppose that $B$ has QR factorization $B = Q R$, where $Q \in \bbC^{k \in n}$ and $R \in \bbC^{n \times n}$. Then, this orthonormal basis is given by
\be{
\label{ONB-irregular-domain-CS}
\upsilon_{\iota_i}(y) = \sum^{i}_{j=1} (R^{-{\top}})_{ji} \phi_{\iota_j}(y),\quad i \in [n].
}
In what follows, we compare two sampling strategies for the basis $\Upsilon$. First, Monte Carlo sampling from the underlying measure $\tau$, i.e.
\be{
\label{MC-CS-ONB}
\bbP ( y = z_i ) = \frac1k,\quad i \in [k].
}
Second, the `optimal' sampling measure identified in \S \ref{ss:CS-opt-disc}. In this case, due to the orthogonalization, this is given by 
\be{
\label{Optimal-CS-ONB}
\bbP(y = z_i) = \frac{\max_{j \in [n]} | Q_{ij} |^2 }{\sum^{k}_{i=1} \max_{j \in [n]} | Q_{ij} |^2},\quad i \in [k].
}
As shown in the previous section, the sample complexity bounds for these two strategies depend on their respective constants $\Theta$ and $\theta$. Because of the previous definition of the basis, these are given by
\bes{
\Theta = \Theta(\Upsilon) = \sqrt{k} \max_{\substack{i \in [k] \\ j \in [n]}} |Q_{ij}|,
}
for the former, and
\bes{
\theta = \theta(\Upsilon) = \sqrt{\sum_{i \in [k]} \max_{j \in [n]} |Q_{ij} |^2 },
}
for the latter.

\subsection{Orderings}

\begin{figure}
	\begin{center}
	\begin{small}
    \begin{tabular}{ccc}
    \includegraphics[width=0.28\textwidth]{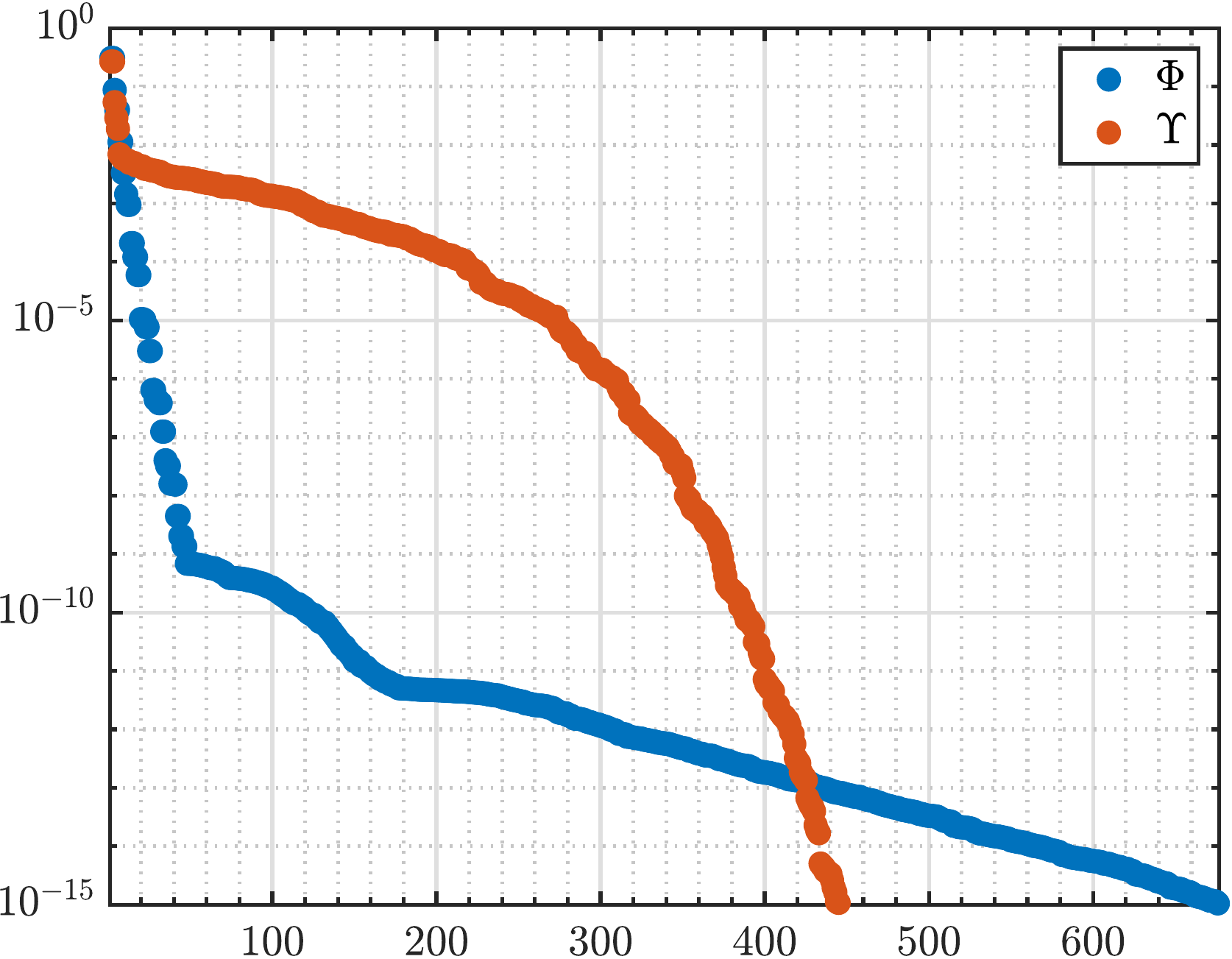} &
    \includegraphics[width=0.28\textwidth]{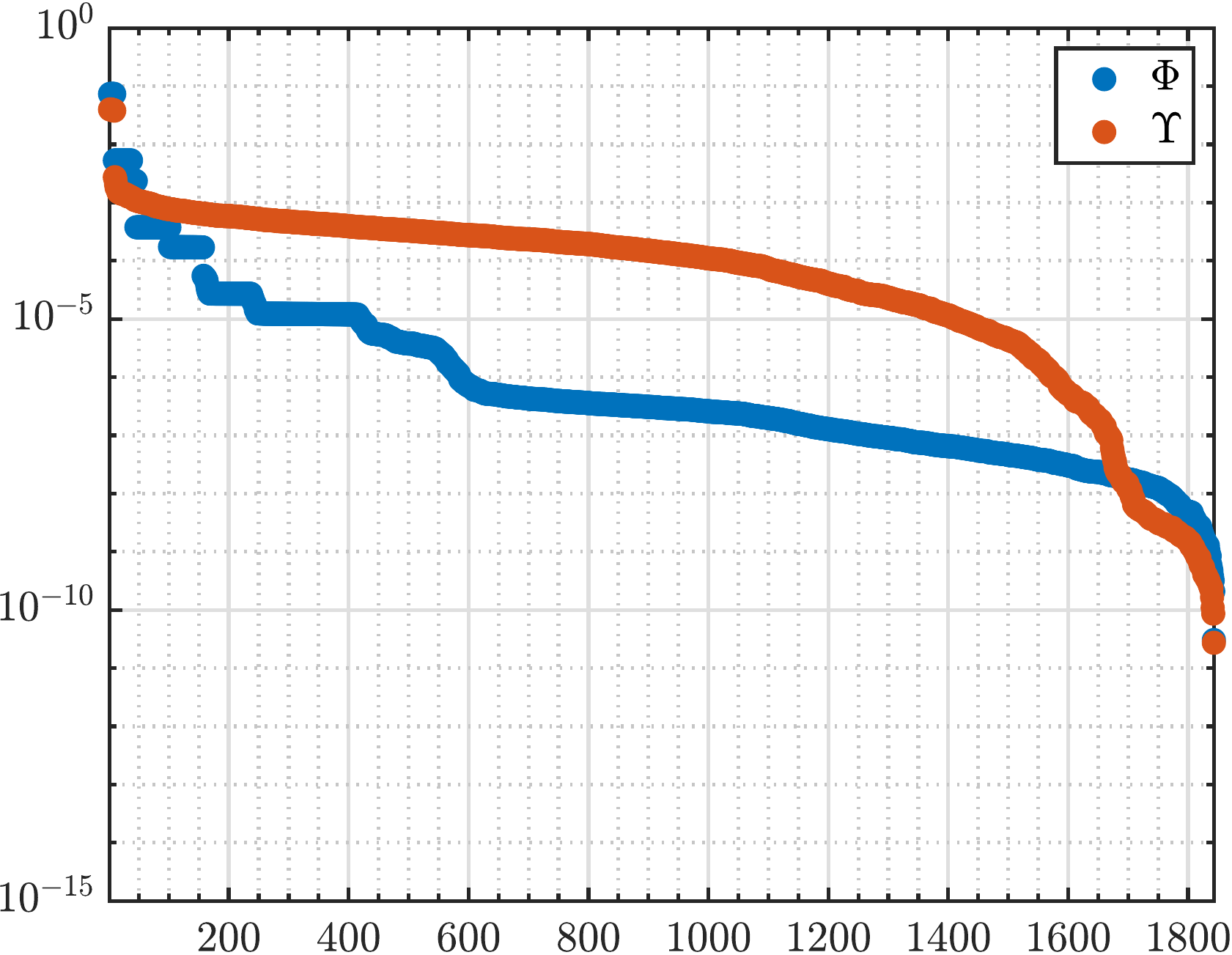} & 
    \includegraphics[width=0.28\textwidth]{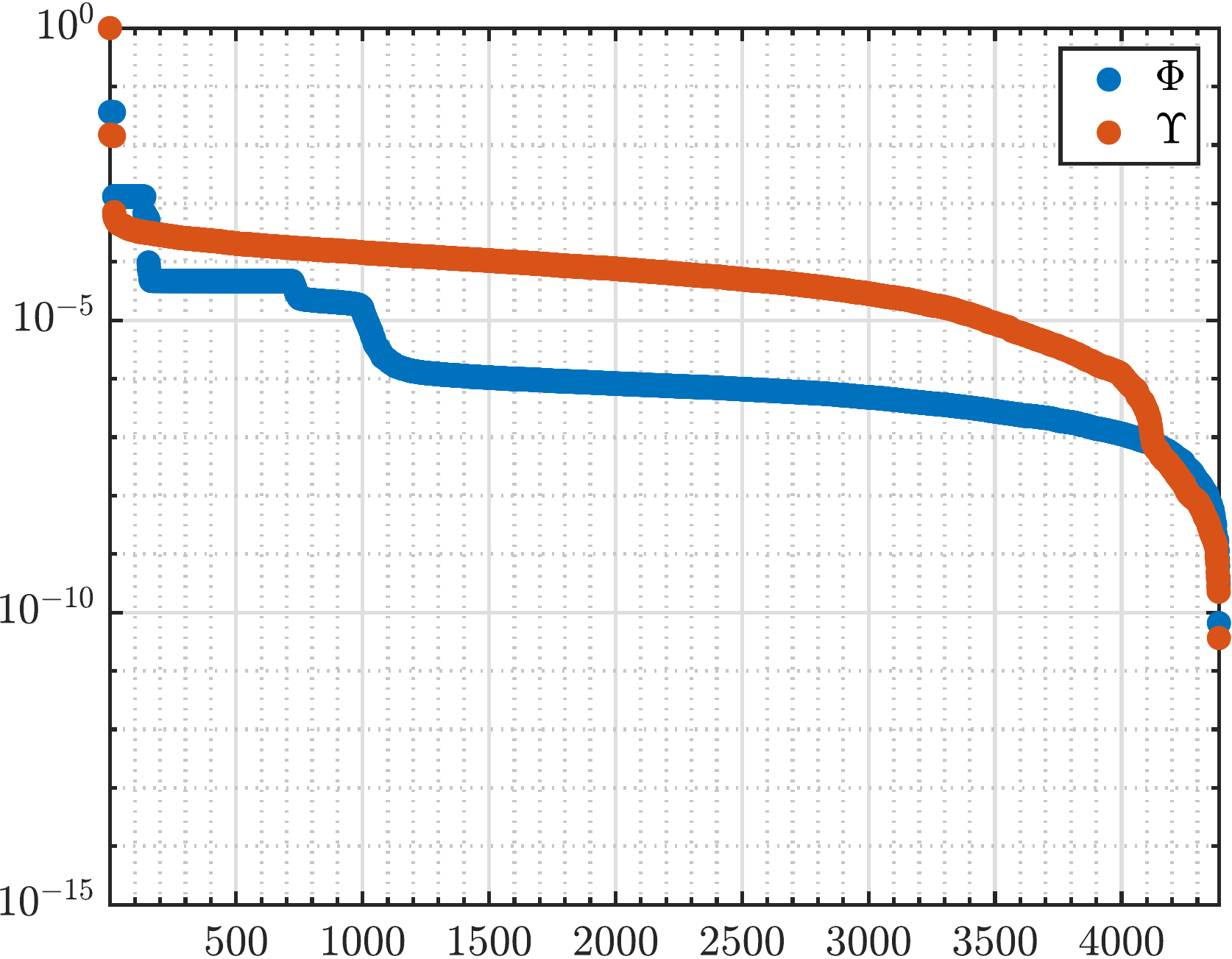} 
    \\
    \includegraphics[width=0.28\textwidth]{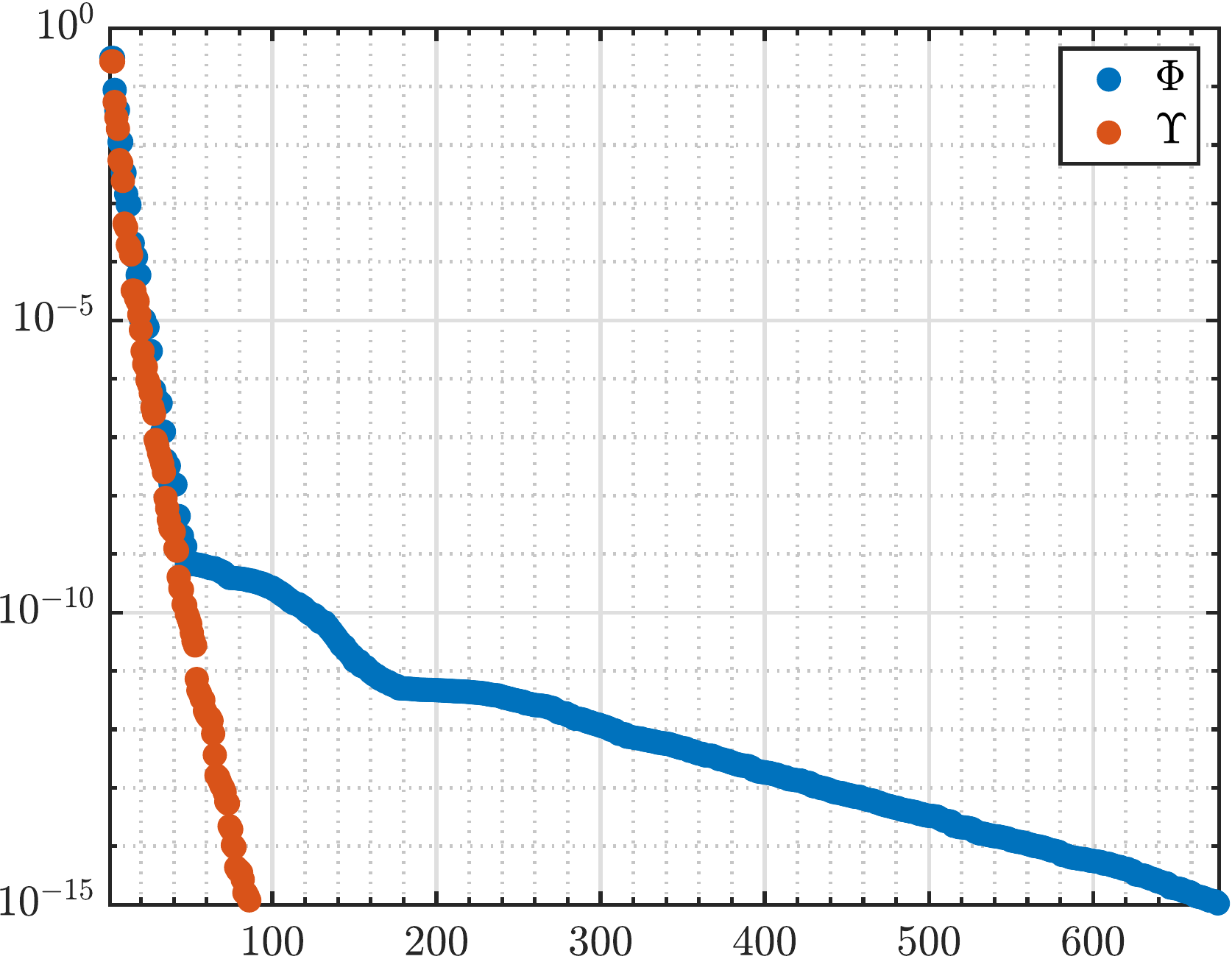} &
    \includegraphics[width=0.28\textwidth]{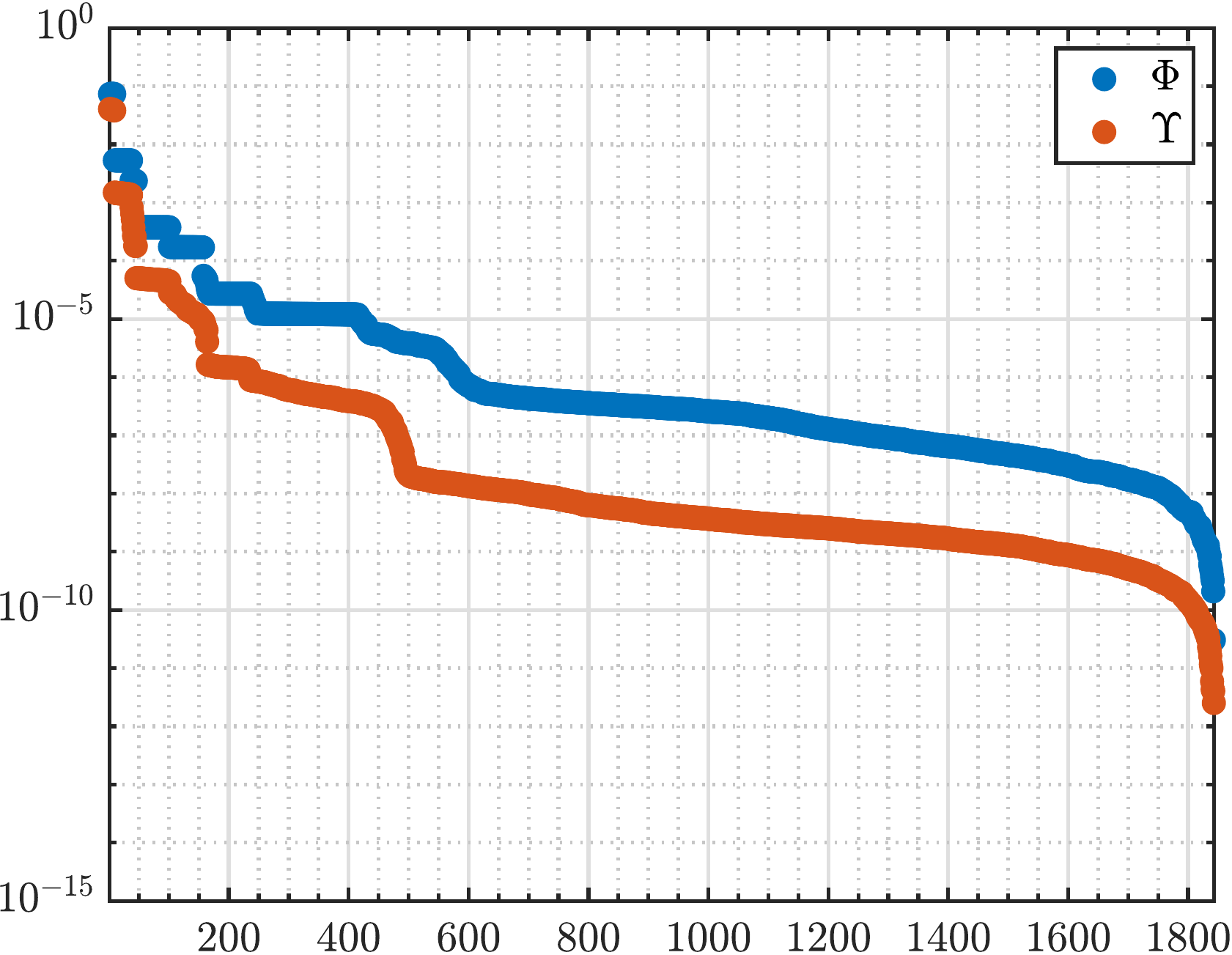} &
    \includegraphics[width=0.28\textwidth]{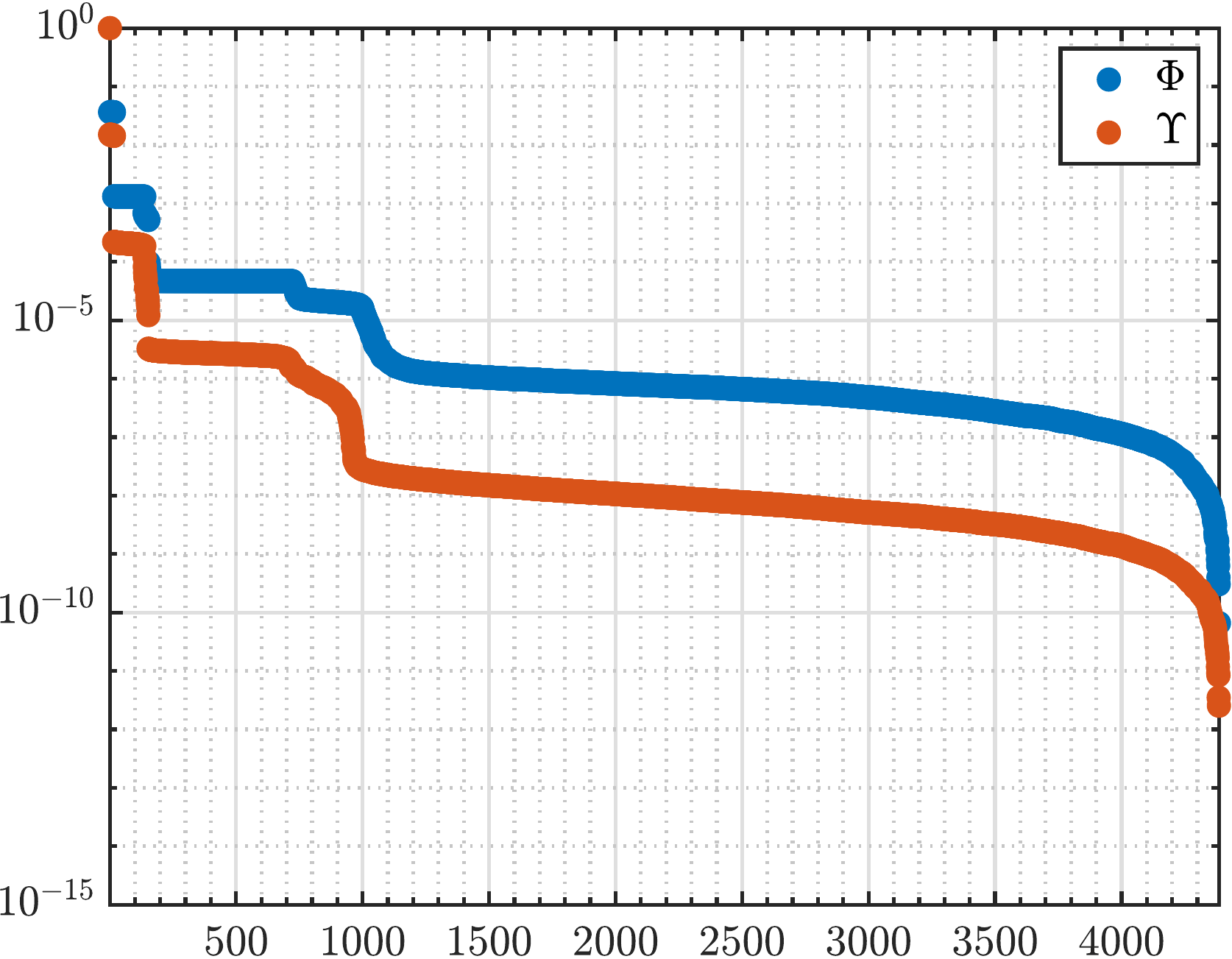} 
    \\
    $(d,t,n)=(2,152,796)$ & $(d,t,n) = (8,22,1843)$ & $(d,t,n) = (16,14,4385)$\\
	\end{tabular}
	\end{small}
	\end{center}
\caption{The absolute values of the coefficients of the function $f = f_1$ over the domain $D = D_2$ with respect to the bases $\Phi$ and $\Upsilon$, where $\cI = \cI^{\mathrm{HC}}_{t-1}$ is the hyperbolic cross index set. The coefficients are sorted from largest in absolute value to smallest.
In the top row, the multi-indices in $\cI$ are sorted lexicographically. In the bottom row, they are sorted according to increasing total degree (i.e.\ the value $\iota_1 + \ldots + \iota_d$ for $\iota = (\iota_k)^{d}_{k=1}$).} 
\label{fig:irregular-coefficients-ex1}
\end{figure}

\begin{figure}[t]
	\begin{center}
	\begin{small}
 \begin{tabular}{ccc}
    \includegraphics[width=0.3\textwidth]{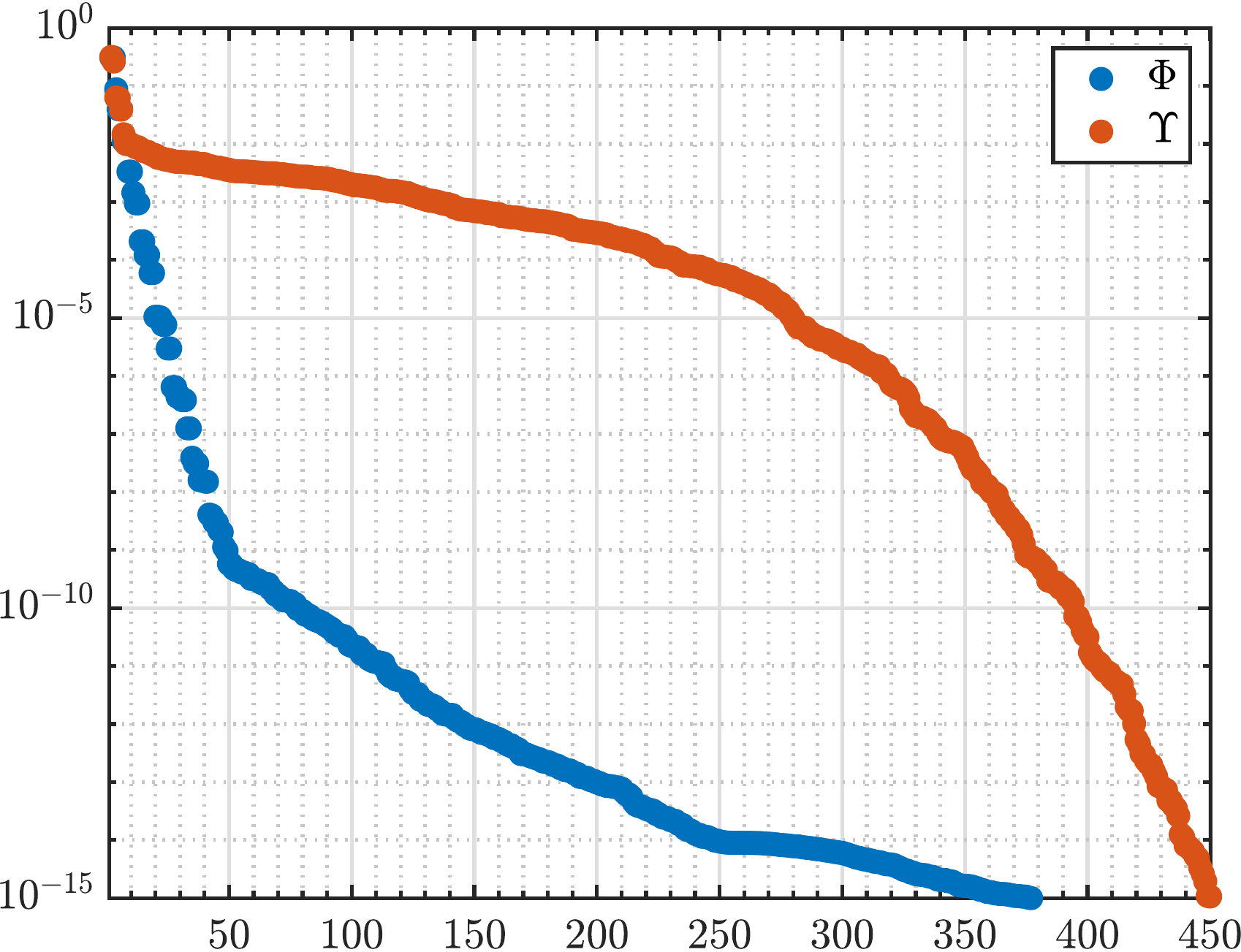} &
    \includegraphics[width=0.3\textwidth]{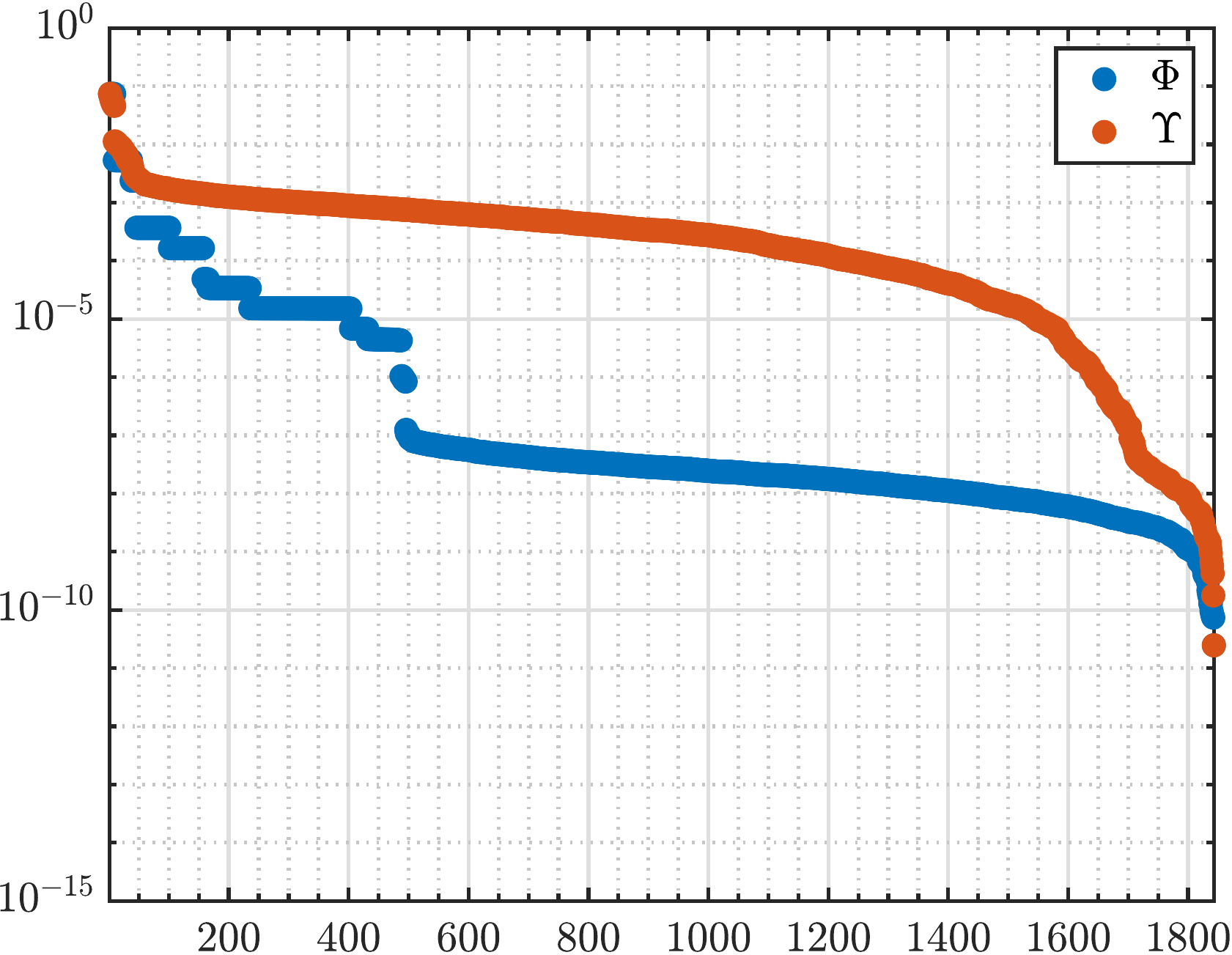} &
    \includegraphics[width=0.3\textwidth]{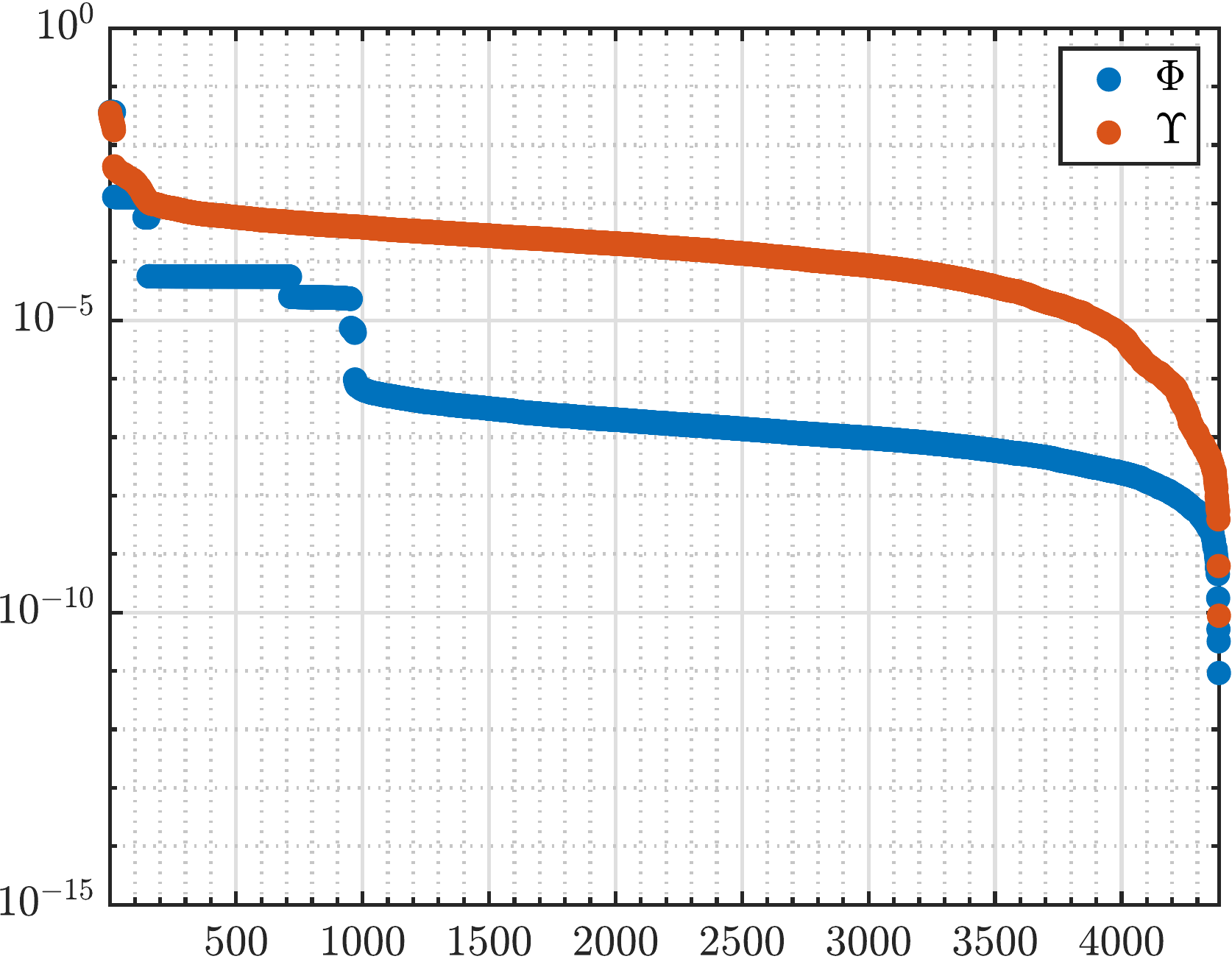} \\
    \includegraphics[width=0.3\textwidth]{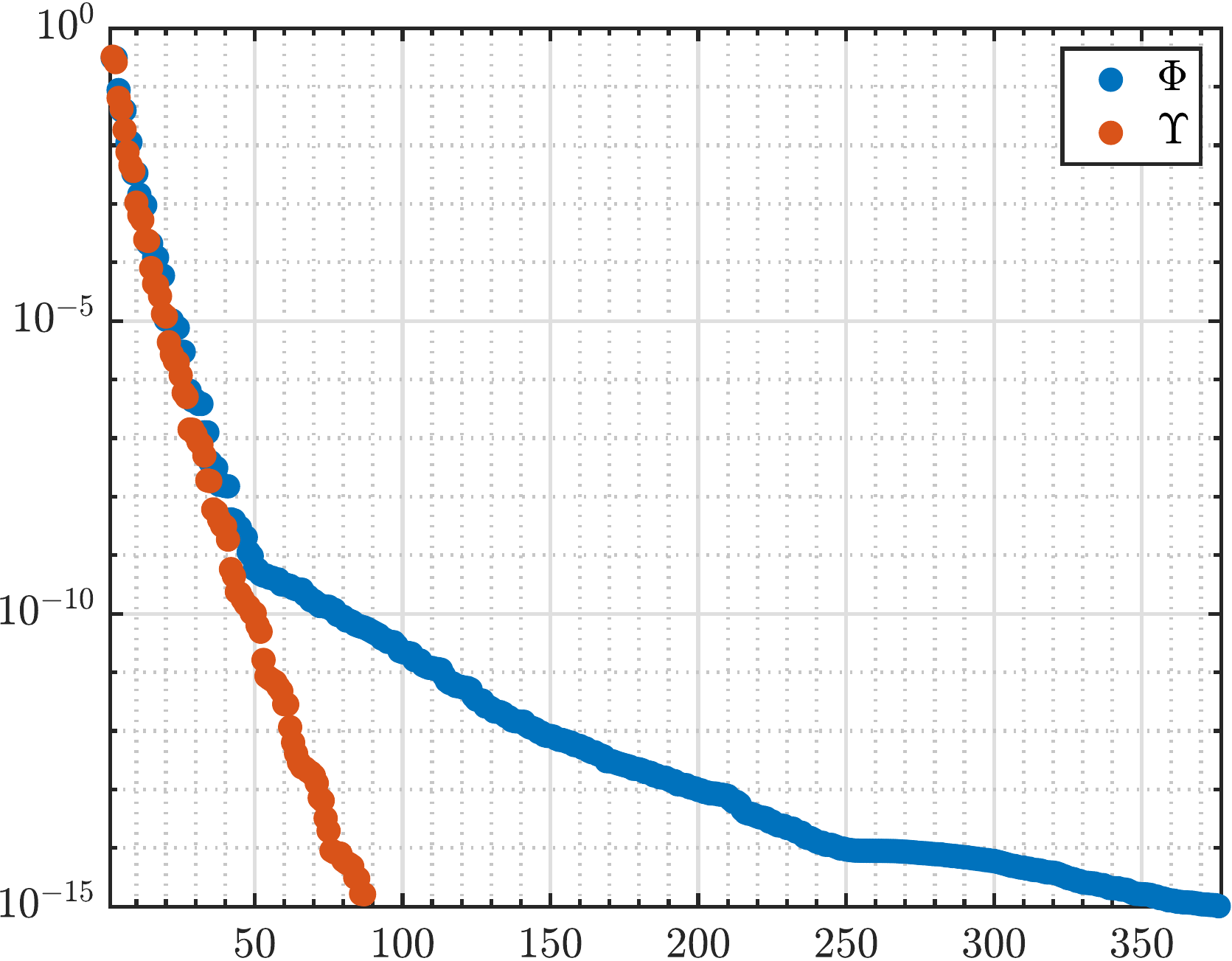} &
    \includegraphics[width=0.3\textwidth]{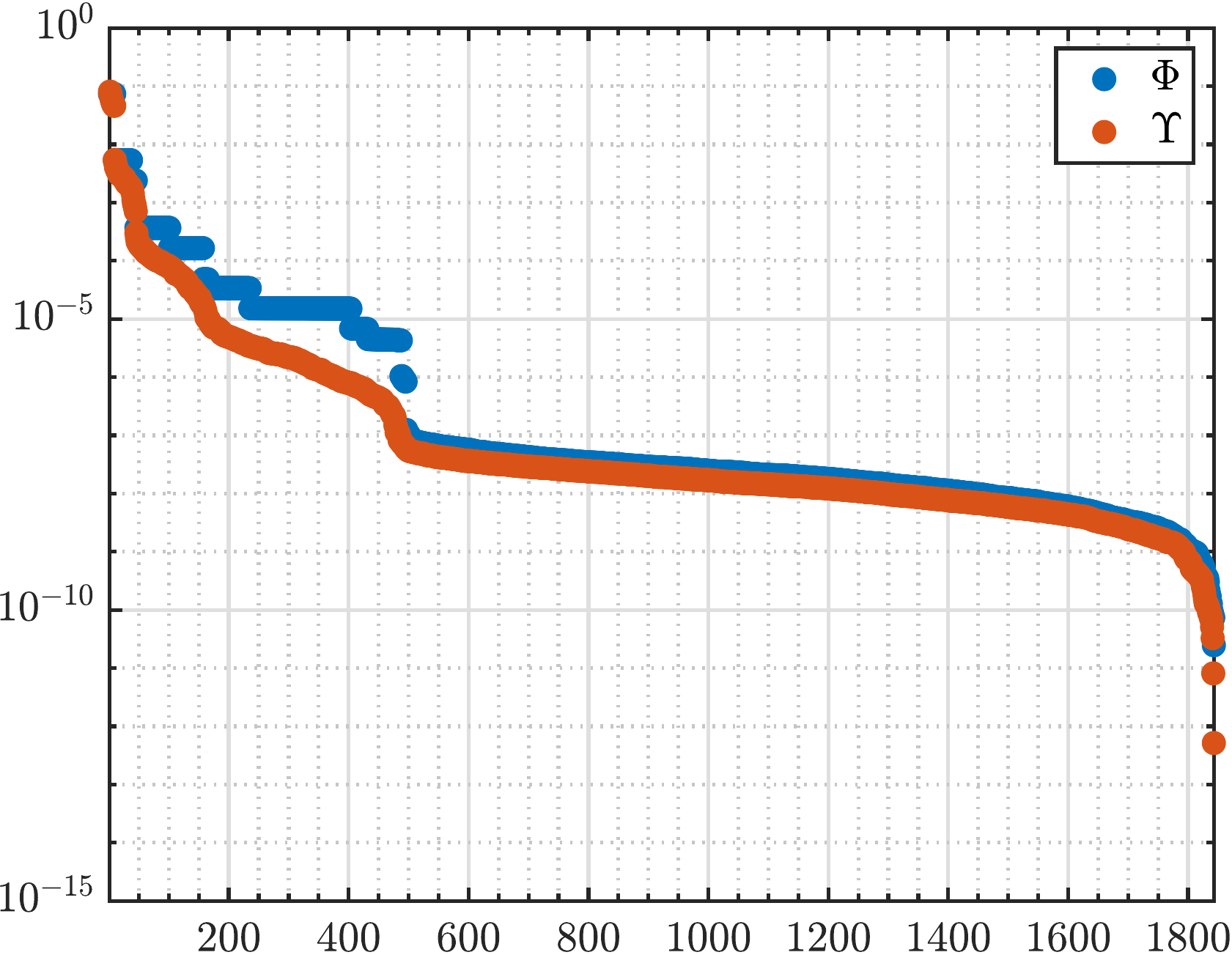} &
    \includegraphics[width=0.3\textwidth]{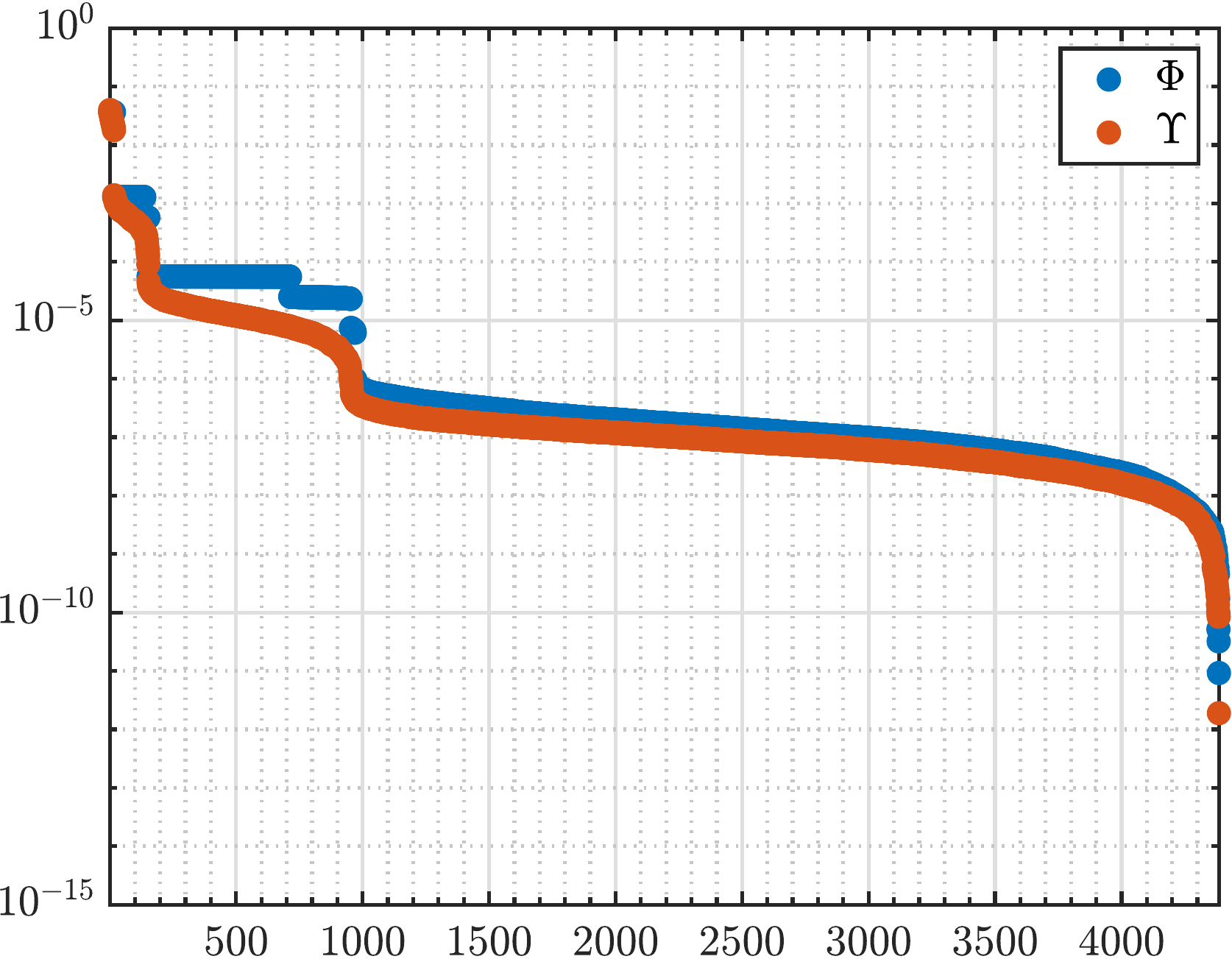} \\    
    $(d,t,n)=(2,152,796)$ & $(d,t,n) = (8,22,1843)$ & $(d,t,n) = (16,14,4385)$\\
	\end{tabular}
	\end{small}
	\end{center}
\caption{The same as in Fig.\ \ref{fig:irregular-coefficients-ex1} except for $f = f_1$ and $D = D_3$.} 
\label{fig:irregular-coefficients-ex2}
\end{figure}

\begin{figure}[t]
	\begin{center}
	\begin{small}
 \begin{tabular}{ccc}
    \includegraphics[width=0.3\textwidth]{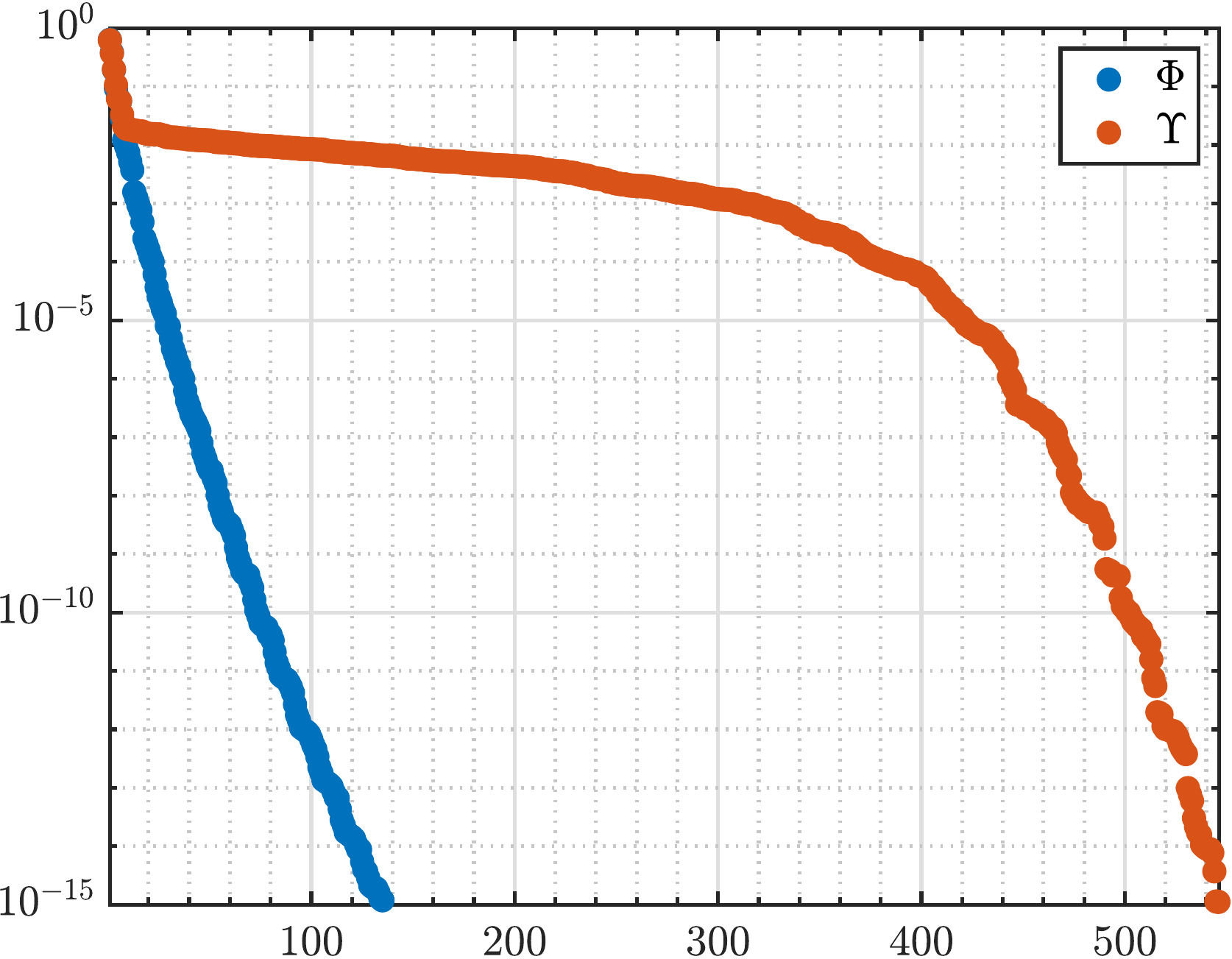} &		
    \includegraphics[width=0.3\textwidth]{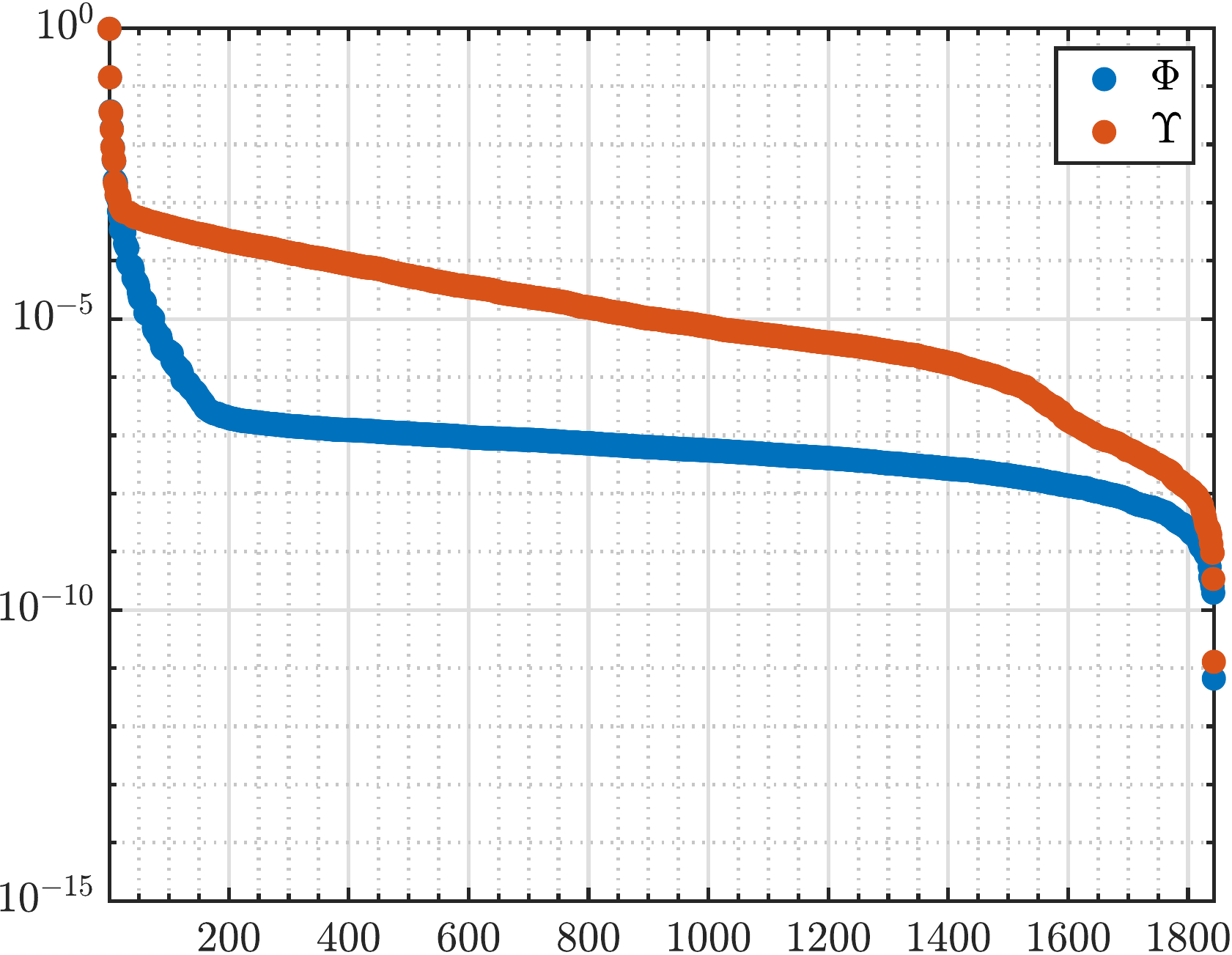} & 
    \includegraphics[width=0.3\textwidth]{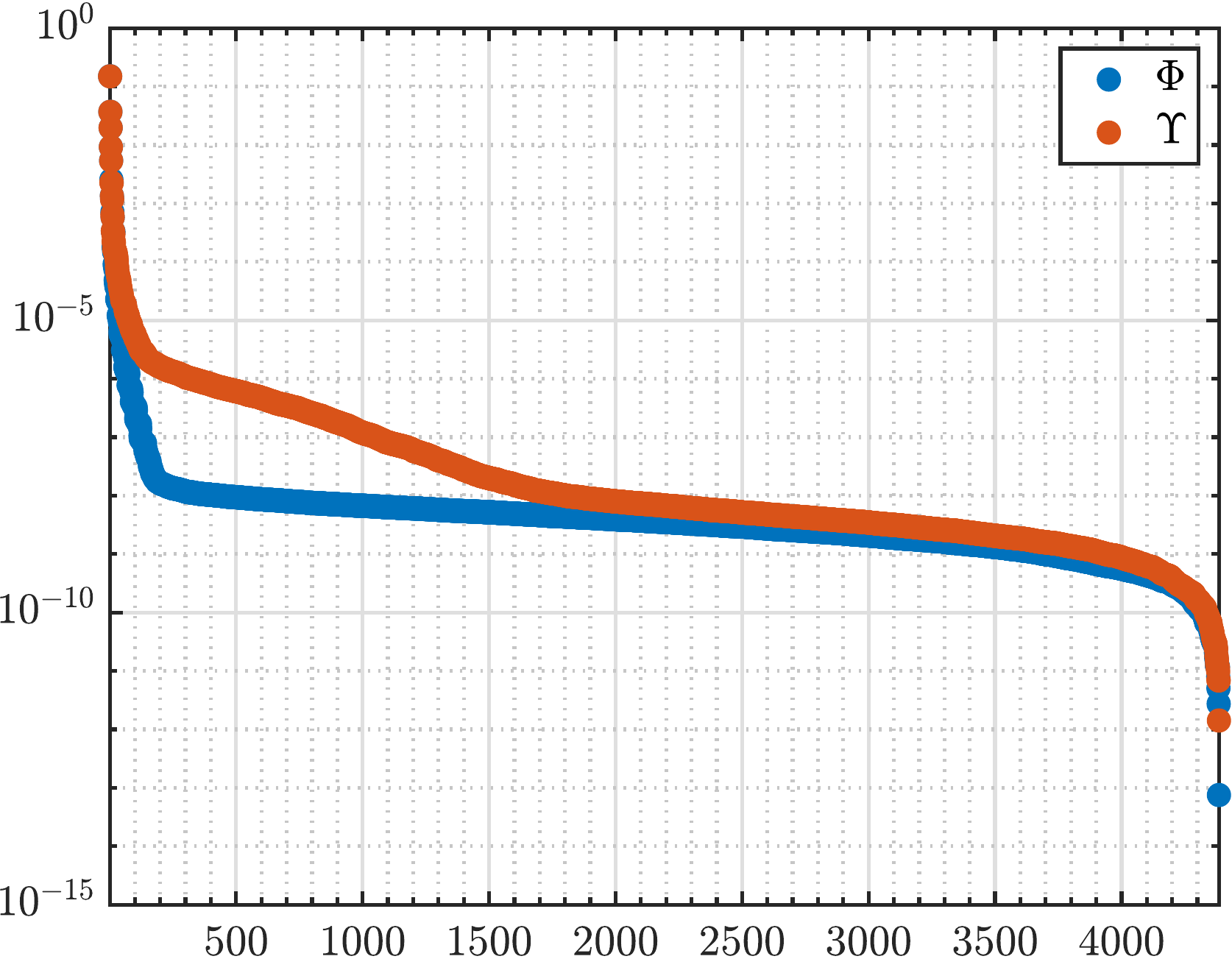} \\
    \includegraphics[width=0.3\textwidth]{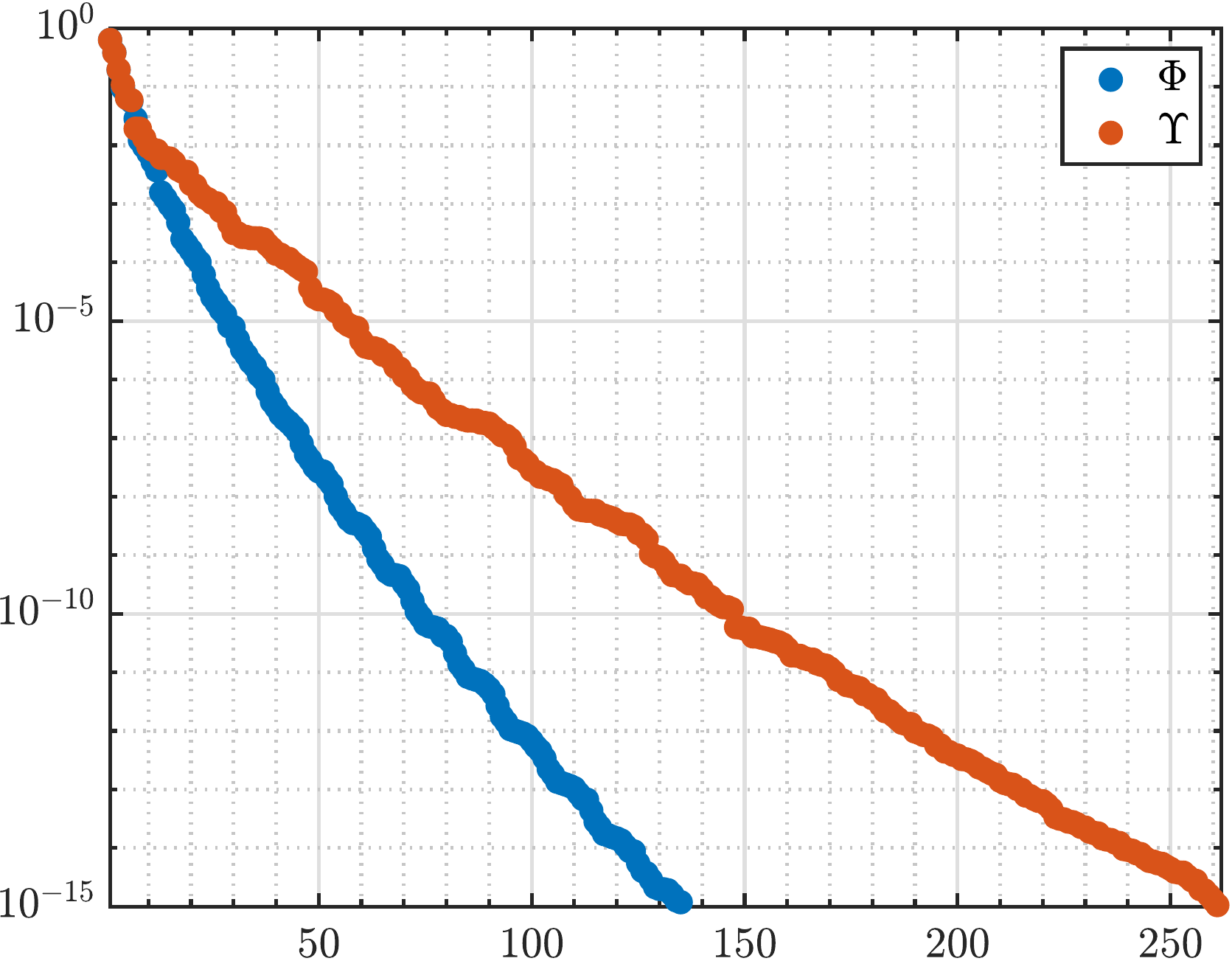} &		
    \includegraphics[width=0.3\textwidth]{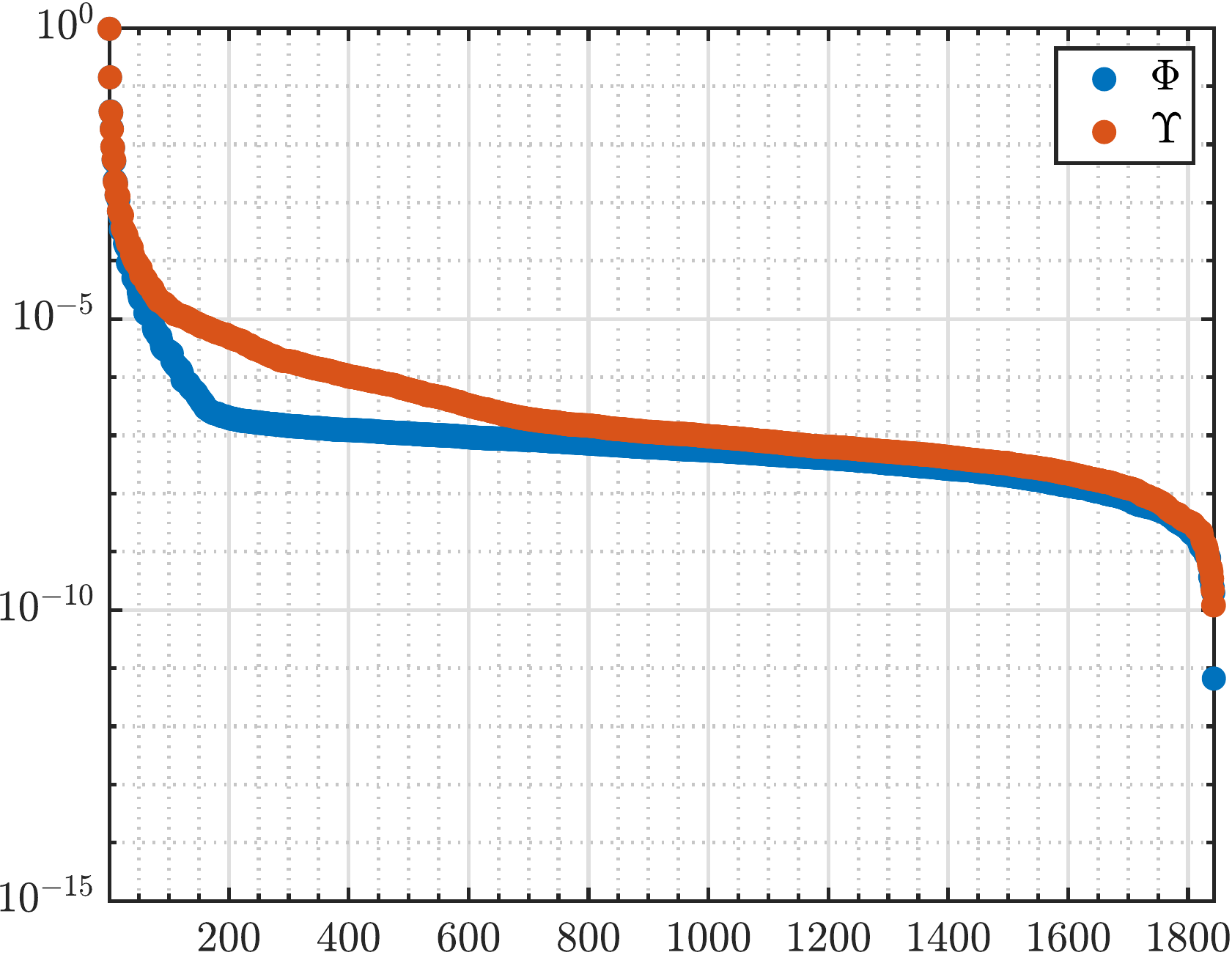} & 
    \includegraphics[width=0.3\textwidth]{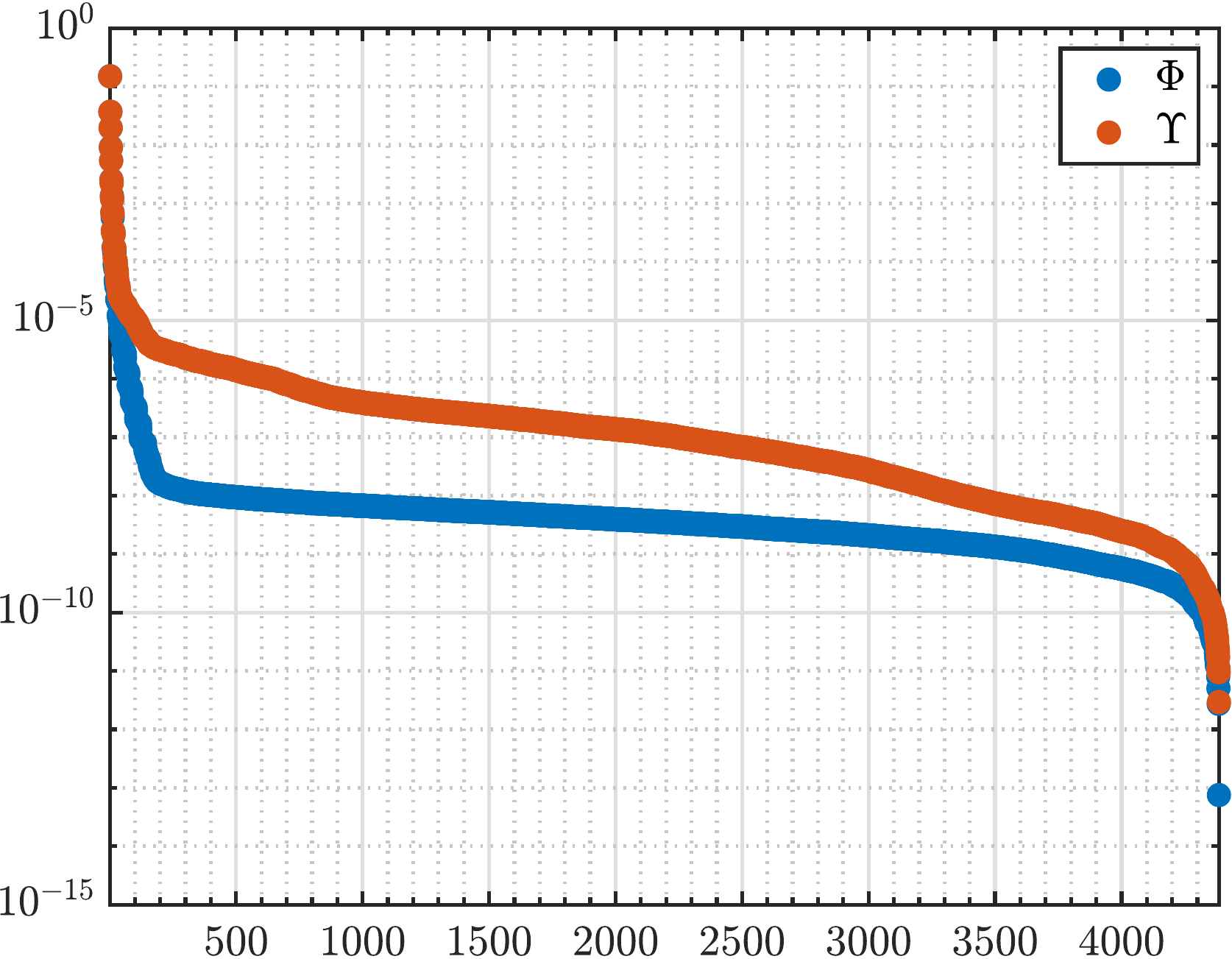} \\   
    $(d,t,n)=(2,152,796)$ & $(d,t,n) = (8,22,1843)$ & $(d,t,n) = (16,14,4385)$\\
	\end{tabular}
	\end{small}
	\end{center}
\caption{The same as in Fig.\ \ref{fig:irregular-coefficients-ex1} except for $f = f_2$ and $D = D_1$.
} 
\label{fig:irregular-coefficients-ex3}
\end{figure}

Suppose that $f$ has a sparse representation in the system $\Phi$. Then it is not guaranteed to possess a sparse representation in the basis $\Upsilon$, since the $j$th basis function $\upsilon_{\iota_j}$ is a linear combination of the functions $\phi_{\iota_1},\ldots,\phi_{\iota_j}$. In particular, the sparsity of $f$ in the representation $\Upsilon$ will be heavily influenced by the ordering of the basis $\{ \iota_1,\ldots,\iota_n \}$ of the indices in $\cI$.

To illustrate this, in Figs.\ \ref{fig:irregular-coefficients-ex1}--\ref{fig:irregular-coefficients-ex3} we compare ordering the multi-indices in $\cI$ lexicographically to ordering them according to increasing total degree, i.e.\ the value $\iota_1 + \ldots + \iota_d$ for $\iota = (\iota_k)^{d}_{k=1}$. We remark in passing that ordering according maximum, i.e.\ the value $\max_{k=1,\ldots,d} \{ \iota_k \}$ produces similar results.
To examine the sparsity in either basis, we plot the coefficients of a function $f$ sorted from largest to smallest in absolute value. In particular, the more rapidly the coefficients decrease, the better $f$ is approximated by a sparse representation in the given basis.

As is evident, lexicographic ordering always leads to a deterioration in sparsity when switching from the basis $\Phi$ to the basis $\Upsilon$. This is of little surprise. On the other hand, using the total degree ordering can substantially improve the situation. In Fig.\ \ref{fig:irregular-coefficients-ex1} it actually leads to better sparsity and in Fig.\ \ref{fig:irregular-coefficients-ex2} it yields better sparsity in $d = 2$ dimensions, and similar sparsity in $d = 8,16$ dimensions. Finally, in Fig.\ \ref{fig:irregular-coefficients-ex3}, while still leading to worse sparsity than in the original basis $\Phi$, it is still generally better than lexicographic ordering.

For this reason, in our subsequent experiments, we employ the total degree ordering. Naturally, this discussion leads to the question of the optimal ordering. We anticipate this to be function dependent, and it is outside the scope of this work to discuss it further. In practice, we expect a good ordering could be estimated from a set of candidate orderings via cross validation.

\subsection{Numerical examples}

\begin{figure}[t]
	\begin{center}
	\begin{small}
 \begin{tabular}{@{\hspace{0pt}}c@{\hspace{-0.5pc}}c@{\hspace{0pt}}}
		
 		\includegraphics[width=0.5\textwidth]{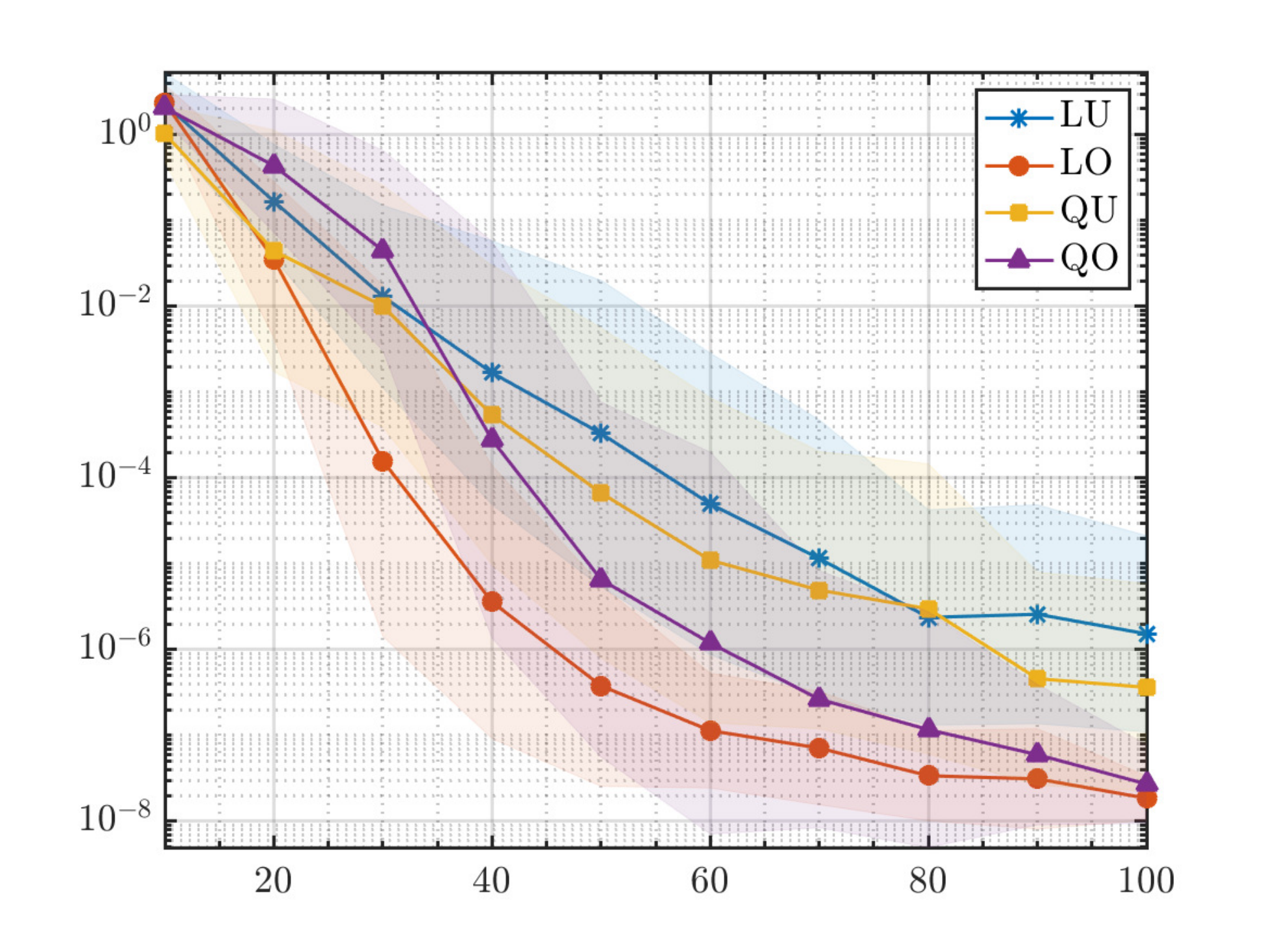} & 
 		\includegraphics[width=0.5\textwidth]{fig_11_1_1-eps-converted-to.pdf} \\
 		$(d,t,n) = (1,399,400)$ & $(d,t,n)=(2,152,796)$ \\
 		$(\theta_L^2,\Theta_L^2) = (2.25,303.73)$ &
 		$(\theta_L^2,\Theta_L^2) = (2.51,30.26)$ \\
 		$(\theta_Q^2,\Theta_Q^2) = (5.19,768.17)$ &  
 		$(\theta_Q^2,\Theta_Q^2) = (9.39,508.06)$\\ 
 		$(a_L,b_L) = (5.29\times 10^{-37},2.55)$ & 
 		$(a_L,b_L) = (3.69\times 10^{-15},1.86)$ \\
 		$(a_Q,b_Q) = (1,1)$ &
 		$(a_Q,b_Q) = (1,1)$ \\
 		\includegraphics[width=0.5\textwidth]{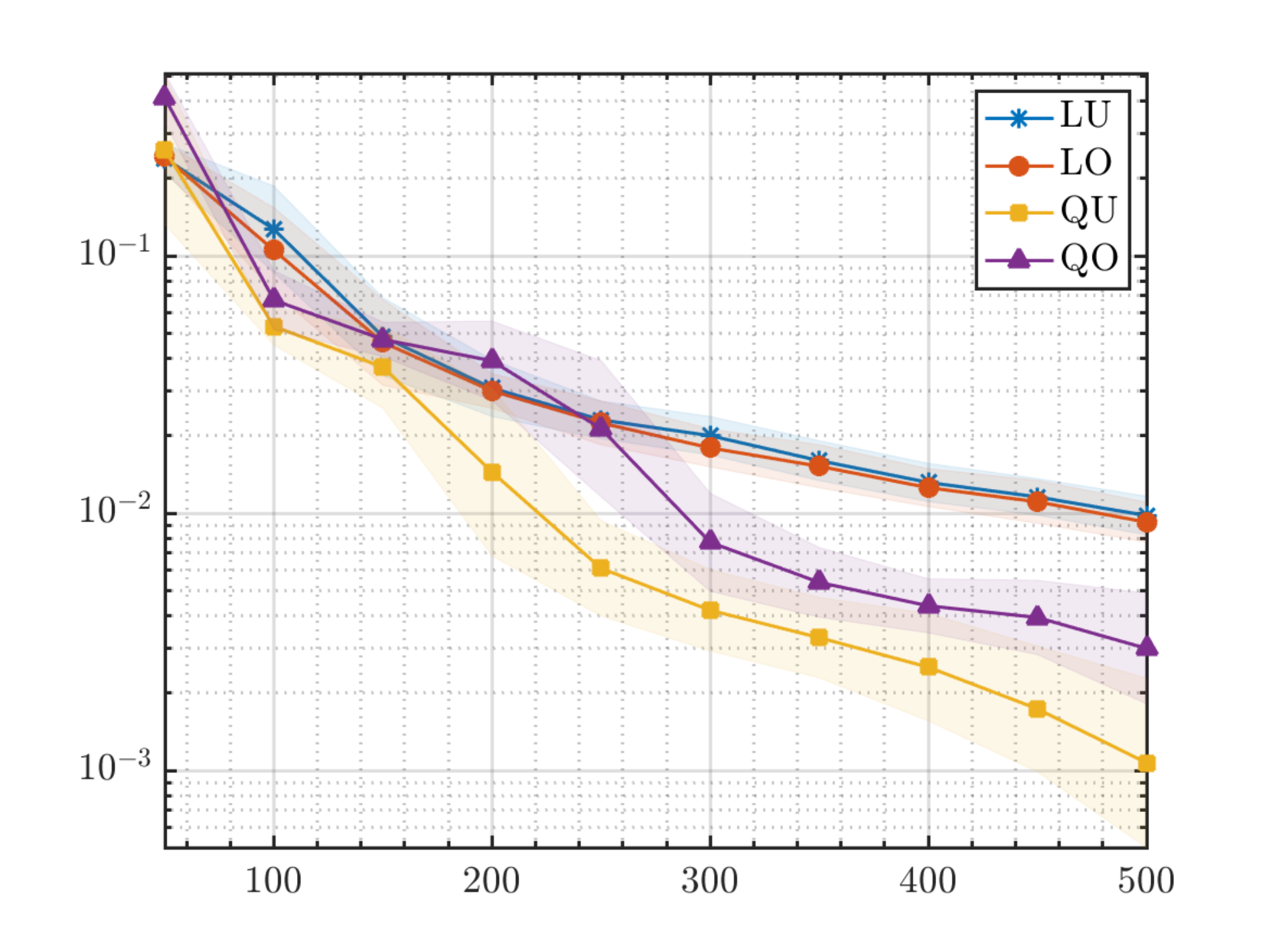} & 
 		\includegraphics[width=0.5\textwidth]{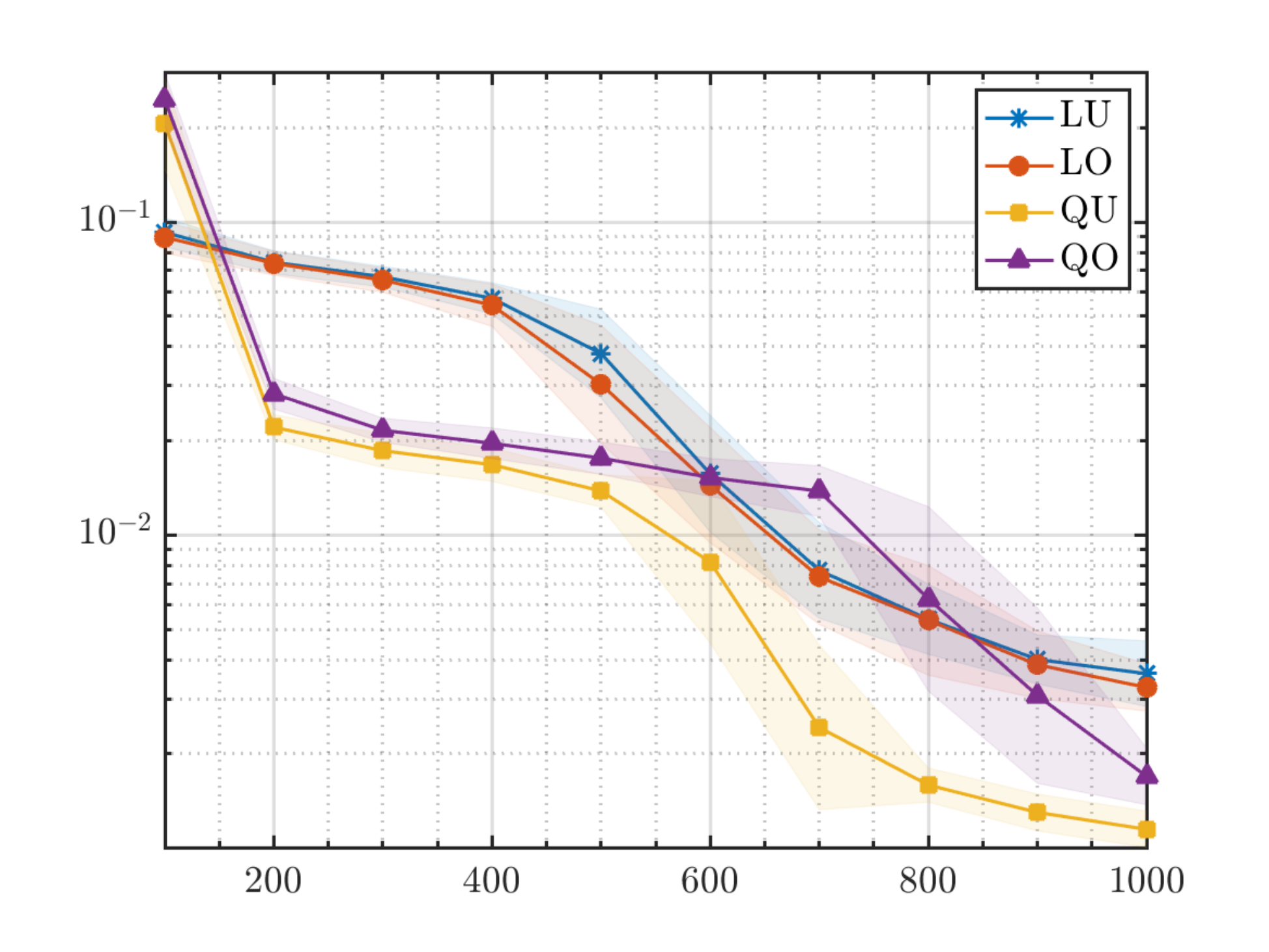} \\
 		$(d,t,n) = (8,22,1843)$ & $(d,t,n)=(16,14,4385)$ \\
 		$(\theta_L^2,\Theta_L^2) = (1.99,5.77)$ &  
 		$(\theta_L^2,\Theta_L^2) = (1.74,2.94)$ \\ 		
 		$(\theta_Q^2,\Theta_Q^2) = (66.99,3659)$ &
 		$(\theta_Q^2,\Theta_Q^2) = (132,11406)$  \\
 		$(a_L,b_L) = (9.31\times 10^{-9},20.86)$ & 
 		$(a_L,b_L) = (2.88\times 10^{-9},188.55)$ \\
 		$(a_Q,b_Q) = (1,1)$ &
 		$(a_Q,b_Q) = (1,1)$ \\ 		 		
	\end{tabular}
	\end{small}
\end{center}
\caption{The relative error \R{relative-error} versus $m$ for $\ell^1$-minimization in the case of Example \ref{ex:alg-poly-general}, where $\cI = \cI^{\mathrm{HC}}_{t-1}$ is the hyperbolic cross index set \R{HC-index} and $f = f_1$ and $D = D_2$ are as in \R{functions-define} and \R{domains-define}, respectively. This figure compares the Legendre basis on $[-1,1]^d$ restricted to $d$ and the orthonormal basis on $D$ constructed via \R{ONB-irregular-domain-CS}. In the former case, the sampling strategies are Monte Carlo sampling from the continuous uniform measure on $D$ (labelled `LU') and sampling from the discrete `optimal' measure \R{mu-optimal-CS} (`LO'). In the latter cases, the sampling strategies are Monte Carlo sampling \R{MC-CS-ONB} from the discrete uniform measure (`QU') and sampling from the discrete `optimal' measure \R{Optimal-CS-ONB} (`QO').  We also report the values of the corresponding constants $\theta^2$ and $\Theta^2$ for both bases (labelled `L' and `Q', respectively), as well as the Riesz basis constants $a,b$ with respect to the discrete measure $\tau$.} 
\label{fig:irregular-exp-1}
\end{figure}

\begin{figure}[t]
	\begin{center}
	\begin{small}
 \begin{tabular}{@{\hspace{0pt}}c@{\hspace{-0.5pc}}c@{\hspace{0pt}}}
  		\includegraphics[width=0.5\textwidth]{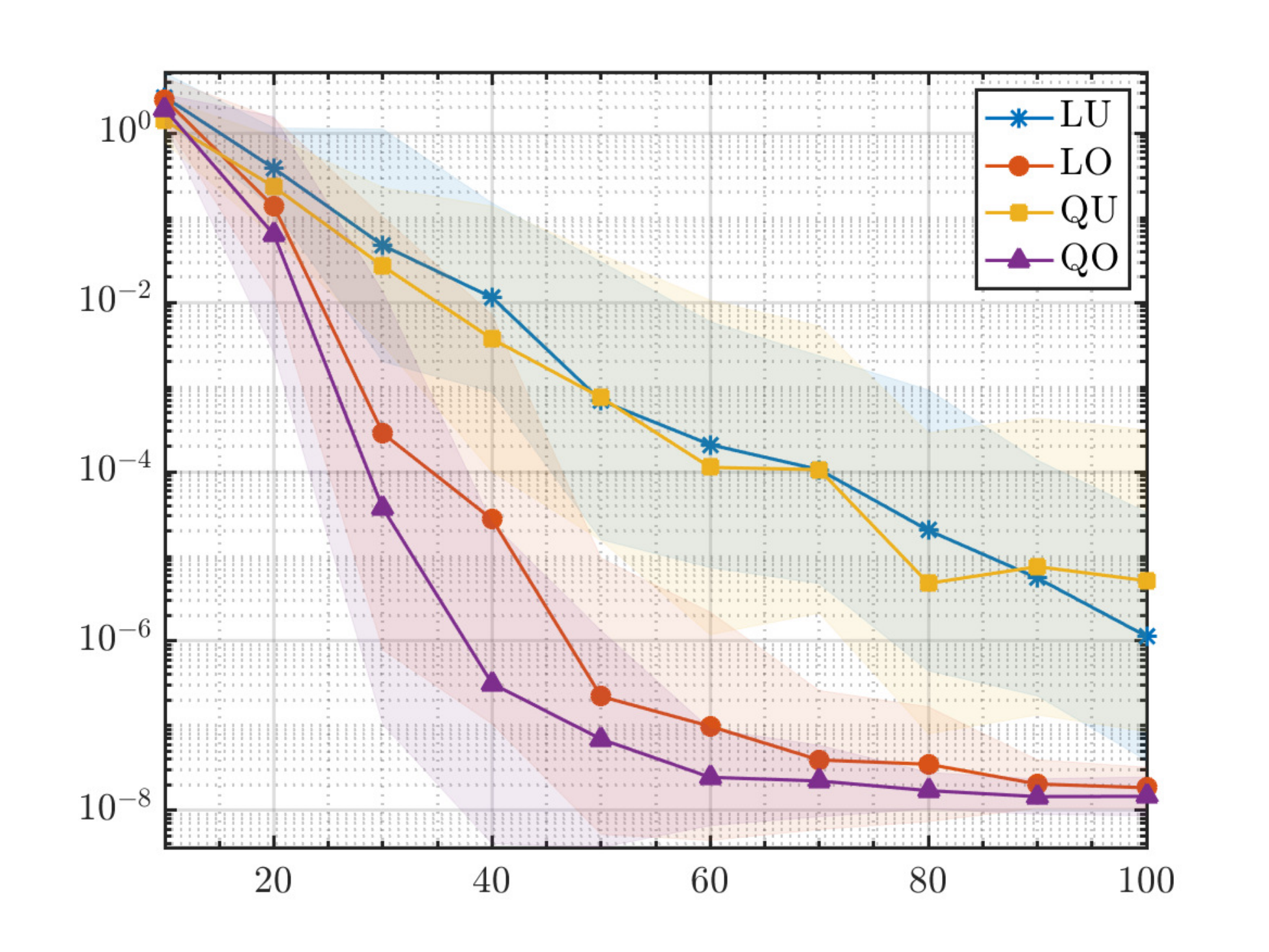} & 
 		\includegraphics[width=0.5\textwidth]{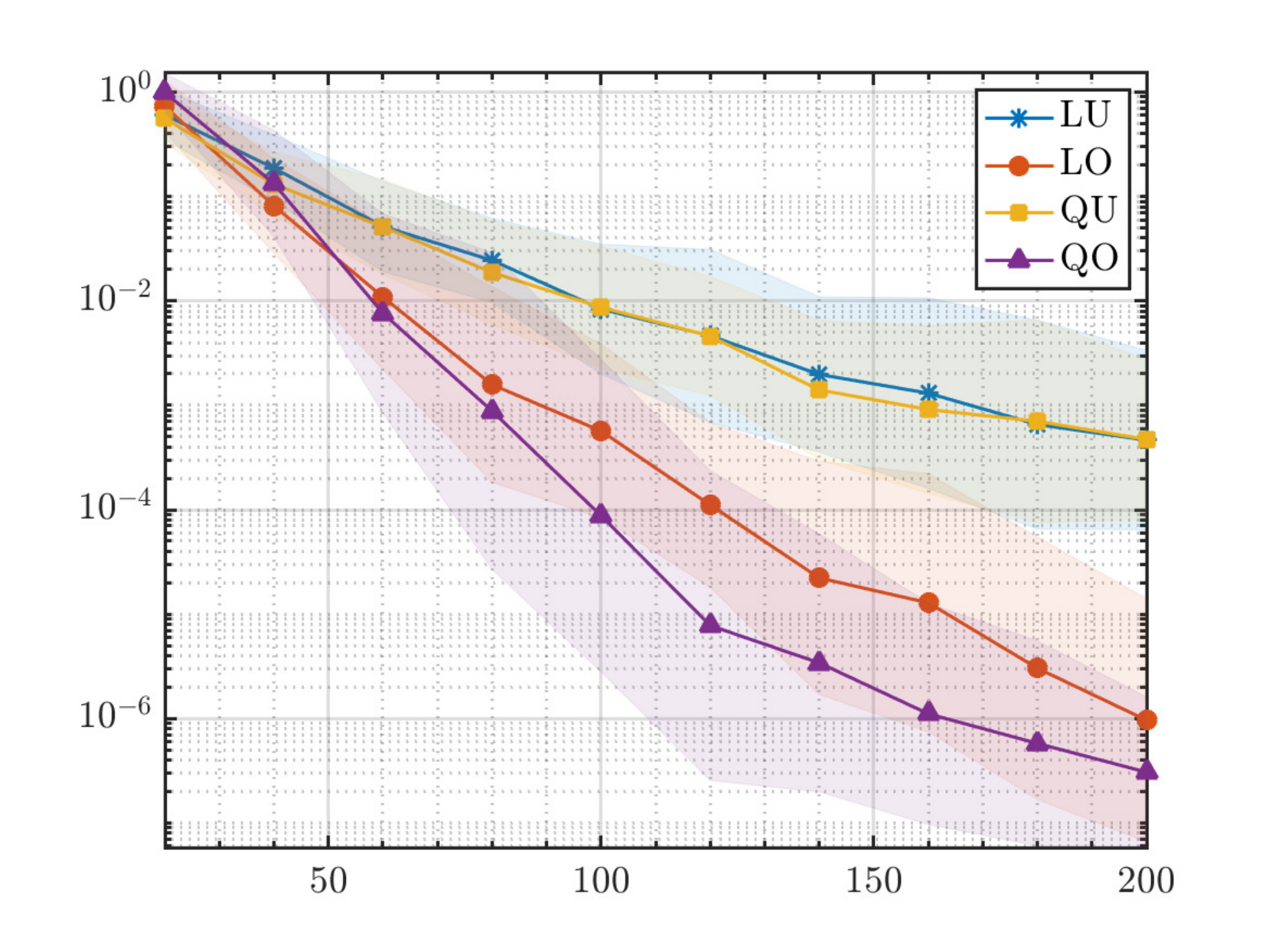} \\
 		$(d,n,N) = (1,399,400)$ & $(d,n,N)=(2,152,796)$ \\
 		$(\theta_L^2,\Theta_L^2) = (1.99,82.24)$ &
 		$(\theta_L^2,\Theta_L^2) = (3.52,281.51)$ \\
 		$(\theta_Q^2,\Theta_Q^2) = (2.62,120.57)$ &  
 		$(\theta_Q^2,\Theta_Q^2) = (7.38,322.03)$\\ 
 		$(a_L,b_l) = (6.47\times 10^{-7},3.39)$ & 
 		$(a_L,b_L) = (9.46\times^{-19},4.19)$ \\
 		$(a_Q,b_Q) = (1,1)$ &
 		$(a_Q,b_Q) = (1,1)$ \\ 	 		 		
 		\includegraphics[width=0.5\textwidth]{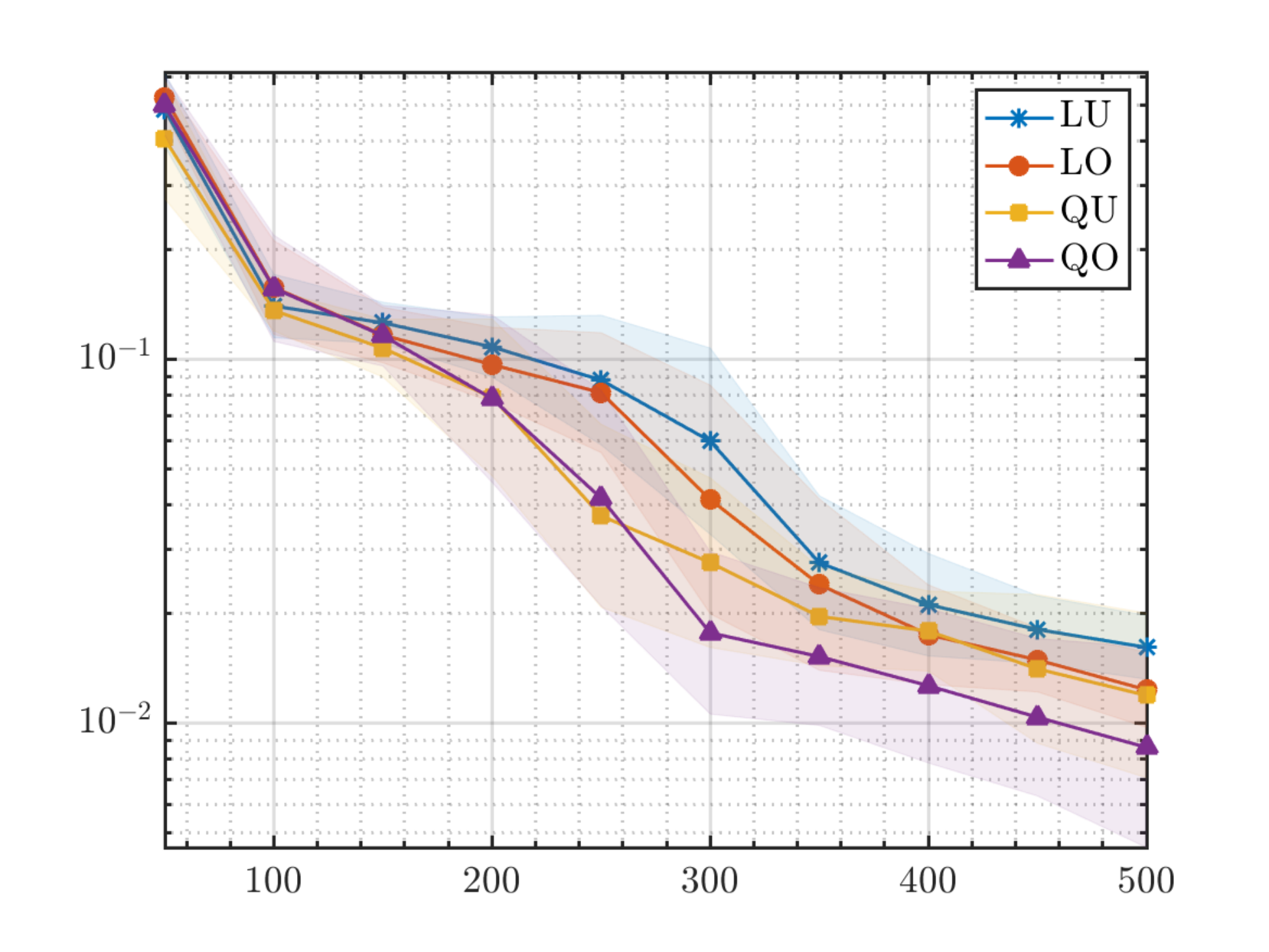} & 
 		\includegraphics[width=0.5\textwidth]{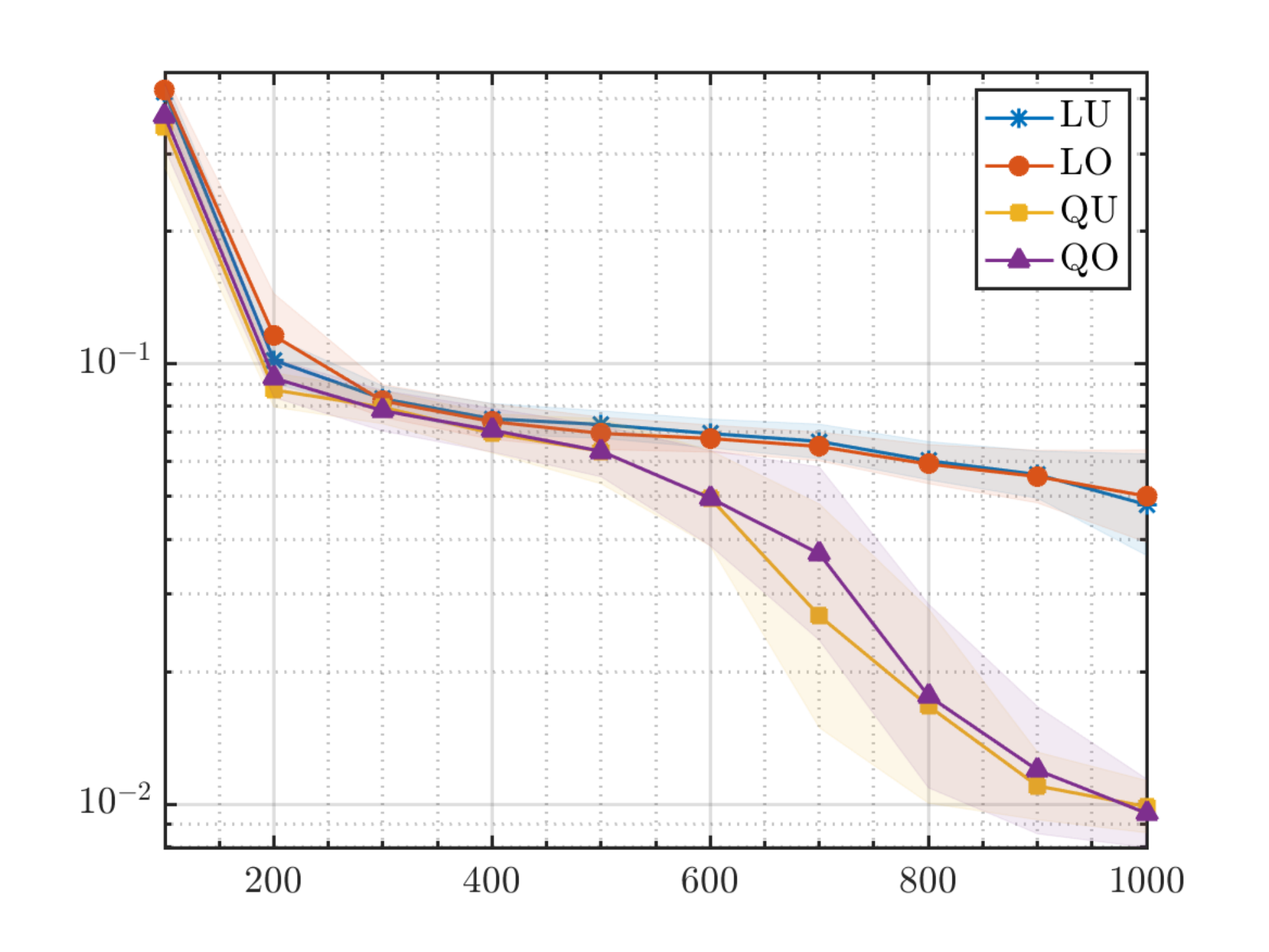} \\
 		$(d,n,N) = (8,22,1843)$ & $(d,n,N)=(16,14,4385)$ \\
 		$(\theta_L^2,\Theta_L^2) = (13.57,79.63)$ &  
 		$(\theta_L^2,\Theta_L^2) = (19.37,44.92)$ \\ 		
 		$(\theta_Q^2,\Theta_Q^2) = (19.16,340.84)$ &
 		$(\theta_Q^2,\Theta_Q^2) = (25.04,192.52)$  \\
 		$(a_L,b_l) = (0.00102,2.24)$ & 
 		$(a_L,b_L) = (0.00513,1.97)$ \\
 		$(a_Q,b_Q) = (1,1)$ &
 		$(a_Q,b_Q) = (1,1)$ \\ 	 				
	\end{tabular}
	\end{small}
	\end{center}
\caption{
The same as in Fig.\ \ref{fig:irregular-exp-1} except for $D = D_3$.
} 

\label{fig:irregular-exp-2}
\end{figure}

\begin{figure}[t]
	\begin{center}
	\begin{small}
 \begin{tabular}{@{\hspace{0pt}}c@{\hspace{-0.5pc}}c@{\hspace{0pt}}}
		\includegraphics[width=0.5\textwidth]{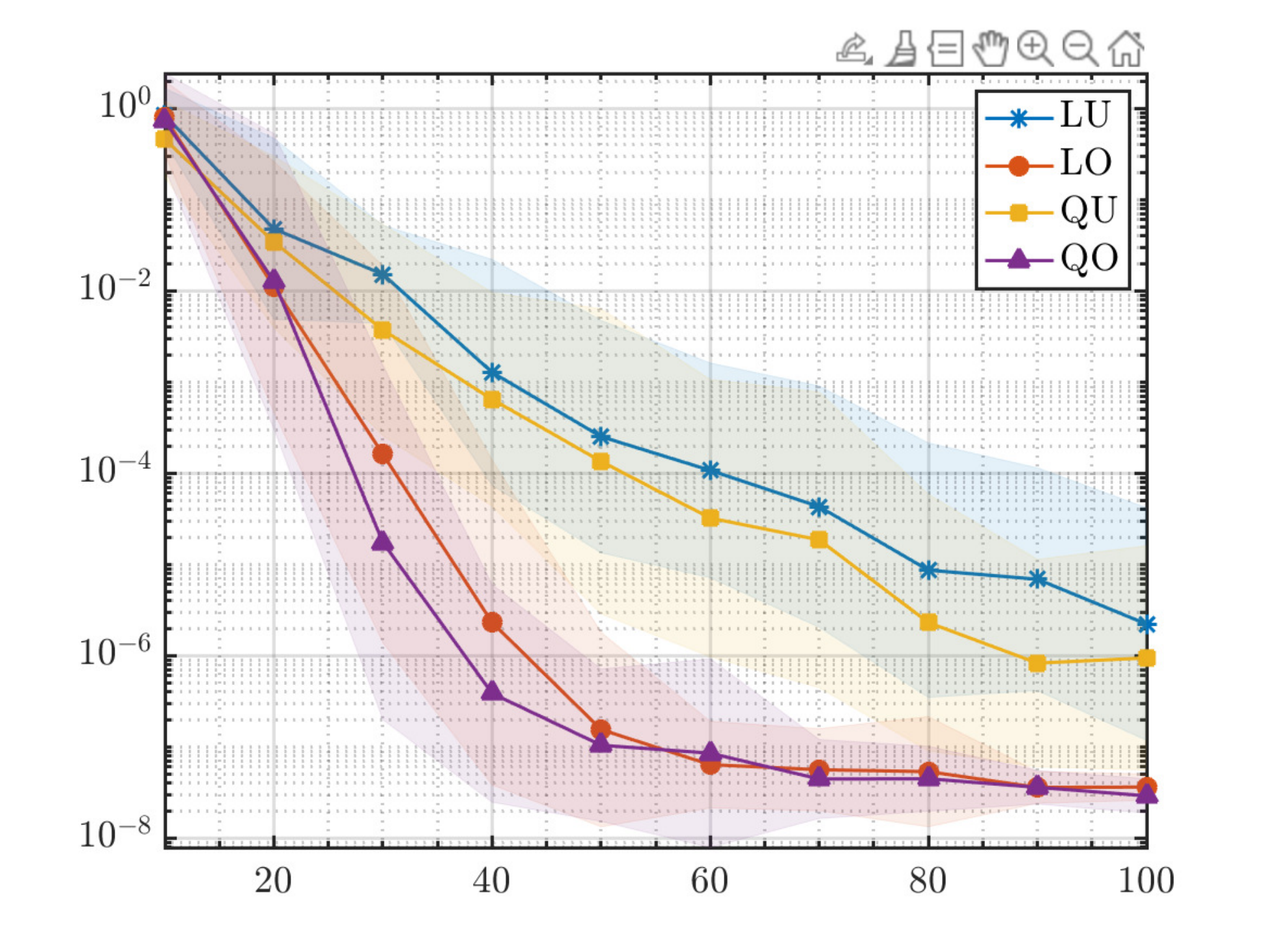} & 
 		\includegraphics[width=0.5\textwidth]{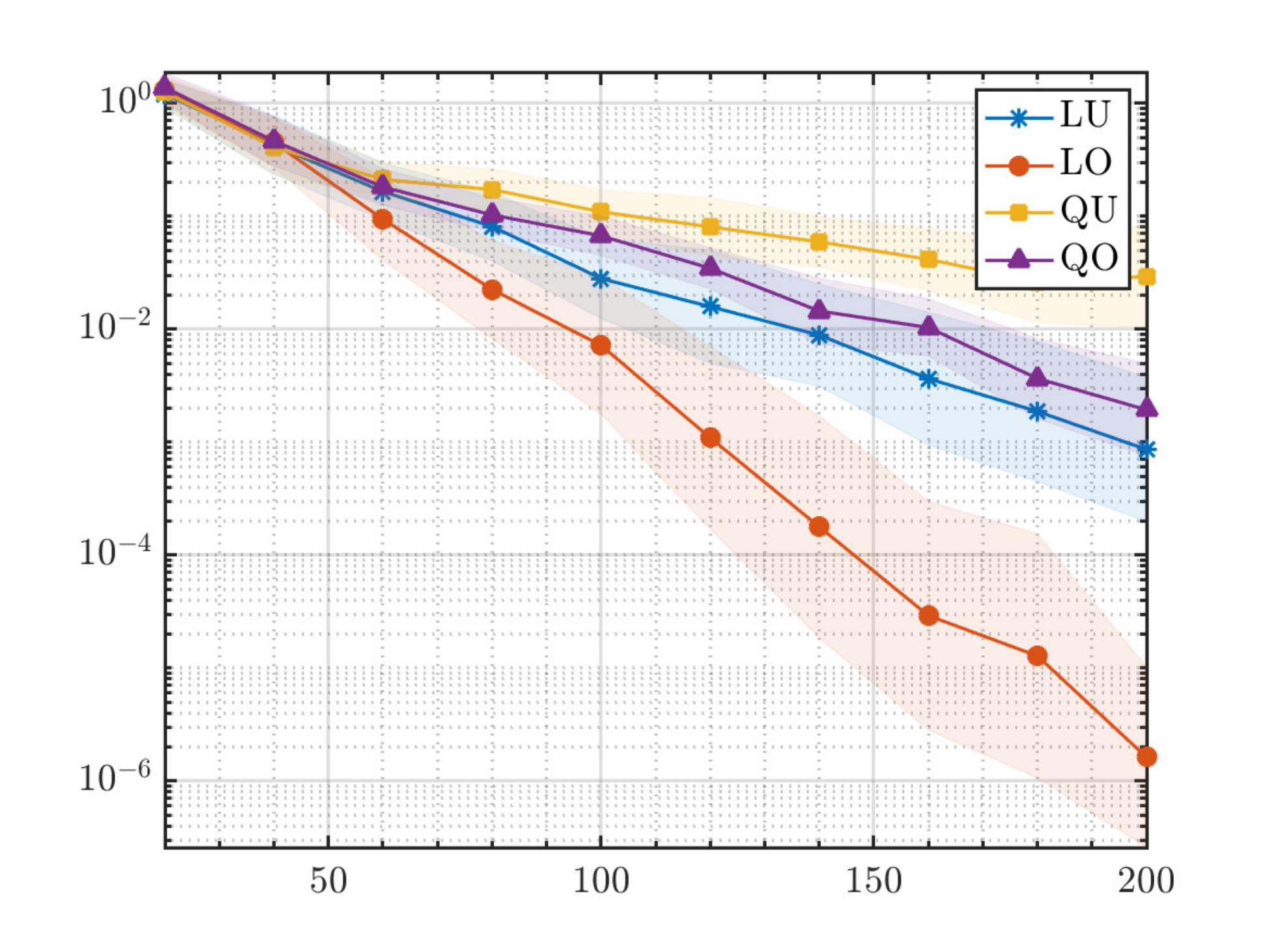} \\
 		$(d,n,N) = (1,399,400)$ & $(d,n,N)=(2,152,796)$     \\
 		$(\theta_L^2,\Theta_L^2) = (1.97,799)$ &
 		$(\theta_L^2,\Theta_L^2) = (3.73,561)$ \\
 		$(\theta_Q^2,\Theta_Q^2) = (2.56,153.84)$ &  
 		$(\theta_Q^2,\Theta_Q^2) = (5.36,205.43)$\\  	
 		$(a_L,b_L) = (0.0027,17.18)$ & 
 		$(a_L,b_L) = (0.012,4.69)$ \\
 		$(a_Q,b_Q) = (1,1)$ &
 		$(a_Q,b_Q) = (1,1)$ \\ 			
 		\includegraphics[width=0.5\textwidth]{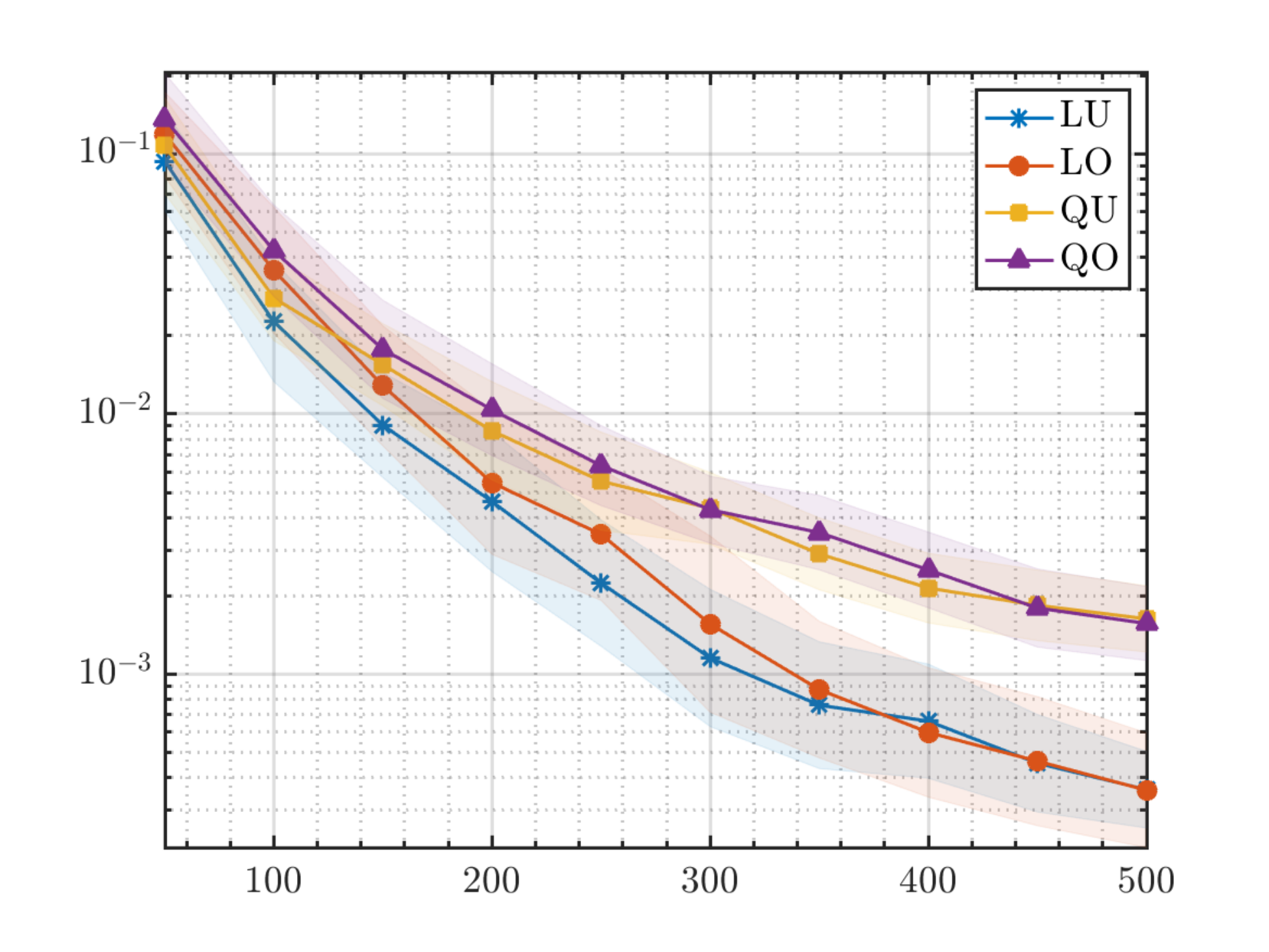} & 
 		\includegraphics[width=0.5\textwidth]{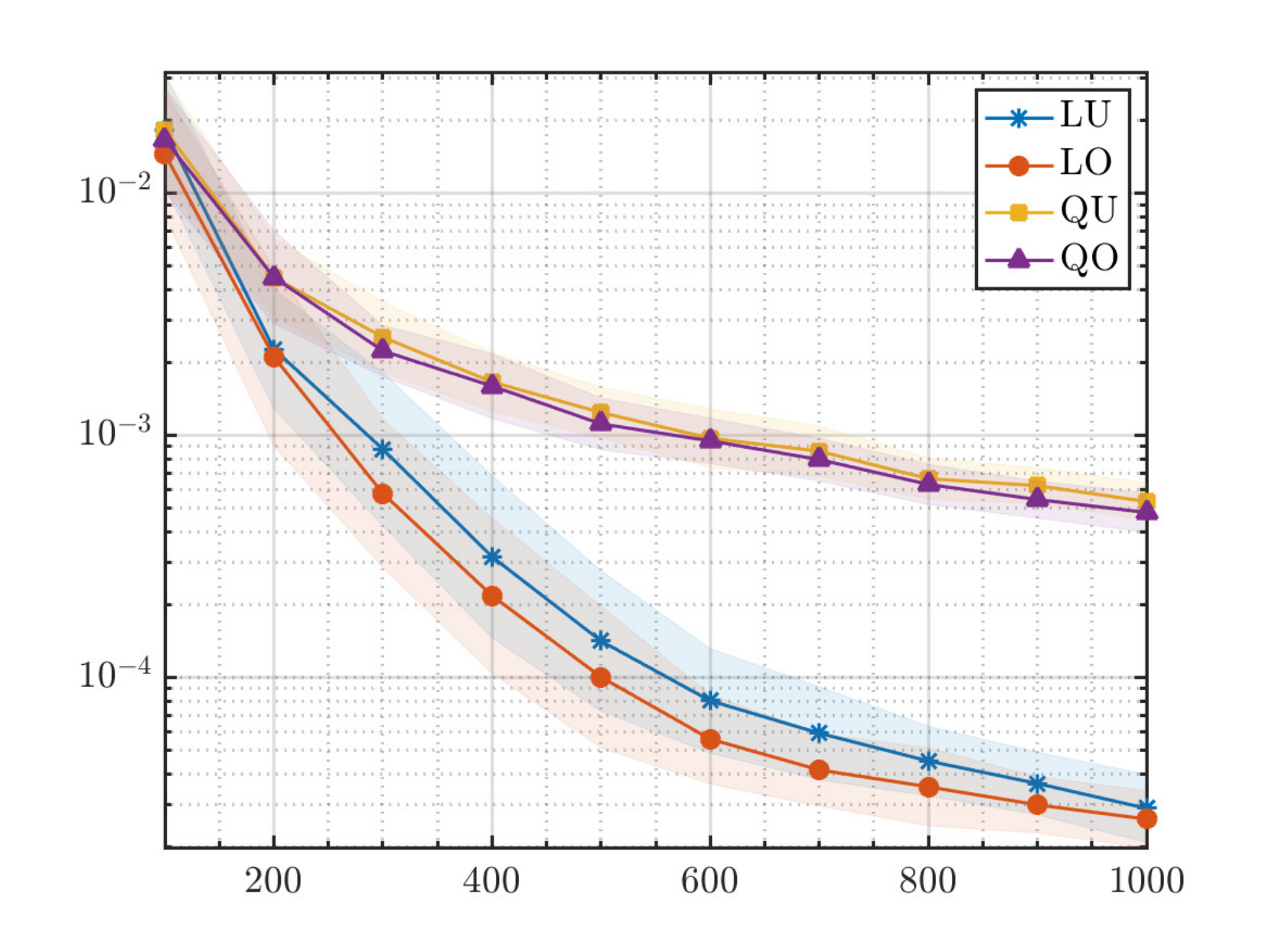} \\
 		$(d,n,N) = (8,22,1843)$ & $(d,n,N)=(16,14,4385)$    \\
 		$(\theta_L^2,\Theta_L^2) = (14.08,81)$ &
 		$(\theta_L^2,\Theta_L^2) = (19.57,45)$ \\
 		$(\theta_Q^2,\Theta_Q^2) = (14.93,91.14)$ &  
 		$(\theta_Q^2,\Theta_Q^2) = (20.18,62.26)$\\  	
 		$(a_L,b_l) = (0.48,1.87)$ & 
 		$(a_L,b_L) = (0.55,1.64)$ \\
 		$(a_Q,b_Q) = (1,1)$ &
 		$(a_Q,b_Q) = (1,1)$ \\ 	 			
	\end{tabular}
	\end{small}
	\end{center}
\caption{The same as in Fig.\ \ref{fig:irregular-exp-1} except for $D = D_1$ and $f = f_2$.} 
\label{fig:irregular-exp-3}
\end{figure}

In Fig.\ \ref{fig:irregular-exp-1} and \ref{fig:irregular-exp-2} we compare the orthogonalization strategy $\Upsilon$ against the original Legendre basis restricted to the irregular domain $D$ (labelled $\Phi$). Several effects are notable. First, orthogonalizing the basis generally leads to better performance than using the original Legendre basis. This is consistent with the observation that the various sample complexity bounds depend on the Riesz basis constants $a,b > 0$ which, as noted and as shown numerically in these figure, can behave wildly for irregular domains. On the other hand, it is clear that the values constants are extremely pessimistic when it comes to predicting the actual performance. Even when the constant $a = a_L$ is exceedingly small, the approximation based on the Legendre basis still offers a reasonable error in most cases. Further, as seen in Figs.\ \ref{fig:irregular-coefficients-ex1} and \ref{fig:irregular-coefficients-ex2} the orthogonalization strategy leads to slightly improved sparsity. Therefore, it is unclear what property of orthogonalization is driving the better approximation, whether it be the smaller Riesz basis constants or the improvement in sparsity. On the other hand, in Fig.\ \ref{fig:irregular-exp-3} we present an example where orthogonalization worsens the approximation. This we expect is due to the worse sparsity in the $\Upsilon$ basis, as shown in Fig.\ \ref{fig:irregular-coefficients-ex3}.

Second, we observe that the `optimal' sampling procedures generally outperform Monte Carlo sampling in lower dimensions, while this improvement lessens in higher dimensions, or may actually lead to slightly worse performance. This is consistent with the observation made previously in \S \ref{ss:examples-CS}. We also report the values of the constants $\theta$ and $\Theta$ in all cases. It is notable that, in the case of the orthogonalized basis, the corresponding constant $\Theta = \Theta_Q$ is much larger than $\theta = \theta_Q$, even in high dimensions. However, this is not reflected in the approximation errors for the two sampling strategies, which are similar, thus suggesting a gap between the theoretical guarantees and performance on actual function approximation problems.

\section{Structured sparse approximation}\label{s:lower-sets}

Our main assumption throughout this chapter has been that $f$ admits an approximately sparse representation in a dictionary $\Phi = \{ \phi_{\iota} : \iota \in \cI \}$. We conclude this chapter with a brief discussion on several types of \textit{structured} sparsity models. Our focus throughout is on the setting of Problem \ref{ass:sparsity_unknown}. Hence, as in \S \ref{s:CS}, we assume that $\Phi$ is a finite set of $n$ linearly-independent functions.

\subsection{Weighted sparsity and weighted $\ell^1$-minimization}

Let $v = (v_{\iota})_{\iota \in \cI}$ be a vector of positive weights. For a set $S \subseteq \cI$, we define its \textit{weighted} cardinality as
\bes{
|S|_{v} = \sum_{\iota \in S} v^2_{\iota}.
}
In the \textit{weighted sparsity} model, given a weights $v$ and a weighted sparsity $k > 0$, we assume that a function $f \in L^2_{\rho}(D;\bbV)$ has a sparse representation $f_S$ of the form \R{sparse_rep} for some set $S \subseteq \cI$ with $|S|_{v} \leq k$. Note that, unlike the case of standard sparsity, the weighted sparsity parameter can take any positive value in this setting -- hence our reason for using the notation $k$ instead of $s$.

Fortunately, promoting weighted sparsity structure is straightforward. Rather than $\ell^1$-minimization, i.e.\ \R{qcbp}, \R{lasso} or \R{sr-lasso}, we consider a \textit{weighted} $\ell^1$-minimization problem. For example, we may replace \R{sr-lasso} by
\be{
\label{sr-lasso-weighted}
\hat{f} \in \argmin{p \in P_{\cI ; \bbV_h}} \left \{ \lambda \nm{c}_{\ell^1_v(\cI ; \bbV)} + \sqrt{\frac1m \sum^{m}_{i=1} w(y_i) \nm{f(y_i) + n_i - p(y_i) }^2_{\bbV}} \right \}.
}
Here, $\nm{c}_{\ell^1_v(\cI ; \bbV)} = \sum_{\iota \in \cI} v_{\iota} \nm{c_{\iota}}_{\bbV}$ is the weighted $\ell^1$-norm of a Hilbert-valued vector $c = (c_{\iota})_{\iota \in \cI}$.

\rem{
Much like the lower set assumption (see \S \ref{MC-sampling-LS}), weighted sparsity is a natural assumption to consider when one expects the most significant coefficients of $f$ to correspond to lower-order terms. We discuss the relation between the two models later. This is typically the case for smooth function approximation using algebraic or trigonometric polynomials. As we see later, incorporating slowly-growing weights can lead to a significant improvement over unweighted $\ell^1$-minimization. Note that weighted sparsity and weighted $\ell^1$-minimization were first elaborated in \cite{rauhut2016interpolation}, before further developments in \cite{adcock2017infinite,adcock2018infinite,adcock2018compressed,chkifa2018polynomial}. 
Other works on incorporating weights into sparse polynomial approximation include \cite{yang2013reweighted,peng2014weighted,adcock2020sparse}.
}

As in the case of standard sparsity, successful recovery via \R{sr-lasso-weighted} follows from a norm equivalence similar to \R{RIP_over_all_S}: namely,
\be{
\label{wRIP_over_all_S}
\alpha \nm{p}^2_{L^2_{\varrho}(D)} \leq \frac{1}{m} \sum^{m}_{i=1} w(y_i) | p(y_i)|^2 \leq \beta \nm{p}^2_{L^2_{\varrho}(D)},\quad \forall p \in P_{T},\ T \subseteq \cI,\ |T|_v\leq t.
}
Under this condition, one obtains an error bound identical to Theorem \ref{t:CS_acc_stab}, except with $\sigma_{s}(c)_{\ell^1(\cI;\bbV)}$ replaced by the weighted term
\bes{
\sigma_{k,v}(x)_{\ell^1_v(\cI;\bbV)} = \inf \left \{ \nm{x-z}_{\ell^1_v(\cI;\bbV)} : \mbox{$z \in \bbV^{n}$ is weighted $(k,v)$-sparse} \right \}, \quad x  \in \bbV^{n},
}
and $s$ replaced by $k$. For the sake of succinctness, we omit the details, and refer to \cite{rauhut2016interpolation,adcock2018compressed} (the results therein are given in the scalar-valued case, but readily extend to the Hilbert-valued case). 

Since our primary focus is on the question of sample complexity, we now state a variant of Theorem \ref{t:BOS_RIP} for \R{wRIP_over_all_S} for the weighted sparse model:

\thm{[Sample complexity for \R{wRIP_over_all_S}]
\label{t:BOS_wRIP}
Let $\Phi = \{ \phi_{\iota} : \iota \in \cI \} \subset L^2_{\rho}(D)$ be a finite dictionary consisting of $n$ linearly-independent functions, with bounds $a,b > 0$ as in \R{Riesz-bounds}. Let $\mu$ be a probability measure satisfying Assumption \ref{ass:abs-cont-pos}, $v = (v_{\iota})_{\iota \in \cI}$ be weights with
\be{
\label{weights-suffic-size}
v_{\iota} \geq \nmu{\sqrt{w(\cdot)}\phi_{\iota}(\cdot)}_{L^{\infty}_{\rho}(\cU)},\quad \forall \iota \in \cI,
}
where $w$ is as in \eqref{mu_weight_fn_single},
$1 \leq t \leq n$, $0 <  \delta < \delta^*$ for some universal constant $0 < \delta^*< 1$, $0 < \epsilon < 1$, $0 < \alpha \leq \beta <\infty$,  and $y_1,\ldots,y_m$ be independent with $y_i \sim \mu$ for $i = 1,\ldots,m$. Suppose that
\bes{
m \geq C \cdot \delta^{-2} \cdot a^{-1} \cdot t \cdot \left ( \log(\E n) \cdot \log^2(\E a^{-1} t / \delta )  + \log(2/\epsilon) \right ),
}
for some universal constant $C >0$.
Then \R{wRIP_over_all_S} holds with $1-\delta \leq \alpha \leq \beta \leq 1+\delta$, with probability at least $1-\epsilon$.
}
We omit the proof of this result, since it is similar to that of Theorem \ref{t:brugiapaglia2021sparse}, the main difference being the use of \cite[Thm.\ 2.13]{brugiapaglia2021sparse} instead of \cite[Thm.\ 1.1]{brugiapaglia2021sparse}.

\subsection{Sparsity in lower sets}\label{ss:lower-set-sparsity}

The lower set sparsity model differs from the weighted sparsity model in that it imposes a lower set structure, as opposed to a weighted sparsity structure. But in practice it can also be effected via weighted $\ell^1$-minimization with specific choices of weights. In the lower set sparsity model, we suppose that $f \in L^2_{\rho}(D;\bbV)$ has a sparse representation $f_S$ of the form \R{sparse_rep} for some set $S \subseteq \cI$ with $|S| \leq s$ that is also lower. In terms of sufficient conditions for lower set recovery, one's first thought may be to consider a variant of \R{RIP_over_all_S} with the additional assumption that the sets $T$ be lower. Unfortunately, it is not known whether such an approach can work. The difficulty lies with the fact that it is unclear how to promote lower set structure directly via a convex penalty term such as the $\ell^1$-norm \cite{adcock2018compressed}.

Instead, the approach originally proposed in \cite{chkifa2018polynomial,adcock2018infinite} 
is to use weighted sparsity as a surrogate for lower set sparsity. This is done by choosing weights $u = (u_{\iota})_{\iota \in \cI}$ as small as possible so that Theorem \ref{t:BOS_wRIP} applies, namely,
\bes{
u_{\iota} =  \nmu{\sqrt{w(\cdot)}\phi_{\iota}(\cdot)}_{L^{\infty}_{\rho}(\cU)},\quad \forall \iota \in \cI,
}
and defining the weighted sparsity as
\bes{
k = k(s;w) = \max \left \{ |S|_u : \mbox{$|S| \leq s$, $S$ lower} \right \}.
}
Note that this ensures that every lower set of size $s$ has weighted cardinality at most $k$, i.e.
\bes{
\left \{ S : \mbox{$|S| \leq s$, $S$ lower}  \right \} \subseteq \left \{ S : |S|_u \leq k(s;w) \right \}.
}
As a result, we can promote lower set sparsity by solving the weighted $\ell^1$-minimization problem \R{sr-lasso-weighted} with weights $v = u$.

\rem{[The choice of $\cI$]
\label{r:choice-I}
Working with lower sets also yields a strategy for choosing the large truncated set $\cI$ (recall the discussion at the beginning of \S \ref{s:CS}) \cite{adcock2018compressed,chkifa2018polynomial}. Indeed, it is a straightforward exercise to show that the union of all lower sets of size at most $s$ is the hyperbolic cross $\cI^{\mathrm{HC}}_{s-1}$. Hence, the target lower set $S$ in the sparse representation \R{sparse_rep} is guaranteed to lie within this index set, thus giving a clear rationale for choosing this set.
}

\subsection{Sampling and numerical experiments}

We now discuss the matter of sampling. Notice that, unlike in the case of Theorem \ref{t:BOS_RIP}, the sample complexity bound in Theorem \ref{t:BOS_wRIP} does not involve a constant $\Gamma$ depending on the basis $\Phi$ and weight function $w$. This dependence only arises in the minimum size condition \R{weights-suffic-size} on the weights $v$. For Monte Carlo sampling ($w \equiv 1$), this condition may be quite stringent if the $L^{\infty}$-norms of the basis functions grow rapidly. Hence, this conditions suggests choosing $w$ to minimize the right-hand side of \R{weights-suffic-size}. This leads to the same choice \R{w-theta-opt-CS} and \R{mu-optimal-CS} as in the standard sparsity setting considered previously.

The case of lower set sparsity allows for a more concrete discussion. Choosing weights $v = u$ as discussed above, and invoking Theorem \ref{t:BOS_wRIP} leads to a measurement condition of the form
\bes{
m \gtrsim a^{-1} \cdot k(s;w) \cdot \left ( \log(\E n) \cdot \log^2(\E a^{-1} k(s;w) / \delta )  + \log(2/\epsilon) \right ),
}
for recovering functions with sparse representations in lower sets. Hence, the objective is to minimize $k(s;w)$ with respect to $w$. In the case of Monte Carlo sampling, we have
\bes{
k(s;w) = k(s;1) = \max \left \{ |S|_u : \mbox{$|S| \leq s$, $S$ lower} \right \}.
}
Consider, for illustration, Example \ref{ex:alg-poly}. In this case, since the Legendre polynomials all attain their maximum at the same point $y = (1,\ldots,1)^{\top}$, we have
\bes{
|S|_u = (\cN(P_S))^2,
} 
where $\cN(P_S)$ is the unweighted Nikolskii constant \R{unweighted-Nikolskii}. Recall from the discussion in \S \ref{MC-sampling-LS} that $\cN(P_S)$ satisfies the sharp bound $(\cN(P_S))^2 \leq s^2$ for lower sets. In other words, $k(s;1) = s^2$ in this case, leading to a sample complexity bound of the form
\bes{
m \gtrsim s^2 \cdot  \left ( \log(\E n) \cdot \log^2(\E s)  + \log(2/\epsilon) \right )
}
In other words, the worst-case sample complexity for lower set recovery via Monte Carlo sampling is the same (up to constants and log factors) as that of least squares in the setting of Problem \ref{ass:sparsity_known}. See \cite{adcock2018compressed} for further discussion.

Having considered Monte Carlo sampling for lower set recovery, one may also consider how to choose the weight function $w$ and corresponding sampling measure $\mu$ to improve the sample complexity. The best solution in this case involves choosing $w$ to minimize
\bes{
k(s;w) = \max \left \{ \sum_{\iota \in S} \nmu{\sqrt{w(\cdot)}\phi_{\iota}(\cdot)}^2_{L^{\infty}_{\rho}(\cU)} :  \mbox{$|S| \leq s$, $S$ lower} \right \},
}
over all strictly positive and finite almost everywhere weight functions on $\supp(\rho)$ for which \R{w_normalization} holds. Unfortunately, even after resorting to a discrete measure as in \S \ref{ss:CS-opt-disc}, it is unclear how to compute such a $w$, since it seemingly involves enumerating all lower sets. As shown in \cite{cohen2017discrete}, there are many lower sets in high dimensions (for example, at least ${d \choose s-1}$ when $s \leq d + 1$).

Since the optimal choice of $w$ (in the sense of minimizing $k(s;w)$) may not be available, it is natural to consider how one might choose a good $w$. One option involves the choice \R{w-theta-opt-CS}. This leads to the bound $k(s;w) \leq \theta^2 s$. This has the benefit of scaling linearly in $s$. But, it gives a sample complexity bound that is no better than the case of standard sparse recovery studied previously. Once more, the question of whether one can choose $w$ in such a way to ensure optimal recovery (scaling linearly in $s$ and at most logarithmically in $d$ and $n$) is currently unresolved.

\begin{figure}[t]
	\begin{center}
	\begin{small}
 \begin{tabular}{@{\hspace{0pt}}c@{\hspace{-0.5pc}}c@{\hspace{-0.5pc}}c@{\hspace{0pt}}}
		
 		\includegraphics[width=0.5\textwidth]{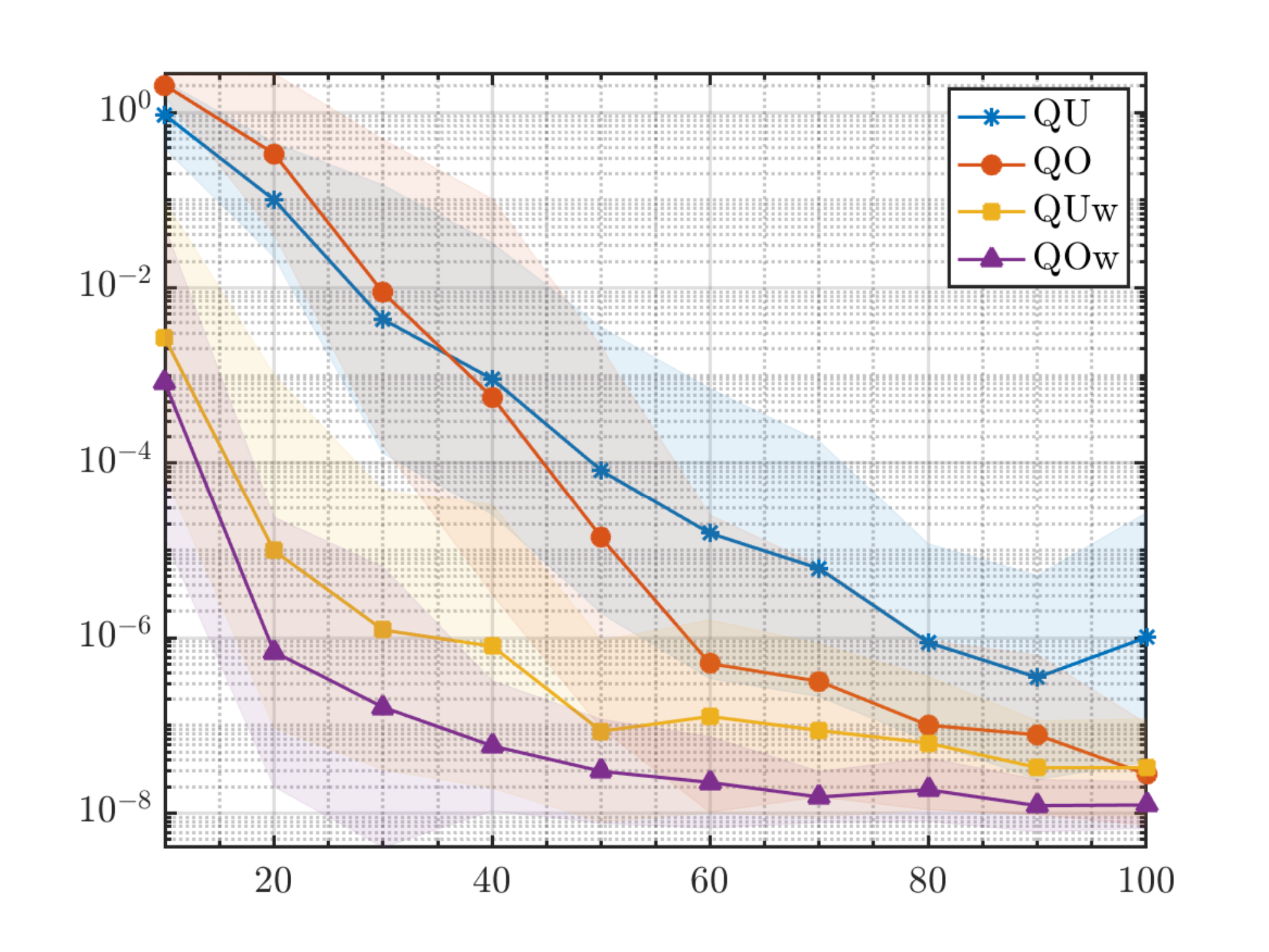} & 
 		\includegraphics[width=0.5\textwidth]{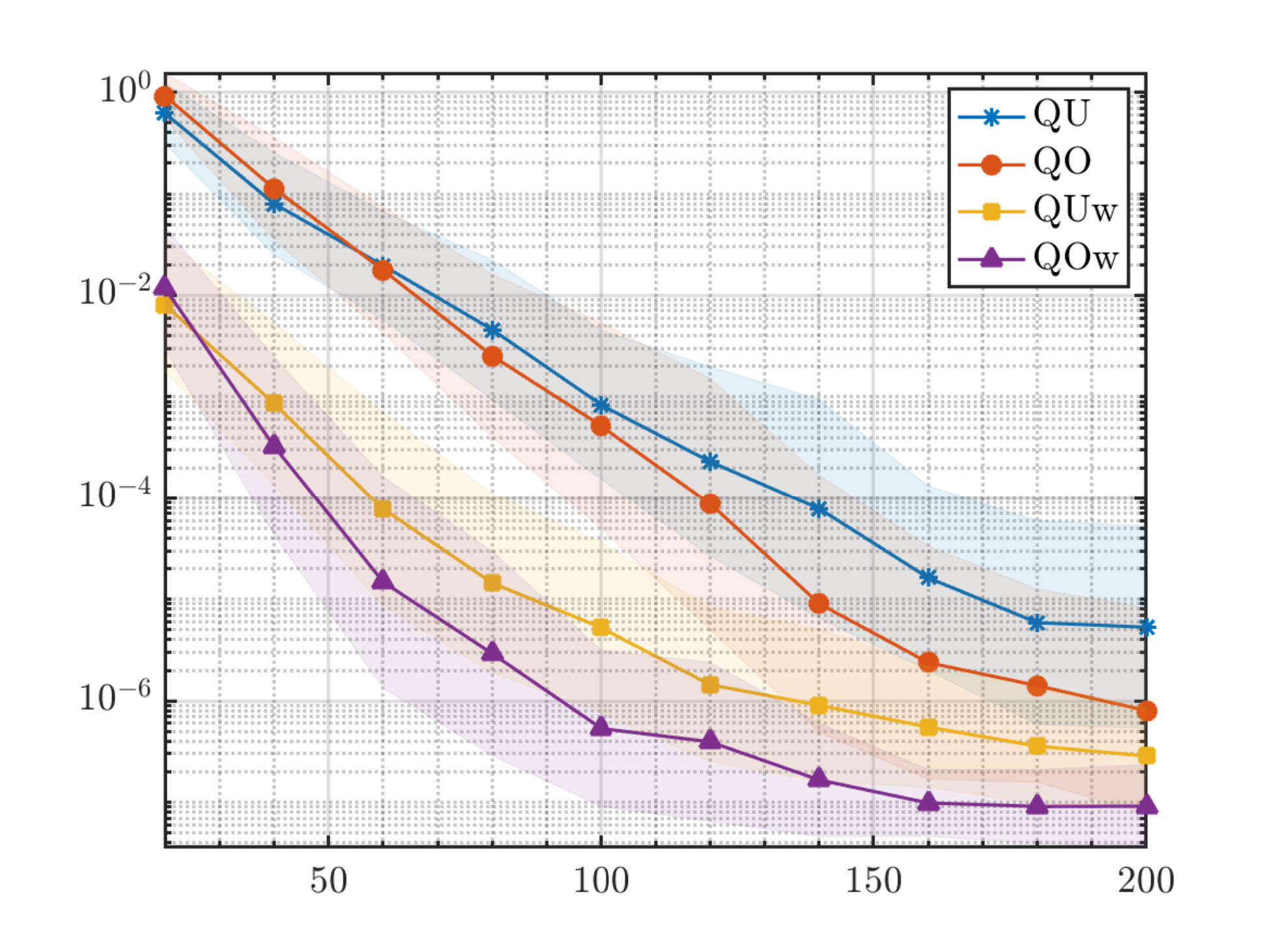} \\
 		$(d,n,N) = (1,399,400)$ & $(d,n,N)=(2,152,796)$ \\
 		\includegraphics[width=0.5\textwidth]{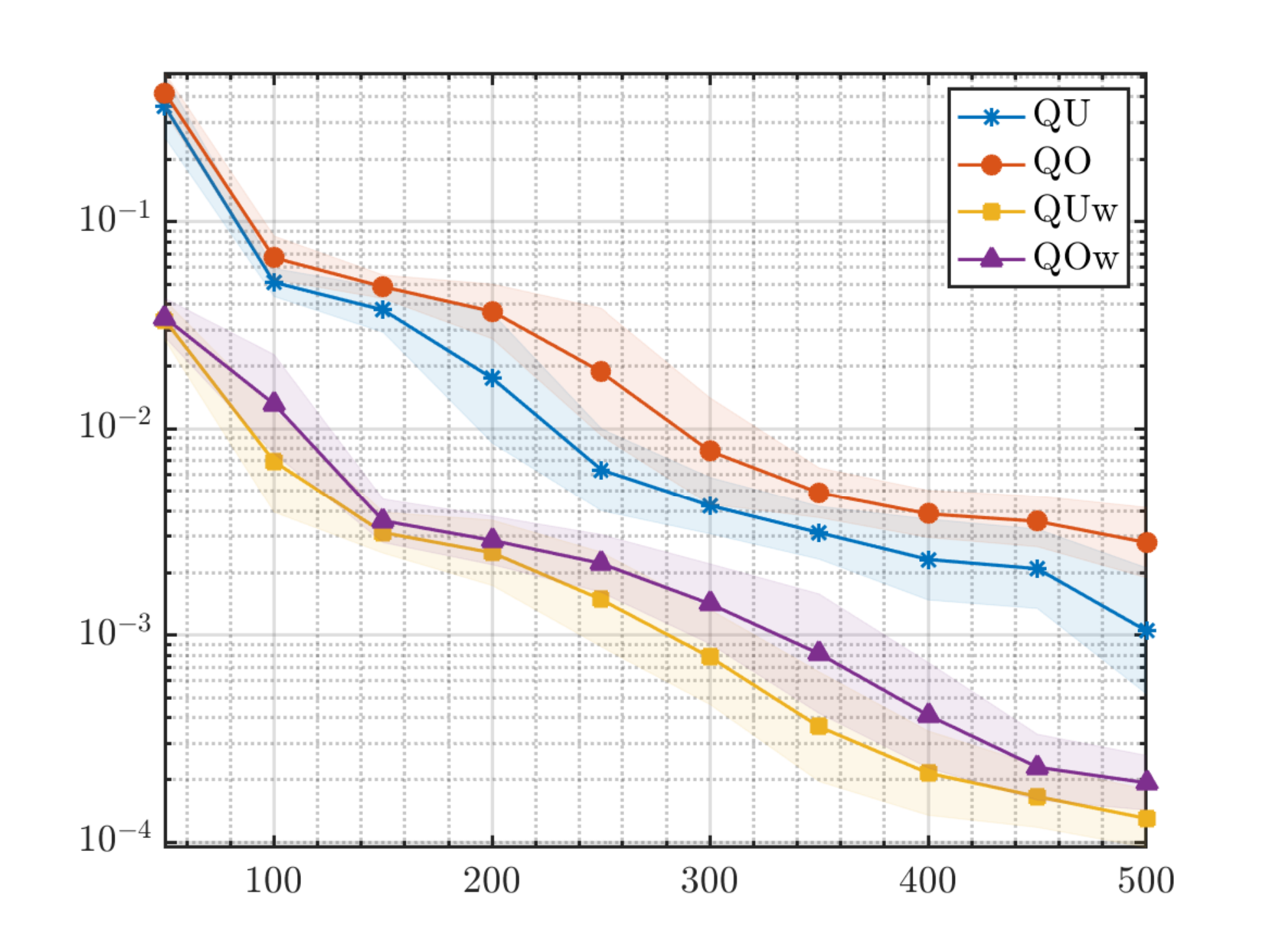} & 
 		\includegraphics[width=0.5\textwidth]{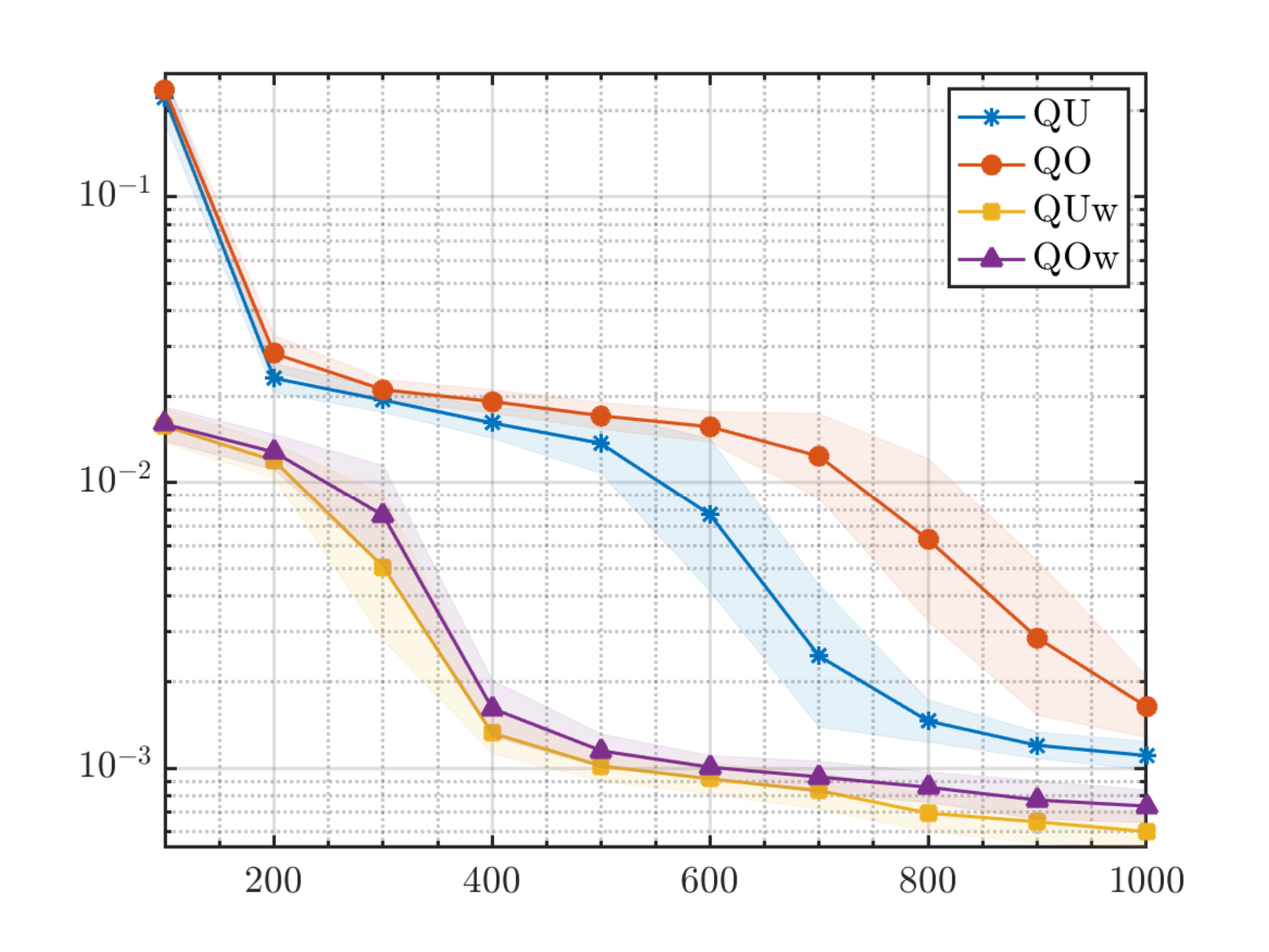} \\
 		$(d,n,N) = (8,22,1843)$ & $(d,n,N)=(16,14,4385)$ \\
	\end{tabular}
	\end{small}
\end{center}
\caption{ The relative error \R{relative-error} versus $m$ for $\ell^1$-minimization and weighted $\ell^1$-minimization in the case of Example \ref{ex:alg-poly-general}, where $\cI = \cI^{\mathrm{HC}}_{t-1}$ is the hyperbolic cross index set \R{HC-index} and $f = f_1$ and $D = D_2$ are as in \R{functions-define} and \R{domains-define}, respectively. This figure compares Monte Carlo sampling \R{MC-CS-ONB} from the discrete uniform measure (labelled `QU') and sampling from the discrete `optimal' measure \R{Optimal-CS-ONB} (`QO').} 
\label{fig:irregular-exp-1-weighted}
\end{figure}

\begin{figure}[t]
	\begin{center}
	\begin{small}
 \begin{tabular}{@{\hspace{0pt}}c@{\hspace{-0.5pc}}c@{\hspace{-0.5pc}}c@{\hspace{0pt}}}
  		\includegraphics[width=0.5\textwidth]{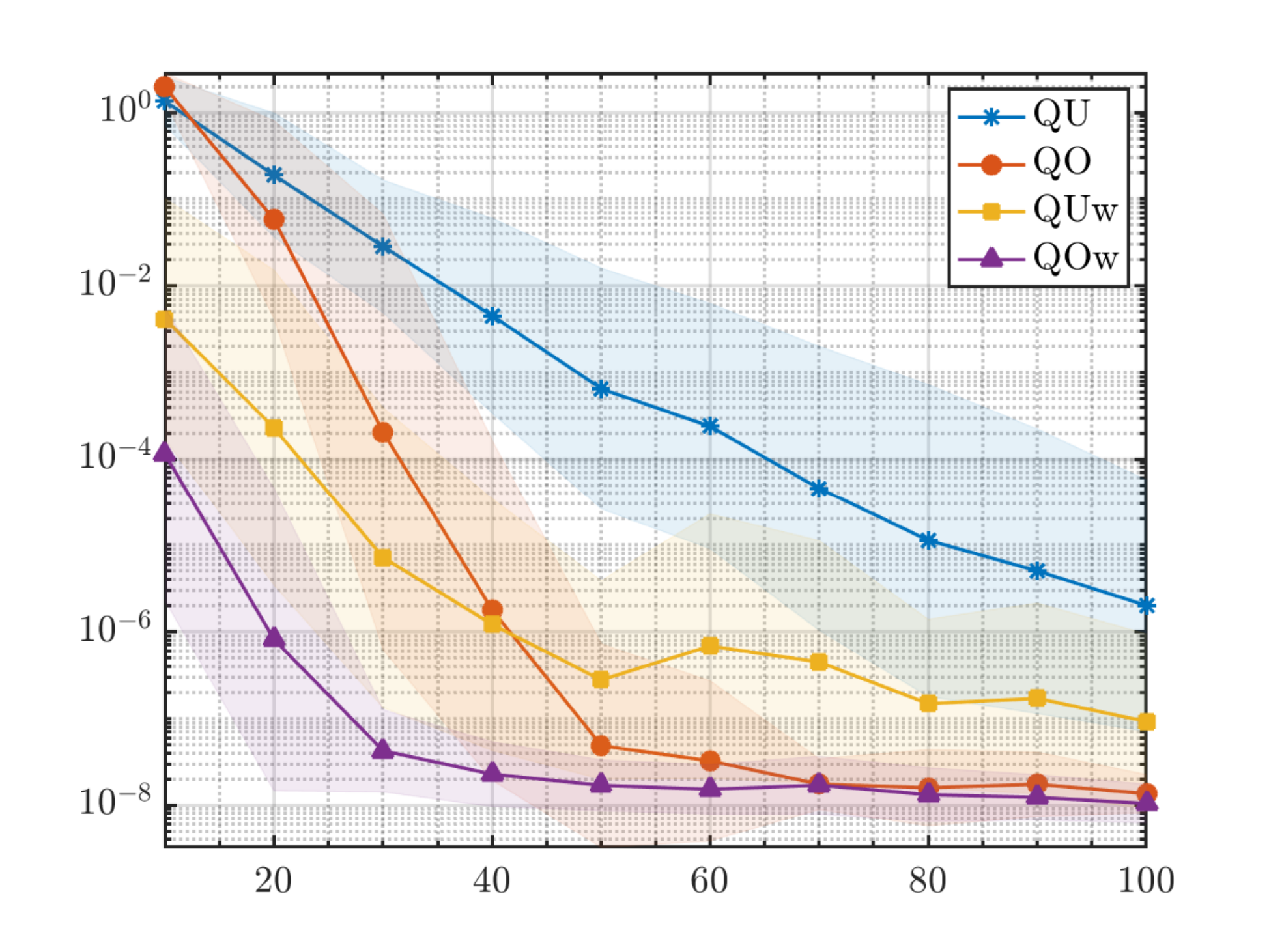} & 
 		\includegraphics[width=0.5\textwidth]{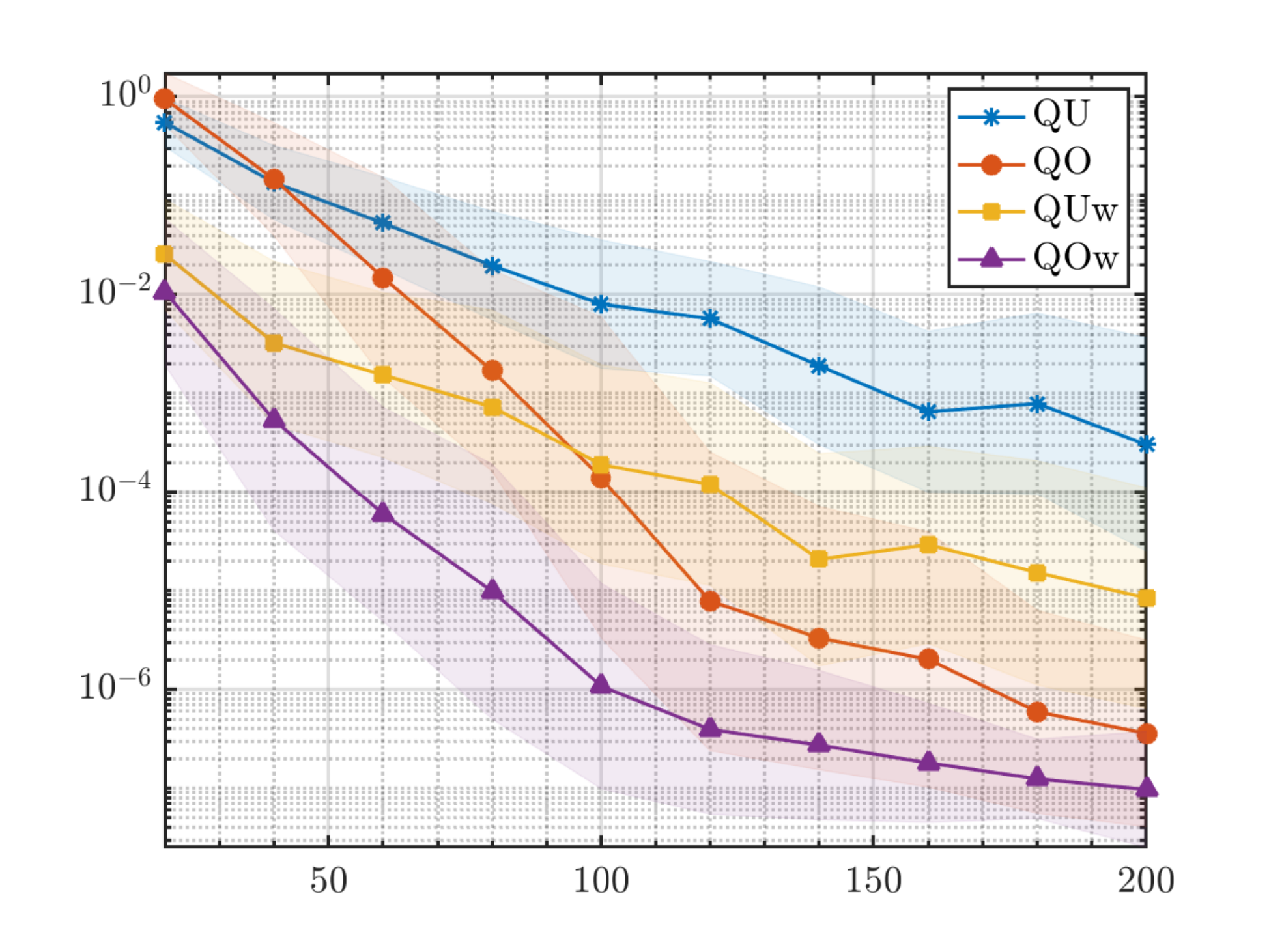} \\
 		$(d,n,N) = (1,399,400)$ & $(d,n,N)=(2,152,796)$ \\
 		\includegraphics[width=0.5\textwidth]{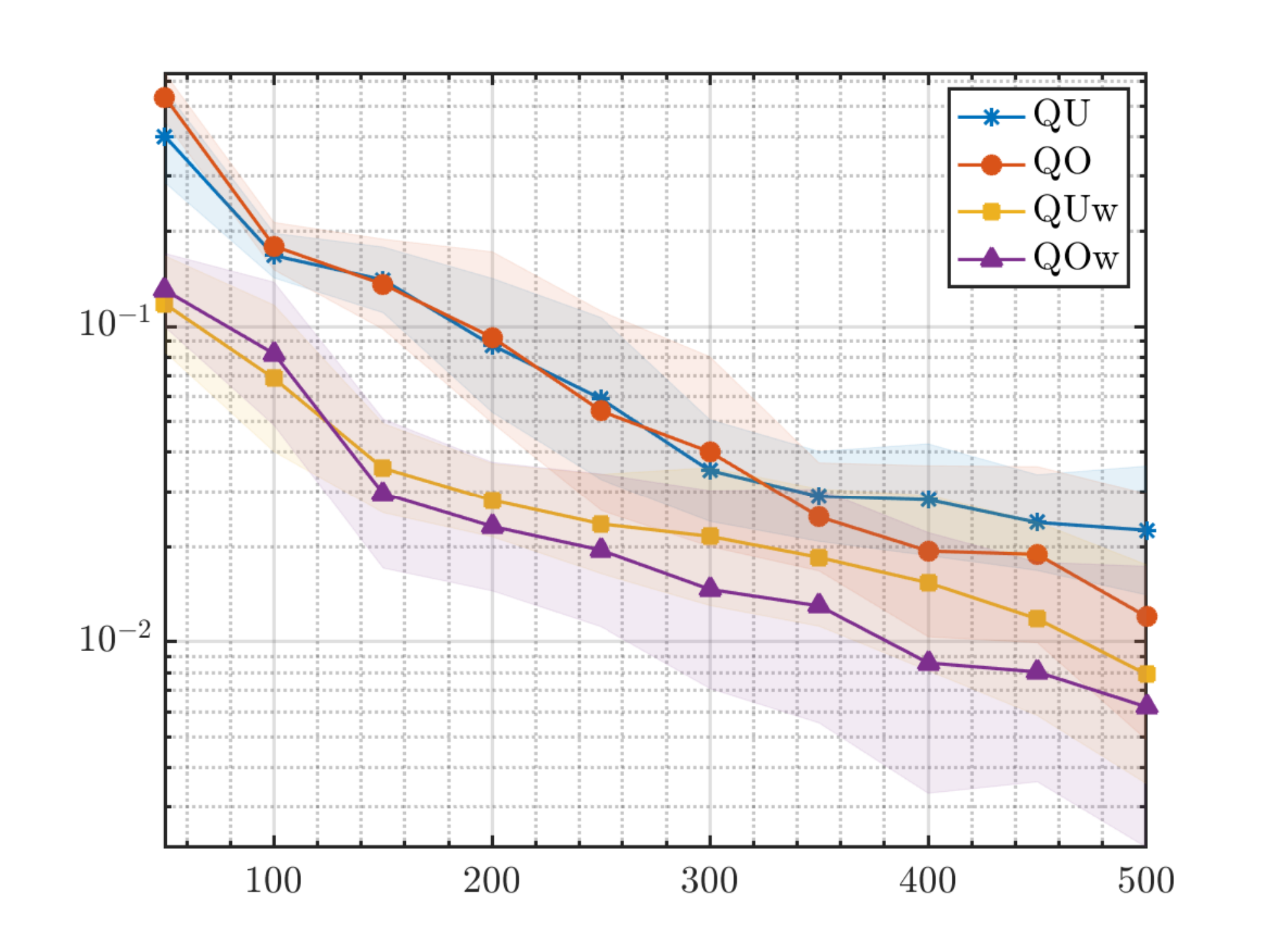} & 
 		\includegraphics[width=0.5\textwidth]{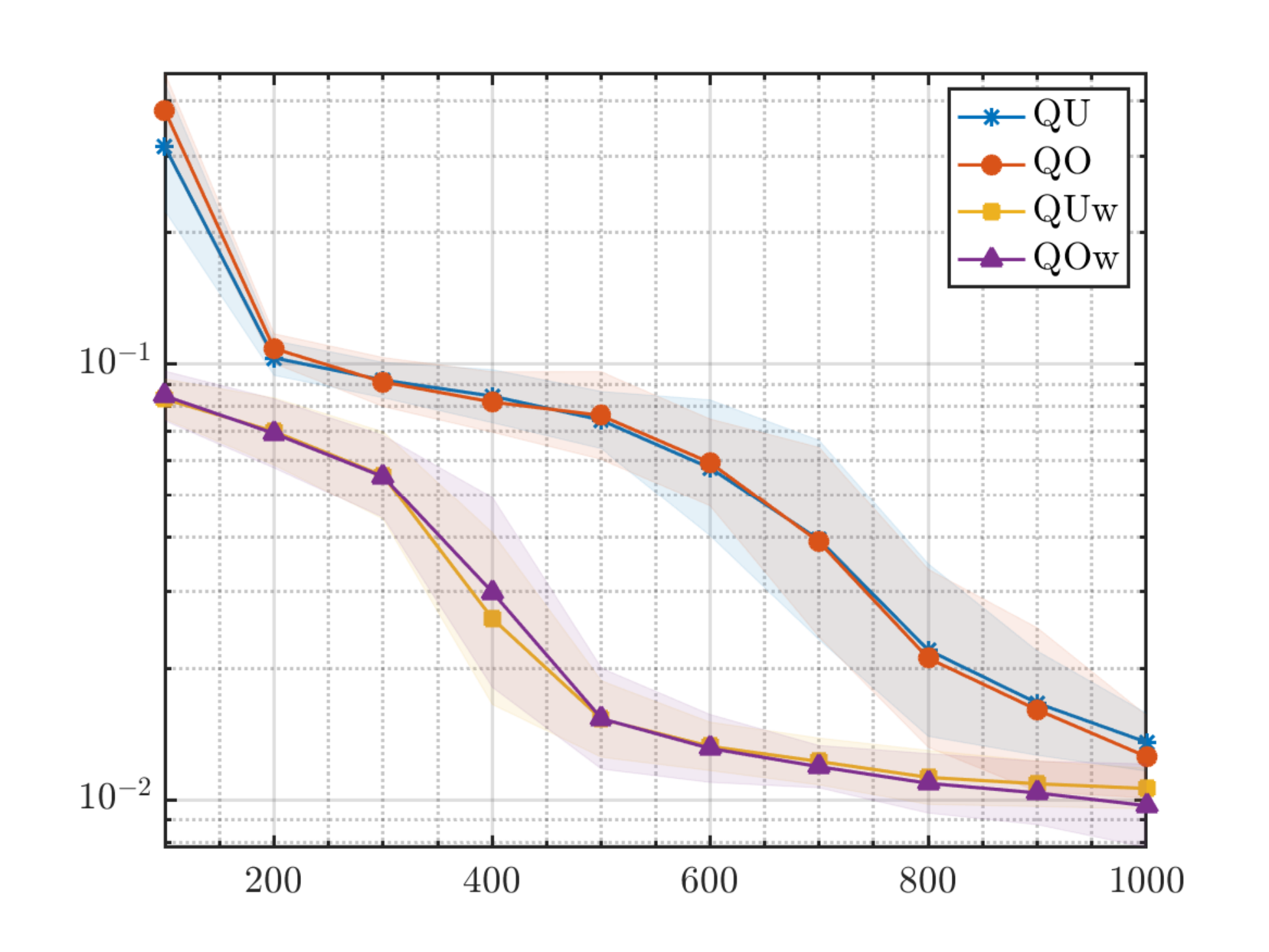} \\
 		$(d,n,N) = (8,22,1843)$ & $(d,n,N)=(16,14,4385)$ \\
	\end{tabular}
	\end{small}
	\end{center}
\caption{The same as Fig,\ \ref{fig:irregular-exp-1-weighted}, except for $D = D_3$.} 
\label{fig:irregular-exp-2-weighted}
\end{figure}

\begin{figure}[t]
	\begin{center}
	\begin{small}
 \begin{tabular}{@{\hspace{0pt}}c@{\hspace{-0.5pc}}c@{\hspace{-0.5pc}}c@{\hspace{0pt}}}
		\includegraphics[width=0.5\textwidth]{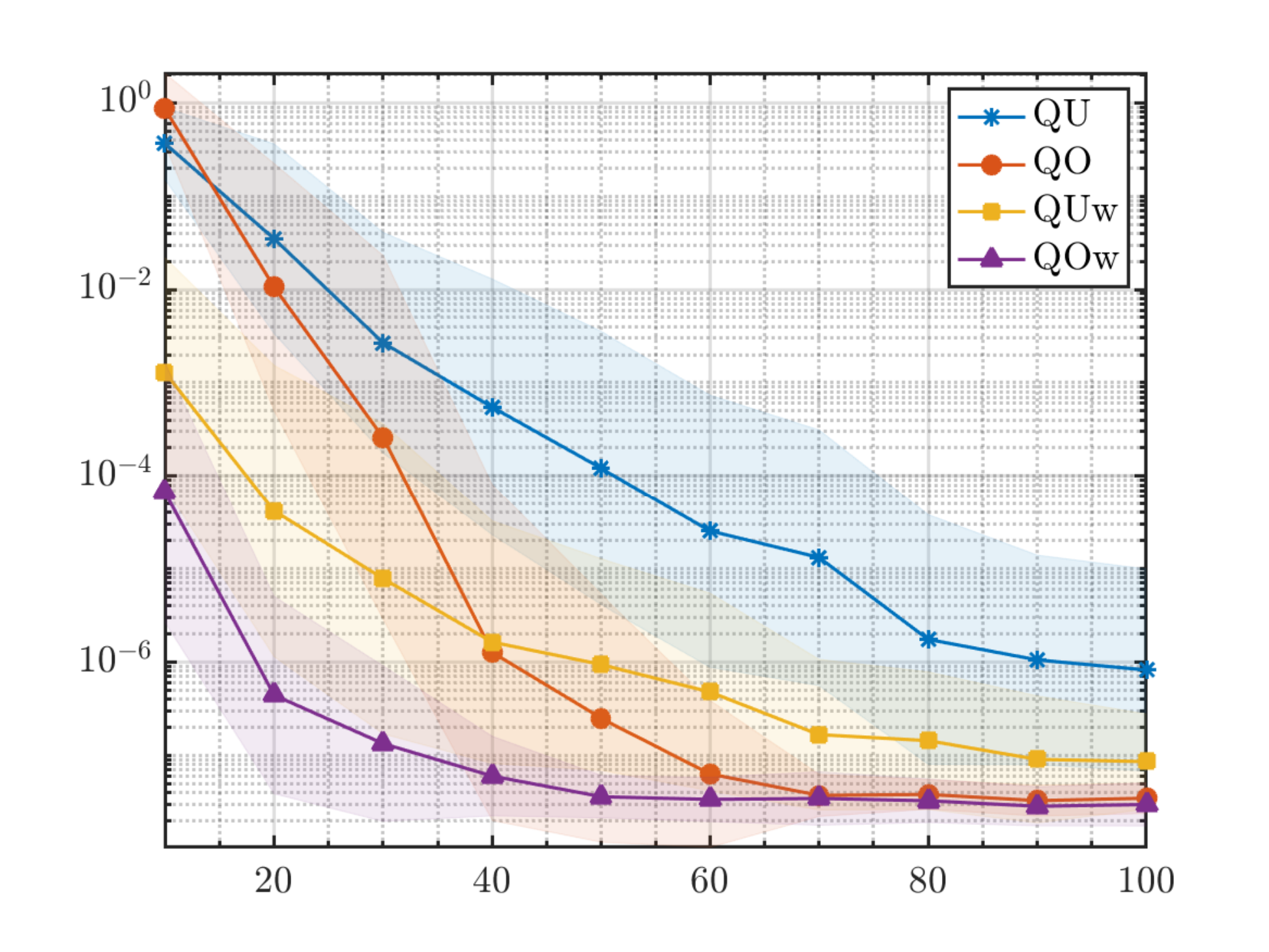} & 
 		\includegraphics[width=0.5\textwidth]{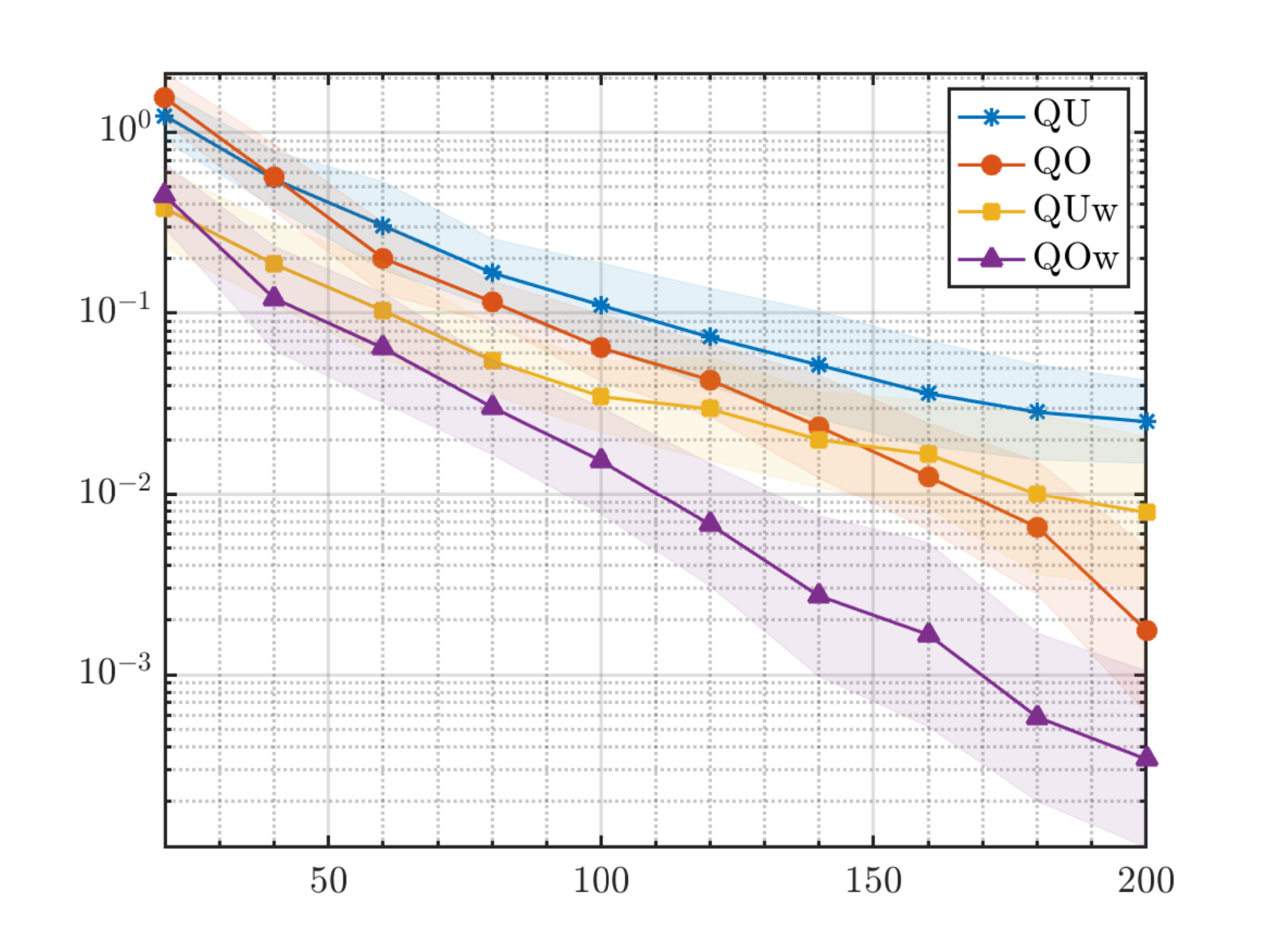} \\
 		$(d,n,N) = (1,399,400)$ & $(d,n,N)=(2,152,796)$     \\
 		\includegraphics[width=0.5\textwidth]{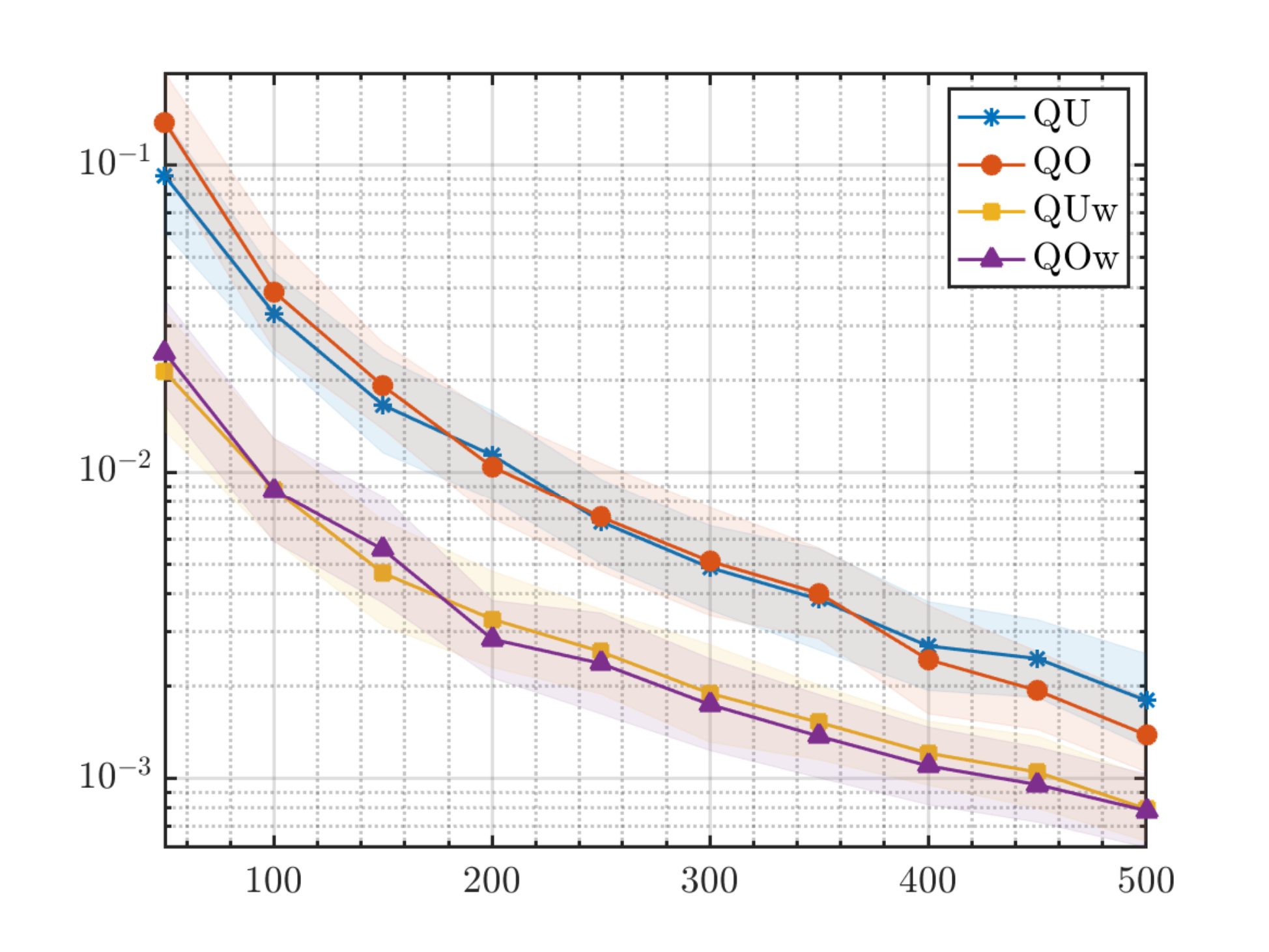} & 
 		\includegraphics[width=0.5\textwidth]{fig_16_2_1-eps-converted-to.pdf} \\
 		$(d,n,N) = (8,22,1843)$ & $(d,n,N)=(16,14,4385)$    \\
	\end{tabular}
	\end{small}
	\end{center}
\caption{The same as Fig,\ \ref{fig:irregular-exp-1-weighted}, except for $f = f_2$ and $D = D_1$.} 
\label{fig:irregular-exp-3-weighted}
\end{figure}

We conclude with a number of numerical experiments. In Figs.\ \ref{fig:irregular-exp-1-weighted}-\ref{fig:irregular-exp-3-weighted} we consider Example \ref{ex:alg-poly-general} and employ the orthogonalization strategy of \S \ref{s:sparse-irregular-new}. We compare unweighted and weighted $\ell^1$-minimization, where in the latter we set the weights $u = (u_{\iota})_{\iota \in \cI}$ to be
\be{
\label{weights-numerics}
u_{\iota} = \nm{\upsilon_{\iota}}_{L^{\infty}_{\tau}(\cU)},
}
where $\Upsilon$ is the basis constructed via the approach of \S \ref{s:sparse-irregular-new}. In other words, these weights follow the approach discussed in \S \ref{ss:lower-set-sparsity} for promoting lower set sparsity. 

For sampling, we consider Monte Carlo sampling \R{MC-CS-ONB} and the discrete `optimal’ measure \R{Optimal-CS-ONB}. In all examples, we see that weighted $\ell^1$-minimization substantially outperforms unweighted $\ell^1$-minimization. This is consistent with the observation that the larger polynomial coefficients tend to occur at smaller multi-indices – a property that the weights \R{weights-numerics} promote by assigning larger weights to higher multi-indices. In terms of sampling, we observe the sampling measure \R{Optimal-CS-ONB} outperforming Monte Carlo sampling \R{MC-CS-ONB}, where, as per usual, the benefit tends to lessen in higher dimensions. As discussed above, we do not claim that \R{Optimal-CS-ONB} is an optimal sampling measure in the weighted case: in fact, unlike in the unweighted case, it does not necessarily minimize the corresponding term $k(s;w)$ in the sample complexity bound. Yet, these experiments appear to suggest that it is a useful strategy when combined with weights to further enhance recovery of smooth functions via polynomials.

\section{Conclusions and challenges}\label{s:conclusions}

The purpose of this chapter has been to explore the question of optimal sampling for learning sparse approximations in high dimensions. In the more straightforward setting of Problem \ref{ass:sparsity_known}, we showed how this can be almost entirely resolved by defining a sampling measure (or measures) in terms of the Christoffel function of the corresponding subspace. We remark in passing recent work that strives to go even further, by removing the log factor in the sample complexity bound. See \cite{cohen2021optimal} 
and references therein. We note, however, that such procedures may not be feasible in practice, or may not guarantee quasi-optimal error bounds.

In the more challenging setting of Problem \ref{ass:sparsity_unknown}, we explored the limitations of Monte Carlo sampling, and showed how to obtain a sampling measure that optimized the sufficient condition of the number of measurements. Empirically, this leads to improved approximation, especially in lower-dimensional problems. Finally, we discussed structured sparsity, via either weighted or lower set sparsity, both of which can be promoted by using weights. Although here we were not even able to find a sampling measure to optimize the sample complexity bound, we found empirically that the same sampling measure used previously worked well in practice when combined with weights.

The major open problem raised by this work is therefore: is it possible to design sampling measures for sparse approximation in dictionaries that are theoretically optimal, with sample complexity bounds that scale log-linearly in $s$ and logarithmically in $n$? Currently, we have no answer to this question. We note in passing that it may be important to take into account more refined structured of the dictionary. See \cite{tran2018analysis} for recent work that uses the envelope bound \R{envelope-bound} to derive improved sample complexity bounds in the case of Example \ref{ex:alg-poly} with Monte Carlo sampling. It is also notable that the various constants $\theta$, $\Theta$ and $a,b$ (in the case of irregular domains) often very poorly explain the observed performance. This is particularly notable in the case of $a,b$. This raises the question of a more refined analysis that avoids these terms.

Let us also mention several extensions. First, while this work has focused on standard dictionaries consisting of algebraic or trigonometric polynomials, it is perfectly applicable to much more general dictionaries. This includes dictionaries now arising commonly in machine learning settings, such as random feature models or learned dictionaries obtained from deep neural network training. We note recent work on learning sparse representations in random feature models \cite{hashemi2021generalization}. An interesting question for future work involves applying the techniques considered herein to these models, to obtain better sampling strategies for such dictionaries. This may be highly relevant for applications using machine learning techniques that are data-starved. There is also the problem {of combining} sampling, via the strategies discussed herein, with learning the dictionary in an adaptive way to boost performance.

Second, we note that the sampling model explored in this work is simple pointwise evaluations. It is possible to extend this to much more general sampling models, many of which occur in practical settings. An example is the discrete-in-space-continuous-in-time model, which can occur when sensors in physical space take continuous recordings a time-dependent function $f(y,t)$. Another problem, which arises commonly in uncertainty quantification (see \cite{adcock2019compressive,guo2017gradient,peng2016polynomial} and references therein), is the problem where one measures both the function $f(y)$ and its gradient $\nabla f(y)$ simultaneously at a sample point $y$. We anticipate that many of the key results of this work can be extended to substantially more general sampling models.

\section*{Acknowledgements}

The authors wish to thank Simone Brugiapaglia for useful feedback. The authors acknowledges support from the Natural Sciences and Engineering Research Council of Canada (NSERC) through grant 611675 and the Pacific Institute for the Mathematical Sciences (PIMS).

\bibliographystyle{abbrv}
\bibliography{SamplingSpringerChptRefs}

\end{document}